# Topics in Normal Bases of Finite Fields

N. A. Carella

# Table of Contents

## Chapter 5 Periods Polynomials

## Chapter 6 Periods Normal Bases

## Chapter 7 Period Normal Bases For Extensions Of Low Degrees

## Chapter 8 Asymptotic Proofs

# Chapter 1

# Bases of Finite Fields

# 1.1 Introduction

This chapter introduces various fundamental ideas and terminologies essential for the understanding of vector representations of finite fields. The study of bases of vector space representations of finite fields and the corresponding computational algorithms is an extensive and important subject.

There are various methods of representing finite fields. The most common are vector spaces, cyclic representations, polynomial quotient rings, quotients of number fields, matrix representations, and binary representations respectively. These are listed here in order.

(1) $\mathbf{F}_{q^n} \cong \{ x = x_0\alpha_0 + x_1\alpha_1 + \cdots + x_{n-1}\alpha_{n-1} : x_i \in \mathbf{F}_q \}$, where $\{\alpha_0, \alpha_1, ..., \alpha_{n-1}\}$ is a basis.

(2) $\mathbf{F}_{q^n} \cong \langle \xi \rangle \cup \{0\}$, where $\xi$ is a generator of the multiplicative group of $\mathbf{F}_{q^n}$.

(3) $\mathbf{F}_{q^n} \cong \mathbf{F}_q[x]/(f(x))$, where $f(x)$ is an irreducible polynomial of degree $n$.

(4) $\mathbf{F}_{q^n} \cong \boldsymbol{O}_K/(\mathcal{J})$, where $\mathcal{J}$ is a maximal ideal and $\boldsymbol{O}_K$ is the ring of integers in a numbers field $\boldsymbol{K}$.

These representations are widely used in algebraic number theory.

(5) $\mathbf{F}_{q^n} \cong \{$ Subset of Nonsingular Matrices $\}$.

(6) $\mathbf{F}_{q^n} \cong \{$ *l–adic* Vectors $\}$, the vectors are defined by a function $\phi : \mathbf{F}_{q^n} \to \mathbf{F}_l^n$. Two instances are the 2–*adic* representation (binary):

$$\phi(x) = [(x^{q^n-1} + x^{(q^n-1)/2} + \alpha_{n-1})/2, ..., (x^{q^n-1} + x^{(q^n-1)/2} + \alpha_0)/2],$$

and the 3–*adic* representation:

$$\phi(x) = [(x+\alpha_{n-1})^{(q^n-1)/2}, ..., (x+\alpha_1)^{(q^n-1)/2}, (x+\alpha_0)^{(q^n-1)/2}]$$

where $\alpha_0, \alpha_1, \ldots, \alpha_{n-1} \in \mathbf{F}_{q^n}$ are fixed.

The fastest methods for addition and subtraction are implemented with vector space representations. And the fastest method for multiplications, divisions, discrete exponentiations, and certain root extractions are implemented with cyclic representations. In the other hand, the fastest algorithms for computing discrete logarithms in finite fields are implemented in polynomial quotient rings and quotient of number fields, see [1, Adleman and DeMarrais], [1, ElGammal] etc. Matrix representations have applications in the construction of hash functions, pseudo numbers generators, and others, see [3, Gieselmann]. Normal bases, (which are vector space representations), and Binary representations are useful in polynomial factorizations, see [1, Nieterreiter], [1, Camion], and [1, Ganz]. Since additions and subtractions are highly efficient operations with respect to most bases, the main focus is on multiplication and multiplicative inverse algorithms with respect to the various bases of the finite fields $\mathbf{F}_{q^n}$ over $\mathbf{F}_q$.

## 1.2 Definitions and Elementary Concepts

Several methods for identifying the bases of the vector space $\mathbf{F}_{q^n}$ over $\mathbf{F}_q$ will be considered in this section. The notion of *basis* of a vector space has already appeared in this text. The concept of basis and the related idea of linear independence, (also algebraic independence), are recurrent themes throughout mathematics.

***Definition 1.1.*** A subset of elements $\{\alpha_0, \alpha_1, ..., \alpha_{n-1}\} \subset \mathbf{F}_{q^n}$ is said to be a *basis* of the vector space $\mathbf{F}_{q^n}$ over $\mathbf{F}_q$ if and only if every element $\alpha \in \mathbf{F}_{q^n}$ can be uniquely written as a linear combination

$\alpha = a_0\alpha_0 + a_1\alpha_1 + \cdots + a_{n-1}\alpha_{n-1}$,

where $a_i \in \mathbf{F}_q$.

A *redundant basis* is a basis such that every element has a representation as linear combination but not necessarily unique. Some algorithms based on redundant bases are more efficient than those based on nonredundant bases. A redundant normal basis $\{\alpha_0, \alpha_1, ..., \alpha_{n-1}, 1\}$ of $\mathbf{F}_{q^n}$ over $\mathbf{F}_q$ is employed in [2, Gao, et al] to improve the exponentiation algorithm. A redundant normal basis permits multiple representations of the elements, e.g., $0 = \alpha_0 + \alpha_1, + \cdots + \alpha_{-n1} + 1$ among others if $q > 2$. Redundant bases are also used to represent integers in fast exponentiation algorithms, and real/complex numbers in numerical algorithms which implement carry free arithmetic operations.

***Definition 1.2.*** Let $\{\alpha_0, \alpha_1, ..., \alpha_{n-1}\}$ be a subset of $\mathbf{F}_{q^n}$. The *regular matrix representation* of the set $\{\alpha_0, \alpha_1, ..., \alpha_{n-1}\}$ is defined by the n×n matrix $A = \left( \sigma^i(\alpha_j) \right) = \left( \alpha_j^{q^i} \right)$.
The matrix

$$A = \begin{bmatrix} \alpha_0 & \alpha_1 & \alpha_2 & \cdot & \cdot & \cdot & \alpha_{n-1} \\ \alpha_0^q & \alpha_1^q & \alpha_2^q & \cdot & \cdot & \cdot & \alpha_{n-1}^q \\ \alpha_0^{q^2} & \alpha_1^{q^2} & \alpha_2^{q^2} & \cdot & \cdot & \cdot & \alpha_{n-1}^{q^2} \\ \cdot & \cdot & \cdot & \cdot & \cdot & \cdot & \cdot \\ \alpha_0^{q^{n-1}} & \alpha_1^{q^{n-1}} & \alpha_2^{q^{n-1}} & \cdot & \cdot & \cdot & \alpha_{n-1}^{q^{n-1}} \end{bmatrix}$$

occurs very frequently in the analysis of bases and matrix representations of linear functionals in finite fields.

***Definition 1.3.*** A pair of bases $\{\alpha_0, \alpha_1, ..., \alpha_{n-1}\}$ and $\{\beta_0, \beta_1, ..., \beta_{n-1}\}$ are equivalent if each $\beta_i = c\alpha_i$ for some constant $c \in \mathbf{F}_q$.

The equivalence class of each basis $\{\alpha_0, \alpha_1, ..., \alpha_{n-1}\}$ can be viewed as a point in (n−1)-dimensional projective space $P^{n-1}(\mathbf{F}_{q^n})$.

***Lemma 1.4.*** A subset $\{\alpha_0, \alpha_1, ..., \alpha_{n-1}\}$ of elements of $\mathbf{F}_{q^n}$ is a basis of the vector space $\mathbf{F}_{q^n}$ over $\mathbf{F}_q$ if and only if the n×n regular matrix representation A associated to $\{\alpha_0, \alpha_1, ..., \alpha_{n-1}\}$ is nonsingular.

Proof: Suppose that $\{\alpha_0, \alpha_1, ..., \alpha_{n-1}\}$ is a basis, and let $\beta \in \mathbf{F}_{q^n}$. Now consider the system of equations

$$\beta = \sum_{i=0}^{n-1} b_i\alpha_i, \quad \beta^q = \sum_{i=0}^{n-1} b_i\alpha_i^q, \quad ..., \quad \beta^{q^{n-1}} = \sum_{i=0}^{n-1} b_i\alpha_i^{q^{n-1}}.$$

Since the subset of elements $\{\alpha_0, \alpha_1, ..., \alpha_{n-1}\}$ is a basis, the system of equations, rewritten as a vector equation

$$\left(\beta, \beta^q, ..., \beta^{q^{n-1}}\right) = A\mathbf{b}.$$

has a unique solution $\mathbf{b} = (b_0, b_1, ..., b_{n-1})$. This implies that the matrix A is nonsingular. Conversely, if the matrix A is nonsingular, then the above system of equations has a unique solution. This in turn implies that each $\beta \in \mathbf{F}_{q^n}$ has a unique representation as a linear combination $\beta = b_0\alpha_0 + b_1\alpha_1 + \cdots + b_{n-1}\alpha_{n-1}$, $b_i \in \mathbf{F}_q$, so $\{\alpha_0, \alpha_1, ..., \alpha_{n-1}\}$ is a basis. ■

***Lemma 1.5.*** *(Basis lifting lemma)* A basis $\{\alpha_0, \alpha_1, ..., \alpha_{n-1}\}$ of $\mathbf{F}_{q^n}$ over $\mathbf{F}_q$ is also a basis of $\mathbf{F}_{q^{nk}}$ over $\mathbf{F}_{q^k}$ for all integers k such that gcd(k, n) = 1.

Proof: Let $a_0, a_1, ..., a_{n-1} \in \mathbf{F}_{q^k}$ and consider the system of equations

$$\begin{aligned}
a_0\alpha_0 + a_1\alpha_1 + a_2\alpha_2 + \cdots + a_{n-1}\alpha_{n-1} &= 0 \\
a_0\alpha_0^{q^k} + a_1\alpha_1^{q^k} + a_2\alpha_2^{q^k} + \cdots + a_{n-1}\alpha_{n-1}^{q^k} &= 0 \\
a_0\alpha_0^{q^{2k}} + a_1\alpha_1^{q^{2k}} + a_2\alpha_2^{q^{2k}} + \cdots + a_{n-1}\alpha_{n-1}^{q^{2k}} &= 0 \\
\cdots & \\
a_0\alpha_0^{q^{(n-1)k}} + a_1\alpha_1^{q^{(n-1)k}} + a_2\alpha_2^{q^{(n-1)k}} + \cdots + a_{n-1}\alpha_{n-1}^{q^{(n-1)k}} &= 0
\end{aligned}$$

Since gcd(k, n) = 1, the map $i \to ik$ is a permutation of {0, 1, 2, ..., n–1} and the matrices

$$A = \left( {}_j^{q^i} \right) \text{ and } A_k = \left( \alpha_j^{q^{ik}} \right)$$

are just rows permutations of each other. Moreover, because the regular matrix representation A attached to this basis is nonsingular, it follows that the system of equations has only a trivial solution $\mathbf{a} = (a_0,a_1,...,a_{n-1}) = (0,0,...,0)$. This proves the linear independence of $\{\alpha_0, \alpha_1, ..., \alpha_{n-1}\}$ over $\mathbf{F}_{q^k}$. ■

The constraint gcd(k, n) = 1 ensures that the subset $\{\alpha_0, \alpha_1, ..., \alpha_{n-1}\}$ remains linear independent over the larger ground field $\mathbf{F}_{q^k}$ and that $\alpha_i \notin \mathbf{F}_{q^k}$, $0 \le i < n$.

***Example 1.6.*** Let n = 5, q = 2 and let $\mathbf{F}_{q^5} = \{ a_0\alpha_0 + a_1\alpha_1 + a_2\alpha_2 + a_3\alpha_3 + a_4\alpha_4 : a_i \in \mathbf{F}_q \}$. Since gcd(5, k) = 1, for k = 1, 2, 3, and 4, the basis $\{\alpha_0, \alpha_1, \alpha_2, \alpha_3, \alpha_4\}$ can (be lifted to) a basis of the finite fields

$\mathbf{F}_{q^5}$ over $\mathbf{F}_q$, $\mathbf{F}_{q^{5\cdot 2}}$ over $\mathbf{F}_{q^2}$, $\mathbf{F}_{q^{5\cdot 3}}$ over $\mathbf{F}_{q^3}$, and $\mathbf{F}_{q^{5\cdot 4}}$ over $\mathbf{F}_{q^4}$, ...

et cetera. For k = 2, the permutation $\pi$ is given by $(0, 1, 2, 3, 4) \to (\pi(0), \pi(1), \pi(2), \pi(3), \pi(4)) = (0, 2, 4, 1, 3)$, and the matrices are

$$A = \begin{bmatrix} \alpha_0 & \alpha_1 & \alpha_2 & \alpha_3 & \alpha_4 \\ \alpha_0^q & \alpha_1^q & \alpha_2^q & \alpha_1^q & \alpha_4^q \\ \alpha_0^{q^2} & \alpha_1^{q^2} & \alpha_2^{q^2} & \alpha_3^{q^2} & \alpha_4^{q^2} \\ \alpha_0^{q^3} & \alpha_1^{q^3} & \alpha_2^{q^3} & \alpha_3^{q^3} & \alpha_4^{q^3} \\ \alpha_0^{q^4} & \alpha_1^{q^4} & \alpha_2^{q^4} & \alpha_3^{q^4} & \alpha_4^{q^4} \end{bmatrix} \text{ and } A_2 = \begin{bmatrix} \alpha_0 & \alpha_1 & \alpha_2 & \alpha_3 & \alpha_4 \\ \alpha_0^{q^2} & \alpha_1^{q^2} & \alpha_2^{q^2} & \alpha_3^{q^2} & \alpha_4^{q^2} \\ \alpha_0^{q^4} & \alpha_1^{q^4} & \alpha_2^{q^4} & \alpha_3^{q^4} & \alpha_4^{q^4} \\ \alpha_0^q & \alpha_1^q & \alpha_2^q & \alpha_1^q & \alpha_4^q \\ \alpha_0^{q^3} & \alpha_1^{q^3} & \alpha_2^{q^3} & \alpha_3^{q^3} & \alpha_4^{q^3} \end{bmatrix}.$$

Lifting a basis changes both the field extension from $\mathbf{F}_{q^n}$ to $\mathbf{F}_{q^{dn}}$ and the ground field from $\mathbf{F}_q$ to $\mathbf{F}_{q^d}$, $d \geq 1$. A related notion is that of a *complete basis*. In this situation the field extension $\mathbf{F}_{q^n}$ remains fixed but the ground field $\mathbf{F}_{q^d}$ varies as d varies over the divisors d of n.

***Definition 1.7.*** A basis $\{\alpha_0, \alpha_1, ..., \alpha_{r-1}\}$ of $\mathbf{F}_{q^n}$ over $\mathbf{F}_q$ is said to be a *complete basis* if is it a basis for $\mathbf{F}_{q^n}$ over $\mathbf{F}_{q^d}$ for all divisors d of n.

Current research in complete bases is limited to *completely* normal bases, confer the chapter on normal bases for more details.

***Theorem 1.8.*** (*Iterated basis theorem*) Let $\{\alpha_0, \alpha_1, ..., \alpha_{r-1}\}$ and $\{\beta_0, \beta_1, ..., \beta_{s-1}\}$ be a pair of bases of $\mathbf{F}_{q^r}$ and $\mathbf{F}_{q^s}$ over $\mathbf{F}_q$ respectively. Then $\{\alpha_0\beta_0, \alpha_0\beta_1, .., \alpha_{r-1}\beta_{s-1}\}$ is a basis of $\mathbf{F}_{q^{rs}}$ over $\mathbf{F}_q$.

Proof: Suppose the set $\{\alpha_0\beta_0, \alpha_0\beta_1, ..., \alpha_{r-1}\beta_{s-1}\}$ of n = rs elements is linearly dependent over $\mathbf{F}_q$ (not a basis), then there exists a nontrivial vector $(a_{i,j}) \neq (0, ..., 0)$, $0 < i, j < n$, such that

$$\sum_{i=0}^{r-1}\sum_{j=0}^{s-1} a_{i,j}\alpha_i\beta_j = \sum_{j=0}^{s-1}\left(\sum_{i=0}^{r-1} a_{i,j}\alpha_i\right)\beta_j$$

$$= \sum_{i=0}^{r-1}\left(\sum_{j=0}^{s-1} a_{i,j}\beta_j\right)\alpha_i = 0.$$

But $\{\alpha_0, \alpha_1, ..., \alpha_{r-1}\}$ and $\{\beta_0, \beta_1, ..., \beta_{s-1}\}$ are both linearly independent over $\mathbf{F}_q$, so the inner sums satisfy both

$$\sum_{i=0}^{r-1} a_{i,j}\alpha_i = 0 \quad \text{and} \quad \sum_{j=0}^{s-1} a_{i,j}\beta_j = 0$$

simultaneously; which implies that $(a_{i,j}) = 0$ for all pairs (i, j). This contradict the existence of the nontrivial vector $(a_{i,j}) \neq (0,0,..,0)$. The converse is also easy to verify. ■

***Lemma 1.9.*** Let $\{\alpha_0, \alpha_1, ..., \alpha_{n-1}\}$ be a basis of $\mathbf{F}_{q^n}$ over $\mathbf{F}_q$. Then there exists an element $\alpha_i \in \{\alpha_0, \alpha_1, ..., \alpha_{n-1}\}$ such that $Tr(\alpha_i) \neq 0$.

Proof: Suppose $\{\alpha_0, \alpha_1, ..., \alpha_{n-1}\}$ is a basis, and $Tr(\alpha_i) = 0$ for all $i \in \{0, 1, 2, ..., n-1\}$. Then $Tr(\alpha) = a_0Tr(\alpha_0) + a_1Tr(\alpha_1) + \cdots + a_{n-1}Tr(\alpha_{n-1}) = 0$ for all $\alpha = a_0\alpha_0 + a_1\alpha_1 + \cdots + a_{n-1}\alpha_{n-1} \in \mathbf{F}_{q^n}$. But this contradicts the fact that any finite field has elements of arbitrary traces. In fact, for any fixed $a \in \mathbf{F}_q$, the equation $Tr(\alpha) = a$ has $q^{n-1}$ distinct solutions in $\mathbf{F}_{q^n}$. ■

# 1.3 The Discriminants of Bases

The theory of discriminant of bases of finite fields is essentially the same as its counterpart in Algebraic Number Theory.

***Definition 1.10.*** Let $\{\alpha_0, \alpha_1, ..., \alpha_{n-1}\}$ be a subset of $\mathbf{F}_{q^n}$ and let T be the n×n matrix ( $Tr(\alpha_i\alpha_j)$ ), called the trace matrix representation of the subset $\{ \alpha_i \}$. The discriminant of this subset is defined by the determinant $disc(\{\alpha_i\}) = det(T)$ of the matrix T.

***Theorem 1.11.*** The subset of elements $\{\alpha_0, \alpha_1, ..., \alpha_{n-1}\}$ constitutes a basis of $\mathbf{F}_{q^n}$ over $\mathbf{F}_q$ if and only if $det(T) \neq 0$.

The proof of this result is fairly standard linear algebra, see [1, Lidl et al, p. 61.]. Under certain condition the regular matrix representation $A = \left(\sigma^i(\alpha_j)\right) = \left(\alpha_j^{q^i}\right)$, and the trace matrix representation $T = ( Tr(\alpha_i\alpha_j) )$ satisfy the relation $det(A)^2 = det(T)$. Accordingly, either of the inequalities $det(A) \neq 0$ or $det(T) \neq 0$ implies that $\{\alpha_0, \alpha_1, ..., \alpha_{n-1}\}$ is a basis.

***Lemma 1.12.*** Suppose $\{\alpha_0, \alpha_1, ..., \alpha_{n-1}\}$ is a subset of conjugate elements, then $det(T) = det(A)^2$.

Proof: Since the elements $\alpha_i$ are conjugates and the maps $\sigma^k$ are automorphisms, the expression $\sigma^k(\alpha_i\alpha_j) = \alpha_{i+k}\sigma^k(\alpha_j)$. Hence the (i, j)th entry in the matrix product $A\cdot A = ( \sigma^i(\alpha_k) )\cdot( \sigma^k(\alpha_j) )$ is

$$\sum_{k=0}^{n-1}\sigma^i(\alpha_k)\sigma^k(\alpha_j) = \sum_{k=0}^{n-1}\alpha_{i+k}\sigma^k(\alpha_j)$$

$$= \sum_{k=0}^{n-1}\sigma^k(\alpha_i\alpha_j) = Tr(\alpha_i\alpha_j).$$

Quod erat demostrandum. ■

***Theorem 1.13.*** Let $\{\alpha_0, \alpha_1, ..., \alpha_{n-1}\}$ and $\{\beta_0, \beta_1, ..., \beta_{n-1}\}$ be a pair of bases of $\mathbf{F}_{q^n}$ over $\mathbf{F}_q$. Then $disc(\{ \beta_i \}) = D^2 disc(\{ \alpha_i \})$, some $0 \neq D \in \mathbf{F}_q$.

Proof: Let n = 2, (it simplifies the argument and there is no loss in generality), and let $Gal(\mathbf{F}_{q^2}/\mathbf{F}_q) = \{ \sigma, \sigma^2 \}$ be the group of automorphisms of $\mathbf{F}_{q^2}$. Since each pair is a basis, there is a 2×2 nonsingular matrix such that

$\beta_0 = a\alpha_0 + b\alpha_1$,
$\beta_1 = c\alpha_0 + d\alpha_1$.

Therefore the matrix $B = ( \sigma^i(\beta_j) )$ attached to the set $\{ \beta_1, \beta_0 \}$ is

$$\begin{bmatrix} \sigma(\beta_0) & \sigma^2(\beta_0) \\ \sigma(\beta_1) & \sigma^2(\beta_1) \end{bmatrix} = \begin{bmatrix} a\sigma(\alpha_0)+b\sigma(\alpha_1) & a\sigma^2(\alpha_0)+b\sigma^2(\alpha_1) \\ c\sigma(\alpha_0)+d\sigma(\alpha_1) & c\sigma^2(\alpha_0)+d\sigma^2(\alpha_1) \end{bmatrix}$$
$$= \begin{bmatrix} a & b \\ c & d \end{bmatrix} \begin{bmatrix} \sigma(\alpha_0) & \sigma^2(\alpha_0) \\ \sigma(\alpha_1) & \sigma^2(\alpha_1) \end{bmatrix}.$$

Taking determinants in both sides yields the claim, with $D = ad - bc \neq 0$. ■

The discriminant function induces an equivalence relation on the set $\mathcal{B}$ = { Bases of $\mathbf{F}_{q^n}$ over $\mathbf{F}_q$ }. One of the equivalence classes consists of all bases which have equal discriminant or are the same up to a linear transformation L of determinant det(L) = 1; these bases are referred to as *integral bases* in Algebraic Number Theory.

The discriminant $disc(\{\alpha^i\})$ of the power basis $\{1, \alpha, \alpha^2, ..., \alpha^{n-1}\}$ coincides with the discriminant of the minimal (characteristic) polynomial of α, . Two well known formulae are presented now.

***Theorem 1.14.*** Let α be a root of $f(x) \in \mathbf{F}_q[x]$. Then the discriminant $disc(\{\alpha^i\})$ of the power basis $\{1, \alpha, \alpha^2, ..., \alpha^{n-1}\}$ is given by

( 1 ) $disc(\{\alpha^i\}) = \prod_{o \le i < j < n} \left(\alpha^{q^i} - \alpha^{q^j}\right)^2$ .

( 2 ) $disc(\{\alpha^i\}) = (-1)^{n(n-1)/2} N(f'(\alpha))$,

where $N : \mathbf{F}_{q^n} \to \mathbf{F}_q$ is the norm, and $f'(x)$ is the derivative of f(x).

.

Proof: Let $\alpha_0, \alpha_1, ..., \alpha_{n-1}$ be the roots of f(x). Then the matrix $A = ( \sigma^i(\alpha^j) )$ is a van der Monde matrix, and the matrix $T = ( Tr(\alpha_i\alpha_j) ) = AA^T$, so the discriminant $disc(\{\alpha^i\}) = \det^2(A)$ is formula (1) above. Since there are n(n − 1)/2 ways of choosing a pair $\alpha_i$, $\alpha_j$ , and $f'(\alpha_i) = \prod_{i \neq j = 0}^{n-1} (\alpha_i - \alpha_j)$, the product of the indexed derivatives leads to

$$\prod_{i=0}^{n-1} f'(\alpha_i) = \prod_{i=0}^{n-1} \prod_{i \neq j=0}^{n-1} (\alpha_i - \alpha_j) = (-1)^{n(n-1)/2} \prod_{0 \le i < j < n} (\alpha_i - \alpha_j)^2 .$$

Similarly

$$\prod_{i=0}^{n-1} f'(\alpha_i) = \prod_{i=0}^{n-1} \sigma^i(f'(\alpha)) = N(f'(\alpha)),$$

where $\sigma^i$ is an automorphism. These complete the proofs. ■

It is easy to verify that the matrix T = ( $Tr(\alpha_i\alpha_j)$ ) is a bilinear form over $\mathbf{F}_q$, and for odd prime powers $q = p^t$, it can be classified as one of two types: It can be either a matrix whose determinant det(T) is a (square) quadratic residue or it can be a matrix whose determinant det(T) is a (nonsquare) nonquadratic residue, both in $\mathbf{F}_q$.

***Lemma 1.15.*** (*Artin 1957*) Let q be an odd prime power, then there are exactly two equivalence classes of nondegenerate symmetric linear forms on an n-dimensional vector space over $\mathbf{F}_q$ represented by the identity matrix $I_n$ = diag(111,…,1) and the diagonal matrix V = diag(1,1,1,...,v), where $v \in \mathbf{F}_q$ is a nonsquare.

This property of the matrix T is employed in [2, Jungnickel *et al.*] to simplify the proof of the self-dual basis theorem. The more general T = ( $Tr(\alpha_i\beta_j)$ ) attached to a pair of subsets $\{\alpha_0, \alpha_1, ..., \alpha_{n-1}\}$ and $\{\beta_0, \beta_1, ..., \beta_{n-1}\}$ also lends itself for the investigation of dual properties of these subsets.

Abel's formula in finite characteristic is a relationship between the determinant of the regular matrix representation *A* attached to a subset of the elements$\{\alpha_{n-1}, \ldots, \alpha_1, \alpha_0\}$and the coefficients of the minimal linear polynomial of the subset of the elements. The case $n = 2$ involves a second order linear differential equation

$$a_2 x^{q^2} + a_1 x^q + a_0 x = 0, \qquad a_2 a_0 \neq 0.$$

The solution space is spanned by any two roots $\{\alpha_1, \alpha_0\}$ of the linear polynomial. Thus

$$\alpha_0^{q^2} + b_1\alpha_0^q + b_0\alpha_0 = 0$$
$$\alpha_1^{q^2} + b_1\alpha_1^q + b_0\alpha_1 = 0,$$

where $b_1 = a_1 a_2$, $b_0 = a_0 a_2$. Adding multiples of the two previous equations result in

$$\left(\alpha_1^{q^2}\alpha_0^q - \alpha_0^{q^2}\alpha_1^q\right) + b_0\left(\alpha_1\alpha_0^q - \alpha_0\alpha_1^q\right) = 0.$$

This is equivalent to $a_1 D^q - a_0 D = 0$, where

$$D = \det\begin{bmatrix} \alpha_0 & \alpha_1 \\ \alpha_0^q & \alpha_1^q \end{bmatrix}.$$

Additional details appear in [1, Goss, p. 22], and the general case is as follows.

***Lemma 1.16.*** Let $\{\alpha_{n-1}, \ldots, \alpha_1, \alpha_0\}$ be any *n* roots of the linear polynomial $f(x) = \sum_{i=0}^{n} a_i x^{q^i}$, and let $D = det(A)$. Then $a_n D^q + (-1)^{n+1} a_0 D = 0$.

Proof: Let $a_n a_0 \neq 0$, otherwise the subset is linearly dependent, and let

$$A = \begin{bmatrix} \alpha_0 & \alpha_1 & \cdots & \alpha_{n-1} \\ \alpha_0^q & \alpha_1^q & \cdots & \alpha_{n-1}^q \\ \vdots & \vdots & \cdots & \vdots \\ \alpha_0^{q^{n-1}} & \alpha_1^{q^{n-1}} & \cdots & \alpha_{n-1}^{q^{n-1}} \end{bmatrix} \quad \text{and} \quad A^q = \begin{bmatrix} \alpha_0^q & \alpha_1^q & \cdots & \alpha_{n-1}^q \\ \alpha_0^{q^2} & \alpha_1^{q^2} & \cdots & \alpha_{n-1}^{q^2} \\ \vdots & \vdots & \cdots & \vdots \\ \alpha_0^{q^n} & \alpha_1^{q^n} & \cdots & \alpha_{n-1}^{q^n} \end{bmatrix}.$$

Replacing $\alpha_i^{q^n} = -\sum_{j=0}^{n-1} a_j \alpha_i^{q^j}$ in the matrix $A^q$, and simplifying the determinant proves the claim.

■

## 1.4 Distribution of Bases

There is a close link between the numbers of various types of bases of $\mathbf{F}_{q^n}$ over $\mathbf{F}_q$ and the size of various subgroups of the general linear group $GL_n(\mathbf{F}_q)$ of n×n nonsingular matrices over $\mathbf{F}_q$. The $GL_n(\mathbf{F}_q)$-orbit of any basis contains every bases of the vector space $\mathbf{F}_{q^n}$. As the matrix A = ( $a_{i,j}$ ) varies in $GL_n(\mathbf{F}_q)$, the change of basis $\alpha_i \rightarrow \beta_i$ given by

$\beta_0 = a_{0,0}\alpha_0 + a_{0,1}\alpha_1 + a_{0,2}\alpha_2 + \cdots + a_{0,n-1}\alpha_{n-1}$,
$\beta_1 = a_{1,0}\alpha_0 + a_{1,1}\alpha_1 + a_{1,2}\alpha_2 + \cdots + a_{1,n-1}\alpha_{n-1}$,
$\beta_2 = a_{2,0}\alpha_0 + a_{2,1}\alpha_1 + a_{2,2}\alpha_2 + \cdots + a_{2,n-1}\alpha_{n-1}$,
…
$\beta_{n-1} = a_{n-1,0}\alpha_0 + a_{n-1,1}\alpha_1 + a_{n-1,2}\alpha_2 + \cdots + a_{n-1,n-1}\alpha_{n-1}$,

varies over all the bases of $\mathbf{F}_{q^n}$. Other relationships similar to this will appear in the investigation of bases. The precise statement about the size of the collection of bases is given below.

***Lemma 1.17.*** The vector space $\mathbf{F}_{q^n}$ over $\mathbf{F}_q$ has a total number of

$$B_q(n) = \prod_{i=0}^{n-1} \left(q^n - q^i\right) = q^{n(n-1)/2} \prod_{i=1}^{n} \left(q^i - 1\right)$$

distinct ordered bases.

Proof: The $GL_n(\mathbf{F}_q)$-orbit of a single basis consists of all the bases of $\mathbf{F}_{q^n}$ over $\mathbf{F}_q$. Thus $B_q(n) = \#GL_n(\mathbf{F}_q)$. ■

The number of bases $\mathbf{F}_{q^n}$ over $\mathbf{F}_q$ grows exponentially as a function of $n$, and very quickly becomes unmanageable. From $B_q(n) = q^{n^2}\prod_{i=1}^{n}\left(1-q^{-i}\right)$, it is clear that the estimates

$$q^{n(n-1)} < B_q(n) < q^{n^2}$$

holds for all pairs $(n, q)$. For tabulation purpose, it is more convenient to consider the number of unordered bases, which is given by $B_q(n)/n!$.

## 1.5 Dual Bases

***Definition 1.18.*** Let x, y $\in \mathbf{F}_{q^n}$ be a pair of vectors and let $\mu$ be a fixed element. The *trace inner product* of x and y is defined as the trace $\mathrm{Tr}(\mu xy)$ of the triple product of x, y, and $\mu$. Two distinct elements are said to be *trace orthogonal* if $\mathrm{Tr}(\mu xy) = 0$ in $\mathbf{F}_q$.

***Definition 1.19.*** A pair of bases $\{\alpha_0, \alpha_1, ..., \alpha_{n-1}\}$ and $\{\beta_0, \beta_1, ..., \beta_{n-1}\}$ of $\mathbf{F}_{q^n}$ over $\mathbf{F}_q$ are said to be dual bases if the trace orthogonal relation

$$Tr(\mu\alpha_i\beta_j) = \begin{cases} c_i & \text{if } i = j, \\ 0 & \text{if } i \neq j, \end{cases}$$

holds for some fixed $\mu \in \mathbf{F}_{q^n}$, (note that $\mathrm{Tr}(\mu\alpha_i ß_i) = ( c_i\delta_{ij} )$ is a diagonal matrix).

The special case $\mathrm{Tr}(\alpha_i ß_i) = \delta_{i,j}$ is referred to as *dual bases*. Moreover, if $\alpha_i = \beta_i$ the basis is called a *self-dual basis*. This corresponds to the simplest linear functional $\mathrm{Tr}((\mu\alpha_i ß_i) = c_i\ \delta_{i,j}$, with the parameters $\mu = 1$, and $c_0 = c_1 = \cdots = c_{n-1} = 1$.

The notion of dual bases is analogous to the notion of orthonormal bases in vector spaces over complex numbers $\mathbb{C}$. If the two elements x, y $\in \mathbf{F}_{q^n}$ are written in terms of a pair of dual bases $\{\alpha_0, \alpha_1, ..., \alpha_{n-1}\}$ and $\{\beta_0, \beta_1, ..., \beta_{n-1}\}$ as $x = x_0\alpha_0 + x_1\alpha_1 + \cdots + x_{n-1}\alpha_{n-1}$ and $y = y_0\beta_0 + y_1\beta_1 + \cdots + y_{n-1}\beta_{n-1}$, $x_i, y_i \in \mathbf{F}_q$, then the trace inner product is given by

$$(x, y) = \mathrm{Tr}(xy) = x_0y_0 + x_1y_1 + \cdots + x_{n-1}y_{n-1}.$$

This is analogous to the standard inner product in vector spaces over the complex numbers $\mathbb{C}$ with respect to orthonormal bases.

It is a simple matter to demonstrate that any finite field has a pair of dual bases. If $x \in \mathbf{F}_{q^n}$ and $\{\alpha_0, \alpha_1, ..., \alpha_{n-1}\}$ is a basis, write

$$\begin{aligned} x &= x_0\alpha_0 + x_1\alpha_1 + \cdots + x_{n-1}\alpha_{n-1} \\ &= c_0(x)\alpha_0 + c_1(x)\alpha_1 + \cdots + c_{n-1}(x)\alpha_{n-1}, \end{aligned}$$

where $c_i(x) = x_i$ is the coefficient function with respect to this basis. The existence of the dual basis $\{\beta_0, \beta_1, ..., \beta_{n-1}\}$ is a consequence of the fact that $c_i : \mathbf{F}_{q^n} \rightarrow \mathbf{F}_q$ is a linear functional on $\mathbf{F}_{q^n}$ ; and the fact that every linear functional on $\mathbf{F}_{q^n}$ has a unique trace representation. It immediately follows that $c_i(x) = Tr(ß_i x) = a_i$ for some unique $\beta_i \in \mathbf{F}_{q^n}$. Specifically $Tr(\alpha_i\beta_j) = \delta_{i,j}$. Hence, the subset $\{\beta_0, \beta_1, ..., \beta_{n-1}\} \subset \mathbf{F}_{q^n}$ is a unique dual basis of $\{\alpha_0, \alpha_1, ..., \alpha_{n-1}\}$. The previous observations proves the following.

***Theorem 1.20.*** (*Dual basis theorem*) Every basis $\{\alpha_0, \alpha_1, ..., \alpha_{n-1}\}$ of $\mathbf{F}_{q^n}$ over $\mathbf{F}_q$ has a unique dual basis $\{\beta_0, \beta_1, ..., \beta_{n-1}\}$.

Let $n = rs$, $\gcd(r, s) = 1$, and let $\eta \in \mathbf{F}_{q^r}$ and $\theta \in \mathbf{F}_{q^s}$. The projected traces of $\eta\theta$ in $\mathbf{F}_{q^r}$ and $\mathbf{F}_{q^s}$ are the elements

$$Tr_{n:r}(\eta\theta) = \eta Tr_{n:r}(\theta) = \eta Tr_{s:1}(\theta),$$

and

$$Tr_{n:s}(\eta\theta) = \theta Tr_{n:s}(\eta) = \theta Tr_{r:1}(\eta),$$

where $Tr_{ab:b} : \mathbf{F}_{q^{ab}} \rightarrow \mathbf{F}_{q^b}$ is the relative trace function defined by

$$Tr_{ab:b}(x) = x + x^{q^b} + x^{q^{2b}} + \cdots + x^{q^{(a-1)b}}.$$

These relationships will be employed to establish the iterated trace orthogonal dual basis theorem.

***Theorem 1.21.*** (*Iterated trace orthogonal dual basis theorem*) Let $\{\alpha_0, \alpha_1, ..., \alpha_{r-1}\}$, $\{\beta_0, \beta_1, ..., \beta_{r-1}\}$ and $\{\eta_0, \eta_1, ..., \eta_{s-1}\}$, $\{\theta_0, \theta_1, ..., \theta_{s-1}\}$ be two pairs of dual bases of $\mathbf{F}_{q^r}$ and $\mathbf{F}_{q^s}$ over $\mathbf{F}_q$ respectively. Then $\{\alpha_i\eta_j\}$ and $\{\beta_i\theta_j\}$, $0 \le i < r$, $0 \le j < s$, is a pair of dual bases of $\mathbf{F}_{q^{rs}}$ over $\mathbf{F}_q$ if and only if $\gcd(r, s) = 1$.

Proof: By the iterated bases theorem it is known that $\{\alpha_i\eta_j\}$ and $\{\beta_i\theta_j\}$ are both bases of $\mathbf{F}_{q^{rs}}$ over $\mathbf{F}_q$. Moreover, by hypothesis $\{\alpha_i\}$ and $\{\beta_i\}$ are dual bases of $\mathbf{F}_{q^r}$ over $\mathbf{F}_q$, or equivalently

(1) $Tr_{r:1}(\alpha_i\beta_j) = \delta_{i,j}$, and $Tr_{n:s}(\alpha_i\beta_j) = Tr_{r:1}(\alpha_i\beta_j)$,
if and only if gcd(r, s) = 1. Similarly

(2) $Tr_{s:1}(\eta_u\theta_v) = \delta_{u,v}$, and $Tr_{n:r}(\eta_u\theta_v) = Tr_{s:1}(\eta_u\theta_v)$.
Now proceed to compute the dual basis relation:

$$\begin{aligned} Tr_{rs:1}(\alpha_i\eta_u\beta_j\theta_v) &= Tr_{s:1}(Tr_{rs:s}(\alpha_i\eta_u\beta_j\theta_v)) \\ &= Tr_{s:1}(\eta_u\theta_v Tr_{rs:s}(\alpha_i\beta_j)) \\ &= \delta_{i,j} Tr_{s:1}(\eta_u\theta_v) = \delta_{i,j}\delta_{u,v}. \end{aligned}$$

Thus the matrix $Tr_{rs:1}(\alpha_i\eta_u\beta_j\theta_v) = \delta_{ur+i,vr+j}$ is the identity matrix. ■

The power basis $\{1, \alpha, \alpha^2, ..., \alpha^{n-1}\}$ of $\mathbf{F}_{q^n}$ over $\mathbf{F}_q$, which corresponds to the standard basis of the vector space $\mathbf{F}_{q^n}$, is the most common and widely used basis.

***Theorem 1.22.*** *(Power Dual Basis Theorem)* Let $\{1, \alpha, \alpha^2, ..., \alpha^{n-1}\}$ be the power basis of $\mathbf{F}_{q^n}$ over $\mathbf{F}_q$ and let $f(x) = (x - \alpha)(b_{n-1}x^{n-1} + \cdots + b_1x + b_0) \in \mathbf{F}_q[x]$ be the minimum polynomial of α. Then the list of elements

$$\beta_0 = \frac{b_0}{f'(\alpha)}, \quad \beta_1 = \frac{b_1}{f'(\alpha)}, \quad ..., \quad \beta_{n-1} = \frac{b_{n-1}}{f'(\alpha)}$$

forms a unique dual basis of the power basis over $\mathbf{F}_q$.

Proof: The polynomial

$$g(x) = x^i - \sum_{j=0}^{n-1} \frac{\alpha^{iq^j}}{f'(\alpha^{q^j})} \frac{f(x)}{(x-\alpha^{q^j})},$$

$0 \le i < n$, of degree deg(g(x)) < n has the same number of zeros as f(x), so g(x) = 0. Rewritten in term of the trace function this becomes

$$\begin{aligned} Tr\left(\frac{\alpha^i}{f'(\alpha)}\frac{f(x)}{x-\alpha}\right) &= Tr\left((b_{n-1}x^{n-1} + \cdots + b_1x + b_0)\frac{\alpha^i}{f'(\alpha)}\right) \\ &= \sum_{j=0}^{n-1} Tr\left(\frac{\alpha^i}{f'(\alpha)}\right) b_j x^j = x^i. \end{aligned}$$

Matching coefficients in both sides of the last equation returns

$$Tr\left(\alpha^i \frac{b_j}{f'(\alpha)}\right) = Tr(\alpha^i\beta_j) = \delta_{i,j}.$$ ■

This proof uses the fact that the ring $\mathbf{F}_q[x^{-1}]$ is σ-invariant, namely, $\sigma^i(x^k) = x^k$ for all i, k ≥ 0; the other possible action of the automorphism $\sigma^i(x) = x^{q^i}$ on $\mathbf{F}_q[x^{-1}]$ also works. The proof presented here is an adaptation to the structure of the roots of polynomials over finite fields. The more general version for separable extensions $\mathbf{F}(\alpha)$ of arbitrary fields **F** is given in [1, Lange, p. 322]. A variation of this proof appears in [1, Menezes, p. 6].

***Example 1.23.*** Let q = 2, n + 1 = 5, and let α be a root of $f(x) = x^4 + x^3 + x^2 + x + 1 \in \mathbf{F}_2[x]$. Since 2 has order 4 modulo 5, f(x) is irreducible over $\mathbf{F}_2$, and

$$\begin{aligned} f(x) &= x^4 + x^3 + x^2 + x + 1 \\ &= (x - \alpha)(b_3x^3 + b_2x^2 + b_1x + b_0) \\ &= b_3x^4 + (\alpha b_3+b_2)x^3 + (\alpha b_2+b_1)x^2 + (\alpha b_1+b_0)x + \alpha b_0. \end{aligned}$$

Matching coefficients returns

$b_0 = \alpha^{-1}$,
$b_1 = \alpha^{-2} + \alpha^{-1}$,
$b_2 = \alpha^{-3} + \alpha^{-2} + \alpha^{-1}$,
$b_3 = \alpha^{-4} + \alpha^{-3} + \alpha^{-2} + \alpha^{-1}$.

And $f'(x) = x^2 + 1$, so $1/f'(\alpha) = \alpha^2 + \alpha$. Hence the dual of the power basis is

$\delta_0 = b_0/f'(\alpha) = \alpha^{-1} (\alpha^2 + \alpha) = \alpha + 1$,
$\delta_1 = b_1/f'(\alpha) = (\alpha^{-2} + \alpha^{-1})(\alpha^2 + \alpha) = \alpha + \alpha^{-1}$,
$\delta_2 = b_2/f'(\alpha) = (\alpha^{-3} + \alpha^{-2} + \alpha^{-1})(\alpha^2 + \alpha) = \alpha + \alpha^{-2}$,
$\delta_3 = b_3/f'(\alpha) = (\alpha^{-4} + \alpha^{-3} + \alpha^{-2} + \alpha^{-1})(\alpha^2 + \alpha) = \alpha + \alpha^{-3}$.

The calculations of the dual basis $\delta_0, \delta_1, \ldots, \delta_{n-1}$ of the power basis $1, \alpha, \alpha^2, \ldots, \alpha^{n-1}$ for the parameters n + 1 = prime, and 2 of order n modulo n + 1, (f(x) is irreducible over $\mathbf{F}_2$), are similar to these.

***Theorem 1.24.*** The finite field $\mathbf{F}_{q^n}$ has a pair of self-dual bases for the following parameters.

(1) q is an even prime power.
(2) q is an odd prime power and n = 2k + 1.

***Theorem 1.25.*** (*Imamura 1983*) The finite field $\mathbf{F}_{q^n}$ has no self-dual power bases.

Proof: A self-dual power basis $\{ 1, \alpha, \alpha^2, \ldots, \alpha^{n-1} \} \Rightarrow Tr(1\cdot\alpha^2) = 0$ and $Tr(\alpha\cdot\alpha) = 1$ simultaneously, which is a contradiction. ■

***Theorem 1.26.*** (*Geiselmann, Gollmann 1993*) The dual basis of the power (polynomial) basis $\{1, \alpha, \alpha^2, ..., \alpha^{n-1}\}$ is a polynomial basis if and only if $n \not\equiv 0 \bmod p$, p being the characteristic of $\mathbf{F}_{q^n}$, and the minimal polynomial of α has the form $f(x) = x^n - c \in \mathbf{F}_q[x]$.

**Weakly Self-Dual Bases**

The matrix $T = ( Tr(\alpha_i\beta_j) )$ attached to a pair of dual bases $\{\alpha_0, \alpha_1, ..., \alpha_{n-1}\}$ and $\{\beta_0, \beta_1, ..., \beta_{n-1}\}$ of $\mathbf{F}_{q^n}$ over $\mathbf{F}_q$ is the identity matrix $( Tr(\alpha_i\beta_j) ) = I_n$ in the general linear group $GL_n(\mathbf{F}_q)$. But the matrix $T = ( Tr(\alpha_i\beta_j) )$ attached to a pair of trace orthogonal bases is a nonsingular diagonal matrix $( Tr(\alpha_i\beta_j) ) = diag(c_0,c_1,...,c_{n-1})$. The result below extends the idea of trace orthogonal bases to allow monomial matrices, which are permutations of diagonal matrices.

***Definition 1.27.*** A pair of bases $\{\alpha_0, \alpha_1, ..., \alpha_{n-1}\}$ and $\{\beta_0, \beta_1, ..., \beta_{n-1}\}$ of $\mathbf{F}_{q^n}$ over $\mathbf{F}_q$ are called weakly self-dual bases if there is an element $t \in \mathbf{F}_{q^n}$, $c_0, c_1, ..., c_{n-1} \in \mathbf{F}_q$, and a permutation π(i) of {0, 1, 2, 3, .., n–1} such that $\beta_i = tc_i\alpha_{\pi(i)}$, i = 0, 1, .., n – 1.

***Theorem 1.28.*** (*Morgan et al, 1997*) The finite field $\mathbf{F}_{q^n}$ of $q^n$ elements has a weakly self-dual basis if and only if there exists an irreducible polynomial $f(x) = x^n - cx^k - d \in \mathbf{F}_q[x]$, $0 < k < n$. The permutation π of {0, 1, 2, 3, ..., n–1} is given by $\pi(i) \equiv k - 1 - i \bmod n$ if $c \neq 0$, and $\pi(i) \equiv -1 - i \bmod n$ if $c = 0$.

The last result generalizes the work of [ Geiselmann ], which establishes it for the polynomial $f(x) = x^n - cx^k - 1$ or $x^n - d$.

A pair of dual bases permit a dual representations of the elements of $\mathbf{F}_{q^n}$ as

$$\alpha = \sum_{i=0}^{n-1} Tr(\alpha\beta_i)\alpha_i = \sum_{i=0}^{n-1} Tr(\alpha\alpha_i)\beta_i .$$

Thus a conversion from one basis to the other basis involves the calculations of all the traces $Tr(\alpha\beta_i)$ or $Tr(\alpha\alpha_i)$. On the other hand, a pair of weakly dual bases permit representations of the elements of $\mathbf{F}_{q^n}$ as

$$\alpha = \sum_{i=0}^{n-1} x_i\alpha_i = \sum_{i=0}^{n-1} x_{\pi(i)}\alpha_{\pi(i)} = t^{-1}\sum_{i=0}^{n-1} c_i x_{\pi(i)}\beta_i .$$

Accordingly, a conversion from one basis to the other basis involves n multiplications in $\mathbf{F}_q$, n permutations of the coefficients, and one inversion, (worst case).
The multiplication complexity of weakly self dual bases is discussed in [1, Shparlinski, p.100].

## 1.6 Distribution of Dual Bases

For a pair of bases $\{\alpha_0, \alpha_1, ..., \alpha_{n-1}\}$ and $\{ß_0, ß_1, ..., ß_{n-1}\}$ and a nonsingular matrix $A = (a_{i,j}) \in GL_n(\mathbf{F}_q)$, consider the linear expressions

$$\alpha_i = \sum_{j=0}^{n-1} a_{i,j}\beta_j$$

and

$$\beta_i = \sum_{j=0}^{n-1} a_{i,j}\alpha_j \ .$$

The trace inner product of the pair $\alpha_i$ and $ß_j$ is given by

$$Tr(\alpha_i\beta_j) = Tr\left(\sum_{s=0}^{n-1} a_{i,s}\alpha_s\right)\left(\sum_{t=0}^{n-1} a_{j,t}\beta_t\right)$$
$$= \sum_{s,t=0}^{n-1} a_{i,s}a_{j,t}Tr(\alpha_s\beta_t).$$

Specializing this equation to dual basis, namely $Tr(\alpha_i ß_j) = \delta_{i,j}$, leads to the dual bases equation

$$Tr(\alpha_i\beta_j) = \sum_{s=0}^{n-1} a_{i,s}a_{j,s} .$$

This equation classifies the set of dual bases as the orbit of the orthogonal group $O_n(\mathbf{F}_q) = \{ A \in GL_n(\mathbf{F}_q) : AA^T = I_n \}$. In particular, if a finite field $\mathbf{F}_{q^n}$ has a self-dual basis $\{\alpha_i\}$ (it occurs only for certain combinations of n and q), then the orbit $Orb(\{\alpha_i\}) = \{ A\{\alpha_i\} : A \in O_n(\mathbf{F}_q) \}$ under the group of all n×n orthogonal matrices contains all the self-dual bases of $\mathbf{F}_{q^n}$ over $\mathbf{F}_q$.

***Theorem 1.29.*** Let n > 1, and let q be a prime power, then there are

$$SD_q(n) = \begin{cases} \prod_{i=1}^{n-1}\left(q^i - \left(\frac{1+(-1)^i}{2}\right)\right) & \text{if } q \equiv 0 \mod 2 \\ 2\prod_{i=1}^{n-1}\left(q^i - \left(\frac{1+(-1)^i}{2}\right)\right) & \text{if } n, q \equiv 1 \mod 2, \\ 0 & \text{otherwise,} \end{cases}$$

distinct ordered self-dual bases of $\mathbf{F}_{q^n}$ over $\mathbf{F}_q$.

A derivation of this equation appears in [2, 3, 4, Jungnickel et al.]. The integer $SD_q(n)$ is precisely the size of the orthogonal group $O_n(\mathbf{F}_q)$ whenever $SD_q(n) \neq 0$, see the chapter on Structured Matrices.

In contrast to ordinary bases, self-dual bases are very rare. The ratio of the numbers of the self-dual bases to the numbers of ordinary bases is given by either

$$\frac{SD_q(n)}{B_q(n)} \approx \left( q^{n(n-2)/4} \prod_{i=0}^{n/2} \left( q^{2i+1} - 1 \right) \right)^{-1}$$

or $SD_q(n)/B_q(n) = 0$. This ratio rapidly vanishes as either n or q increases. Thus almost every basis of a finite field is a nonself-dual basis.

## 1.7 Polynomials Bases

Polynomial bases of finite fields are constructed from polynomials with coefficients in the ground field. These bases are quite easy to find and very common in many applications. Let $f(x) = x^n + \cdots + f_1x + f_0 \in \mathbf{F}_q[x]$ be irreducible, and $\mathbf{F}_{q^n} \cong \mathbf{F}_q[x]/(f(x))$.

***Definition 1.30.*** A subset of polynomials $\{ p_{n-1}(x), \ldots, p_1(x), p_0(x) \}$ is a basis of $\mathbf{F}_q[x]/(f(x))$ if every $a(x) = a_{n-1}x^{n-1} + \cdots + a_1x + a_0 \in \mathbf{F}_q[x]/(f(x))$ has a unique representation as

$$a(x) = b_{n-1}p_{n-1}(x) + \cdots + b_1p_1(x) + b_0p_0(x),$$

$b_i \in \mathbf{F}_q$. The basis is an *ordered polynomial basis* if $\deg(p_i(x)) = i$ for $i = 0, 1, 2, ..., n-1$.
A large portion of the polynomial bases of $\mathbf{F}_{q^n}$ are ordered polynomial bases up to a relabeling of the indices, for instance, if $\deg(p_i) \neq i$, then a permutation $\pi$ produces $\deg(p_{\pi(i)}) = \pi(i)$. Another important class of polynomial bases is the class of *equal degree polynomial bases*, for these bases $\deg(p_i(x)) = n - 1$ for $i = 0, 1, 2, ..., n-1$.

***Example 1.31.*** (1) (*Binomial Basis*) The subset $\{ p_i(x) = (x - a)^i : i = 0,1,2,...,n-1 \}$, $a \in \mathbf{F}_q$ is an ordered polynomial basis. The ubiquitous polynomial basis $\{ x^{n-1}, ..., x, 1 \}$ is the best known, and perhaps the most important ordered polynomial basis. Further, since

$$x^i = (x - a + a)^i = \sum_{j=0}^{i} \binom{i}{j} a^{i-j} (x-a)^j,$$

the matrix

$$A = \left(a_{i,j}\right) = \left(\binom{i}{j} a^{i-j}\right)$$

is the corresponding change of basis matrix, $(x - a)^i \rightarrow x^i$, provided $n < p$.

(2) Let $a(x) = x^n + a_{n-1}x^{n-1} + \cdots + a_1x + a_0 \in \mathbf{F}_q[x]$, and let $p_i(x) \equiv a(x) \bmod x^{i+1}$, $(0 \leq i < n)$. Then $p_0(x), p_1(x), \ldots, p_{n-1}(x)$ is an ordered polynomial basis of $\mathbf{F}_q[x]/(f(x))$, $f(x)$ irreducible of degree $\deg(f) = n$ if and only if the condition $a_n \cdots a_0a_1 \neq 0$ holds.

(3) (*Newton's Basis*) Let $p_i(x) = (x - a_i)\cdots(x - a_1)(x - a_0)$, $\in \mathbf{F}_q[x]$, $a_i \neq a_j$, for $i \neq j$. Then $p_{n-1}(x), \ldots, p_1(x), p_0(x)$, is an ordered polynomial basis of $\mathbf{F}_q[x]/(f(x))$.

***Lemma 1.32.*** The total number of ordered monic polynomial bases of $\mathbf{F}_{q^n}$ over $\mathbf{F}_q$ is given by

$$P_q(n) = \frac{q^{n(n-1)/2}}{n}\sum_{d\,|\,n}\mu(d)\,q^{n/d}$$

Proof: Clearly $p_0(x) = 1$. Since for each $i > 0$, each $p_i(x) = x^i + a_{i-1}x^{i-1} + \cdots + a_1x + a_0$ in an ordered monic polynomial basis $\{ p_0(x), p_1(x), \ldots, p_{n-1}(x) \}$ can be chosen in $q^i$ different ways, there is a total of $qq^2q^3 \cdots q^{n-1} = q^{n(n-1)/2}$ subsets of the form $\{p_i(x) : \deg(p_i) = i\} \subset \mathbf{F}_q[x]$. Moreover, in order to have the operations of multiplication and division uniquely defined, a unique irreducible polynomial $f(x)$ of degree $n$ is required. And since there are precisely

$$I_n(q) = \frac{1}{n}\sum_{d\,|\,n}\mu(d)\,q^{n/d}$$

such polynomials, the claim is proved. ■

**Remark**: Each n-degree irreducible polynomial generates $q^{n(n-1)/2}$ unique ordered monic polynomial bases of $\mathbf{F}_{q^n}$ over $\mathbf{F}_q$, including the power basis. For all fixed pair $(n, q)$ there are precisely $I_n(q)$ distinct power bases $\{ p_i(x) = x^i : i = 0,1,2,\ldots,n-1 \}$.

***Lemma 1.33.*** (*Lagrange basis*) Let $x_{n-1}, \ldots, x_1, x_0 \in \mathbf{F}_q$ be distinct elements, and let $f(x) \in \mathbf{F}_q[x]$ be irreducible of degree n. Then

( 1 ) The list of polynomials $L_i(x) = \prod_{i \neq j} \frac{x - x_j}{x_i - x_j}$ forms a unique equal degree polynomial basis of $\mathbf{F}_q[x]/(f(x))$, $(\deg(L_i(x)) = n - 1)$.
( 2 ) The matrix for the change of basis $x^i \rightarrow L_i(x)$ is the Vandermonde matrix $V = V(x_{n-1}\ldots x_1x_0)$.
( 3 ) The number of distinct Lagrange bases is

$$L_q(n) = \frac{q(q-1)\cdots(q-1+n)}{n}\sum_{d\,|\,n}\mu(d)q^{n/d}\ .$$

Proof: (2) Let $a(x) = a_{n-1}x^{n-1} + \cdots + a_1x + a_0 = b_{n-1}L_{n-1}(x) + \cdots + b_1L_1(x) + b_0L_0(x)$. Since $L_i(x_j) = \delta_{i,j}$, the substitution $x = x_i$ returns $a(x_i) = b_i$, which leads to the vector equation $\mathbf{b} = V\mathbf{a}$, where $\mathbf{b} = (b_{n-1}, \ldots, b_1, b_0)$, and $\mathbf{a} = \cdot(a_{n-1}, \ldots, a_1, a_0)$. ■

Other varieties of polynomial bases of $\mathbf{F}_q[x]/(f(x))$ besides the ordered and equal degree are also possible. For example, the list $p_{n-1}(x) = x^{n-1}$, and

$$p_i(x) = \prod_{i\neq j}\frac{x - x_j}{x_i - x_j}, \qquad (0 \le i < n-1),$$

where $x_{n-2}, \ldots, x_1, x_0 \in \mathbf{F}_q$ are distinct elements, is a polynomial basis that is neither ordered nor equal degree.

**The $\sigma^k$-Bases and the Matrices $Q_k$**

Let n, k > 0 be positive integers, and $\sigma \in Gal(\mathbf{F}_{q^n}/\mathbf{F}_q)$ be an automorphism of the n degree extension $\mathbf{F}_{q^n}$ of $\mathbf{F}_q$. Consider the image $\sigma^k(\{\alpha_i\}) = \{ß_0, ß_1, \ldots, ß_{n-1}\} : i \in \mathbf{N}\}$ of the basis $\{\alpha_0, \alpha_1, \ldots, \alpha_{n-1}\}$ of $\mathbf{F}_{q^n}$ over $\mathbf{F}_q$, called a σ-set of the basis $\{\alpha_0, \alpha_1, \ldots, \alpha_{n-1}\}$. For example, fix a map $\sigma^k \in Gal(\mathbf{F}_{q^n}/\mathbf{F}_q)$, then the $\sigma^k$-set of the basis $\{\alpha_0, \alpha_1, \ldots, \alpha_{n-1}\}$ has the form

$$\alpha_0^{q^k},\quad \alpha_1^{q^k},\quad \alpha_2^{q^k},\quad \ldots,\quad \alpha_{n-1}^{q^k}$$

where

$$\sigma^k(x) = x^{q^k},\ 0 \le k < n\,.$$

***Lemma 1.34.*** The σ-set $\{\ \sigma^k(\alpha_0), \sigma^k(\alpha_1), \ldots, \sigma^k(\alpha_{n-1})\ \}$ of a basis $\{\ \alpha_0, \alpha_1, \ldots, \alpha_{n-1}\ \}$ is a basis of $\mathbf{F}_{q^n}$ over $\mathbf{F}_q$ for all $\sigma \in Gal(\mathbf{F}_{q^n}/\mathbf{F}_q)$.

Proof: Let $\{\ \alpha_0, \alpha_1, \ldots, \alpha_{n-1}\ \}$ be a basis of $\mathbf{F}_{q^n}$ over $\mathbf{F}_q$, and choose a map $\sigma^k \in Gal(\mathbf{F}_{q^n}/\mathbf{F}_q)$.To show that the $\sigma^k$-set is a basis, it is sufficient to prove that if the linear combination

$$a_0\ {}_0^{q^k} + a_1\ {}_1^{q^k} + a_2\ {}_2^{q^k} + \cdots + a_{n-1}\ {}_{n-1}^{q^k} = 0\,,$$

then the coefficient vector $\mathbf{a} = (a_0,a_1,\ldots,a_{n-1}) = (0,0,\ldots,0)$. Suppose there exists $\mathbf{a} \neq (0,0,\ldots,0)$, and

$$0 = a_0\alpha_0^{q^k} + a_1\alpha_1^{q^k} + a_2\alpha_2^{q^k} + \cdots + a_{n-1}\alpha_{n-1}^{q^k}$$
$$= \left( a_0\alpha_0 + a_1\alpha_1 + a_2\alpha_2 + \cdots + a_{n-1}\alpha_{n-1} \right)^{q^k}.$$

Then it follows that the set { $\alpha_0, \alpha_1, \ldots, \alpha_{n-1}$ } is linearly dependent over $\mathbf{F}_q$, in contradiction of the fact that it is a basis. Therefore, the $\sigma^k$-set is a basis of $\mathbf{F}_{q^n}$ over $\mathbf{F}_q$. ■

The change of basis matrix $Q_k = ( a_{i,j} )$ defined by

$$x^{iq^j} = \sum_{j=0}^{n-1} a_{i,j} x^j \ mod \ \ f(x)$$

is associated with the change of $\sigma^k$-basis

$$\left\{1, \quad x, \quad x^2, \quad \ldots, \quad x^{n-1}\right\} \rightarrow \left\{x, \quad x^{q^k}, \quad x^{2q^k}, \quad \ldots, \quad x^{(n-1)q^k}\right\}.$$

for the ring $\mathbf{R} = \mathbf{F}_q[x]/(f(x))$, where $f(x) \in \mathbf{F}_q[x]$ is a polynomial of degree n
For an arbitrary f(x) the matrix $Q_k$ is not always invertible. However, if the polynomial f(x) is irreducible, the matrix $Q_k$ is nonsingular. Conversely if the matrix $Q_k$ is nonsingular then the polynomial f(x) is irreducible.

**Remark:** The matrix $Q_1 = ( a_{i,j} )$ is the one utilized in the Berlekamp polynomial factorization algorithm. This procedure takes advantage of the singularity/nonsingularity of the matrix $Q_1$ to factor the polynomial f(x) or declares it irreducible, see [1, 2, 3 Berlekamp].

The system of equations

$\det(Q_1) \neq 0$,
$\det(Q) \neq 0$,

where the matrix $Q = ( b_{i,j} )$ is defined by

$$x^{q^i} = \sum_{j=0}^{n-1} b_{i,j} x^j \ mod \ \ f(x)$$

associated with the change of basis

$$\{1, \quad x, \quad x^2, \quad \ldots, \quad x^{n-1}\} \rightarrow \{x, \quad x^q, \quad x^{q^2}, \quad \ldots, \quad x^{q^{n-1}}\} ,$$

provides a deterministic test for normal polynomial f(x) but the calculation is extensive.

# Chapter 2

# Structured Matrices

## 2.1 Basic Concepts

The emphasis will be on the applications of linear algebra to bases of finite fields. Several important classes of matrices that are widely used in analysis of bases are introduced in this chapter. These matrices are employed in change of bases, and in determining the distributions of various types bases.

The set of all n×n matrices is denoted by $\mathcal{M}_n(\mathbf{F}_q) = \{ A = (a_{i,j}) : 0 \le i, j < n, a_{i,j} \in \mathbf{F}_q \}$.

Let $A = (a_{i,j})$, $B = (b_{i,j}) \in \mathcal{M}_n(\mathbf{F}_q)$. The sum and product of a pair of n×n matrices are defined by

( 1 ) $A + B = (a_{i,j} + b_{i,j})$,
( 2 ) $AB = C = (c_{i,j})$, where $c_{i,j} = \sum_{k=0}^{n-1} a_{i,k} b_{k,j}$ .

The trace and determinant are the functions tr, det : $\mathcal{M}_n(\mathbf{F}_q) \to \mathbf{F}_q$ defined by

$\mathrm{tr}(A) = a_{n-1,n-1} + \cdots + a_{1,1} + a_{0,0}$

and

$$\det(A) = \sum_{\pi \in S_n} \mathrm{sgn}(\pi) a_{n-1,\pi(n-1)} \cdots a_{1,\pi(1)} a_{0,\pi(0)} .$$

The summation index $\pi$ runs through all the permutations of { 0, 1, 2, …, n–1 }, and $\mathrm{sgn}(\pi) = \pm 1$ is the sign of the even/odd permutation. The trace and determinant satisfy the following properties.

( 1 ) $\mathrm{tr}(aA + bB) = a\mathrm{tr}(A) + b\mathrm{tr}(B)$, linearity,
( 2 ) $\mathrm{tr}(AB) = \mathrm{tr}(BA)$, commutative,
( 3 ) $\mathrm{tr}(B^{-1}AB) = \mathrm{tr}(A)$, the trace is a class function,
( 4 ) $\det(aA) = a^n\det(A)$, isogenous,
( 5 ) $\det(AB) = \det(BA)$, commutative,
( 6 ) $\det(B^{-1}AB) = \det(A)$, the determinant is a class function.

Another important map is the transpose function $t : \mathcal{M}_n(\mathbf{F}_q) \to \mathcal{M}_n(\mathbf{F}_q)$ given by $A \to A^T$. Transposition is an involution on $\mathcal{M}_n(\mathbf{F}_q)$. The subset of all symmetric matrices $\mathrm{Sym}\mathcal{M}_n(\mathbf{F}_q) = \{ A \in \mathcal{M}_n(\mathbf{F}_q) : A = A^T \}$ is the fixed subset of the transpose function.

The most common and widely used matrices in finite fields analysis are:

( 1 ) The Vandermonde matrix V = ( $x_i^j$ ), $0 \le i, j < n$.

( 2 ) The regular matrix representation A = ( $x_i^{q^j}$ ), $0 \le i, j < n$.

( 3 ) The trace matrix representation T = ( $Tr(x_i x_j)$ ), $0 \le i, j < n$.

( 4 ) The circulant matrix C = $circ[c_{n-1},\ldots,c_1,c_0]$, $c_i \in \mathbf{F}_q$.

The Vandermonde matrix V = ( $x_i^j$ ) arises in the polynomial $f(x) = a_{n-1}x^{n-1} + \cdots + a_1x + a_0$ evaluation/interpolation problem at the n points $x_n, \ldots, x_1, x_0$. The regular matrix representation A = ( $x_i^{q^j}$ ) arises in the determination of the linear independence and other properties of a subset of elements $\{ x_n, \ldots, x_1, x_0 \} \subset \mathbf{F}_{q^n}$. And the trace matrix representation T arises in the determination of the discriminants of bases, also linear independence. Last but not least the circulant matrices are essential in the analysis of normal bases. Most of these matrices are of the form $A_f = ( f(x_i,x_j) )$, where $f : \mathbf{F}_{q^n} \to \mathbf{F}_q$ is a function and $\{ x_n, \ldots, x_1, x_0 \} \subset \mathbf{F}_{q^n}$ is a subset of n points.

The regular matrix representation A is the σ-image of the Vandermonde matrix V. This is accomplished with the assignment

$$x_i^j \quad \to \quad \sigma(x_i^j) = x_i^{q^j}.$$

The determinant of the first matrix V is quite simple. However, the image does not preserve the determinant. In fact the determinant of the matrix A is not so simple.

***Theorem 2.1.*** **Let** $\{ x_n, \ldots, x_1, x_0 \} \subset \mathbf{F}_{q^n}$. Then

( 1 ) $\det(V) = \prod_{0 \le i < j < n} (x_i - x_j)$

( 2 ) $\det(A) = x_1 \prod_{j=1}^{n-1} \prod_{a_1 \cdots a_j} \left( x_{j+1} - \sum_{i=0}^{j} a_i x_i \right).$

Proof: The first appears in many publications, and for the second see [1, Lidl et al, p. 109]. ■

The subset of nonsingular matrices $GL_n(\mathbf{F}_q) = \{ A \in \mathcal{M}_n(\mathbf{F}_q) : \det(A) \ne 0 \}$ coincides with the multiplicative group of $\mathcal{M}_n(\mathbf{F}_q)$, and $SL_n(\mathbf{F}_q) = \{ A \in GL_n(\mathbf{F}_q) : \det(A) = 1 \}$ is one of the most important subgroup of $GL_n(\mathbf{F}_q)$.

***Theorem 2.2.*** The cardinalities $\#GL_n(\mathbf{F}_q)$ and $\#SL_n(\mathbf{F}_q)$ of the sets of nonsingular matrices $GL_n(\mathbf{F}_q)$ and $SL_n(\mathbf{F}_q)$ are given by

( 1) $\#GL_n(\mathbf{F}_q) = \prod_{i=0}^{n-1}(q^n - q^i) = q^{n(n-1)/2}\prod_{i=1}^{n}(q^i - 1)$,

( 2) $\#SL_n(\mathbf{F}_q) = \prod_{i=2}^{n}(q^i - 1)$.

## 2.2 Circulant Matrices

The collection of circulant matrices is important and appear frequently in the analysis of normal bases in all characteristic. Only the essential details for application in finite characteristic will be considered here.

***Definition 2.3.*** A *circulant* matrix, denoted by circ[$c_{n-1}\ldots c_1c_0$], is a matrix of the form

$$circ[c_{n-1}...c_1c_0] = \begin{bmatrix} c_0 & c_1 & c_2 & \cdots & c_{n-1} \\ c_{n-1} & c_0 & c_1 & \cdots & c_{n-2} \\ c_{n-2} & c_{n-1} & c_0 & \cdots & c_{n-3} \\ \vdots & \vdots & \vdots & \vdots & \vdots \\ c_1 & c_2 & c_3 & \cdots & c_0 \end{bmatrix}.$$

**Properties Of Circulant Matrices**

Let C = circ[$c_{n-1},\ldots,c_1,c_0$], and D = circ[$d_{n-1},\ldots,d_1,d_0$].

( 1 ) C + D = circ[$c_{n-1}+d_{n-1},\ldots,c_1+d_1,c_0+d_0$].

( 2 ) CD = DC, multiplication is commutative.

( 3 ) A matrix C is circulant if and only if CP = PC, where P = circ[0,1,…,0,0] = ( $p_{i,j}$ ), is a permutation matrix with

$$p_{i,j} = \begin{cases} 1 & if \;\; j - i \equiv 1 \bmod n, \\ 0 & otherwise. \end{cases}$$

( 4 ) The inverse $C^{-1}$ of a circulant C is a circulant matrix.

( 5 ) Let C ∈ $C_n(\mathbf{F}_q)$ = { nonsingular circulant matrix }. Then the inverse of C is $C^{-1}$ = circ[$b_0,b_{n-1},\ldots,b_2,b_1$], where $a(x)(x^n - 1) + b(x)c(x) = 1$, and $b(x) = b_{n-1}x^{n-1} + \cdots + b_1x + b_0$, $c(x) = c_{n-1}x^{n-1} + \cdots + c_1x + c_0 \in \mathbf{F}_q[x]/(x^n-1)$ are the polynomial representations of C and $C^{-1}$ respectively.

A wide range of algorithms for computing the inverse are available in the literature, see [1, Bini and Pan], p. ?, and check [2, Bini et al] for recent developments.

Let $\mathfrak{C}_n(\mathbf{F}_q)$ = { Circulant Matrices } denotes the set of circulant matrices over $\mathbf{F}_q$, and

consider the maps

( 1 ) $\mu : \mathfrak{C}_n(\mathbf{F}_q) \rightarrow \mathbf{F}_q[x]/(x^n-1)$, defined by $\rho(cir[c_{n-1},\ldots,c_1,c_0]) = c(x) = c_{n-1}x^{n-1} + \cdots + c_1x + c_0$, circulant matrix to polynomial,

( 2 ) $\kappa : \mathbf{F}_{q^n} \rightarrow \mathbf{F}_q[x]/(x^n-1)$, defined by $\kappa(\eta) = c_\eta(x) = c_{n-1}x^{n-1} + \cdots + c_1x + c_0$, where $c_i = Tr(\eta^{1+q^i})$.

These maps are depicted in the commutative diagram below.

$$\begin{array}{ccccc} \mathbf{F}_{q^n} & \xrightarrow{\kappa} & \dfrac{\mathbf{F_q}[x]}{(x^n-1)} & \xrightarrow{\mu} & \mathfrak{C}_n(\mathbf{F}_q) \\ | & & | & & | \\ \mathfrak{C}_n(\mathbf{F}_q) & \xrightarrow{\mu^{-1}} & \dfrac{\mathbf{F_q}[x]}{(x^n-1)} & \xrightarrow{\kappa^{-1}} & \mathbf{F}_{q^n} \end{array}$$

Clearly $\mu$, and $\kappa$ are one to one and since all these subsets are finite, $\mu^{-1}$, and $\kappa^{-1}$ exist. The map $\rho = \mu \circ \kappa$, defined by $\rho(\alpha) = circ[Tr(\alpha^{1+q^{n-1}}),\ldots,Tr(\alpha^{1+q}),Tr(\alpha^2)]$, is the representation $\rho : \mathbf{F}_{q^n} \rightarrow \mathfrak{C}_n(\mathbf{F}_q)$ of the additive group $\mathbf{F}_{q^n}$ in the matrix group $\mathfrak{C}_n(\mathbf{F}_q)$.

***Theorem 2.4.*** (*Convolution Theorem*) The map $circ[c_{n-1},\ldots,c_1,c_0] \rightarrow c(x) = c_{n-1}x^{n-1} + \cdots + c_1x + c_0$ is a ring isomorphism, viz, $\mathfrak{C}_n(\mathbf{F}_q) \cong \mathbf{F}_q[x]/(x^n-1)$.

Proof: The additive property is clear. To verify the multiplicative property, let $A = circ[a_{n-1},\ldots,a_1,a_0]$, $B = circ[b_{n-1},\ldots,b_1,b_0]$, and $\mu(A) = a(x)$, $\mu(B) = b(x)$, now use direct calculations to check $\mu(AB) = \mu(A)*\mu(B)$, where the product $a(x)*b(x) \equiv a(x)b(x) \bmod (x^n - 1)$ is the circular convolution of the polynomials $a(x)$ and $b(x) \in \mathbf{F}_q[x]/(x^n-1)$. ■

Another proof is via the correspondence $x^i \rightarrow P^i$, where P is a permutation matrix. A generalization to c-circulant matrices and the polynomials ring $\mathbf{F}_q[x]/(x^n-c)$, $0 \neq c \in \mathbf{F}_q$ is described in [2, Pan, p. 134]. If $c = -1$, a c-circulant matrix is called anti-circulant. If n is even, then a circulant matrix decomposes as a direct sum $n/2\times n/2$ circulant matrix and a $n/2\times n/2$ anti-circulant matrix. An anti-circulant matrix has the form

$$a\,circ[c_{n-1}\ldots c_1c_0] = \begin{bmatrix} c_0 & -c_1 & -c_2 & \cdots & -c_{n-1} \\ c_{n-1} & c_0 & -c_1 & \cdots & -c_{n-2} \\ c_{n-2} & c_{n-1} & c_0 & \cdots & -c_{n-3} \\ \vdots & \vdots & \vdots & \vdots & \vdots \\ c_1 & c_2 & c_3 & \cdots & c_0 \end{bmatrix}$$

Thus the set of circulant matrices $\mathfrak{C}_n(\mathbf{F}_q) \cong \frac{\mathbf{F_q}[x]}{(x^{n/2}-1)} \otimes \frac{\mathbf{F_q}[x]}{(x^{n/2}+1)}$, see [id bid, p. 215].

In [2, MacWilliams], and other authors, the isomorphism is continued to the direct product

$$\frac{\mathbf{F_q}[x]}{(x^n-1)} \cong \frac{\mathbf{F_q}[x]}{(x\pm 1)^v} \otimes \frac{\mathbf{F_q}[x]}{(f_1(x)^v)} \otimes \frac{\mathbf{F_q}[x]}{(f_2(x)^v)} \otimes \cdots \otimes \frac{\mathbf{F_q}[x]}{(f_d(x)^v)}$$

where gcd(n, q) = v, and $x^n - 1 = (x\pm 1)^v f_1(x)^v f_2(x)^v \cdots f_d(x)^v$. This technique readily leads to the classification of two important subgroups of $C_n(\mathbf{F}_q)$:

( 1 ) The subset of *symmetric* circulant matrices $SC_n(\mathbf{F}_q) = \{ C \in C_n(\mathbf{F}_q) : C = C^T \}$.

( 2 ) The subset of *orthogonal* circulant matrices $OC_n(\mathbf{F}_q) = \{ C \in C_n(\mathbf{F}_q) : CC^T = I_n \}$.

The transpose $\text{circ}[c_{n-1},c_{n-2},\ldots,c_1,c_0]^T = \text{circ}[c_0,c_{n-1},\ldots,c_2,c_1]$ of a circulant matrix corresponds to the reciprocal polynomial

$c^*(x) = x^n c(1/x) = c_0x^{n-1} + c_1x^{n-2} + \quad + c_{n-2}x + c_{n-1}$

of $c(x) = c_{n-1}x^{n-1} + \cdots + c_1x + c_0$ in $\mathbf{F}_q[x]/(x^n-1)$. Similarly, the symmetric circulant matrix

$\text{circ}[c_{n-1},\ldots,c_1,c_0] = \text{circ}[c_{n-1},\ldots,c_1,c_0]^T$

corresponds to the self-reciprocal $c(x) = c^*(x)$, and the orthogonal circulant matrix $\text{cir}[c_{n-1},\ldots,c_1,c_0]$ such that

$\text{cir}[c_{n-1},\ldots,c_1,c_0]\text{cir}[c_{n-1},\ldots,c_1,c_0]^T = I_n$

corresponds to the polynomial $c(x) \in \mathbf{F}_q[x]/(x^n-1)$ such that $c(x)c^*(x) \equiv 1 \bmod (x^n - 1)$.

These concepts are utilized to complete the enumeration of dual and self-dual normal bases.

The size $\#C_n(\mathbf{F}_q)$ of $C_n(\mathbf{F}_q)$ determines the total number of normal bases up to a multiplicative factor 1/n. The size $\#SC_n(\mathbf{F}_q)$ of the subset of nonsingular circulant orthogonal matrices $SC_n(\mathbf{F}_q)$ determines the totality of dual normal bases up to a multiplicative factor. And the size $\#OC_n(\mathbf{F}_q)$ of the subset $OC_n(\mathbf{F}_q)$ ) = $O_n(\mathbf{F}_q) \cap C_n(\mathbf{F}_q)$ of nonsingular circulant orthogonal matrices over $\mathbf{F}_q$ determines the totality of self-dual normal bases up to a multiplicative factor. All these statistics are conditional on the existence of self-dual normal bases. The order of the subgroup $OC_n(\mathbf{F}_q)$ = { Orthogonal Circulant Matrices } of $GL_n(\mathbf{F}_q)$ is computed in [2, MacWilliams], and other.

The *reciprocal ordered factorization* of $x^n - 1$ is the arrangement $x^n - 1 = (x-1)(x+1)^e f_1(x)\cdots f_t(x)g_{t+1}(x)\cdots g_u(x)$, where e = 1 if n is even, otherwise e = 0, $f_i(x)$ is self-reciprocal of degree $deg(f_i(x)) = 2c_i$, and $g_j(x)$ is not self-reciprocal of degree $deg(g_j(x)) = c_j$.·

***Theorem 2.5.*** The cardinality $\#C_n(\mathbf{F}_q)$ of the collection $C_n(\mathbf{F}_q)$ of nonsingular circulant matrices over $\mathbf{F}_q$ is as follows.

$$\#C_n(\mathbf{F}_q) = \begin{cases} (q-1)\prod_{i=1}^{t}(q^{2c_i}-1)\prod_{j=t+1}^{u}(q^{d_j}-1)^2 & \text{if } \gcd(n,q)=1, \\ q^{(k-1)q}\#C_k(\mathbf{F}_q) & \text{if } n=kq, \gcd(k,q)=1, \text{and } q\neq 2, \\ q^{(k+1)/2}\#C_k(\mathbf{F}_q) & \text{if } n=kq=2k, \text{ and } k=\text{odd}, \\ q^{k/2}\#C_k(\mathbf{F}_q) & \text{if } n=kq=2k, \text{ and } k=\text{even}. \end{cases}$$

The formulae for the cardinalities of the various subgroups are simpler to derived whenever gcd(n, q) = 1. In this case the polynomial $x^n - 1$ is separable:

$x^n - 1 = (x - 1)f_1(x)f_2(x)\cdots f_d(x)$ or $(x - 1)(x + 1)f_1(x)f_2(x)\cdots f_d(x)$,

depending on whether n is odd or even, with $f_i(x)$ irreducible.

***Theorem 2.6.*** The cardinality $\#SC_n(\mathbf{F}_q)$ of the collection $SC_n(\mathbf{F}_q)$ of nonsingular symmetric circulant matrices over $\mathbf{F}_q$ is as follows.

$$\#SC_n(\mathbf{F}_q) = \begin{cases} (q-1)^{\delta}\prod_{i=1}^{t}(q^{c_i}-1)\prod_{j=t+1}^{u}(q^{d_j}-1) & \text{where } \delta = \begin{cases} 1 & \text{if } q=2, \\ 1 & \text{if } n=\text{odd}, \\ 2 & \text{if } n \text{ even}, \end{cases} \\ q^{(q-1)k/2}\#SC_k(\mathbf{F}_q) & \text{if } n=kq, \text{and } \gcd(k,q)=1. \end{cases},$$

***Theorem 2.7.*** The cardinality $\#OC_n(\mathbf{F}_q)$ of the collection $OC_n(\mathbf{F}_q)$ of nonsingular orthogonal circulant matrices over $\mathbf{F}_q$ is as follows.

$$\#OC_n(\mathbf{F}_q) = \begin{cases} 2^{\delta}\prod_{i=1}^{t}(q^{c_i}+1)\prod_{j=t+1}^{u}(q^{d_j}-1) & \text{where } \delta = \begin{cases} 0 & \text{if } q=2, \\ 1 & \text{if } n=\text{odd}, \\ 2 & \text{if } n \text{ even}, \end{cases} \\ q^{(q-1)k/2}\#OC_k(\mathbf{F}_q) & \text{if } n=kq, \gcd(k,q)=1. \end{cases}$$

Another approach to the calculation of the sizes of the various subgroups of $C_n(\mathbf{F}_q)$ uses only information about the degrees of the factors of $x^n - 1$. The resulting formulae uses only integers arithmetics, confer [1, Byrd et al] for more details.

***Theorem 2-8.*** Let $F = (\omega^{ij})$ be the *Fourier matrix* of order n, where $\omega$ is a nth root of unity, $0 \le i, j < n$. Then $F^*CF = \text{diag}(c_{n-1},\ldots,c_1,c_0)$.

Proof: The matrix $F^*$ is the complex conjugate of F, and $F^*F = I_n$. ■

***Theorem 2-9.*** (*Konig and Rados*) The polynomial $f(x) = a_{q-2}x^{q-2} + \cdots + a_1x + a_0 \in \mathbf{F}_q[x]$ has $q - 1 - r$ distinct nonzero roots in $\mathbf{F}_q$, where r is the rank of the circulant matrix $C = \text{circ}[a_{q-2},\ldots,a_1,a_0]$.

Proof: Let ω be of order $n = q - 1$ in $\mathbf{F}_q$, and $F = ( \omega^{ij} )$, $0 \le i, j < n$. The triple product $F^*CF = \text{diag}(f(\omega^{q-2}),f(\omega^{q-3}),\ldots,f(\omega),f(1))$. Hence The rank of C is $r = \#\{ 0 \ne x \in \mathbf{F}_q : f(x) \ne 0 \}$. ■

This result, which is well known, counts the distinct nonzero roots of a polynomial f(x) of degree $\le q - 2$, or equivalently a polynomial from $\mathbf{F}_q[x]/(x^{q-1}-1)$. The count does not include the multiplicities.

**Corollary 2.10.** Let $n \ge 1$, and $k < m = q^n - 1$. If the polynomial $f(x) = a_kx^k + \cdots + a_1x + a_0 \in \mathbf{F}_q[x]$ has no root in $\mathbf{F}_q$, then the m×m circulant matrix $C = \text{circ}[a_{m-1},\ldots,a_1,a_0]$ is nonsingular for all n for which $\gcd(k, n) = 1$.

Proof: The constraints $f(x) \ne 0$ for all $0 \ne x \in \mathbf{F}_q$, and $\gcd(k, n) = 1 \Rightarrow$ the polynomial f(x) has no root in $\mathbf{F}_{q^n}$ for all such $n \ge 1$, so the rank of the matrix is $q^n - 1$. ■

***Example 2.11.*** Let $x^2 + ax + b \in \mathbf{F}_q[x]$ be irreducible, let $n = 2v + 1$, and $2 < m = q^n - 1$. Then $C = \text{circ}[a_{m-1},\ldots,a_1,a_0]$ is nonsingular for all n. This produces nonsingular circulant matrices of with 3 nonzero entries per rows and arbitrary dimensions $m = q^n - 1$.

***Theorem 2.12.*** The determinant of the circulant matrix $C = \text{circ}[c_{n-1},\ldots,c_1,c_0]$ is given by the formulae

$$(1)\ \det(C) = \begin{cases} c_{n-1}^n + c_{n-2}^n + \cdots + c_0^n - nc_{n-1}c_{n-2}\cdots c_0 & \text{if } n = 2d+1, \\ c_{n-1}^n + c_{n-2}^n + \cdots + c_0^n - \frac{n}{2}\left(c_{2d}c_{2d-2}\cdots c_0 + c_{2d+1}c_{2d-1}\cdots c_1\right) & \text{if } n = 2(d+1). \end{cases}$$

$$(2)\ \det(C) = c_{n-1}^n \prod_{i=0}^{n-1} (\alpha_i^n - 1),$$

where $\alpha_i$ are the roots of $c(x) = c_{n-1}x^{n-1} + \cdots + c_1x + c_0 \in \mathbf{F}_q[x]/(x^n-1)$,

$$(3)\ \det(C) = \prod_{i=0}^{n-1}\left( \sum_{j=0}^{n-1} c_j\theta^{ij} \right),$$

where θ is a primitive nth root of unity.

Proof: Lines (2) and (3) are derived from the resultants, viz, equation (2) $\det(C) = R(c(x), x^n-1)$, et cetera. ■

Further simplification of det(C) is achieved if $n \equiv 0$ in $\mathbf{F}_q$.

The regular matrix representation $A = ( \sigma^i(\alpha_j) )$ of a subset of conjugate elements $\{ \alpha_j \}$ is a circulant matrix. In addition, if these elements are periods, then the determinant of this matrix can be evaluated in terms of gamma sums. The beta and gamma sums satisfy the properties

( 1 ) $G(\theta^s) = -1$ if n divides s,

( 2 ) $|G(\theta^s)|^2 = p$ if n does not divide s,

( 3 ) $B(s,t) = -\dfrac{G(\theta^s)G(\theta^t)}{G(\theta^{s+t})}$,

where $\eta_{n-1}, \ldots, \eta_1, \eta_0$ are periods of degree n, $G(x) = \eta_{n-1}x^{n-1} + \cdots + \eta_1 x + \eta_0$, $B(s,t) = G(s)G(t)/G(st)$, and θ is a primitive root of unity.

***Corollary 2.13.*** If $p = kn + 1$ is a prime, and $\eta_{n-1}, \ldots, \eta_1, \eta_0$ are periods of degree n, then the determinant of $C = \text{circ}[\eta_{n-1},\ldots,\eta_1,\eta_0]$ is $\det(C) = \pm p^{(n-1)/2}$.
Proof: By the previous theorem

$$\det(C) = \prod_{s=0}^{n-1} G(\theta^s).$$

And substituting the values of the gamma sum returns $|\det(C)|^2 = p^{n-1}$. ■

The concept of circulant matrix is extendable in several directions. One possibility is to consider a function $f : \mathbf{Z}_n \to \mathbf{F}_q$, and let the circulant matrix $C_f = \text{circ}[f(n-1), \ldots, f(1), f(0)]$. Some of the functions widely used are characters of $\mathbf{F}_r$, r a prime power. For an additive character ψ, the circulant matrix is given by

$Q = (\psi(ij)) = \text{circ}[[\psi(n-1),\ldots, \psi(1), \psi(0)]$.

And for a multiplicative character χ, the circulant matrix is given by

$R = ( \chi(j-i) ) = \text{circ}[[\chi(n-1),\ldots,\chi(1),\chi(0)]$.

***Lemma 1.14.*** Let χ be the quadratic symbol in $\mathbf{F}_r$, and let n = r. Then
( 1 ) R is symmetric if $r \equiv 1 \bmod 4$, other wise it is skew symmetric, viz, $R = -R^T$.
( 2 ) R is not invertible.

Proof: To see that R is singular, consider the assignment

$R = \text{circ}[[\chi(n-1),\ldots,\chi(1),\chi(0)] \to R(x) = \chi(n-1)x^{n-1} + \cdots + \chi(1)x + \chi(0)$.

Then $R(1) = 0 \Rightarrow$ gcd($x^n - 1$, $R(x) = (x - 1)f(x)$, some f(x). Thus R(x) is not invertible in $\mathbf{F}_q[x]/(x^n-1)$. ■

The class of circulant matrices R are used in the construction of other classes of structured matrices.

***Definition 2.15.*** An n×n matrix $A = (a_{i,j})$, $a_{i,j} \in \{-1, 0, 1\}$ is a *conference* matrix provided that $a_{i,i} = 0$, $AA^T = (n-1)I_n$ and n is even.

***Lemma 2.116.*** Let $r \equiv 1 \bmod 4$, and let $\chi$ be the quadratic symbol in $\mathbf{F}_r$. Then
( 1 ) The n×n matrix $S = (s_{i,j})$, where $n = r + 1$, and

$$s_{i,j} = \begin{cases} 1 & \text{if } i = 0 \text{ or } j = 0, \text{and } i \neq j, \\ \chi(j-i) & \text{otherwise,} \end{cases}$$

is a nonsingular symmetric conference matrix.
( 2 ) The matrix $\frac{1}{\sqrt{r}}S$ is orthonormal.
( 3 ) If the basis $\{\alpha_{n-1}, \ldots, \alpha_1, \alpha_0\}$ is self-dual, then the change of basis

$$\alpha_i \to \beta_i = s_{i,0}\alpha_0 + s_{i,1}\alpha_1 + \cdots + s_{i,n-1}\alpha_{n-1}$$

is also self-dual.

Proof: (1) The entry $\chi(i-j) = \chi(j-i)$ since $\chi(-1) = 1$ in $\mathbf{F}_r$, so the matrix S is symmetric and $S^2 = rI_{r+1}$. ■

***Example 2.17.*** Let $r = 4e + 1$ and $2n + 1 = 2(r + 1) + 1$ be primes. If $q = 2$ is of order 2n modulo 2n + 1, then the element $\eta = \omega + \omega^{-1}$ generates a self-dual normal basis of $\mathbf{F}_{2^n}$ over $\mathbf{F}_2$, where $1 \neq \omega \in \mathbf{F}_{2^{2n}}$ and $\omega^r = 1$. Moreover since the matrix S is orthogonal in characteristic p = 2, the change of basis

$$\begin{aligned} \eta \quad \to \quad \tau &= s_{0,0}\eta + s_{0,1}\eta^q + s_{0,2}\eta^{q^2} + \cdots + s_{0,n-1}\eta^{q^{n-1}} \\ &= \omega^2 + \omega^{-2} + \omega^{2^2} + \omega^{-2^2} + \cdots + \omega^{2^{n-1}} + \omega^{-2^{n-1}} \end{aligned}$$

is also a self-dual normal basis of $\mathbf{F}_{2^n}$ over $\mathbf{F}_2$. The first instance is for the parameters r = 5, n = r + 1 =6. The corresponding matrix is

$$S_6 = \begin{bmatrix} 0 & 1 & 1 & 1 & 1 & 1 \\ 1 & 0 & 1 & -1 & -1 & 1 \\ 1 & 1 & 0 & 1 & -1 & -1 \\ 1 & -1 & 1 & 0 & 1 & -1 \\ 1 & -1 & -1 & 1 & 0 & 1 \\ 1 & 1 & -1 & -1 & 1 & 0 \end{bmatrix}.$$

Hence $\eta = \omega + \omega^{-1}$ and $\tau = \omega^2 + \omega^{-2} + \omega^{2^2} + \omega^{-2^2} + \cdots + \omega^{2^5} + \omega^{-2^5}$ generate self-dual normal bases of $\mathbf{F}_{2^5}$ over $\mathbf{F}_2$.

## 2.3 Triangular Matrices

***Lemma 2.18.*** Let n > 1, and let q be a prime power. Then
( 1 ) The set of n×n upper triangular matrices $T_n = \{ A = ( a_{i,j} ) : a_{i,j} = 0 \text{ for } i > j \}$ is a group.
( 2 ) The subset of n×n upper triangular matrices $T_n(1) = \{ A \in T_n : a_{i,j} = 1 \}$ is a Sylow subgroup of $GL_n(\mathbf{F}_q)$.

Proof: (2) The cardinality of the set $T_n$ of falling diagonal matrices is $\# T_n = (q-1)^n q^{n(n-1)/2}$, but the cardinality of the set $T_n(1)$ of falling constant diagonal $a_{i,i} = 1$ matrices is $\# T_n(1) = q^{n(n-1)/2}$. This proves the claim. ■

The modification $A \rightarrow B = J_nA$, where $J_n$ is the anti identity matrix, from a falling diagonal matrix to a rising diagonal matrix is used in finite fields multipliers. For n = 3, this has the form

$$\begin{bmatrix} a_3 & a_2 & a_1 \\ 0 & a_3 & a_2 \\ 0 & 0 & a_3 \end{bmatrix} \rightarrow \begin{bmatrix} a_3 & a_2 & a_1 \\ 0 & a_3 & a_2 \\ 0 & 0 & a_3 \end{bmatrix} \begin{bmatrix} 0 & 0 & 1 \\ 0 & 1 & 0 \\ 1 & 0 & 0 \end{bmatrix} = \begin{bmatrix} a_1 & a_2 & a_3 \\ a_2 & a_3 & 0 \\ a_3 & 0 & 0 \end{bmatrix}.$$

The coefficients of a monic polynomial $f(x) = a_nx^n + a_{n-1}x^{n-1} + \cdots + a_1x + a_0$ are utilized to construct the triangular matrix $C_f = ( c_{ij} )$ and its inverse $D = ( d_{ij} )$. The entries of these matrices are defined by

$$c_{i,j} = \begin{cases} a_{i+j+1} & \text{if } 0 \le i+j \le n-1, \\ 0 & \text{if } n \le i+j \le 2n-2, \end{cases} \quad \text{and} \quad d_{i,j} = \begin{cases} 0 & \text{if } 0 \le i+j \le n-1, \\ b_{i+j+1-n} & \text{if } n \le i+j \le 2n-2, \end{cases}$$

where

$$b_k = \sum_{i=0}^{n} a_{n-k+i}\, b_i .$$

The subset of the matrices

$$C_f = \begin{bmatrix} a_1 & a_2 & a_3 & \cdot & \cdot & \cdot & a_n \\ a_2 & a_3 & a_4 & \cdot & \cdot & \cdot & 0 \\ . & . & . & . & . & . & . \\ a_{n-1} & a_n & 0 & . & . & . & 0 \\ a_n & 0 & 0 & . & . & . & 0 \end{bmatrix}$$

does not have a group structure. The matrix $C_f$ is a special case of a Bezout matrix, see [2, Pan, p. 156]. The determinant of a triangular matrix is a very simple expression, namely, $\det(C_f) = -(a_n)^n = -1$.
The triangular change of basis is given by the formula

$$\beta_i = \sum_{j=0}^{n-1} c_{i,j}\alpha_j = \sum_{j=0}^{n-1-i} a_{i+j+1}\alpha_j \ .$$

Specifically, these equations are:

$\beta_0 = a_1\alpha_0 + a_2\alpha_1 + a_3\alpha_2 + \cdots + a_{n-1}\alpha_{n-2} + a_n\alpha_{n-1}$, $a_n = 1$,
$\beta_1 = a_2\alpha_0 + a_3\alpha_1 + a_3\alpha_2 + \cdots + a_n\alpha_{n-2}$,
$\beta_2 = a_3\alpha_0 + a_4\alpha_1 + \cdots + a_n\alpha_{n-3}$,
...
$\beta_{n-2} = a_{n-1}\alpha_0 + a_n\alpha_1$,
$\beta_{n-1} = a_n\alpha_0$.

The collection of rising diagonal triangular matrices $T_n(2)$defined by these polynomials is closely linked to the collection $T_n(1)$. However it does not have a group structure, and since the polynomials are irreducible there are fewer matrices in $T_n(2)$.

***Example 2.19.*** The trinomial $f(x) = x^n + x^k + 1$ in $\mathbf{F}_q[x]$, the triangular change of basis is nearly a permutation of the power basis. As an instance, there is the sequence of irreducible polynomials $f_2(x) = x^2 + x + 1$, $f_3(x) = x^3 + x + 1$, $f_4(x) = x^4 + x + 1$, $f_6(x) = x^6 + x + 1$, ..., in $\mathbf{F}_2[x]$, and the corresponding sequence of matrices

$$C_1 = \begin{bmatrix} 1 & 1 \\ 1 & 0 \end{bmatrix}, \quad C_3 = \begin{bmatrix} 1 & 0 & 1 \\ 0 & 1 & 0 \\ 1 & 0 & 0 \end{bmatrix}, \quad C_4 = \begin{bmatrix} 1 & 0 & 0 & 1 \\ 0 & 0 & 1 & 0 \\ 0 & 1 & 0 & 0 \\ 1 & 0 & 0 & 0 \end{bmatrix}, \quad C_6 = \begin{bmatrix} 1 & 0 & 0 & 0 & 1 \\ 0 & 0 & 0 & 1 & 0 \\ 0 & 0 & 1 & 0 & 0 \\ 0 & 1 & 0 & 0 & 0 \\ 1 & 0 & 0 & 0 & 0 \end{bmatrix}, \ldots.$$

In this case the change of basis from the power basis $\{1, \alpha, \alpha^2, ..., \alpha^{n-1}\}$ to triangular basis $\{\beta_0, \beta_1, \beta_2, ..., \beta_{n-1}\}$ is simply

$\beta_0 = \alpha^{n-1} + 1 = \alpha^{-1}$, since $\alpha^n + \alpha + 1 = 0$,
$\beta_1 = \alpha^{n-2}$,
$\beta_2 = \alpha^{n-3}$,
...
$\beta_{n-2} = \alpha$,
$\beta_{n-1} = 1$.

Comprehensive tables of the irreducible/primitive polynomials $f_n(x) = x^n + x + 1 \in \mathbf{F}_2[x]$ appear in [1, Zieler et al.], [1, 2, Zivkovic], and [2, Menezes, p. 158]. A list of primitive trinomials $x^n + x^k + 1$ in $\mathbf{F}_2[x]$ appears [2, Menezes, p. 161], and a table of trinomials whose degrees $n = 2^p - 1$ are primes is given in [1, Kurita et al].

The next lemma indicates that if $\alpha$ is a root of an irreducible polynomial $f_{2k}(x) = x^{2k} + x + 1$ of even degree $n = 2k$ in the sequence { irreducible $f_n(x) = x^n + x + 1 \in \mathbf{F}_2[x]$ }, then $\{1, \alpha, \alpha^2, ..., \alpha^{n-1}\}$ and $\{\alpha^{n-1}+1, \alpha^{n-1}, \alpha^{n-2}, ..., \alpha, 1\}$ are dual bases.

***Lemma 2.20.*** The power basis and the triangular basis are dual bases if and only if n is even and the coefficients of the irreducible polynomial $f(x) = a_n x^n + a_{n-1}x^{n-1} + \cdots + a_1 x + a_0 \in \mathbf{F}_2[x]$ satisfy

$$a_i = \begin{cases} 1 & \text{if } i = 1, \\ 0 & \text{if } i \text{ is odd and } 0 < i < n. \end{cases}$$

Proof: See [4, Hasan and Bhargava]. ■

## 2.4 Hadamard Matrices

In finite fields analysis the collection of Hadamard Matrices has many applications, one of them is the generation of dual bases in characteristic $\text{char}(\mathbf{F}_q) > 2$, but not self-dual.

***Definition 2.21.*** Let A and B be n×n and k×k matrices respectively. The tensor product $A \otimes B$ of A and B is defined by the kn×kn matrix $A \otimes B = ( a_{i,j}B )$.

***Theorem 2.22*** If $A = ( a_{i,j} )$,and $B = ( b_{i,j} )$ are matrices such that $a_{i,j}, b_{i,j} \in \{ -1, 1 \}$, then the tensor product $A \otimes B = ( c_{i,j} )$, $c_{i,j} \in \{ -1, 1 \}$.

***Definition 2.23*** A n×n matrix $H_n = ( h_{i,j} )$, $h_{i,j} \in \{ -1, 1 \}$, such that $HH^t = nI_n$ is called a Hadamard matrix of order n.

The tensor product of two matrix $A_k = (a_{i,j})$, and $B_n = (b_{i,j})$ of orders k and n respectively yields a Hadamard matrix $H_{kn} = (h_{i,j}) = A \otimes B$ of order kn. This simple mechanism is very useful in the construction of larger matrices and infinite sequences of these matrices.

The first is taken to be $H_1 = [\ 1\ ]$. Repeatedly tensoring the trivial matrix $H_1$ reproduces an infinite sequence of Hadamard matrices:

$$H_{2^k} = H_{2^{k-1}} \otimes H_{2^{k-1}} = \begin{bmatrix} H_{2^{k-1}} & H_{2^{k-1}} \\ H_{2^{k-1}} & -H_{2^{k-1}} \end{bmatrix},$$

of order $n = 2^k$, $k \geq 0$. The first two are

$$H_2 = H_1 \otimes H_1 = \begin{bmatrix} 1 & 1 \\ 1 & -1 \end{bmatrix}, \quad \text{and} \quad H_{2^2} = H_2 \otimes H_2 = \begin{bmatrix} 1 & 1 & 1 & 1 \\ 1 & -1 & 1 & -1 \\ 1 & 1 & -1 & -1 \\ 1 & -1 & -1 & 1 \end{bmatrix}.$$

The apparatus for producing other well known infinite sequences of Hadamard matrices is stated below.

***Theorem 2.24.*** (*Paley 1933*) If q is a prime power, then there is a Hadamard matrix of order $n = 2^e(q + 1)$, where $e \geq 0$ such that $n \equiv 0 \bmod 4$.

Proof: See [1, Beth et al., p. 55] or similar text.

The entries of the matrices specified above are as follow:
For $q \equiv 3 \bmod 4$

$$h_{i,j} = \begin{cases} 1 & \text{if } i = 0 \text{ or } j = 0, \\ \chi(j-1) & \text{if } i \neq j \text{ and } 1 \leq i,j \leq q, \\ -1 & \text{if } i = j \text{ and } 1 \leq i \leq q. \end{cases}$$

And for $q \equiv 1 \bmod 4$,

$$h_{i,j} = \begin{cases} 1 & \text{if } i = 0 \text{ or } j = 0, \text{ and } i \neq j, \\ \chi(j-1) & \text{otherwise.} \end{cases}$$

***Theorem 2.25.*** If $H_n = (h_{i,j})$ is a Hadamard matrix of order n, then
( 1 ) The integer n = 1, 2, 12, or 4k, $k \geq 1$.
( 2 ) The integer $n = 2^k$, $k \geq 0$.
( 3 ) The integer $n = 2^e(q + 1)$, where q is a prime power, and $e \geq 0$ such that $n \equiv 0 \bmod 4$.

The form of the integer k in statement (1) above is not known. In other words, the converse n = 4k $\Rightarrow$ there exists a Hadamard matrix of order n is not known. This is one of the main open problem in the theory of orthogonal arrays.
***Lemma 2.26.*** ***(****Hadamard 1893****)*** The determinant of the $H_n$ is $|\det(H_n)| = n^{n/2}$.

Proof: The matrix $H_n$ is the maximum (equality) of the determinant inequality

$$|\det(A)| \le \prod_{i=0}^{n-1}\sqrt{\sum_{j=0}^{n-1} a_{i,j}^2}\,. \qquad \blacksquare$$

The only known circulant Hadamard matrix is

$$H_4 = \begin{bmatrix} 1 & 1 & 1 & -1 \\ -1 & 1 & 1 & 1 \\ 1 & -1 & 1 & 1 \\ 1 & -1 & -1 & 1 \end{bmatrix}.$$

The existence of other circulant Hadamard matrix is an open problem. However, it is known that the order of the next one is n ≥ 1898884.

## 2.5 Multiplication Tables

Let $x = x_{n-1}\alpha_{n-1} + \cdots + x_1\alpha_1 + x_0\alpha_0, y \in \mathbf{F}_{q^n}$. The multiplication matrix

$$M = \begin{bmatrix} \alpha_0\alpha_0 & \alpha_0\alpha_1 & \cdots & \alpha_0\alpha_{n-1} \\ \alpha_1\alpha_0 & \alpha_1\alpha_1 & \cdots & \alpha_1\alpha_{n-1} \\ \vdots & \vdots & \vdots & \vdots \\ \alpha_{n-1}\alpha_0 & \alpha_{n-1}\alpha_0 & \cdots & \alpha_{n-1}\alpha_{n-1} \end{bmatrix}$$

extracted from the product

$$xy = (x_{n-1},\ldots,x_1,x_0)\begin{pmatrix} \alpha_0 \\ \alpha_1 \\ \vdots \\ \alpha''_{n-1} \end{pmatrix}(\alpha_0,\alpha_1,\ldots,\alpha_{n-1})\begin{pmatrix} y_{n-1} \\ \vdots \\ y_1 \\ y_0 \end{pmatrix}.$$

The equivalent description as an $n \times n^2$ matrix $T = (t_{i,j,k})$, where

$$\alpha_i \alpha_j = \sum_{k=0}^{n-1} t_{i,j,k} \alpha_k ,$$

and $t_{i,j,k} \in \mathbf{F}_q$, is more useful.

There is practical interest in determining bases $\{ \alpha_{n-1},\ldots,\alpha_1,\alpha_0 \}$ with very sparse multiplication matrices. A sparse multiplication matrix $M$ has low discrete weight $w(M) = \#\{ t_{i,j,k} \neq 0 : 0 \leq i,j,k < n \}$.

***Definition 2.27.*** The complexity of multiplication with respect to the basis $\{ \alpha_{n-1},\ldots,\alpha_1, \alpha_0 \}$ of $\mathbf{F}_{q^n}$ over $\mathbf{F}_q$, is defined by

$$C(n) = \frac{1}{n} \#\{ t_{i,j,k} \neq 0 : 0 \leq i, j,k < n \} .$$

The complexity of some highly structured $n \times n^2$ matrices $M = [ M_0\ M_1 \cdots M_{n-1} ]$, where $M_i$ is an $n \times n$ matrix, reduces to $w(M_0) = \#\{ t_{0,j,k} \neq 0 : 0 \leq j,k < n \}$. Specifically, for a normal basis $\{ \eta, \eta^{q}, \eta^{q^2}, ..., \eta^{q^{n-1}} \}$, the product

$$\begin{aligned} \eta^{q^i} \eta^{q^j} &= \sum_{k=0}^{n-1} t_{i,j,k} \eta^{qk} \\ &= \sigma^i (\eta\eta^{q^{j-i}}) = \sum_{k=0}^{n-1} t_{0,j-i,k} \eta^{qk} . \end{aligned}$$

Thus a single submatrix $T_0 = (t_{0,j-i,k})$ generates the entire $n \times n^2$ multiplication matrices $T = [ T_0\ T_1 \cdots T_{n-1} ]$.

***Definition 2.28.*** Let $A = ( a_{i,j} )$ be an $n \times n$ matrix. The weight and excess of $A$ are defined by $w(A) = \#\{a_{i,j} \neq 0\}$ and $e(A) = w(A) - n$ respectively.

Clearly, $n \leq w(A) \leq n^2$, and $0 \leq e(A) \leq n^2 - n$. Efficient algorithms for self-dual bases multipliers demand change of bases (between the dual bases) matrices of nearly zero excess, see [4, Berlekamp], [1, Stinson], and [2, Morgan et al]. As far as multiplication algorithm based on self-dual bases is concerned, the ideal change of bases (between the dual bases) matrix $M : A \to B$ will be monomial matrix; permutations matrices are monomial matrices. Monomial matrices have excess $e(M) = 0$. Obviously, not all monomial matrices change of bases correspond to dual bases.

Let $\{ \delta_0, \delta_1, \ldots, \delta_{n-1} \}$ be the dual basis of the power basis $\{ 1, \alpha, \alpha^2, \ldots, \alpha^{n-1} \}$, and put $\beta_i = \dfrac{\delta_i}{\alpha^{n_s} f'(\alpha)}$. Here the element $\alpha \in \mathbf{F}_{q^n}$ is a root $f(x) = x^n + \sum_{i=0}^{t} a_{n_i} x^{n_i}$, and $s = [t/2]$.

***Theorem 2.29***. (*Morgan et al, 1997*) Suppose that the minimal polynomial of $\alpha \in \mathbf{F}_{q^n}$ is $f(x) = x^n + \sum_{i=0}^{t} a_{n_i} x^{n_i}$ , $0 < n_0 < n_1 < \cdots , n_t$. Then the change of basis matrix $A : \{ \alpha_i \} \rightarrow \{ \beta_i \}$ has excess equal to

$$e(A) = \begin{cases} \sum_{s+1=i}^{t} n_i - \sum_{1=i}^{s-1} n_i & t \text{ odd,} \\ \sum_{s+1=i}^{t} n_i - \sum_{1=i}^{s} n_i & t \text{ even.} \end{cases}$$

***Example 2.30***. (1) The change of basis matrix $A$ has excess $e(A) = 0$ if and only if either $f(x) = x^n + a_k x^k + a_0$ or $x^n + a_0 \in \mathbf{F}_q[x]$, $0 < k < n$.
(2) The excess $e(A) = 1$ if and only if $f(x) = x^n + a_k x^k + a_{k-1} x^{k-1} + a_0 \in \mathbf{F}_q[x]$, $1 < k < n$.
(3) The excess $e(A) = 2$ if and only if either $f(x) = x^n + a_k x^k + a_{k-1} x^{k-1} + a_{k-2} x^{k-2} + a_0$ or $x^n + a_k x^k + a_{k-2} x^{k-2} + a_0 \in \mathbf{F}_q[x]$, $2 < k < n$.

# Chapter 3

# Normal Bases

## 3.1 Basic Concepts

An important collection of bases of finite fields is known as normal bases. Theses bases are of interest in both theory and applications. Normal bases are employed in factorization algorithms, [1, Niederreiter], cyclic codes, testing/signature analysis, [1, Hoffner et al.]. This chapter presents a variety of basic ideas and techniques applicable to normal bases.

Let G be a group and let S be a set. The G-orbit orb(x) of $x \in S$ is defined by orb(x) = { ρ(x) : ρ ∈ G }. The structure and size of the orbit of a point or a subset of points depends on the group G. In particular, if G = **Z** is the set of integers, then the **Z**-orbit orb(α) = { $\alpha^n : n \in \mathbf{Z}$ } of the fixed element $\alpha \in \mathbf{F}_{q^n}$ coincides with the cyclic group generated by α. But if G = Gal($\mathbf{F}_{q^n}/\mathbf{F}_q$) = { $\sigma, \sigma^2, \ldots, \sigma^{n-1}$ } is the set of automorphisms of $\mathbf{F}_{q^n}$ over $\mathbf{F}_q$, then the G-orbit orb(α) = {σ(α) : σ ∈ Gal($\mathbf{F}_{q^n}/\mathbf{F}_q$)} of the element $\alpha \in \mathbf{F}_{q^n}$ coincides with the set

$$\{ \alpha^{q^{n-1}}, \ldots, \alpha^q, \alpha \}$$

of conjugates of α.

***Definition 3.1.*** A basis N = { $\eta_{n-1},\ldots, \eta_1, \eta_0$ } of $\mathbf{F}_{q^n}$ over $\mathbf{F}_q$ is said to be a *normal basis* if it is the orbit of some element η in $\mathbf{F}_{q^n}$, viz, N = orb(η) = { σ(η) : σ ∈ Gal($\mathbf{F}_{q^n}/\mathbf{F}_q$) }. A *primitive normal basis* is a basis generated by a primitive element η of $\mathbf{F}_{q^n}$.

If the element $\eta \in \mathbf{F}_{q^n}$ is normal, then each element $\alpha \in \mathbf{F}_{q^n}$ has an expansion of the form

$$\alpha = a_{n-1}\eta^{q^{n-1}} + \cdots + a_1\eta^q + a_0\eta$$

$a_i \in \mathbf{F}_q$.

**Computational Properties Normal Bases**

Let $\eta \in \mathbf{F}_{q^n}$ be a normal element, and let $\alpha, \beta \in \mathbf{F}_{q^n}$.

( 1 ) $\alpha + \beta = \sum_{i=0}^{n-1} (a_i + b_i)\eta^{q^i}$ ,

( 2 ) $\alpha\beta = \sum_{i=0}^{n-1} c_i(\alpha, \beta, \eta)\eta^{q^i}$ ,

( 3 ) $\alpha^q = a_{n-2}\eta^{q^{n-1}} + \cdots + a_0\eta^q + a_{n-1}\eta$ ,

( 4 ) $\alpha^{1/q} = a_0\eta^{q^{n-1}} + \cdots + a_2\eta^q + a_1\eta$ .

The coefficient function $c_i(\alpha,\beta,\eta)$ for multiplication is completely and uniquely determined by the element $\eta$. Property (3) above states that computing the qth-power of $\alpha \in \mathbf{F}_{q^n}$ is equivalent to a left cyclic shift

$(a_{n-1}, a_{n-2}, \ldots, a_1, a_0) \rightarrow (a_{n-2}, a_{n-3}, \ldots, a_0, a_{n-1})$

of the coordinates of $\alpha$ with respect to a normal basis. The next property can be interpreted as a statement about the existence of qth root of elements of $\mathbf{F}_{q^n}$. A right cyclic shift of the coefficients

$(a_{n-1}, a_{n-2}, \ldots, a_1, a_0) \rightarrow (a_0, a_{n-1}, \ldots, a_2, a_1)$

computes the $q^{th}$ root of any element of $\mathbf{F}_{q^n}$. In particular, in the field $\mathbf{F}_{2^n}$, the square $\alpha^2$ of $\alpha$ is computed by a left cyclic shift of the coefficients of $\alpha$ with respect to a normal basis. And the square root is computed by a right cyclic shift of the coefficients.

The computational properties normal bases are the cornerstones of finite field multipliers, see [1, Geiselmann et al], [1, Hasan et al], [1, Itoh et al], et cetera.

**Remark**: A normal basis representation of $\mathbf{F}_{2^n}$ provides a trivial proof of a result in fields of characteristic char(**F**) = 2, which states that every element in the field is a square, so it has a square root.

## 3.2 Existence of Normal Bases

The vector space $\mathbf{F}_{q^n}$ has the structure of a cyclic vector space. Accordingly all the tools of linear algebra applicable to cyclic vector spaces are available.

***Lemma 3.2.*** Let $A \in GL_n(\mathbf{F}_q)$ be a linear map on the vector space V. Then there exists a vector v such that $\{ v, Av, A^2v, A^3v, ..., A^{n-1}v\}$ is a basis of V if and only if the characteristic polynomial $c(x) = \det(A - xI_n)$ and the minimal polynomial $f(x) = a_dx^d + \cdots + a_1x + a_0$ of A are identical.

Both the minimal polynomial f(x) of A and the characteristic polynomial c(x) annihilate the linear map A, id est, the equation $c(A) = f(A) = a_dA^d + \cdots + a_1A + a_0 = 0$. Furthermore, *deg*(f(x)) ≤ *deg*(c(x)). For application to finite fields, let $v = \eta$ be a normal element and the linear map $A = \sigma$ be the canonical automorphism. In this setting the subset $\{ v, Av, A^2v, ..., A^{n-1}v\} = \{ \eta, \eta^q, ... \}$ is precisely a normal basis.

***Theorem 3.3.*** (*Normal Basis Theorem*). The vector space $\mathbf{F}_{q^n}$ over $\mathbf{F}_q$ has a normal basis for all n ≥ 1.

Proof: This is a consequence of the previous lemma; see also [1 Lidl et al, p.60], [1, Waterhouse], [1, van der Waerden, p. 200-204].

The normal basis theorem was developed by at least three authors between 1850 to 1888, [1 Lidl et al, p.76].

***Lemma 3.4.*** (*Normal basis projection lemma*) Suppose η generates a normal basis for $\mathbf{F}_{q^n}$ over $\mathbf{F}_q$, and n = de. Then $\gamma = Tr_{n:d}(\eta)$ generates a normal basis for $\mathbf{F}_{q^d}$ over $\mathbf{F}_q$, where $Tr_{n:d} : \mathbf{F}_{q^n} \to \mathbf{F}_{q^d}$ is the relative trace.

Proof: Suppose that $\gamma = \sum_{j=0}^{e-1} \eta^{q^{dj}}$ does not generate a normal basis of $\mathbf{F}_{q^d}$ over $\mathbf{F}_q$. Then

$$a_{n-1}\gamma^{q^{d-1}} + \cdots + a_1\gamma^q + a_0\gamma = 0$$

for some nontrivial vector $(a_{d-1},\ldots,a_1,a_0) \neq (0,\ldots,0,0)$. But this implies that

$$\sum_{i=0}^{d-1} a_i \left( \sum_{j=0}^{e-1} \eta^{q^{dj}} \right)^{q^i} = \sum_{i=0}^{d-1}\sum_{j=0}^{e-1} a_i \eta^{q^{dj+i}}$$
$$= \sum_{k=0}^{n-1} b_i \eta^{q^k} = 0$$

in contradiction of the linear independence of the subset $\{ \eta^{q^{n-1}}, ..., \eta^q, \eta \}$. ■

## 3.3 Iterative Construction of Normal Bases

The idea of iterative construction of normal bases is the synthesis of normal bases for $\mathbf{F}_{q^n}$ from the normal bases of the individual subfields of $\mathbf{F}_{q^n}$. Moreover, to construct the multiplication table of $\mathbf{F}_{q^n}$, simply take the tensor product of the individual multiplication tables of the subfields. The result below is the foundation of iterated construction of normal bases of $\mathbf{F}_{q^n}$ over $\mathbf{F}_q$ for composite degree n.

***Lemma 3.5.*** (*Semaev 1988*) Let n = rs, gcd(r, s) = 1, and let { $\alpha_{n-1}, \ldots, \alpha_1, \alpha_0$ } and{ $ß_{n-1}, \ldots, ß_1, ß_0$ } be normal bases of $\mathbf{F}_{q^r}$ and $\mathbf{F}_{q^s}$ over $\mathbf{F}_q$ respectively. Then

( 1 ) The subset { $\alpha_{r-1}\beta_{s-1}, \ldots, \alpha_0\beta_1, \alpha_0\beta_0$ } is a normal basis of $\mathbf{F}_{q^n}$ over $\mathbf{F}_q$.

( 2 ) If R = $(r_{i_1,j_1})$ and S = $(s_{i_2,j_2})$ are the multiplication tables (submatrices) of $\mathbf{F}_{q^r}$ and $\mathbf{F}_{q^s}$, respectively, then T = ( $t_{i,j}$ ) = $(r_{i_1,j_1}s_{i_2,j_2})$, where

$$i \equiv \frac{i_1 s + i_2 r}{r+s} \bmod rs, \quad j \equiv \frac{j_1 s + j_2 r}{r+s} \bmod rs,$$

is the multiplication table (submatrix) of $\mathbf{F}_{q^n}$, with $0 \le i, j < n$.

Proof: (1) Since gcd(r+s, rs) = 1, the composition map $\phi : \mathbf{Z}_r \times \mathbf{Z}_s \to \mathbf{Z}_{rs}$ defined by

$$(x, y) \quad \to \quad \frac{xs + yr}{r+s} \bmod rs$$

is one-to-one. Put $\eta_k = \alpha_i\beta_j$, $k \equiv (is + jr)/(r + s) \bmod rs$, and let $\sigma \in \mathrm{Gal}(\mathbf{F}_{q^n}/\mathbf{F}_q)$. Then

$$\sigma^d(\eta_k) = \sigma^d(\alpha_i\beta_j) = \alpha_{i+d}\beta_{j+d} = \eta_{k+d},$$

since $k + d \equiv [(i+d)s + (j+d)r]/(r + s) \bmod rs$. This implies that $\eta_k$ is normal.

(2) The entries of the matrix T = ( $t_{i,j}$ ) are defined by $\eta_0\eta_i = \sum_{j=0}^{n-1} t_{i,j}\eta_j$. In addition $\eta_0\eta_i = \alpha_0\beta_0\alpha_{i_1}\beta_{i_2} = \alpha_0\alpha_{i_1}\beta_0\beta_{i_2}$. Hence

$$\eta_0\eta_i = \left(\sum_{j_1=0}^{r-1} r_{i_1,j_1}\alpha_{j_1}\right)\left(\sum_{j_2=0}^{s-1} s_{i_2,j_2}\beta_{j_2}\right) = \sum_{j_1=0}^{r-1}\sum_{j_2=0}^{s-1} r_{i_1,j_1}s_{i_2,j_2}\alpha_{j_1}\beta_{j_2} = \sum_{j=0}^{n-1} t_{i,j}\eta_j.$$

These prove the claims about the subset { $\eta_{n-1}, \ldots, \eta_1, \eta_0$ } = { $\alpha_{r-1}\beta_{s-1}, \ldots, \alpha_0\beta_1, \alpha_0\beta_0$ }. ■

## 3.4 Additive Order and Decomposition Theorem

The $\mathbf{F}_q$–linear map $\circ : \mathbf{F}_q[x] \times \mathbf{F}_{q^n} \to \mathbf{F}_{q^n}$ is defined by

$$f(x) \circ \alpha = a_k \alpha^{q^k} + \cdots + a_1 \alpha^q + a_0 \alpha .$$

The pairing (f(x), α) = f(x)°α turns the additive group $\mathbf{F}_{q^n}$ into a cyclic $\mathbf{F}_q[x]$-module. And the $\mathbf{F}_q[x]$-module structure of the additive group $\mathbf{F}_{q^n}$ induces the additive order of elements.

***Definition 3.6.*** The *additive order* Ord(α) of an element $\alpha \in \mathbf{F}_{q^n}$ is the polynomial a(x) of least degree such that a(x)°α = 0 but a(x)°α ≠ 0 for all divisor b(x) of a(x).

The additive order Ord(α) = a(x) is always a divisor of $x^n - 1$. Thus there is a monic polynomial of degree < n for which $\text{Ord}(\alpha) = a_{n-1}\sigma^{n-1}(\alpha) + \cdots + a_1\sigma(\alpha) + a_0(\alpha)$. The following identities are readily verified.

( 1 ) Ord(0) = 1, since a(x) = 1 is the divisor of $x^n - 1$ of least degree such that a(x)°0 = 0,
( 2 ) Ord(α) = x – 1, since (x – 1)°α = 0 for all $0 \neq \alpha \in \mathbf{F}_q$, and
( 3 ) $\text{Ord}(\alpha) = a_1(x)a_2(x) \cdots a_r(x)$, where $\text{Ord}(\alpha_i) = a_i(x)$ divides $x^n - 1$.
( 4 ) Ord(α+β) = Ord(α)Ord(β)/gcd(Ord(α),Ord(β)).

***Definition 3.7.*** An element $\eta \in \mathbf{F}_{q^n}$ is a normal over $\mathbf{F}_q$ if and only if $\text{Ord}(\eta) = x^n - 1$.

***Lemma 3.8.*** Let $\alpha \in \mathbf{F}_{q^n}$. Then
( 1 ) The set $V_\alpha = \{ f(x)\circ\alpha : f(x) \in \mathbf{F}_q[x] \}$ is a vector space of dimension *deg*(Ord(α)).
( 2 ) The set $\{ \alpha, \sigma(\alpha), \ldots, \sigma^{n-1}(\alpha) \}$ spans a vector space of dimension *deg*(Ord(α)).

The additive order characterization of elements leads straight to the *standard normal test*. The standard normal, which consists of a system of inequalities, is the additive version of the Lucas test for primitive elements in a cyclic group; all the intricacies of the multiplicative version come through mutatis mutandis..

***Theorem 3.9.*** *(Standard Normal Test)* An element $\eta \in \mathbf{F}_{q^n}$ is a normal element over $\mathbf{F}_q$ if and only if the system of inequalities

$$\frac{x^n - 1}{a(x)} \circ \eta \neq 0$$

holds for all irreducible factors a(x) of $x^n - 1 \in \mathbf{F}_q[x]$ .

Proof: Let a(x) be an irreducible factor of $x^n - 1$, and let $a(x)^v \mid x^n - 1$, but $a(x)^{v+1}$ not a divisor of $x^n - 1$, $v \geq 0$. Then the hypothesis $[(x^n - 1)/a(x)]\circ\eta \neq 0$ implies that the polynomial $(x^n - 1)/a(x)$ cannot annihilate $\eta$, so $a(x)^v$ must be a factor of the annihilator $\mathrm{Ord}(\eta)$ of $\eta$. Since the factor $a(x)$ is arbitrary, the order $\mathrm{Ord}(\eta)$ of $\eta$ is divisible by all divisors $a(x)^v$ of $x^n - 1$. Specifically $x^n - 1 \mid \mathrm{Ord}(\eta)$. But by definition $\mathrm{Ord}(\eta) \mid x^n - 1$, so $\mathrm{Ord}(\eta) = x^n - 1$. ■

***Lemma 3.10.*** Let $\eta$ be a normal element in $\mathbf{F}_{q^n}$ over $\mathbf{F}_q$ and let $x^n - 1 = a(x)b(x)$. Then $\gamma = b(x)\circ\eta$ is an element of order $\mathrm{Ord}(\gamma) = a(x)$ in $\mathbf{F}_{q^n}$ over $\mathbf{F}_q$.

Proof: Same technique as above. ■

The techniques employed in these proofs are extensively used in the literature, see [1, Schwarz], [1, von sur Gathen et al.], [1, Lenstra et al.], etc.

***Example 3.11.*** For n = 3, there are two possibilities. ( 1 ) The case $x^3 - 1 = (x - 1)(x^2 + x + 1) \in \mathbf{F}_q[x]$ if $q \equiv 1 \bmod 3$, (or $q = 2^v$, v odd). The test has two inequalities:

( i ) $\dfrac{x^n - 1}{x - 1}\circ\eta = Tr(\eta) \neq 0$,

( ii ) $\dfrac{x^n - 1}{x^2 + x + 1}\circ\eta = \eta^q - \eta \neq 0$.

( 2 ) And the case $x^3 - 1 = (x - 1)(x - a)(x - b)\,) \in \mathbf{F}_q[x]$, if $q \equiv 2 \bmod 3$, (or $q = 2^v$, v even). The test has three inequalities:

( i ) $\dfrac{x^n - 1}{x - 1}\circ\eta = Tr(\eta) \neq 0$,

( ii ) $\dfrac{x^n - 1}{x - a}\circ\eta = \eta^{q^2} - (b + 1)\eta^q + b\eta \neq 0$,

( iii ) $\dfrac{x^n - 1}{x - b}\circ\eta = \eta^{q^2} - (a + 1)\eta^q + a\eta \neq 0$.

If any one of these systems of inequalities ( 1 ) or ( 2 ) holds, then the element $\eta$ is normal in $\mathbf{F}_{q^3}$ over $\mathbf{F}_q$.

Let $x^n - 1 = (x - 1)f_1(x)f_2(x) \cdots f_k(x)$, $\gcd(f_i(x)f_j(x)) = 1$, and let $V_i = \{ \alpha \in \mathbf{F}_{q^n} : f_i(x)\circ\alpha = 0 \}$ be a cyclic submodule of the $\mathbf{F}_q[x]$-module $\mathbf{F}_{q^n}$. Each $V_i$ is an $\sigma$-invariant subspace of the vector space. Since $x - 1$ divides $x^n - 1$, the first submodule $V_0 = \{ \alpha \in \mathbf{F}_{q^n} : (x - 1)\circ\alpha = 0 \} = \mathbf{F}_q$ is common to all cases.

***Theorem 3.12.*** (*Decomposition theorem*) (1) The vector space $\mathbf{F}_{q^n} = V_0 \oplus V_1 \oplus \cdots \oplus V_k$.
(2) Every $\alpha \in \mathbf{F}_{q^n}$ has a unique decomposition as $\alpha = \alpha_0 + \alpha_1 + \cdots + \alpha_k$, where $\alpha_i \in V_i$.

This is due to [1, Pencin],, and {1, Semaev]; see also [1, Menezes et al] for a proof. An element $\eta \in \mathbf{F}_{q^n}$ is a normal element over $\mathbf{F}_q$ if and only if the decomposition $\eta = \eta_0 + \eta_1 + \cdots + \eta_k$ satisfies $\eta_i \neq 0$ for all $i \geq 0$. The sum $\eta = \eta_0 + \eta_1 + \cdots + \eta_k$ of the generators of the cyclic submodules $V_i$'s is a generator of the module $\mathbf{F}_{q^n}$, that is, $\eta$ is a normal element in $\mathbf{F}_{q^n}$ over $\mathbf{F}_q$.

***Lemma 3.13.*** Let $\mathrm{ord}_n(q) = n - 1$, n prime, and let $\alpha \in \mathbf{F}_{q^n} - \mathbf{F}_q$. Suppose that $a\mathrm{Tr}(\alpha) + bn \neq 0$, $a, b \in \mathbf{F}_q$, $a \neq 0$. Then $\eta = a\alpha + b$ is a normal element in $\mathbf{F}_{q^n}$.

Proof: The multiplicative order $\mathrm{ord}_n(q) = n - 1$ implies that $x^n - 1 = (x - 1)f_1(x)$, where $f_1(x)$ is irreducible of degree $n - 1$, and the system of inequalities of the standard normal test has only two lines:

( i ) $(x - 1)\circ\eta = \eta^q - \eta = (a\alpha + b)^q - (a\alpha + b) = a(\alpha^q - \alpha) \neq 0$,

( ii ) $f_1(x)\circ\eta = \mathrm{Tr}(\eta) = a\mathrm{Tr}(\alpha) + bn \neq 0$.

By hypothesis $\alpha \notin \mathbf{F}_q$ and $a\mathrm{Tr}(\alpha) + bn \neq 0$; so the test is valid. ■

The situation describes in the previous lemma is the simplest decomposition possible, namely, $\mathbf{F}_{q^n} = V_0 \oplus V_1 = \mathbf{F}_q \oplus \mathbf{F}_{q^n} - \mathbf{F}_q$. Accordingly, a normal element is a sum of two elements $\eta = \eta_0 + \eta_1$, such that $0 \neq \eta_0 \in \mathbf{F}_q$ and $\eta_1 \in \mathbf{F}_{q^n} - \mathbf{F}_q$.

## 3.5 Normal Tests

Other techniques for identifying normal elements are discussed in this section.

Let $\{ \eta_{n-1}, \ldots, \eta_1, \eta_0 \} \subset \mathbf{F}_{q^n}$ be a subset of conjugate elements. The regular matrix representation attached to this subset is given by

$$N = \left(\eta^{q^{i+j}}\right) = \begin{bmatrix} \eta & \eta^{q} & \eta^{q^2} & \cdot & \cdot & \cdot & \eta^{q^{n-1}} \\ \eta^{q} & \eta^{q^2} & \eta^{q^3} & \cdot & \cdot & \cdot & \eta \\ \eta^{q^2} & \eta^{q^3} & \eta^{q^4} & \cdot & \cdot & \cdot & \eta^{q} \\ \cdot & \cdot & \cdot & \cdot & \cdot & \cdot & \cdot \\ \eta^{q^{n-1}} & \eta & \eta^{q} & \cdot & \cdot & \cdot & \eta^{q^{n-2}} \end{bmatrix}$$

and the trace matrix representation is given by

$$T = \left(Tr(\eta^{q^i + q^j})\right)$$

A cursory observation reveals that these two matrices are symmetric.

***Theorem 3.14.*** (*Matrix Normal Test*) A subset of conjugates { $\eta_{n-1}, \ldots, \eta_1, \eta_0$ } elements is a normal basis of $\mathbf{F}_{q^n}$ over $\mathbf{F}_q$ if and only if the matrices N and T are nonsingular.

A few other normality tests are based on polynomials computations. The resultant of a pair of polynomials and the formula for the determinant of a circulant matrix will be employed to prove the next result, see the chapter on circulant matrices for more details.

***Theorem 3.15.*** (*Gcd Normal Test*) Let $\eta \in \mathbf{F}_{q^n}$. Then the conjugate set { $\eta^{q^{n-1}}, \ldots, \eta^q, \eta$ } of $\eta$ forms a normal basis of $\mathbf{F}_{q^n}$ over $\mathbf{F}_q$ if and only if the polynomials

$f(x) = x^n - 1$ and $g(x) = \eta^{q^{n-1}} x^{n-1} + \cdots + \eta^{q^2} x^2 + \eta^q x + \eta$

are relatively prime in $\mathbf{F}_{q^n}[x]$, i.e., $0 \neq \gcd(f(x), g(x)) \in \mathbf{F}_q$.

Proof: Clearly N is a circulant matrix, so the determinant det(N) is equal to the resultant

$$R(f, g) = \pm \prod_{i=0}^{n-1} g(\theta^i),$$

where $\theta$ is a root of $f(x) = x^n - 1$. Moreover, $R(f(x), g(x)) \neq 0$ if and only if f(x) and g(x) shares do not any roots in common, id est, $\gcd(f(x), g(x)) \neq$ nonzero constant. This proves the claim. ■

The core of this proof appears in [1, Davenport]; see also [1, Lidl and Niederreiter].

Each term $g(\omega)$ is an eigenvalue of the circulant matrix N. The first term g(1) in the product R(f,

g) = $\Pi_i g(\theta^i)$ is just the trace Tr(η) = g(1) of the element η. If n is even then −1 is a root of $x^n - 1$, and the term g(−1) is equal to the alternating trace

$$Tr^*(\eta) = (-1)^{n-1}\eta^{q^{n-1}} + \cdots + \eta^{q^2} - \eta^q + \eta$$

which is also a factor of R(f, g) = $\Pi_i g(\theta^i)$, has the value $Tr^*(\eta) = g(-1) \neq 0$. If the element η is a normal element and n is even then both the trace Tr(η) = g(1) and the alternating trace $Tr^*(\eta)$ = g(−1) have nonzero values.

## 3.6 Polynomial Representations

The set of units in (the multiplicative group) the polynomials ring $\mathbf{F}_q[x]/(x^n-1)$ has many properties in common with the set of normal bases of $\mathbf{F}_{q^n}$ over $\mathbf{F}_q$. This analogy provides a vehicle for deriving various properties of the set of normal bases.

The map $\kappa : \mathbf{F}_{q^n} \rightarrow \mathbf{F}_q[x]/(x^n-1)$ defined by

$$\alpha \quad \rightarrow \quad c_\alpha(x) = \sum_{i=0}^{n-1} Tr(\alpha^{1+q^i})x^i = \sum_{i=0}^{n-1} c_i x^i \ .$$

This map is of considerable importance in the investigation of dual normal bases and normal polynomials of degree n in $\mathbf{F}_q[x]$. The polynomial $c_\alpha(x) \in \mathbf{F}_q[x]/(x^n-1)$ is the *polynomial representation* of the element $\alpha \in \mathbf{F}_{q^n}$.

***Lemma 3.16.*** The additive groups $\mathbf{F}_{q^n}$ and $\mathbf{F}_q[x]/(x^n-1)$ are isomorphic.

The relationship between the subset of normal elements of $\mathbf{F}_{q^n}$ and the multiplicative group $\mathbf{U}_n(q)$ of $\mathbf{F}_q[x]/(x^n-1)$ is significant and intricate. The first idea consider the into property of the map under consideration.

***Lemma 3.17.*** The correspondence $\alpha \rightarrow c_\alpha(x)$ is an - to - one, a ≥ 1.

Since the correspondence η → c(x) is not one-to-one, the image of a subset of normal elements is not necessarily a subgroup of the group of units of $\mathbf{F}_q[x]/(x^n-1)$.

Let the map $\tau : \mathbf{F}_q(x) \rightarrow \mathbf{F}_q(x)$ be the involution $\tau(x) = x^{-1}$ in the polynomials ring $\mathbf{F}_q[x]/(x^n-1)$. The self-reciprocal polynomials $f(x) = f^*(x)$, (1-cycles), are the fixed points of the involution τ, and the non self-reciprocal polynomials $f(x) \neq f^*(x)$, (2-cycles), are mapped to their counterparts: $f(x) \rightarrow f^*(x) \rightarrow f(x)$. The map τ is very closely associated with the matrix transpose map t : $GL_n(\mathbf{F}_q) \rightarrow GL_n(\mathbf{F}_q)$, given by $A \rightarrow A^T$.

The polynomial $c_\alpha(x) = a(x)\tau(a(x))$, where $\tau(a(x)) = a^*(x)$ is the reciprocal polynomial of a(x) in $\mathbf{F}_q[x]/(x^n-1)$, establishes a reversed correspondence between the normal element η and $\gamma = a^{-1}(x)\circ\eta$. The coefficients of $a(x) = a_{n-1}x^{n-1} + \cdots + a_1x + a_0 \in \mathbf{F}_q[x]/(x^n-1)$ are the coordinates of the expansion of η with respect to the normal element γ, that is,

$$\eta = a_{n-1}\gamma^{q^{n-1}} + \cdots + a_2\gamma^{q^2} + a_1\gamma^q + a_0\gamma\,,$$

and the coefficients of $a^{-1}(x) = b_{n-1}x^{n-1} + \cdots + b_1x + b_0 \in \mathbf{F}_q[x]/(x^n-1)$ are the coordinates of the expansion of γ with respect to the normal element η, that is,

$$\gamma = b_{n-1}\eta^{q^{n-1}} + \cdots + b_2\eta^{q^2} + b_1\eta^q + b_0\eta\,.$$

The set of polynomials $\mathfrak{I} = \{ c_\alpha(x) : \alpha \in \mathbf{F}_{q^n} \}$ is not a subset of multiplicative group $\mathbf{U}_n(q)$ of units in $\mathbf{F}_q[x]/(x^n-1)$, however the subset of polynomials $\wp = \{ c_\eta(x) : \eta \in \mathbf{F}_{q^n}$ is normal $\}$ is subset of $\mathbf{U}_n(q)$.

***Lemma 3.18.*** *(Lacan et al. )* Suppose $\mathbf{F}_{q^n}$ has at least one self-dual normal basis over $\mathbf{F}_q$. Then the followings hold.

( 1 ) The sets of polynomials $\wp = \{ c_\eta(x) : \eta \in \mathbf{F}_{q^n}$ is normal $\}$, and $\mathfrak{R} = \{ c_\eta(x) = a(x)\tau(a(x)) \}$ are equal.

( 2 ) The set of polynomials $\wp = \{ c_\eta(x) : \eta \in \mathbf{F}_{q^n}$ is normal $\}$ is a multiplicative subgroup of the group of units $\mathbf{U}_n(q)$.

Proof: (1) Assume there exists self-dual normal basis $\{ \eta^{q^{n-1}}, ..., \eta^q, \eta \}$ of $\mathbf{F}_{q^n}$ over $\mathbf{F}_q$, and let

$$\alpha = a_{n-1}\eta^{q^{n-1}} + \cdots + a_1\eta^q + a_0\eta$$

Then using the self-dual basis relation $Tr(\eta_i\eta_j) = \delta_{i,j}$, it follows that

$$\begin{aligned} c_\alpha(x) &= \sum_{i=0}^{n-1} Tr(\alpha^{1+q^i})x^i \\ &= \sum_{i=0}^{n-1} Tr\left(\sum_{j=0}^{n-1} a_j\eta^{q^j}\right)\left(\sum_{j=0}^{n-1} a_j\eta^{q^{j+i}}\right)x^i \\ &= \sum_{i=0}^{n-1}\left(\sum_{j=0}^{n-1} a_j a_{j-i}\right)x^i = a(x)\tau(a(x)), \end{aligned}$$

where $a(x) = a_{n-1}x^{n-1} + \cdots + a_1x + a_0$. Ergo each $\alpha \in \mathbf{F}_{q^n}$ is mapped to a product $a(x)\tau(a(x))$. On the contrary, if there is no self-dual normal basis $\{ \eta^{q^{n-1}}, ..., \eta^q, \eta \}$ of $\mathbf{F}_{q^n}$ over $\mathbf{F}_q$ and sets of polynomials $\wp$ and $\mathfrak{R}$ are equal, then the polynomial $a(x) = 1 \in \mathfrak{R}$, $a(x)\tau(a(x)) = 1$, implies that $c_\alpha(x) = 1 \in \wp$. This in turns implies the existence of a self-dual normal basis of $\mathbf{F}_{q^n}$ of over $\mathbf{F}_q$ in contradiction of the hypothesis. (2) Again let assume that there is a self-dual normal basis generated by η. Then

$$c_\eta(x) = \sum_{i=0}^{n-1} Tr(\eta^{q^i+1})x^i = 1$$

But since $1 \in \wp = \{ c_\eta(x) : \eta \in \mathbf{F}_{q^n}$ is normal $\}$ and $c_\eta(x) = a(x)\tau(a(x))$, $a(x)$ is invertible in $\mathbf{F}_q[x]/(x^n-1)$, and $1 \in \wp$, it is also clear that every normal element $\eta \in \mathbf{F}_{q^n}$ is mapped to an invertible polynomial $c_\eta(x)$ in $\mathbf{F}_q[x]/(x^n-1)$. Thus the set $\wp$ is a group. Conversely, if there is no self-dual normal basis, then $1 \notin \wp$, so it is not a group. ■

Let $\boldsymbol{N} = \{$ normal bases of $\mathbf{F}_{q^n}$ over $\mathbf{F}_q \}$, and let $\ker(\kappa) = \{ c_\eta(x) \in \wp : c_\alpha(x) = 1 \}$ be the kernel of the map $\kappa(\eta) = c_\eta(x)$. The inverse image $\kappa^{-1}(1)$ is the set of all self-dual normal bases of $\mathbf{F}_{q^n}$ over $\mathbf{F}_q$, and the quotient $G = \boldsymbol{N}/\ker(\kappa)$ is a group if $\ker(\kappa) \neq \emptyset$. Note that the subset $\{\alpha, \alpha^q, ..., \alpha^{q^{n-1}}\} \subset \kappa^{-1}(1)$ if α generates a non self-dual normal basis.

***Lemma 3.19.*** Let n = 2k + 1, and let q = 2. Then the quotient group G is a group of order

$$\#G = \prod_{i=1}^{t}(2^{c_i}-1)\prod_{i=t+1}^{u}(2^{d_i}-1).$$

Proof: n odd implies the existence a self-dual normal basis, so $\ker(\kappa) \neq \emptyset$. Moreover, the ratio of the cardinalities of the multiplicative group $\mathbf{U}_n$ and the kernel $\ker(\kappa)$ is

$$\frac{\Phi(x^n-1)}{SDN_n(q)} = (q-1)\prod_{i=1}^{t}(q^{c_i}-1)\prod_{i=t+1}^{u}(q^{d_i}-1),$$

see the sections on distributions of normal bases, and self-dual normal bases below for specific on the integers $c_i$ and $d_i$. ■

***Lemma 3.20.*** Let n be an odd integer, $x^n - 1 = (x-1)f_1(x)f_2(x) \cdots f_{r-1}(x)$, and $q = 2^v$. Then
( 1 ) The polynomial $\kappa(\eta) = c_\eta(x) = c(x)$ is a square in $\mathbf{F}_q[x]/(x^n-1) = A_r \oplus \cdots \oplus A_2 \oplus A_1$, and $c(x) = c_r(x) \oplus \cdots \oplus c_2(x) \oplus c_1(x)$, where each $c_i(x) \equiv c(x) \bmod f_i(x)$ is a square in $A_i = \mathbf{F}_q[x]/(f_i(x))$.
( 2 ) The inverse b(x) of the polynomial c(x) is given by

$$b(x) = \sum_{i=0}^{n-1} b_i x^i = c(x^2)^{2^{uv-1}-1}$$

where u is the order of q modulo n.

Proof: See [3, Poli]. ■

## 3.7 DUAL NORMAL BASES

All pairs of dual bases are linked via nonsingular matrices. Moreover if the bases are normal, then the normal bases change of bases matrices are circulants. To investigate the structure of the matrices in normal bases change of bases and the dual bases equation, let C = ( $c_{i,j}$ ) $\in GL_n(\mathbf{F}_q)$ be a nonsingular matrix, and let.

$$\delta_i = \sum_{j=0}^{n-1} c_{i,j}\eta^{q^j}$$

define a change of basis.

***Theorem 3.21.*** The basis { $\delta_{n-1}, \ldots, \delta_1, \delta_0$ } is a normal basis if and only if the matrix C is circulant.

Proof: Suppose that the matrix is circulant. Then C = cir[$c_{n-1},\ldots,c_1,c_0$] = ( $c_{j-i}$ ), and

$$\delta_i = \sum_{j=0}^{n-1} c_{j-i}\eta^{q^j} = \left(\sum_{j=0}^{n-1} c_{j-i}\eta^{q^{j-i}}\right)^{q^i} = \delta_0^{q^i} .$$

Conversely, if the elements $\delta_{n-1}, \ldots, \delta_1, \delta_0$ are normal, and C = ( $c_{i,j}$ ) is a nonsingular matrix, let $\delta_0$, = δ, then

$$\delta^{q^i} = \left(\sum_{j=0}^{n-1} c_{0,j}\eta^{q^j}\right)^{q^i} = \sum_{j=0}^{n-1} c_{0,j-i}\eta^{q^j} .$$

Thus C = ( $c_{i,j}$ ) = ( $c_{j-i}$ ) is circulant. ■

***Lemma 3.22.*** The dual basis { $\delta_{n-1}, \ldots, \delta_1, \delta_0$ } of a normal basis { $\eta_{n-1}, \ldots, \eta_1, \eta_0$ } is also a normal basis.

Proof: Consider the dual bases matrix equation $ND^t = I_n$. Since regular representation matrix N of a normal basis { $\eta_1, \eta_2, \ldots, \eta_n$ } is a symmetric matrix, we have $ND^t = N^tD^t = DN^t = D^tN = I_n$. This implies that the regular representation matrix D of dual normal basis { $\delta_1, \delta_2, \ldots, \delta_n$ } is a symmetric matrix. But since the regular representation matrix of a basis is a symmetric matrix if and only if the basis is a normal basis, we conclude that D is also a normal basis. ■

Naturally every element is a linear combination of the basis elements, but to give a specific linear combination corresponding to a given element is a different matter. In some cases this can be done. The coefficients of the inverse $b(x) = b_{n-1}x^{n-1} + \cdots + b_1x + b_0$ of the polynomial $a_\eta(x) = a_{n-1}x^{n-1} + \cdots + a_1x + a_0 \in \mathbf{F}_q[x]/(x^n-1)$ are the coordinates of the dual $\delta$ of $\eta$ with respect to the basis generated by $\eta$.

***Theorem 3.23.*** (*Gao 1994*) Let $\{ \eta_{n-1}, \ldots, \eta_1, \eta_0 \}$ be a normal basis of $\mathbf{F}_{q^n}$ over $\mathbf{F}_q$, and let $a(x) = a_{n-1}x^{n-1} + \cdots + a_1x + a_0 \in \mathbf{F}_q[x]$, where $a_i = Tr(\eta_0\eta_i)$. Then the generator of the dual basis $\{ \delta_{n-1}, \ldots, \delta_1, \delta_0 \}$ is given by $\delta = b_{n-1}\eta_{n-1} + \cdots + b_1\eta_1 + b_0\eta_0$, where $a(x)b(x) \equiv 1 \bmod (x^n-1)$.

Proof: The product a(x)b(x) of a(x) and b(x) in $\mathbf{F}_q[x]/(x^n-1)$ is given by

$$\begin{aligned} a(x)b(x) &\equiv \sum_{i=0}^{n-1}\sum_{j=0}^{n-1} a_i b_j x^{i+j} \quad \bmod (x^n-1) \\ &\equiv \sum_{i=0}^{n-1}\left(\sum_{j=0}^{n-1} a_{j-i} b_j\right) x^i, \end{aligned}$$

where the inner sum is the circular convolution of the coordinates of a(x) and b(x). Moreover, the constraint $a(x)b(x) \equiv 1 \bmod (x^n - 1)$ implies

$$\sum_{j=0}^{n-1} a_{j-i} b_j = \begin{cases} 1 & \text{if } i = 0, \\ 0 & \text{if } i \neq 0. \end{cases}$$

Evaluating the dual bases equation $Tr(\eta_i\delta_j)$ returns

$$\begin{aligned} Tr\left(\eta^{q^i}\delta^{q^j}\right) &= Tr\left(\eta^{q^i}\left(\sum_{k=0}^{n-1} b_k \eta^{q^k}\right)^{q^j}\right) \\ &= \sum_{k=0}^{n-1} b_k Tr\left(\eta^{q^i}\eta^{q^{k+j}}\right) \\ &= \sum_{k=0}^{n-1} b_k Tr\left(\eta\eta^{q^{k+j-i}}\right) = \sum_{k=0}^{n-1} a_{k+j-i} b_k \end{aligned}$$

This is precisely the dual bases equation $Tr(\eta_i\delta_j) = \delta_{i,j}$. ■

From this, it can be deduced that the two generators are equal, $\delta = \eta$, if and only if the polynomial $a(x) = b(x) = 1$; so the basis $\{ \eta_{n-1}, \ldots, \eta_1, \eta_0 \}$ is a self-dual basis whenever this occurs. Let

$$\eta\eta^{q^i} = \sum_{j=0}^{n-1} a_{i,j}\eta^{q^j}.$$

***Theorem 3.24.*** The normal basis of $\mathbf{F}_{q^n}$ over $\mathbf{F}_q$ generated by η is a self-dual normal basis if and only if the matrix A = ( $a_{i,j}$ ) = ( $a_{j-i}$ ) is symmetric and $\mathrm{Tr}(\eta^2) = 1$.

Proof: See [3, Gieselmann and Gollmann]. ■

## 3.8. Distribution Of Normal Bases

***Lemma 3.25.*** (*Perlis 1942*) Let element η be a generator of a normal basis of finite field $\mathbf{F}_{q^n}$ over $\mathbf{F}_q$, and let

$$\gamma = c_{n-1}\eta^{q^{n-1}} + \cdots + c_2\eta^{q^2} + c_1\eta^q + c_0\eta,$$

where $c_i \in \mathbf{F}_q$. Then the element γ is a generator of a normal basis of $\mathbf{F}_{q^n}$ over $\mathbf{F}_q$ if an only if

$x^n - 1$ and $c(x) = c_{n-1}x^{n-1} + \cdots + c_1x + c_0$

are relatively prime.

Proof: First observe that if the circulant matrix C = circ$[c_{n-1},\ldots,c_1,c_0]$ is invertible, then the conjugates of the element γ are linearly independent, so the only solution of the equation

$$a_{n-1}\gamma^{q^{n-1}} + \cdots + a_2\gamma^{q^2} + a_1\gamma^q + a_0\gamma = 0$$

is the trivial solution **a** = (0,...,0,0). Conversely, if the subset $\{ \gamma^{q^{n-1}}, ..., \gamma^q, \gamma \}$ are linearly independent, then **a** = (0,...,0,0) is the only solution of the vector equation

$$\begin{pmatrix} a_0 \\ a_1 \\ a_2 \\ \vdots \\ a_{n-1} \end{pmatrix}^T \begin{bmatrix} c_0 & c_1 & c_2 & \cdots & c_{n-1} \\ c_1 & c_0 & c_{n-1} & \cdots & c_2 \\ c_2 & c_1 & c_0 & \cdots & c_3 \\ \vdots & \vdots & \vdots & \ddots & \vdots \\ c_{n-1} & c_{n-2} & c_{n-3} & \cdots & c_0 \end{bmatrix} \begin{pmatrix} \eta \\ \eta^q \\ \eta^{q^2} \\ \vdots \\ \eta^{q^{n-1}} \end{pmatrix} = \mathbf{0}.$$

This implies that C is invertible. To establish the correspondence $\gcd(c(x), x^n-1) = 1$ $\Longleftrightarrow$ $\{ \gamma^{q^{n-1}}, ..., \gamma^q, \gamma \}$ is a basis, it is sufficient to work with the identification c(x) $\Longleftrightarrow$ circ$[c_{n-1},\ldots,c_1,c_0]$ of the multiplicative groups isomorphism $\mathbf{F}_q[x]/(x^n-1)^* \cong C_n(\mathbf{F}_q)$. ■

This result leads to an enumeration formula for the distribution of normal bases over $\mathbf{F}_q$.

***Theorem 3.26.*** The total number of normal bases of $\mathbf{F}_{q^n}$ over $\mathbf{F}_q$ is given by

$$N_n(q) = \frac{1}{n}\Phi(x^n - 1).$$

Proof: The quotient ring $\mathbf{F}_q[x]/(x^n-1)$ contains $\Phi(x^n-1)$ invertibles, which have a one to one correspondence with the normal elements in $\mathbf{F}_{q^n}$, and each normal basis requires n conjugates normal elements. ■

***Example 3.27.*** Determine the number of ordinary bases of $\mathbf{F}_{2^{15}}$ over $\mathbf{F}_2$.
The factorization of the polynomial $x^{15} - 1$ over $\mathbf{F}_2$ is

$x^{15} - 1 = (x - 1)(x^2 + x + 1)(x^4 + x^3 + x^2 + x + 1)(x^4 + x^3 + 1)(x^4 + x^2 + 1) \in \mathbf{F}_2[x]$.

The formula for the number of normal bases gives

$$N_{15}(q) = \frac{1}{15}\Phi(x^{15} - 1) = 2^{15}(1-2^{-1})(1-2^{-2})(1-2^{-4})^3 = 625$$

normal bases over $\mathbf{F}_2$.

Computational techniques for evaluating $\Phi(x^n-1)$ and other arithmetic functions are covered in the chapter on arithmetic functions.

The special case of $n = p^v$, $v \geq 1$, in characteristic p, has the polynomial $x^n - 1 = (x - 1)^n$, so $\Phi(x^n-1) = q^{n-1}(q - 1)$ independent of $v \geq 1$.The number $q^{n-1}(q - 1)$ is the maximal number units or invertibles in $\mathbf{F}_q[x]/(x^n-1)$ possible. For example, if $\eta \in \mathbf{F}_{q^n}$ is a fixed normal element over $\mathbf{F}_q$, and $c(x) = c_{n-1}x^{n-1} + \cdots + c_1x + c_0$, then the element

$$c_{n-1}\eta^{q^{n-1}} + \cdots + c_2\eta^{q^2} + c_1\eta^q + c_0\eta$$

is a normal element in $\mathbf{F}_{q^n}$ over $\mathbf{F}_q$ whenever $c(0) \neq 0$, since $c(x) = c_{n-1}x^{n-1} + \cdots + c_1x + c_0 \in \mathbf{F}_q[x]/(x^n-1)$ is invertible whenever $c(0) \neq 0$.

## 3.9 Distribution Of Self-Dual Normal Bases

The existence of self-dual normal bases of $\mathbf{F}_{2^n}$ over $\mathbf{F}_2$ for $q = 2^n$, and n odd, was recognized about two decades ago. This follows from the uniqueness of dual bases pairs and the fact that the number $\Phi(x^n-1)$ of normal bases of $\mathbf{F}_{2^n}$ over $\mathbf{F}_2$ for odd n is an odd integer. In addition to extension of degree n = odd, some extension of degree n = even and $q = 2^v$ also have self-dual normal bases.

***Theorem 3.28.*** The extension $\mathbf{F}_{q^n}$ of $\mathbf{F}_q$ has a self-dual normal basis for the following parameters.

( 1 ) If and only if n is odd and $q = p^v$ is an odd prime power.
( 2 ) If and only if $n \neq 4k$ and $q = 2^v$, $v \geq 1$

Proof: The subfields of a field which has a self-dual normal basis also have self-dual normal bases, so it is sufficient to prove the nonexistence of self-dual normal basis of the two smaller subfields:
( 1 ) $\mathbf{F}_{q^2}$ over $\mathbf{F}_q$ in characteristic char($\mathbf{F}_q$) ≠ 2, and
( 2 ) $\mathbf{F}_{q^4}$ over $\mathbf{F}_q$ in characteristic char($\mathbf{F}_q$) = 2.
( 1 ) To verify this case, suppose that $\eta$, $\sigma(\eta)$ is a self-dual normal basis of $\mathbf{F}_{q^2}$ over $\mathbf{F}_q$, then

$$0 = Tr(\eta\sigma(\eta)) = 2\eta\sigma(\eta),$$

which is a contradiction in characteristic char($\mathbf{F}_q$) ≠ 2.

( 2 ) And to verify the other case, suppose $\eta$, $\sigma(\eta)$, $\sigma^2(\eta)$, $\sigma^3(\eta)$ is a self-dual normal basis of $\mathbf{F}_{q^4}$ over $\mathbf{F}_q$. Then

$$0 = Tr(\eta\sigma(\eta) = [\eta + \sigma^2(\eta)][\sigma(\eta) +\sigma^3(\eta)].$$

But $\eta + \sigma^2(\eta) = 0$ or $\sigma(\eta) +\sigma^3(\eta) = 0$ contradict the fact that $\eta$, $\sigma(\eta)$, $\sigma^2(\eta)$,$\sigma^3(\eta)$ is basis. ■

The approach taken here is treated in full generality in [1, Bayer et al, p.369.].

***Example 3.29.*** Self-dual normal bases of extensions of $\mathbf{F}_2$. The first two extensions $\mathbf{F}_4 = \mathbf{F}_2[x]/(x^2+x+1) = \{ a_1\eta^2 + a_0\eta : a_i \in \mathbf{F}_2 \}$ and $\mathbf{F}_8 = \mathbf{F}_2[x]/(x^3+x^2+1) = \{ a_2\eta^2 + a_1\eta^2 + a_0\eta : a_i \in \mathbf{F}_2 \}$ of $\mathbf{F}_2$, where $\eta$ is a root of the normal polynomials $x^2 + x + 1$ and $x^3 + x^2 + 1$ respectively, have self-dual normal bases, but the next extension $\mathbf{F}_{16} = \mathbf{F}_2[x]/(x^4+x^3+1)$ of $\mathbf{F}_2$ does not have self-dual normal basis since n = 4. More generally, the sequence of finite fields

$$\mathbf{F}_{2^4},\quad \mathbf{F}_{2^8},\quad \mathbf{F}_{2^{12}},\quad \mathbf{F}_{2^{16}},\quad \ldots,\quad \mathbf{F}_{2^{4k}},\quad \ldots,$$

$k \geq 1$, do not have self-dual normal bases.

***Example 3.30.*** Consider a root $\alpha$ of the irreducible polynomial $f(x) = x^p - x - a \in \mathbf{F}_q[x]$, $q = p^v$,

and $\beta \in \mathbf{F}_q$. Let $\eta = (\beta - \alpha)^{p-1} - 1$, then the list of elements

$$\eta^{q^{p-1}},\ \ldots, \eta^{q^2}, \eta^q, \eta$$

is a self-dual normal basis of $\mathbf{F}_{q^p}$ over $\mathbf{F}_q$.

Now take an odd prime power $q = p^v$ and an odd integer $n = 2k + 1$ or $q = 2^v$ and an odd integer $n \neq 4k$, for $\gcd(n, q) = 1$, let

$$x^n - 1 = (x-1)f_1(x)f_2(x) \cdots f_t(x)g_{t+1}(x)g_{t+1}^*(x)g_{t+2}(x)g_{t+2}^*(x) \cdots g_u(x)g_u^*(x)$$

be the *ordered reciprocal polynomial factorization*. The irreducible factors $f_i(x)$ of even degree $\deg(f_i(x)) = 2c_i$ are self-reciprocals, id est., $f_i(x) = f_i^*(x)$, but the irreducible factors $g_i(x)$ of degree $\deg(g_i(x)) = d_i$ are not self-reciprocals, id est, $g_i(x) \neq g_i^*(x)$.

***Theorem 3.31.*** The number of distinct self-dual normal bases of $\mathbf{F}_{q^n}$ over $\mathbf{F}_q$ is given by

$$SDN_n(q) = \begin{cases} \dfrac{2^a}{n}\displaystyle\prod_{i=1}^{t}(q^{c_i}+1)\prod_{i=t+1}^{u}(q^{d_i}-1) & \text{if } \gcd(n,q)=1, \\ q^{-1+(q-1)(k+b)/2}SDN_k(q) & \text{if } n = kq. \end{cases}$$

where

$$a = \begin{cases} 0 & \text{if } q = 2 \text{ and } n \neq 0 \bmod 4, \\ 1 & \text{if } q \equiv 1 \bmod 2, \text{and } n \equiv 1 \bmod 2, \end{cases} \quad \text{and} \quad b = \begin{cases} 1 & \text{if } q = 2 \text{ and } k \equiv 1 \bmod 2, \\ 0 & \text{if } q \equiv 1 \bmod 2, \text{and } n \equiv 1 \bmod 2. \end{cases}$$

These statistic are derived from cardinalities of the subset $OC_n(\mathbf{F}_q) = O_n(\mathbf{F}_q) \cap C_n(\mathbf{F}_q)$ of nonsingular circulant orthogonal matrices over $\mathbf{F}_q$, see the chapter on structured matrices. A proof of the theorem is given in [2, Jungnickel et al].

***Example 3.32.*** Determine the number of self-dual bases of finite field $\mathbf{F}_{2^{15}}$ over $\mathbf{F}_2$.
The *ordered self-reciprocal factorization* of the polynomial $x^{15} - 1$ over $\mathbf{F}_2$ is

$$\begin{aligned} x^{15} - 1 &= (x-1)f_1(x)f_2(x)g_2(x)g_3^*(x) \\ &= (x-1)(x^2+x+1)(x^4+x^3+x^2+x+1)(x^4+x^3+1)(x^4+x^2+1) \in \mathbf{F}_2[x]. \end{aligned}$$

Using the formula for the number of self-dual bases, to find that $\mathbf{F}_{2^{15}}$ has

$$SDN_{15}(q) = \frac{1}{15}(2+1)(2^2+1)(2^4-1) = 15$$

self-dual normal bases over $\mathbf{F}_2$.

Self-dual bases are of considerable interest in practical in finite field multipliers, see [4, Berlekamp], [3, Geiselmann et al], [1, Stinson] etc.

## 3.10. Formulae For Normal Elements/Polynomials

The general methods for constructing arbitrary normal elements and polynomials are not very practical. Any moderately large parameters n and q requires a significant amount of computations. However there are many specific cases, which are relatively easy to construct for specific parameters n and q.

***Definition 3.33.*** A polynomial with linearly independent roots is called a *normal* polynomial.

The normal polynomial

$$N(x) = (x - g^{q^{n-1}}) \cdots (x - g^q)(x - g),$$

of degree n and parametized by $g = g(r_{n-1}, \ldots, r_1, r_0)$, $r_i \in \mathbf{F}_q$, ranges though all the normal polynomials over $\mathbf{F}_q$ as the vector $r = (r_{n-1}, \ldots, r_1, r_0)$ ranges over certain subset of vectors. This formula is the additive counterpart of the formula $g^v$, $\gcd(v, p-1) = 1$, for reproducing all the primitive roots modulo p from a single primitive root g.

This formula is practical for very small n, and it is illustrated here for n = 2. The derivations are given in [1, Nemoga et al, p.89].

***Example 3.34.*** ( 1 ) If $q = 2^v$, and $x^2 + ax + b$ is the initial normal polynomial, then the entire collection of $q(q-1)^2/2$ normal polynomials over $\mathbf{F}_q$ is reproduced by the formula

$$N(x) = x^2 + a(r_0 + r_1)x + b(r_0 + r_1)^2 + r_0 r_1 a^2,$$

where $r_i \in \mathbf{F}_q$, and $r_0 \neq r_1$. The generator $g = r_0(\alpha + \alpha^q) + r_1(\alpha - \alpha^q)$, $\alpha^2 + a\alpha + b = 0$.

( 2 ) If $q \neq 2^v$, and $x^2 + ax + b$ is the initial normal polynomial, then the entire collection of $(q - 1)^2/2$ normal polynomials over $\mathbf{F}_q$ is reproduced by the formula

$$N(x) = x^2 + 2ar_0 x + b^2 r_0^2 - r_1^2(a^2 - 4b),$$

where $r_i \in \mathbf{F}_q$, and $r_0 r_1 \neq 0$.

An *effective* deterministic algorithm for constructing normal elements is developed in [2, Poli]. This algorithm is also based on the factorization of the polynomial $x^n - 1$, and the standard normal test.

**Other Generators Of Normal Bases**
Several techniques for generating normal elements are explored in this section.

The list of polynomials

$$e_i(x) = \frac{f(x)}{(x-\alpha^{q^i})f(\alpha^{q^i})}$$

is an idempotent basis of the quotient ring $\mathbf{F}_q[x]/(f(x))$. Specifically,
( 1 ) $(e_i(x))^2 \equiv e_i(x) \bmod f(x)$
( 2 ) $e_i(x)e_j(x) \equiv 0 \bmod f(x)$, $i \neq j$, and
( 3 ) $e_{n-1}(x) + \cdots + e_1(x) + e_0(x) \equiv 1$.

***Theorem 3.35.*** *(Artin* 1966) Let α be a root of the irreducible polynomial f(x) of degree n over $\mathbf{F}_q$, and let $\beta \in \mathbf{F}_{q^n}$ .Then the element $e(ß) = e_0(\beta)$ generates a normal basis of $\mathbf{F}_{q^n}$ over $\mathbf{F}_q$ for at least $q - n(n-1)$ distinct $\beta \in \mathbf{F}_q$, $q > n^2$.

Proof: Consider the subset $\{ e_i(x) : 0 \le i < n \}$ of polynomials of degree $\deg(e_i) = i < n$, and the matrix $E(x) = ( e_i(x)e_j(x) )$. From the properties of the polynomials $e_i(x)$ it follows that the matrix $E(x) = ( e_i(x)e_j(x) ) = \mathrm{diag}(e_{n-1}(x),\ldots,e_1(x), e_0(x) )$, (use from the properties (1) and (2) of the set $\{ e_i(x) \}$ just stated), and its determinant $\det(E(x)) = \det((e_{i+j}(x))) = e_{n-1}(x)\cdots e_1(x)e_0(x)$ is a polynomial of degree $n(n-1)$. Further, the elements in the set $\{ e_i(\beta) : \text{fixed } \beta \in \mathbf{F}_q \}$ is a basis if and only if the matrix $E(\beta) = ( e_i(\beta)e_j(\beta) ) = ( e_{i+j}(\beta) )$ is nonsingular. But since the degree of $\det(E(x))$ is $n(n-1)$, the matrix degree is singular for at most $n(n-1)$ distinct $\beta \in \mathbf{F}_q$. ■

**Cyclic Convolution**
The *cyclic convolution* of α and β is defined by

$$\alpha * \beta = \sum_{i=0}^{n-1} \alpha^{q^i} \beta^{q^{-i}} \ ,$$

for all $\alpha, \beta \in \mathbf{F}_{q^n}$.

***Lemma 3.36.*** *(Lacan et al. )* Let $\eta \in \mathbf{F}_{q^n}$ be normal over $\mathbf{F}_q$. Then $L(x) = \eta * x$ is a one-to-one linear map $L : \mathbf{F}_{q^n} \to \mathbf{F}_{q^n}$.

Proof: The linearity $L(ax + by) = \eta*(ax + by) = \eta*ax + \eta*by = aL(x) + bL(y)$ is clear. The verification of the one-to-one property appears in the stated source. ■
The cyclic convolution of normal elements is a method of generating new normal elements from a given normal element.

***Lemma 3.37.*** *(Lacan et al. )* The cyclic convolution $\eta*\gamma$ of a pair of normal elements η, γ in $\mathbf{F}_{q^n}$ is again a normal element in $\mathbf{F}_{q^n}$ over $\mathbf{F}_q$ if and only if n is odd. Furthermore, for every such pair $\eta^*\gamma \neq \eta$ and $\eta^*\gamma \neq \gamma$ if and only if q = 2 and n is odd.

**Algorithm for Constructing Normal Elements**

An algorithm for constructing normal elements of $\mathbf{F}_{q^n}$ over $\mathbf{F}_q$ based on the polynomial $c_\eta(x) = c_{n-1}x^{n-1} + \cdots + c_1x + c_0 = Tr(\eta_0\eta_{n-1})x^{n-1} + \cdots + Tr(\eta_0\eta_1)x + Tr(\eta_0\eta_0) \in \mathbf{F}_q[x]$ is outlined here.

*Step 1.* Choose a random element $\alpha \in \mathbf{F}_{q^n}$ and compute $Tr(\alpha)$. If $Tr(\alpha) = 0$, repeat step 1.

*Step 2.* Compute

$$Tr(\alpha^2),\quad Tr(\alpha^{q+1}),\quad Tr(\alpha^{q^2+1}),\quad ...,\quad Tr(\alpha^{q^{n-1}+1}),$$

and $\gcd(c(x), x^n - 1)$. If $\gcd(c(x), x^n - 1) = 1$, then the element $\alpha$ is normal $\mathbf{F}_{q^n}$ over $\mathbf{F}_q$, in addition if $c_{n-1} = \cdots = c_1 = 0$, and $c_0 = 1$, then $\alpha$ generates a self-dual normal basis. Otherwise repeat step 1.

## 3.11 Completely Normal Bases and Testing Methods

***Definition 3.38.*** A normal element $\eta \in \mathbf{F}_{q^n}$ over $\mathbf{F}_q$ is said to be completely normal if it generates a normal basis $\mathbf{F}_{q^n}$ over $\mathbf{F}_{q^d}$ for any divisor d of n.

***Theorem 3.39.*** (*Blessenohl and Johnsen 1986*) There is a normal element $\eta \in \mathbf{F}_{q^n}$ over $\mathbf{F}_q$ completely normal over $\mathbf{F}_{q^d}$ for all divisors d of n.

The normal basis tests introduced before are easily extended to the completely normal basis case. The extended versions of three tests most widely used are as follows.

***Theorem 3.40.*** (*Matrix Completely Normal Test*) A subset of conjugates $\{\eta_{n-1}, \ldots, \eta_1, \eta_0\}$ elements is a completely normal basis of $\mathbf{F}_{q^n}$ over $\mathbf{F}_q$ if and only if the $e \times e$ circulant submatrices

$$A_d = \left(\eta^{q^{d(j-i)}}\right)$$

are nonsingular for all $d \mid n = de$, $0 \leq i, j < e$.

***Example 3.41.*** For n = 6 the matrix test consist of three matrices:

$$A_3 = \begin{bmatrix} \eta^{q^3} & \eta \\ \eta & \eta^{q^3} \end{bmatrix}, \quad A_2 = \begin{bmatrix} \eta & \eta^{q^2} & \eta^{q^4} \\ \eta^{q^4} & \eta & \eta^{q^2} \\ \eta^{q^2} & \eta^{q^4} & \eta \end{bmatrix}, \quad A_1 = \begin{bmatrix} \eta & \eta^{q} & \eta^{q^2} & \eta^{q^3} & \eta^{q^4} & \eta^{q^5} \\ \eta^{q} & \eta^{q^2} & \eta^{q^3} & \eta^{q^4} & \eta^{q^5} & \eta \\ \eta^{q^2} & \eta^{q^3} & \eta^{q^4} & \eta^{q^5} & \eta & \eta^{q} \\ \eta^{q^3} & \eta^{q^4} & \eta^{q^5} & \eta & \eta^{q^2} & \eta^{q^3} \\ \eta^{q^4} & \eta^{q^5} & \eta & \eta^{q} & \eta^{q^2} & \eta^{q^3} \\ \eta^{q^5} & \eta & \eta^{q} & \eta^{q^2} & \eta^{q^3} & \eta^{q^4} \end{bmatrix}.$$

The element $\eta \in \mathbf{F}_{q^6}$ is a completely normal element over $\mathbf{F}_q$ if and only if $\det(A_3) \neq 0$, $\det(A_2) \neq 0$, and $\det(A_1) \neq 0$.

***Theorem 3.41.*** (*Gcd Completely Normal Test*) The conjugate set $\{ \eta^{q^{d(e-1)}}, ..., \eta^{q^d}, \eta \}$ forms a normal basis of $\mathbf{F}_{q^n}$ over $\mathbf{F}_{q^d}$ if and only if the polynomials

$$f(x) = x^{n/d} - 1 \quad \text{and} \quad g(x) = \eta^{q^{d(e-1)}} x^{e-1} + \cdots + \eta^{q^{2d}} x^2 + \eta^{q^d} x + \eta$$

are relatively prime in $\mathbf{F}_{q^n}[x]$, i.e., $0 \neq \gcd(f(x), g(x)) \in \mathbf{F}_q$. Further if this holds for all $d \mid n$, then the element $\eta \in \mathbf{F}_{q^n}$ is completely normal.

***Theorem 3.42.*** (*Standard Completely Normal Test*) An element $\eta \in \mathbf{F}_{q^n}$ is a completely normal element over $\mathbf{F}_q$ if and only if the system of inequalities

$$\frac{x^e - 1}{a(x)} \circ \eta \neq 0$$

holds for all irreducible factors $a(x)$ of $x^e - 1 \in \mathbf{F}_{q^d}[x]$, and all divisors e of n = de.

The computational complexity of the standard completely normal test is determined by the factorization of the integer n = de, and the polynomial $x^e - 1$. For each divisor d of n, there are $\Omega(x^e - 1)$ inequalities in the system of inequalities, which is the number of irreducible factors in $x^e - 1$ over $\mathbf{F}_{q^d}$. The test is repeated $\tau(n)$ times, which is the number of divisors of the integer n.

For the parameter $n = p^v$, p being the characteristic of the finite field $\mathbf{F}_q$ and $v \geq 1$, the polynomial $x^n - 1 = (x - 1)^n$. Accordingly the system of inequalities in the standard completely normal test collapses to a single line:

$$\frac{x^e - 1}{x - 1} \circ \eta = \eta^{q^{d(e-1)}} + \eta^{q^{d(e-2)}} + \cdots + \eta^{q^d} + \eta \neq 0,$$

but the integer e ranges over the divisors of n, or n = de.

In the standard normal test the parameter k = n is fixed, so the condition $Tr(\alpha) \neq 0$ is sufficient to test the normality of an element in $\mathbf{F}_{q^n}$, $n = p^v$. The transitive property of the trace function allows an extension of this short test to completely normal elements.

***Theorem 3.43.*** Let $f(x) = x^n + a_{n-1}x^{n-1} + \cdots + a_1x + a_0 \in \mathbf{F}_q[x]$ be irreducible of degree deg(f) = $n = p^v$, $p \mid q$, and d be a divisor of n. Then a root η of f(x) is a completely normal element in $\mathbf{F}_{q^n}$ over $\mathbf{F}_{q^d}$ if and only if $a_{n-1} \neq 0$.

Proof: The relative trace is precisely

$$Tr_{n:d}(\eta) = \eta^{q^{d(e-1)}} + \eta^{q^{d(e-2)}} + \cdots + \eta^{q^d} + \eta$$

for all $d \mid n$. Now to confirm that $Tr_{n:d}(\eta) \neq 0$, use the transitivity of the trace function and the data $a_{n-1} \neq 0$, to arrive at $Tr_{n:1}(\eta) = Tr_{d:1}(Tr_{n:d}(\eta)) = a_{n-1} \neq 0$. Ergo $Tr_{n:d}(\eta) \neq 0$. This proves that the element η is completely normal in $\mathbf{F}_{q^n}$ over $\mathbf{F}_{q^d}$. ■

***Theorem 3.441.*** (*Blake et al., 1997*) Suppose $x^n - a \in \mathbf{F}_q[x]$ is irreducible. Then a root $\eta \in \mathbf{F}_{q^n}$ of $ax^n - (x-1)^n$ is a completely normal element over $\mathbf{F}_q$.

The pure equation $x^n - a$ is irreducible whenever the two conditions below hold.
( 1 ) The prime $p_i$ divides m, but not (q – 1)/m, and
( 2 ) If 4 divides n, then 4 divides q – 1,
where the $p_i$ are the prime divisors of n, and m is the order of q modulo n, see [1, Lidl et al, p.124]. The two specific cases:

( i ) The constant *a* is a nonsquare in $\mathbf{F}_q$, and $n = 2^k$,
( ii ) The constant *a* is primitive in $\mathbf{F}_q$, and $n = p_1^{u_1} \cdots p_s^{u_s}$, $q - 1 = p_1^{v_1} \cdots p_s^{v_s}$, $v_i \leq u_i$, are easy to verify.

A table of completely normal polynomials for $p^n < 10^{50}$, $p \leq 97$, of minimal weight is provided in [2. Morgan/Mullin].

**Example 3.45.**

| $\mathbf{F}_2[x]$ | $\mathbf{F}_3[x]$ | $\mathbf{F}_5[x]$ |
|---|---|---|
| $x^2 + x + 1$ | $x^2 + x + 2$ | $x^2 + x + 2$ |
| $x^3 + x^2 + 1$ | $x^3 + 2x^2 + 1$ | $x^3 + x^2 + 2$ |
| $x^4 + x^3 + 1$ | $x^4 + 2x^3 + 2$ | $x^4 + x^2 + 2$ |
| $x^5 + x^4 + x^3 + x^2 + 1 + 1$ | $x^5 + 2x^4 + 1$ | $x^5 + 2x^4 + 3$ |

**Distribution**
The distribution of completely normal elements in for arbitrary parameters n and q is not known. However, for a few specific parameters n , and q it is known. For small pairs (n, q) the machine calculations of the cardinalities of completely normal elements $CN_n(q)$ and completely normal primitive elements $CNP_n(q)$ are feasible.

**Iterated Construction**
The iterated construction of completely normal elements has many similarities to the iterated construction of normal elements. The product $\eta = \eta_1\eta_2$ of two completely normal elements $\eta_1 \in \mathbf{F}_{q^r}$ and $\eta_2 \in \mathbf{F}_{q^s}$ over $\mathbf{F}_q$ such that gcd(r, s) = 1 is again a completely normal element in $\mathbf{F}_{q^n}$ over $\mathbf{F}_q$, n = rs. This is the basic building block employed to assemble arbitrary completely normal elements in $\mathbf{F}_{q^n}$.

# 3.12 Infinite Sequences of Normal Elements/Polynomials

An iterative technique for constructing infinite sequences of irreducible polynomials rests on the polynomials transformation

$$f^Q(x) = (ax)^r f(Q(x))$$

where $f(x) = x^r + a_{r-1}x^{r-1} + a_{r-2}x^{r-2} + \cdots + a_1x + a_0$, Q(x) is a rational function, and $0 \neq a \in \mathbf{F}_q$.

***Lemma 3.46.*** (*Meyn 1990*) Let $f^Q(x) = x^n f(x+1/x)$, $f(x) \in \mathbf{F}_q[x]$ irreducible. Then $f^Q(x)$ is self–reciprocal irreducible over $\mathbf{F}_q$ if and only if $f(2)f(-2) \neq$ square if $\mathbf{F}_q$ .

The individual polynomials in the sequences are given recursive formula

$$f_n(x) = (ax)^{r^n} f_{n-1}(Q(x)),$$

where the initial polynomial $f_0(x) = x^r + a_{r-1}x^{r-1} + a_{r-2}x^{r-2} + \cdots + a_1x + a_0 \in \mathbf{F}_q[x]$ is irreducible, and $n \geq 1$.

The infinite sequence of roots $\{ \alpha_n \}$ of the infinite sequence of polynomials $\{ f_n(x) \}$ induces an infinite tower of finite extensions

$$\mathbf{F}_q \subset \mathbf{F}_{q^r} \subset \mathbf{F}_{q^{r^2}} \subset \mathbf{F}_{q^{r^3}} \subset \cdots \subset \mathbf{F}_{q^{r^\infty}} .$$

***Theorem 3.47.*** (*Varshamov 1989*) Let $f_0(x) = x^p + x^{p-1} + \cdots + x + 1 \in \mathbf{F}_p[x]$, $f_1(x) = f_0(x^p - x - 1)$, and let $f_n(x) = f^*_{n-1}(x^p - x - 1)$. Then $f^*_n(x)$ is a completely normal polynomial of degree $\deg(f^*_n(x)) = p^{n+1}$ for all $n \geq 2$.

One of the simplest construction of an infinite sequence of polynomials uses a quadratic polynomial $f_0(x) = x^2 + a_1x + a_0$, and the parameter r = 2. Several infinite sequences have been constructed in this manner. Here are some typical cases.

***Theorem 3.48.*** (*Chapman 1997*) Let $q \equiv 1 \bmod 4$ be a prime power, and define the sequence

$$f_n(x) = (2x)^{2^n} f_{n-1}\left(\frac{x^2+1}{2x}\right),$$

where the initial polynomial $f_0(x) = x^2 + ax + 1 \in \mathbf{F}_q[x]$ is irreducible, c is a not a square in $\mathbf{F}_q$, and $n \geq 1$. Then any root $\alpha_n \in \mathbf{F}_{q^{2^n}}$ of $f_n(x)$ is a completely normal element over $\mathbf{F}_q$.

***Theorem 3.49.*** (*Chapman 1997*) Let $q \equiv 3 \bmod 4$ be a prime power, and define the sequence

$$f_n(x) = (2x)^{2^n} f_{n-1}\left(\frac{x^2+c}{2x}\right),$$

where the initial polynomial $f_0(x) = x^2 + ax + b \in \mathbf{F}_q[x]$ is irreducible, $b \neq 0$, and c is a not a square in $\mathbf{F}_q$, $n \geq 1$. Then any root $\alpha_n \in \mathbf{F}_{q^{2^n}}$ of $f_n(x)$ is a completely normal element over $\mathbf{F}_q$.

***Theorem 3.50.*** (*Blake et al, 1997*) Let $p \equiv 3 \bmod 4$, let $f_0(x) = x^2 - bx - c \in \mathbf{F}_q[x]$ be irreducible such that $b \neq 2$, and c is a square in $\mathbf{F}_q$. Then any root $\alpha_n \in \mathbf{F}_{q^{2^n}}$ of the polynomial

$$f_n(x) = (x-1)^{2^{n+1}} - b(x-1)^{2^n} x^{2^n} - cx^{2^{n+1}},$$

is a completely normal element over $\mathbf{F}_q$, $n \geq 1$, and $q^v$, v odd.

These results are described in fine details in [1,2 Myen], [1, Chapman], and [1, Blake et al. 1997]. Methods for generating infinite sequences of normal polynomials $\{ f_n(x) \}$ of degree $\deg(f_n) = r^n$ whose roots are trace compatible, and the parameters

( 1 ) $r = 2$, $q \equiv 1 \bmod 4$,
( 2 ) $r = p$, the characteristic of $\mathbf{F}_q$, $q = p^k$,
( 3 ) r = odd, and r divides $q^2 - 1$,

are developed in [1, Scheerhorn, 1994].

***Theorem 3.51.*** (*Blake et al, 1997*) Let $\alpha_n$ be a root of the cyclotomic polynomial $\Phi_{p^{n+1}}(x)$, and let 2 be a primitive root modulo $p^2$. Then

( 1 ) $\eta_n = \alpha_n + \alpha_n^p + \alpha_n^{p^2} + \cdots + \alpha_n^{p^n}$ is a normal element in $\mathbf{F}_{2^{(p-1)p^n}}$ over $\mathbf{F}_2$, $n \geq 1$.

( 2 ) $\delta_n = \eta_n^{(p-1)p^n} + \sum_{i=0}^{(p-1)p^n-1} \eta_n^{2^i}$ is the dual normal element.

The proof of (1) is based on the determinant of the trace matrix T = ( $Tr(\alpha_n^{2^i}\alpha_n^{2^j})$ ). And that of (2) uses the polynomial representation

$$N(x) = \sum_{i=0}^{(p-1)p^n-1} Tr(\eta_n^{1=2^i})x^i = x_n^{(p-1)p^n/2} + \sum_{i=0}^{(p-1)p^n-1} x^i$$

of the normal element $\eta_n$. The inverse of N(x) is easy to calculate because this polynomial is an idempotent in $\mathbf{F}_2[x]/(x^{(p-1)p^n}-1)$ > Specifically

$$N^2(x) = \left( x_n^{(p-1)p^n/2} + \sum_{i=0}^{(p-1)p^n-1} x^i \right)^2 = x_n^{(p-1)p^n} + \left( \frac{x^{(p-1)p^n/2} - 1}{x-1} \right)^2 (x^{(p-1)p^n} - 1).$$

Thus $N^{-1}(x) = N(x)$ since $N^2(x) \equiv 1 \bmod (x^{(p-1)p^n} - 1)$.

***Definition 3.52.*** A sequence of polynomials { $f_n(x)$ } of degree $\deg(f_n) = r^n$ is trace compatible if the relative trace $Tr_{n+1:n}(\alpha_{n+1}) = \alpha_n$. for each pair of roots $\alpha_{n+1} \in \mathbf{F}_{q^{r^{n+1}}}$ and $\alpha_n \in \mathbf{F}_{q^{r^n}}$ of $f_{n+1}(x)$ and $f_n(x)$. More precisely $f_{n+1}(\alpha_{n+1}) = 0$, and $f_n(Tr_{n+1:n}(\alpha_{n+1})) = f_n(\alpha_n) = 0$.

# 3.13 Characteristic Functions

A characteristic function encapsulates certain properties of a subset of elements of $\mathbf{F}_{q^n}$. It effectively filters out those elements that do not satisfy the constraints. The equation of a characteristic function is of the form

$$C(\alpha) = \begin{cases} 1 & \text{if the properties are valid,} \\ 0 & \text{otherwise,} \end{cases}$$

for all $\alpha \in \mathbf{F}_{q^n}$.

**The Characteristic Function of Primitive Elements**

The characteristic function of primitive elements of $\mathbf{F}_{q^n}$ is constant on the subset of primitive elements and vanishes otherwise.

Let $\chi$ be a multiplicative character of order d = *ord*($\chi$), d | $q^n - 1$, on $\mathbf{F}_{q^n}$.The characteristic function of primitive elements in $\mathbf{F}_{q^n}$ is defined by

$$C_P(\alpha) = \frac{\varphi(q^n - 1)}{q^n - 1} \sum_{d \mid q^n - 1} \frac{\mu(d)}{\varphi(d)} \sum_{ord(\chi) = d} \chi(\alpha),$$

where $\alpha \in \mathbf{F}_{q^n}$, and the arithmetic functions $\mu$ and $\varphi$ are the Mobius and Euler functions on the ring of integers **Z** respectively.

A product version of this formula

$$C_P(\alpha) = \frac{\varphi(q^n - 1)}{q^n - 1} \prod_{p \mid q^n - 1} \left( 1 - \frac{1}{p - 1} \sum_{ord(\chi) = p} \chi(\alpha) \right) ,$$

where p runs through the prime divisors of $q^n - 1$, is also effective in certain applications. The transformation required is straightforward, see [1, Hua, p. 177], etc.

The function $C_P(\alpha)$ is one of the basic tools used in the investigation of the distribution of primitive elements in finite fields. Typical applications are illustrated in [1, Jungnickel, et al], [1, Moreno], etc.

**The Characteristic Function of Normal Elements**

Let $\psi$ be an additive character of order d(x) = *Ord*($\psi$) | $x^n - 1$ on $\mathbf{F}_{q^n}$. The characteristic function of normal elements in $\mathbf{F}_{q^n}$ is defined by

$$C_N(\alpha) = \frac{\Phi(x^n - 1)}{q^n} \sum_{d(x) \mid x^n - 1} \frac{M(d(x))}{\Phi(d(x))} \sum_{Ord(\psi) = d(x)} \psi([(x^n - 1)/d(x)] \circ \alpha),$$

where $\alpha \in \mathbf{F}_{q^n}$, and the arithmetic functions M and $\Omega$ are the Mobius and Euler functions on the ring of polynomials $\mathbf{F}_q[x]$ respectively.

A product version of this formula

$$C_N(\alpha) = \frac{\Phi(x^n - 1)}{q^n} \prod_{f(x) \mid x^n - 1} \left( 1 - \frac{1}{q^{\deg(f(x))} - 1} \sum_{Ord(\psi) = f(x)} \psi([x^n - 1)/f(x)] \circ \alpha) \right)$$

where f(x) runs through the irreducible factors of $x^n - 1$, is also effective in certain applications. The transformation required is straight forward, see [1, Lenstra et al], etc.

***Example 3.53.*** For the parameter $n = p^u$, $q = p^v$, the polynomial $x^n - 1 = (x-1)^n \in \mathbf{F}_q[x]$ and the expression $[(x^n - 1)/g(x)]\circ\alpha = Tr(\alpha)$ is the trace $Tr : \mathbf{F}_{q^n} \to \mathbf{F}_q$, (since $g(x) = x - 1$ is the only irreducible factor of $x^n - 1$). Thus the characteristic function of normal elements in $\mathbf{F}_{q^n}$ is

$$C_N(\alpha) = (1 - 1/q)\left(1 - \frac{1}{q-1}\sum_{Ord(\psi)=x-1}\psi(Tr(\alpha))\right) = \begin{cases} 0 & if\ Tr(\alpha) = 0, \\ 1 & if\ Tr(\alpha) \neq 0. \end{cases}.$$

**The Characteristic Function of Primitive Normal Elements**

The product of the characteristic function of primitive elements and the characteristic function of normal elements in $\mathbf{F}_{q^n}$ yield the characteristic function of primitive normal elements.

Let $\chi$ be a multiplicative character of order $d = ord(\chi)$, $d \mid q^n - 1$, and let $\psi$ be an additive character of order $e(x) = Ord(\psi) \mid x^n - 1$ on $\mathbf{F}_{q^n}$**.** The characteristic function of primitive normal elements in $\mathbf{F}_{q^n}$ is defined by

$$C_{PN}(\alpha) = \frac{\varphi(q^n - 1)}{q^n - 1}\frac{\Phi(x^n - 1)}{q^n}\sum_{d \mid q^n - 1}\frac{\mu(d)}{\varphi(d)}\sum_{e(x)\mid x^n - 1}\frac{M(e(x))}{\Phi(e(x))}\sum_{Ord(\psi)=e(x)}\sum_{ord(\chi)=d}\chi(\alpha)\psi(\beta)$$

where $\alpha \in \mathbf{F}_{q^n}$, and $\beta = [(x^n - 1)/e(x)]\circ\alpha$.

The function $C(\alpha)$ is one of the basic tools used in the investigation of the distribution of primitive normal elements in finite fields. Typical applications are illustrated in [1, Lenstra et al], and [1, Carella], etc.

**The Characteristic Function of Completely Normal Elements**

The characteristic function of completely normal elements in $\mathbf{F}_{q^n}$ is constructed from a series of characteristic functions of normal elements in $\mathbf{F}_{q^n}$ over $\mathbf{F}_{q^d}$, $d \mid n = de$. These functions are for extensions of degree $e = [\mathbf{F}_{q^n} : \mathbf{F}_{q^d}]$. Since an element $\eta \in \mathbf{F}_{q^n}$ is a completely normal if and only if $\eta$ is a normal element in $\mathbf{F}_{q^n}$ over $\mathbf{F}_{q^d}$ for all $d \mid n$, it follows that the characteristic function of completely normal elements is the product of the individuals functions:

$$C_{CN}(\alpha) = \prod_{e \mid n}\left(\frac{\Phi(x^e - 1)}{q^n}\sum_{f(x)\mid x^e - 1}\frac{M(f(x))}{\Phi(f(x))}\sum_{Ord(\psi)=f(x)}\psi([x^e - 1)/f(x)]\circ\alpha)\right)$$

$$= \prod_{e \mid n}\prod_{g(x)\mid x^e - 1}\frac{\Phi(x^e - 1)}{q^n}\left(1 - \frac{1}{q^{\deg(g(x))} - 1}\sum_{Ord(\psi)=g(x)}\psi([x^e - 1)/g(x)]\circ\alpha)\right)$$

where $\alpha \in \mathbf{F}_{q^n}$, and $g(x)$ runs through the irreducible factors of $x^e - 1 \in \mathbf{F}_{q^d}[x]$.

***Example 3.54.*** For the parameter $n = p^u$, $q = p^v$, the polynomial $x^n - 1 = (x - 1)^n \in \mathbf{F}_{q^{p^i}}[x]$ and the expression $[(x^n - 1)/g(x)]°\alpha = Tr_i(\alpha)$ is the trace $Tr_i : \mathbf{F}_{q^{p^u}} \rightarrow \mathbf{F}_{q^{p^i}}$, $0 \le i < u$, (since $g(x) = x - 1$ is the only irreducible factor of $x^n - 1$). Thus the characteristic function of completely normal elements in $\mathbf{F}_{q^n}$ is

$$C_{CN}(\alpha) = \prod_{i=0}^{u-1}\left(1 - 1/q^{p^i}\right)\left(1 - \frac{1}{q^{p^i} - 1}\sum_{Ord(\psi)=x-1}\psi(Tr(\alpha))\right).$$

# 3.14 Primitive Normal Bases

A primitive normal basis is generated by a primitive element in $\mathbf{F}_{q^n}$. The asymptotic proof of the Primitive Normal Basis Theorem was first established by both [1, Carlitz], and [1, Davenport]. And the final version for all pair n, q was established by [1, Lenstra and Schoof].

***Theorem 3.55.*.** (*Primitive Normal Basis Theorem*) Let $\mathbf{F}_{q^n}$ be an n-degree extension of $\mathbf{F}_q$. Then $\mathbf{F}_{q^n}$ has a primitive normal basis over $\mathbf{F}_q$.

The next result is a refinement of the Primitive Normal Basis Theorem, it calls for primitive normal elements of arbitrary traces.

***Theorem 3.56.*** (*Primitive Normal Basis Theorem Of Arbitrary Trace*) For every $a \neq 0$ in $\mathbf{F}_q$, there exists a primitive normal element in $\mathbf{F}_{q^n}$ of trace $a$.

This modification was proposed as a conjecture by [1, Morgan and Mullin]. The asymptotic proof of the Primitive Normal Basis Theorem of Arbitrary Trace was first completed in [1, Carella]. And about a year later it was extended to all pairs $n$, $q$ by [2, Cohen et al.].

**The Distribution Of Primitive Normal Elements**

The distribution of primitive normal elements is more intricate than either the distribution of primitive elements or the distribution of normal elements An exact closed form formula for the number of primitive normal bases of $\mathbf{F}_{q^n}$ over $\mathbf{F}_q$ appear to be unknown, however there is an asymptotic approximation due to [1, Carlitz]. The approximation is

$$PN_n(q) = \frac{\varphi(q^n - 1)\Phi(x^n - 1)}{q^n} + O(q^{(.5+\varepsilon)n}),$$

for all $\varepsilon > 0$.

***Conjecture 3.57.*** *(Morgan-Mullin* 1996) Let $q$ be a prime power and let $n > 3$ be an integer. Then there exists a completely normal primitive polynomial of degree $n$ over $\mathbf{F}_q$ .

## 3.15 Applications Of Fractional Linear Transformations To Normal Bases

The projective line over the finite field $\mathbf{F}_q$ consists of the set of points $P^1(\mathbf{F}_q) = \mathbf{F}_q$ □ {□},∞ represents the point at infinity. The projective linear group $PGL_2(\mathbf{F}_q) = GL_2(\mathbf{F}_q)/\{\ aI_2\ \}$ acts on $P^1(\mathbf{F}_q)$ via fractional linear transformations. The operation is composition of maps.

Let $\gamma$ be the $2\times2$ matrix $\begin{bmatrix} a & b \\ c & d \end{bmatrix} \in GL_2(F_q)$. A fractional linear transformation is defined by the map $\gamma(z) = (az + b)/(cz + d)$ on $P^1(F_q)$.

More generally, a fractional linear transformation is a permutation on $\mathbf{F}_{q^n}$ □ {□} for all $n \geq 1$.

**Classification of Fractional Linear Transformations**
A map $\gamma(z) = (az + b)/(cz + d)$ on $P^1(\mathbf{F}_q)$. has at most two fixed points. The fixed points of the map $\gamma(z)$ are the solution of the equation $\gamma(z) = z$. The discriminant $disc(\gamma) = (a - d)^2 + 4bc$ is the discriminant of the fixed points equation $cz^2 - (a - d)z - b = 0$.

The order $n = ord(\gamma)$ of the map $\gamma(z)$ is the smallest integer n which satisfies $\gamma^{n-1}(z) = 1$ for all z □ $P^1(\mathbf{F}_q)$. The integer n is a divisor of the order $q(q - 1)(q + 1) = \#\ PGL_2(\mathbf{F}_q)$ of the group $PGL_2(\mathbf{F}_q)$. The discriminant $disc(\gamma)$ of the fixed points equation serves as an indicator of the order of the map $\gamma(z)$.

( 1 ) If discriminant $disc(\gamma) = 0$, then order $ord(\gamma)$ divides q.
( 2 ) If discriminant $disc(\gamma)$ □ 0 is a quadratic residue in $\mathbf{F}_q$, then order $ord(\gamma)$ divides $q - 1$.
( 3 ) If discriminant $disc(\gamma)$ □ 0 is a nonquadratic residue in $\mathbf{F}_q$, then order $ord(\gamma)$ divides $q + 1$.

***Lemma 3.58.*** If $f(x) \in \mathbf{F}_q[x]$ is irreducible, then $f(\gamma(x))$ is again irreducible over $\mathbf{F}_q$.

Proof: The invertible map $\alpha^{q^i} \to \gamma(\alpha^{q^i}) = \dfrac{a\alpha^{q^i} + b}{c\alpha^{q^i} + d}$ preserves the number of conjugates in a set of conjugates, so the factors of $f(x)$ are matched 1-to-1 to the factors of $f(\gamma(x))$. ■

The fractional linear transformations $\gamma(z) = (az + b)/(cz + d)$, $ad - bc \neq 0$ for which the list of elements { $\gamma(\eta)^{q^{n-1}}$ ..., $\gamma(\eta)^q, \gamma(\eta)$ } forms a normal basis of $\mathbf{F}_{q^n}$ over $\mathbf{F}_q$ are stated below.

***Theorem 3.59.*** (*Sidel'nikov 1988*) Let $\gamma(z) = (az + b)/(cz + d)$ be a map of order $n = ord(\gamma)$. Suppose that the element $\eta \in \mathbf{F}_{q^n} - \mathbf{F}_q$ of nonzero trace $\mathrm{Tr}(\eta) \in \mathbf{F}_q$ is a root of an irreducible factor $f(x)$ of $F(x) = (cx + d)x^q - (ax + b) \in \mathbf{F}_q[x]$. Then the lists of elements { $\gamma(\eta)^{q^{n-1}}$ ..., $\gamma(\eta)^q, \gamma(\eta)$ } $\subset \mathbf{F}_{q^n}$ are linearly independent over $\mathbf{F}_q$. Furthermore, the matrix for change of basis $\eta_t \rightarrow \eta\eta_t$ is associated with the multiplication submatrix

$$\eta_0 \begin{pmatrix} \eta_0 \\ \eta_1 \\ \vdots \\ \eta_{n-1} \end{pmatrix} = \begin{bmatrix} \tau - \varepsilon & -e_{n-1} & \cdots & e_1 \\ e_{n-1} & e_{n-1} & \cdots & 0 \\ \vdots & \vdots & \vdots & \vdots \\ e_{n-1} & 0 & \cdots & e_{n-1} \end{bmatrix} \begin{pmatrix} \eta_0 \\ \eta_1 \\ \vdots \\ \eta_{n-1} \end{pmatrix} + \begin{pmatrix} b_0 \\ b_1 \\ \vdots \\ -b_{n-1} \end{pmatrix},$$

where $e_1 = a$, $e_{i+1} = \gamma(e_i)$, $t = Tr(\eta_0)$, and $\tau^* = \tau - e$, with
$e = \notin \mathbf{F}_q$.

By definition of the polynomial f(x) and the map φ(z), the φ-images of a root η of f(x) are

$$\eta^q = \varphi(\eta) = \frac{a\eta + b}{c\eta + d}, \eta^{q^2} = \varphi^2(\eta) = \frac{a\eta^q + b}{c\eta^q + d}, \ldots, \eta^{q^{n-1}} = \varphi^{n-1}(\eta) = \frac{a\eta^{q^{n-2}} + b}{c\eta^{q^{n-2}} + d}$$

conjugates.

The map φ(z) imposes certain structure on the φ–orbit ορβ(α) = {α, φ(α), $\varphi^2(\alpha)$, ..., $\varphi^{\nu-1}(\alpha)$}of an element $\alpha \notin \mathbf{F}_{q^n}$. And by definition of the polynomial $F(x) = (cx + d)x^q - (ax + b) \notin \mathbf{F}_q[x]$, so if the element α is a root of F(x), then we have

$$\alpha^q = (c\alpha + d)/(a\alpha + b) = \varphi(\alpha).$$

In other words, the map φ acts as conjugation on the subset $\{\alpha \notin \mathbf{F}_{q^n} : F(\alpha) = 0\}$, $\mathbf{F}_{q^n}$ being the splitting field of F(x). Thus the φ-orbit orb(α) of any root α of F(x) coincides with the conjugates of α, so the subset orb(α) are the roots of some factor f(x) of F(x).

Since the map φ(z) is of order n = ord(φ), the φ-orbit orb(α) are of length 1 #orb(α) or n #orb(α). For example, if the element $\alpha \in \mathbf{F}_q$ is a root of the polynomial $F(x) = (cx + d)x^q - (ax + b) \notin \mathbf{F}_q[x]$, then the φ-orbit orb(α) = {α} is itself, so it corresponds to a linear factor x - α of F(x); but if the element $\alpha \in \mathbf{F}_{q^n} - \mathbf{F}_q$ is a root of F(x), then the φ-orbit orb(α) = {α, φ (α),$\varphi^2(\alpha)$, ..., $\varphi^{n-1}(\alpha)$} corresponds to a nonlinear factor f(x) of F(x) of degree n = deg(f(x)).

***Theorem 3.59.*** (*Sidel'nikov 1988*) Let _(z) = (az + b)/(cz + d) be a fractional linear transformation and suppose $\eta_0 = \eta$, $\eta_1 = \varphi(\eta)$, $\eta_2 = \varphi^2(\eta)$, ..., $\eta_{v-1} = \varphi^{v-1}(\eta)$ is a normal basis of $GF(q^n)$ over GF(q). Then the product terms are given by

$$\eta_i\eta_j = c_{i-j}\eta_i + c_{j-i}\eta_j + s,\ c_i, s \in \mathbf{F}_q.$$

The factorization of the polynomial falls into four different cases.
(1) c = 0, and d = 1;
(2) $(a - d)^2 + 4bc = 0$;
(3) $(a - d)^2 + 4bc \neq 0$ is a quadratic residue in $\mathbf{F}_q$; and
(4) $(a - d)^2 + 4bc \neq 0$ is not a quadratic residue in $\mathbf{F}_q$.

***Theorem 3.60.*** (Blake et al *1997*) Let _(z) = (az + b)/(cz + d), $c \neq 0$, $ad - bc \neq 0$, $d(\varphi) = (a - d)^2 + 4bc = 0$ and $x_0$ be a root of $\varphi(z) = z$. Then

$$(cx + d)x^q - (ax + b) = (x - x_0)\prod_{\beta \in T}\left((x - x_0)^p + (a/c - x_0)\beta^{-1}(x - x_0)^{p-1} - (a/c - x_0)^p\beta^{-1}\right)$$

where $T = \{ 0 \neq \beta \in \mathbf{F}_q : Tr(\beta) = 1 \}$, and $Tr : \mathbf{F}_q \to \mathbf{F}_p$.

For a proof, see Blake et al., [1], Theorem 3.4, p.505.

***Theorem 3.61.*** (Blake *et al 19*) (1) The roots of the polynomial $x^p - a\beta^{-1}x - a^p\beta^{-1}$ form a normal basis of $\mathbf{F}_{q^p}$ over $\mathbf{F}_q$. Moreover this basis has a complexity of $C_N = 3p - 2$.

(2) If the element $a = \beta$, then the roots of the polynomial $x^p - x - a^{p-1}$ form a self-dual normal basis of $\mathbf{F}_{q^p}$ over $\mathbf{F}_q$.

***Theorem 3.62.*** (Blake *et al 19*) (1) Let $q = p^s$, and let $n \mid (q - 1)(q + 1)$ or $n = p$. Then the element $\delta_i = \delta^{q^i} = \dfrac{\eta^{q^i} + v}{t(t+nv)}$ is the dual of $\eta_i = \dfrac{a\eta^{q^i} + b}{c\eta^{q^i} + d}$.

Proof: For $i \neq j$,

$$Tr(\eta_i\delta_j) = Tr\left(\eta_i\,\frac{\eta_j + v}{t(t+nv)}\right) = \frac{tv + Tr(\eta_i\eta_j)}{t(t+nv)}.$$

Taking the trace of $\eta_i\eta_j = e_{j-i}\eta_i + e_{i-j}\eta_j + u$, yields $Tr(\eta_i\delta_j) = 0$. And for i = j,

$$Tr(\eta_i\delta_i) = Tr\left(\eta_i\,\frac{\eta_i + v}{t(t+nv)}\right) = \frac{1}{t(t+nv)}Tr\left(\eta_i\left(t + v - \sum_{i\neq j=0}^{n-1}\eta_j\right)\right) = 1$$

as claimed. ■

# Chapter 4

# General Periods

## 4.1 Concept of Periods

A period can be viewed as a finite sum $\Sigma_{x\in K}\psi(x)$ over a subset $K \subset \mathbf{R}$ of a ring $\mathbf{R}$. The function $\psi$ is defined on $\mathbf{R}$, and assumed to be fixed. The notion of periods is quite general and can be defined in any arbitrary sets as well as groups. The most common and best known case is cyclotomic period. A cyclotomic period is simply an incomplete exponential sum $\Sigma_{x\in K}\psi(x)$ with respect to a fixed character. The character is the complex-valued function $\psi(x) = e^{i2\pi x/p}$ or the finite field-valued function $\psi(x) = \omega^x$, where $\omega$ is a primitive pth root of unity in a finite extension $\mathbf{F}_{q^n}$ of $\mathbf{F}_q$, and the coset is a residue class $K = K_i$ of $\mathbf{F}_p$. The cyclotomic periods for several other exponential functions have also been investigated in the literature. In [2, McEnliece et al], the coset $K$ is a residue class $K_i$ of an extension $\mathbf{F}_{q^n}$ of $\mathbf{F}_p$, and the cyclotomic hyperperiods using $\psi(x) = \Sigma_{t\neq 0}e^{i2(t+x/t)/p}$ are considered in [2, Lehmer] and other. . Several generalizations of the cyclotomic periods are considered and some of the basic results on this topic will be considered in this chapter.

***Definition 4.1.*** Let $S$ be a nonempty set, put $S = S_0 \cup S_1 \cup\cdots \cup S_{n-1}$ and let $f$ be a function on $S$. The periods are defined by the list of elements

$$\eta_0 = \sum_{x\in S_0} f(x),\quad \eta_1 = \sum_{x\in S_1} f(x),\ \ldots,\ \eta_{n-1} = \sum_{x\in S_{n-1}} f(x)$$

The number of distinct periods in the list is the *degree* of the periods. The simplest nontrivial periods are the quadratic periods for which $S = S_0 \cup S_1$, $k = \#S_i$ $2k = \#S$, and the periods of degree $kn$.

As an instance, consider the function $f$ on the set $\Omega = \{ 1 \neq \omega : \omega^p = 1 \} \subset \mathbf{R}$, $\mathbf{R}$ a ring, $p = 2n + 1$, and let the subsets $S_i = \{ \omega_i, \omega_{n+i} \}$, $0 \le i < n$, be a partition of $\Omega$. Then the periods of degree n are given by

$$\eta_0 = f(\omega_0) + f(\omega_n),\quad \eta_1 = f(\omega_1) + f(\omega_{n+1}),\quad \ldots,\quad \eta_{n-1} = f(\omega_{n-1}) + f(\omega_{2n-1}).$$

Each summation here is over a subset of points $S_i \subset S$, and the function f is fixed.

Another related elements are the *coperiods*. The summation in the definition of a coperiod is over a subset of functions $L_i \subset L \subset \mathbf{L}(S)$ of the set of functions on a set $S$, and the point $x \in S$ is fixed.

***Definition 4.2.*** Let $L = L_0 \cup L_1 \cup \cdots \cup L_{n-1}$ be a partition of the set of functions $L$ on $S$, and let $x \in S$ be a fixed point. The coperiods are defined by the list of elements

$$\gamma_0 = \sum_{f\in L_0} f(x),\quad \gamma_1 = \sum_{f\in L_1} f(x),\quad \ldots,\quad \gamma_{n-1} = \sum_{f\in L_{n-1}} f(x).$$

***Example 4.3.*** Take the prime integer $r = 2n + 1$, and a prime power $q$ of order $ord_r(q) = 2n$ or n modulo $r$. Consider the partition $L_0 = \{ 1, \sigma^n \}$, $L_1 = \{ \sigma, \sigma^{n+1} \}$, .., $L_i = \{ \sigma^i, \sigma^{n+i} \}$, ..., $L_{n-1} = \{ \sigma^{n-1}, \sigma^{2n-1} \}$ of galois group $Gal(\mathbf{F}_{q^n}/\mathbf{F}_q) = \{ 1, \sigma, \sigma^2, ..., \sigma^{2n-1} \}$ of the extension $\mathbf{F}_{q^n}$ of $\mathbf{F}_q$ of degree $2n$. Then the coperiods in $\mathbf{F}_{q^n}$ are given by the list

$$\gamma_0 = \sum_{j \in S_0} \sigma^j(\omega) = \omega + \omega^{-1}, \quad \gamma_1 = \sum_{j \in S_1} \sigma^j(\omega) = \omega^q + \omega^{-q}, \quad ...,$$
$$\gamma_{n-1} = \sum_{j \in S_{n-1}} \sigma^j(\omega) = \omega^{^{n-1}q} + \omega^{-q^{n-1}}$$

In practice, is quite common to use the partition $S_0 = \{ 0, n \}$, $S_1 = \{ 1, n + 1 \}$, ..., $S_i = \{ i, n + i\}$, ..., $S_{n-1} = \{ n - 1, 2n - 1 \}$ of the isomorphic group $\mathbf{Z}_{2n}$ as index in place of $L_i = \{ \sigma^i, \sigma^{n+i} \}$.

In some cases the periods and coperiods are exactly the same elements. This is the case for the function $f_x(\omega) = \sigma^x(\omega)$, where $\sigma^x \in G$, so $\eta_i = \gamma_i$.

**Periods Generating Maps**

The functions appearing in the definition of the periods are called *periods generating maps*. The automorphisms of the group $G = Gal(\mathbf{Q}(\omega),\mathbf{Q})$, where $\omega$ is a primitive pth root of unity in $\mathbf{Q}(\omega)$ or $G = Gal(\mathbf{F}_{q^n}/\mathbf{F}_q) = \{ 1, \sigma, \sigma^2, .., \sigma^{n-1} \}$, where $\mathbf{F}_{q^n}$ is an extension of $\mathbf{F}_q$ of degree n, and characters group $G^\wedge = \{$ characters on a group $G \}$ have very rich structures, and are probably some of the most interesting periods generating maps. Other structured functions are of interest too.

***Definition 4.4.*** Let $k, n, r \in \mathbf{N}$ be integers such that $\varphi(r) = kn$, $r = r_1 r_2$, where $r_1$ is squarefree, and let $v_l = v_l(r_2)$ be the valuation at the prime $l$. Then the function

$$f(x) = x^{r_2} \prod_{\ell | r_2} \sum_{1 \le i \le v_\ell} x^{r\ell^{-i}}$$

is a periods generating map.

For squarefree integers $r$, this reduces to $f(x) = x$, for more details see [1, Feisel et al.].

**Partitions of Sets**

The partitions used in the construction of *cyclotomic periods* and period normal bases require well-structured partitions, and subsets of uniform cardinalities. The nth power residues partitions have the structure sought after, so the sets S are the multiplicative groups $\mathbf{Z}_r^*$ of the residue numbers systems $\mathbf{Z}_r$, $r \in \mathbf{N}$, and $\varphi(r) = kn$.

The concentration here will be on the following groups and sets:

( 1 ) $\mathbf{Z}_r^* = \{ x \in \mathbf{Z}_r : gcd(x, r ) = 1 \}$,

( 2 ) $Gal(\mathbf{F}_{q^n}/\mathbf{F}_q) = \{ 1, \sigma, \sigma^2, ..., \sigma^{n-1} \}$, where $\mathbf{F}_{q^n}$ is a finite extension of $\mathbf{F}_q$ of degree n.

( 3 ) $Gal(\mathbf{Q}(\omega),\mathbf{Q})$, where ω is a primitive $p$th root of unity in $\mathbf{Q}(\omega)$,
( 4 ) $G^\wedge$ = { characters on a group $G$ }, and the set
( 5 ) $\Omega = \{ 1 \neq \omega : \omega$ is an $r$th primitive root of unity in $\mathbf{C}$ or $\mathbf{F}_{q^n}$ $\}$, $r = kn + 1$ an integer.

The multiplicative group of $\mathbf{Z}_r$ comes in three varieties depending on the prime factors of the integer $r$. The three structures are as follows.

**Structure of the Group $\mathbf{Z}_r^*$**
( 1 ) $\mathbf{Z}_r^* = \mathbf{Z}_2 \times \mathbf{Z}_{2^{v-2}} \cong \langle -1, 5 \rangle$.
( 2 ) $\mathbf{Z}_{p^v}^* \cong \langle g \rangle$, prime $p > 2$, $g$ a generator, and $v \geq 1$.
( 3 ) $\mathbf{Z}_r^* = \mathbf{Z}_{p_1^{v_1}} \times \mathbf{Z}_{p_{2a}^{v_2}} \times \cdots \times \mathbf{Z}_{p_a^{v_a}}$, where $r = p_1^{v_1} p_2^{v_2} \cdots p_a^{v_a}$.

A union $S = S_0 \cup S_1 \cup \cdots \cup S_{n-1}$ of disjoint subsets $S_0, S_1, ..., S_{n-1}$ of a set $S$ is called a partition of $S$. The disjoint condition stipulates that $S_i \cap S_j = \varnothing$ for all $i \neq j$.

The number of partitions of a set $S$ of cardinality $m = \#S$ into $n$ nonempty subsets of cardinalities $k_0, k_1, ..., k_{n-1}$ is given by the (Sterling) number

$$S(m,n) = \frac{1}{n!}\sum_{d=0}^{n}(-1)^{n-d}\binom{n}{d}d^m.$$

However, the recursive relation

$$S(m,n) = S(m-1,n-1) + nS(m-1,n),$$

where $S(m,1) = S(m,m) = 1$, $m \geq 2$, is more efficient in numerical calculations.

The total number of partitions of uniform cardinalities $k_0 = k_1 = \cdots = k_{n-1} = k$ is given by

$$\binom{kn}{k}\binom{kn-k}{k}\binom{kn-2k}{k}\cdots\binom{2k}{k}\binom{k}{k}$$

where $kn = \#S$.

## 4.2 Cyclotomic Periods

A cyclotomic period is simply an incomplete exponential sum $\Sigma_{x\in K}\psi(x)$ with respect to a fixed character or exponential function ψ on a coset K of some group G. The best known case, called a Gaussian period, is constructed with the complex-valued function $\psi(x) = e^{i2\pi x/p}$ or the finite field-valued function $\psi(x) = \omega^x$, where $\omega$ is a primitive pth root of unity in a finite extension $\mathbf{F}_{q^n}$ of $\mathbf{F}_q$, and the coset is a residue class $K = K_i$ of $\mathbf{F}_p$.

The periods based on the complex-valued function $\psi(x) = e^{i2\pi x/p}$ have been investigated for quite sometimes, perhaps, over two centuries, confer [1, Gauss, Art. 343-366].

**Gaussian Periods**

An important class of periods is derived from the algebraic structure of the residue numbers systems $\mathbf{Z}_r$. The partitions used in the construction of the gaussian periods consists of the cosets of $\mathbf{Z}_r$.

***Definition 4.5.*** Let $p = kn + 1 \in N$ be prime. An element $0 \neq a \in \mathbf{F}_p^*$ is an *nth power residue* if $a^{(p-1)/n} = 1$ in $\mathbf{F}_p$. In other words, the pure equation $x^n - a = 0$ has at least one solution $0 \neq x$ in $\mathbf{F}_p$. Otherwise $a$ is an *nth power nonresidue*.

***Definition 4.6.*** Let $p = kn + 1$ be prime, and let g be a generator of the multiplicative group of $\mathbf{F}_p$. Then the cosets decomposition of $\mathbf{F}_p^*$ with respect to g is defined by

$$\mathbf{F}_r^* = \prod_{i=0}^{k-1} g^i \{ g^{jn} : 0 \leq j < k \},$$

where the integers $g^{jn+i}$ in each coset $K_i = \{ g^{jn+i} : 0 \leq j < k \} = g^i K_0$ are reduced modulo p.

The first coset $K_0 = \{ g^{jn} : 0 \leq j < k \} = (\mathbf{F}_p^*)^n$ consists of all the nth power residues in $\mathbf{F}_p^*$, and the other cosets consist of subsets of nth power nonresidues. If p is not a prime, then the subgroup $K_0$ of nth power residues is not unique, and the set of periods of degree n is not unique.

If the generator g is not readily available, the cosets can be generated recursively:

$K_0 = \{ x^n : 0 \neq x \in \mathbf{F}_p \}$, $K_1 = x_1 K_0$ for some $x_1 \notin K_0$, $K_2 = x_2 K_0$ for some $x_2 \notin K_0 \cup K_1$, ...

The cosets form a disjoint partition of $\mathbf{F}_p^*$. Nondisjoint unions of $\mathbf{F}_p^*$ are also useful in the construction of more general periods.

***Definition 4.7.*** Let $p = kn + 1$ be prime, and let g be a primitive root in $\mathbf{F}_p$. The gaussian periods of degree n are the elements

$$\eta_0 = \sum_{j=0}^{k-1} e^{i2\pi g^{jn}/p}, \quad \eta_1 = \sum_{j=0}^{k-1} e^{i2\pi g^{jn+1}/p}, \quad \ldots, \quad \eta_{n-1} = \sum_{j=0}^{k-1} e^{i2\pi g^{jn+n-1}/p} .$$

It is convenient to refer to these elements as periods of type (k, n). A type (k, n) period is a sum of k elements, and the minimal polynomial $\psi_p(x) \in \mathbf{F}_q[x]$ of the period $\eta_0$ is a polynomial of degree n. The polynomial $\psi_p(x)$ are investigated in Chapter 5.

***Lemma 4.8.*** If $p \equiv 1 \bmod 2n$, then the followings are equivalent.
( 1 ) The pure equation $x^n - 1 = 0$ has a nonzero solution in $\mathbf{F}_p$.
( 2 ) The periods $\eta_0, \eta_1, ..., \eta_{n-1}$ are real numbers.

Proof: (1) This is equivalent to $-1 \in K_0 = \{ g^{jn} : 0 \le j < k \} = (\mathbf{F}_p^*)^n$. In addition the hypothesis $p \equiv 1 \bmod 2n \Rightarrow -1 = g^{kn/2} \in K_0 = \{ g^{jn} : 0 \le j < k \}$, since k is even.
(2) Since $-1 = g^{kn/2}$, the nth power residue occur in pairs $(g^n, g^{n+kn/2} = -g^n)$, $(g^{2n}, g^{2n+kn/2} = -g^{2n})$, ..., $(g^{(k-1)n}, g^{(k-1)n+kn/2} = -g^{(k-1)n})$, all together there are k/2 pairs. Thus the exponential sums $\eta_i$ are actually sums of cosine functions. ■

A different proof from a different point of view is given in [1, Berndt et al., p. 176.].

## 4.3 Cyclotomic Numbers

The theory of cyclotomic numbers is tightly linked to the arithmetic of cyclotomic fields. A limited amount of background materials on cyclotomic fields will be provided here. For finer analysis consult the literatures.

***Theorem 4.9.*** Let p be a prime integer, g a primitive root modulo p, and let $\omega = e^{i2\pi/p} \in \mathbf{C}$. Then the followings hold.

( 1 ) The cyclotomic field $\mathbf{Q}(\omega)$ is an extension of $\mathbf{Q}$ of degree $k = [\mathbf{Q}(\omega): \mathbf{Q}]$, with galois group $\mathrm{Gal}(\mathbf{Q}(\omega),\mathbf{Q}) = \{ 1, \tau, \tau^2, ..., \tau^{kn-1} \}$, where the generating automorphism is defined by $\tau(\omega) = \omega^g$.

( 2 ) $\mathbf{Q}(\eta_0)$ is an extension of $\mathbf{Q}$ of degree $n = [\mathbf{Q}(\eta_0) : \mathbf{Q}]$, with galois group $\mathrm{Gal}(\mathbf{Q}(\eta_0),\mathbf{Q}) = \{1, \sigma, \sigma^2, ..., \sigma^{n-1} \}$, where $\sigma = \tau^k$. This is a unique subfield of the field $\mathbf{Q}(\omega)$.

( 3 ) $\mathbf{Q}(\omega)$ is an extension of $\mathbf{Q}(\eta_0)$ of degree $k = [\mathbf{Q}(\omega) : \mathbf{Q}(\eta_0)]$,with galois group $\mathrm{Gal}(\mathbf{Q}(\omega),\mathbf{Q}(\eta_0)) = \{1, \tau^n, \tau^{2n}, ..., \tau^{(k-1)n} \}$.

The subfield $\mathbf{Q}(\eta_0)$ of $\mathbf{Q}(\omega)$ is fixed by the subgroup $\mathrm{Gal}(\mathbf{Q}(\omega),\mathbf{Q}(\eta_0)) = \{1, \tau^n, \tau^{2n}, ..., \tau^{(k-1)n} \}$ of $\mathrm{Gal}(\mathbf{Q}(\omega),\mathbf{Q}) = \{ 1, \tau, \tau^2, ..., \tau^{kn-1} \}$, and the coset $K_i$ is identified with the subgroup $H = \{ 1, \sigma, \sigma^2, ..., \sigma^{n-1} \}$ of $\mathrm{Gal}(\mathbf{Q}(\omega),\mathbf{Q})$, see [1, Washington, p. 16].

***Lemma 4.10.*** ( 1 ) The set $\{ 1, \omega, \omega^2, ..., \omega^{kn-2} \}$ is a $\mathbf{Z}$-basis of the ring of algebraic integers $\mathbf{Z}[\omega]$ of the field $\mathbf{Q}(\omega)$.
( 2 ) The periods $\{ \eta_0, \eta_1, ..., \eta_{n-1} \}$ forms a $\mathbf{Z}$-basis of the ring of algebraic integers $\mathbf{Z}[\eta_0,\eta_1,...,\eta_{n-1}]$ of the field $\mathbf{Q}(\eta_0)$.
( 3 ) The set $\{ 1, \omega, \omega^2, ..., \omega^{k-1} \}$ is a $\mathbf{Z}[\eta_0,\eta_1,...,\eta_{n-1}]$-basis of the ring of algebraic integers $\mathbf{Z}[\omega]$ of the field $\mathbf{Q}(\omega)$.

Quite often the subring of algebraic integers $\mathbf{Z}[\eta_0,\eta_1,...,\eta_{n-1}] = \mathbf{Z}[\eta_0] = \mathbf{Z}[\eta_1] = \cdots = \mathbf{Z}[\eta_{n-1}]$, for example, if $2 = k = (p - 1)/n$, then the ring of integers is $\mathbf{Z}[\eta_0,\eta_1,...,\eta_{n-1}] = \mathbf{Z}[\omega+\omega^{-1}]$, and the subfield is $\mathbf{Q}(\omega+\omega^{-1})$, these are the maximal subring and subfield of $\mathbf{Z}[\omega]$ and $\mathbf{Q}(\omega)$ respectively. But in general for $k > 2$, this is not true. This fails to hold because there are primes for which the linear expansions

$$\eta_i = \sum_{j=0}^{n-1} a_{i,j}\, \eta_0^j$$

require some rational coefficients $a_{i,j} \notin \mathbf{Z}$, see [1, Washington, p. 17.].

In the subring of cyclotomic integers $\mathbf{Z}[\eta_0,\eta_1,...,\eta_{n-1}]$ the product of two gaussian periods is expressible as a linear combination of the same periods with integers coefficients.

***Definition 4.11.*** The multiplication table or matrix $T = ( t_{i,j,k} )$ attached to the cyclotomic periods is defined by

$$\eta_i\eta_j = \sum_{k=0}^{n-1} t_{i,j,k}\eta_k \,.$$

The most important linear expansions of the pairwise products are the followings:

$$\eta_0\eta_0 = \left(\sum_{a=0}^{k-1} e^{i2\pi g^{an}/p}\right)\left(\sum_{b=0}^{k-1} e^{i2\pi g^{bn}/p}\right) = \varepsilon_i + \sum_{j=0}^{n-1}(0, j)\eta_j \,,$$

$$\eta_0\eta_1 = \left(\sum_{a=0}^{k-1} e^{i2\pi g^{an}/p}\right)\left(\sum_{b=0}^{k-1} e^{i2\pi g^{bn+1}/p}\right) = \varepsilon_i + \sum_{j=0}^{n-1}(1, j)\eta_j \,,$$

… … …

$$\eta_0\eta_{n-1} = \left(\sum_{a=0}^{k-1} e^{i2\pi g^{an}/p}\right)\left(\sum_{b=0}^{k-1} e^{i2\pi g^{bn+n-1}/p}\right) = \varepsilon_i + \sum_{j=0}^{n-1}(n-1, j)\eta_j \,,$$

where the constant term $\varepsilon_i = k$ if $-1 \in K_i$, otherwise $\varepsilon_i = 0$. Sometimes this is characterized in terms of the impulse function as

$$\varepsilon_i = \begin{cases} \delta_{0,i} & k \text{ even,} \\ \delta_{n/2,i} & k \text{ odd.} \end{cases}$$

The other linear expansions of the pairwise products $\eta_i\eta_j$ are derived by repeated application of the automorphism. In particular, from the identity $\tau^v(\eta_i\eta_j) = \eta_{i+v}\eta_{j+v}$, the pairwise product become

$$\eta_i\eta_j = \tau^i(\eta_0\eta_{j-i}) = \sum_{k=0}^{n-1}\tau^i(t_{0,j-i,k}\eta_k)$$
$$= \sum_{k=0}^{n-1}t_{0,j-i,k-i}\eta_k = \sum_{k=0}^{n-1}t_{j-i,k-i}\eta_k,$$

where the coefficients $t_{a,b} = t_{0,a,b}$.

***Definition 4.12.*** The complexity of multiplication in the subfield $\mathbf{Q}(\eta_0)$ over $\mathbf{Q}$ is defined by $w(T_0) = \#\{ t_{0,j,k} \neq 0 \}$.

The coefficients (i, j) appearing in the linear expansion

$\eta_0\eta_i = (i, 0)\eta_0 + (i, 1)\eta_1 + \cdots + (i, n-1)\eta_{n-1} + \varepsilon_i$

are classically known as *cyclotomic numbers*.

In the literature there are a few equivalent working definitions of these numbers. One of these definition is given below.

***Definition 4.13.*** Let p = kn + 1 be prime, and let g be a primitive root modulo p, the cyclotomic numbers (i, j) counts the number of solutions of the congruence equation

$1 + g^{nx + i} \equiv g^{ny + j} \bmod p$,

where $0 \leq i, j < n$, and $0 \leq x, y < k$. These numbers depend on both the prime p and the primitive root g.

For small prime p, an easy and intuitive method of computing (i, j) is to count the number of times that the coset $K_j = \{ g^{nv + j} \bmod p : 0 \leq v < k \}$ appears in the sequence $< z_d = x + g^i y : x, y \in K_0, d = 0, 1, 2, ..., k^2 - 1 >$.

***Example 4.14.*** Let p = kn + 1 = 11, (n = 2 and k = 5), and let the cosets of quadratic residues and non residues in $(\mathbf{F}_{11})^*$ be $K_0$ = { 1, 4, 9, 5, 3 } and $K_1$ = { 2, 8, 7, 10, 6 }. Then using g = 2, the cyclotomic numbers are

$(0, 0) = \#\{ (x, y) : 1 + 2^{2x} \equiv 2^{2y} \bmod 11, 0 \leq x, y < 5 \} = 2$,

$(0, 1) = \#\{ (x, y) : 1 + 2^{2x+1} \equiv 2^{2y} \bmod 11, 0 \leq x, y < 5 \} = 3$,

$(1, 0) = \#\{ (x, y) : 1 + 2^{2x} \equiv 2^{2y+1} \bmod 11, 0 \leq x, y < 5 \} = 2$, and

$(1, 1) = \#\{ (x, y) : 1 + 2^{2x+1} \equiv 2^{2y+1} \bmod 11, 0 \leq x, y < 5 \} = 2$.

In general an arbitrary partition of a set induces a set of periods. However, the complete set of the numbers (i, j) that arises in the linear combinations

$\eta_0\eta_i = (i, 0)\eta_0 + (i, 1)\eta_1 + \cdots + (i, n-1)\eta_{n-1} + \varepsilon_i$

exist (are defined) only if the set **Z**[$\eta_0,\eta_1,...,\eta_{n-1}$] is equipped with a ring structure over the ring of integers **Z**.

## 4.4 Linear and Algebraic Properties of the Integers (i, j)

The linear properties of the cyclotomic numbers periods are basic tools used in the investigations of the periods.

**Entries Relations**
( 1 ) (i, j) = (i + sn, j + sn), the shift-invariant property for s $\in$ N.

( 2 ) (i, j) = (n − i, n − j), circulant property.

( 3 ) $(i, j) = \begin{cases} (j, i) & \text{if } p = 2 \text{ or k is even,} \\ (j + n/2, i + n/2) & \text{if k and p are both odd.} \end{cases}$

**Sums Relations**
( 4 ) $\sum_{i=0}^{n-1} (i, j) = k - a_j$ where $a_j = \begin{cases} 1 & \text{if } j = 0, \\ 0 & \text{otherwise.} \end{cases}$

( 5 ) $\sum_{j=0}^{n-1} (i, j) = k - b_i$ where $b_i = \begin{cases} 1 & \text{k even and } i = 0, \\ 1 & \text{k odd and } i = n/2, \\ 0 & \text{otherwise.} \end{cases}$

(6) If g is a primitive root modulo p, and the pair (i, j) = $(i, j)_1$ and $(i, j)_s$ are the cyclotomic numbers with respect to the primitive roots respectively g and $g^s$, with s = gcd(s, p − 1) = 1, then $(i, j)_s$ = (si, sj).

Property (3) specifies the symmetries of the integers (i, j) along the diagonals. The last two properties (4) and (5) are the row sums of the (i, j) with respect to a fixed index i and the column sums of the (i, j) with respect to a fixed index j.

**Cyclotomic Matrices**
A cyclotomic matrix of type (k, n) is a n×n array $C_k$(,n) = ( (i, j) ) with entries in **Z**. The linear properties are used in the construction of these matrices. Properties (1), (2), and (3) specify the symmetries along the diagonals of the matrix. The structure of this matrix depends on parity of the parameter k in the prime p = kn + 1. Specifically, property (3) indicates that the matrix is symmetric if k is even and nonsymmetric if k is odd.

As an illustration, let consider the simplest ones of type (k, 2) and type (k, 3). The type (k, 2) matrix for the prime kn + 1 = 2k + 1 comes in two varieties: The matrix is either

$$C_k(2) = \begin{bmatrix} (0,0) & (0,1) \\ (1,0) & (1,1) \end{bmatrix} = \begin{bmatrix} (0,0) & (0,1) \\ (1,1) & (1,1) \end{bmatrix} \text{ or } C_k(2) = \begin{bmatrix} (0,0) & (0,1) \\ (1,0) & (1,1) \end{bmatrix} = \begin{bmatrix} (0,0) & (0,1) \\ (0,0) & (0,0) \end{bmatrix}.$$

for k even or k odd respectively. And the type (k, 3) matrix for p = kn + 1 = 3k + 1 has the form

$$C_k(3) = \begin{bmatrix} (0,0) & (0,1) & (0,2) \\ (1,0) & (1,1) & (1,2) \\ (2,0) & (2,1) & (2,2) \end{bmatrix} = \begin{bmatrix} (0,0) & (0,1) & (0,2) \\ (0,1) & (0,2) & (1,2) \\ (0,2) & (1,2) & (0,1) \end{bmatrix}.$$

## 4.5 Characterization of the Cyclotomic Periods

Several properties of the gaussian periods uniquely identify the subset $\{ \eta_0,\eta_1,\ldots,\eta_{n-1} \}$ of periods, and give a complete characterization.

***Theorem 4.15.*** Let $\theta_0, \theta_1, \ldots, \theta_{n-1} \in \mathbf{Q}(\omega)$, and suppose that
( 1 ) $\theta_{i+nj} = \theta_i$, $j \in \mathbf{N}$.
( 2 ) The subset $\{ \theta_0, \theta_1, \ldots, \theta_{n-1} \}$ is linearly independent over $\mathbf{Q}$.
( 3 ) $\theta_0 + \theta_1 + \cdots + \theta_{n-1} = -1$.
( 4 ) $\sum_{j=0}^{n-1} \theta_j \theta_{j+i} = p\varepsilon_i - n$.
( 5 ) $p^{-1}\left( k^2 + \sum_{d=0}^{n-1} \eta_d \eta_{d+i} \eta_{d+j} \right) \in \mathbf{Z}$.

Then $\theta_0, \theta_1, \ldots, \theta_{n-1}$ are the gaussian periods of degree n. Conversely if $\theta_0, \theta_1, \ldots, \theta_{n-1}$ are the gaussian periods of degree n, then properties (1) to (5) hold.

Proof: Confer [1, Thaine, p.38]. ■

**Inversion Formula**
The individual numbers (i, j) can be recovered via the inversion formula.

***Theorem 4.16.*** Let p kn + 1 be a prime integer, and let $\eta_0, \eta_1, \ldots, \eta_{n-1}$ be the periods of degree n. Then

$$(i, j) = \begin{cases} p^{-1}\left(k^2 + \sum_{d=0}^{n-1} \eta_d \eta_{d+i} \eta_{d+j}\right) & p = 2 \text{ or } p > 2 \text{ and } k \text{ even,} \\ p^{-1}\left(k^2 + \sum_{d=0}^{n-1} \eta_d \eta_{d+i} \eta_{d+i+n/2}\right) & p > 2 \text{ and } k \text{ odd.} \end{cases}$$

Proof: In the sum of triple products $\eta_d\eta_{d+i}\eta_{d+j}$ replace the linear expansion for $\eta_d\eta_{d+i}$, and rescale the index. Then

$$\begin{aligned} \sum_{d=0}^{n-1} \eta_d \eta_{d+i} \eta_{d+j} &= \sum_{u=0}^{n-1} t_{i,u} \sum_{d=0}^{n-1} \eta_{d+u} \eta_{d+j} \\ &= \sum_{u=0}^{n-1} t_{i,u} \sum_{d=0}^{n-1} \eta_d \eta_{d+j-u}. \end{aligned}$$

Again substitute the linear expansion for $\eta_d\eta_{d+j-u}$, and sum over d, to obtain

$$\sum_{d=0}^{n-1} \eta_d \eta_{d+i} \eta_{d+j} = \sum_{u=0}^{n-1} (p\varepsilon_{j-u} - k) t_{i,u} .$$

Now apply the sum identity and simplify. ■

Further materials, identities and more, are available in [1, Storer], [1, Thane, p. 37.] or similar sources.

# 4.6 Cyclotomic Numbers of Short Type

The cyclotomic numbers corresponding to small parameter n will be refer to as cyclotomic numbers of short type (k, n). And those with a large parameter n as cyclotomic numbers of long type (k, n). There is interest in determining closed form formulae for the numbers (i, j) as function of the primes p = kn + 1.

There are several techniques used in the calculations of the cyclotomic numbers. Two of these are the followings.

( 1 ) *Quadratic Partition Method.* This technique uses the solutions of the system of equations involving the quadratic partition of the prime p and other related parameters. The cyclotomic numbers are expressed in terms of these parameters. This is an effective method for small n only, since as n increases, the system of equations describing the quadratic partition of the prime p becomes very complex.

( 2 ) *Exponential Sums Method.* This is also an effective method but for small n only since the evaluations of the exponential sums involved for large n are not known.

**Case n = 2**
Properties (4) and (2), yield the system of equations

$$\begin{array}{lcl} (0,0)+(0,1)=k-1 & \text{or} & (0,0)+(0,1)=k \\ \quad 2(0,1)=k & & \quad 2(0,1)=k-1, \end{array}$$

depending on parity of k even or odd respectively. The solutions of this system yield the next results.

***Lemma 4.17.*** If 2k + 1 is a prime, then the cyclotomic numbers of type (k, 2) are given by

( 1 ) (0, 1) = (1, 0) = (1, 1) = k/2, and (0, 0) = (k – 2)/2, if k is even, or

( 2 ) (0, 0) = (1, 0) = (1, 1) = (k – 1)/2, and (0, 1) = (k + 1)/2, if k is odd.

Thus the type (k, 2) matrix is given by either

$$C_k(2) = \begin{bmatrix} (k-2)/2 & k/2 \\ k/2 & k/2 \end{bmatrix} \text{ or } C_k(2) = \begin{bmatrix} (k-1)/2 & (k+1)/2 \\ (k-1)/2 & (k-1)/2 \end{bmatrix},$$

where the integer k is even or odd respectively.

The numbers (i, j) of type (k, 2) are also related to the numbers, counting multiplicities, of solutions of the equations x + y, x + ay, and –(x + ay), as the variables x and y range over $K_0$, and a fixed nonquadratic residue *a* in $\mathbf{F}_p$. And to the coefficients, up to k depending on the parity of k, of the linear expansions

$\eta_0\eta_0 = a\eta_0 + b\eta_1$ and $\eta_0\eta_1 = c\eta_0 + d\eta_1$,

where a, b, c, d $\in$ **Z** respectively. The coefficients are the followings:

( 1 ) a = (0,0) – k, b = (0,1) – k, c = (0,1), and d = (1,1), if k is even, or
( 2 ) a = (0,0), b = (0,1), c = (0,1) – k, and d = (1,1) – k, if k is odd.

As an instance, for p = 11, (k = 5), the coefficients are a = (0,0) = 2, b = (0,1) = 3, c = (0,1) − k = -3, and d = (1,1) − k = −3, so

$\eta_0\eta_0 = 2\eta_0 + 3\eta_1$, and $\eta_0\eta_1 = -3(\eta_0 + \eta_1)$.

## Case n = 3

Let p = kn + 1 = 3k + 1 be a prime and g be a primitive root modulo p. The prime p splits in the quadratic field $\mathbf{Q}(e^{i2\pi/3})$ in four different ways. But the quadratic partition specified by the relations below gives a unique one.

( 1 ) $4p = A^2 + 3B^2$,
( 2 ) $A \equiv 1 \bmod 3$, and $B \equiv 0 \bmod 3$, and
( 3 ) $3B \equiv (g^{(p-1)/3} - g^{2(p-1)/3})A \bmod p$.

Line (2) above specifies a solution (A, ±B) of equation (1). And condition (3) determine a unique B or −B. This in turns determine a unique cyclotomic matrix of type (k, 3), one of the two possible matrices determines by ±B. However, if the prime is large condition (3) is cumbersome to obtain or very difficult to compute The actual 3×3 symmetric matrix is

$$C_k(3) = \begin{bmatrix} (0,0) & (0,1) & (0,2) \\ (1,0) & (1,1) & (1,2) \\ (2,0) & (2,1) & (2,2) \end{bmatrix} = \begin{bmatrix} (0,0) & (0,1) & (0,2) \\ (0,1) & (0,2) & (1,2) \\ (0,2) & (1,2) & (0,1) \end{bmatrix}.$$

***Lemma 4.18.*** Let p = kn + 1 = 3k + 1 be a prime such that $4p = A^2 + 27B^2$, and $A \equiv 1 \bmod 3$. Then the cyclotomic numbers of type (k, 3) are given by

( 1 ) (0, 0) = (p − 8 + A)/9,
( 2 ) (0, 1) = (1, 0) = (2, 2) = (2p − 4 − A + 9B)/18,
( 3 ) (0, 2) = (2, 0) = (1, 1) = (2p − 4 − A − 9B)/18, and
( 4 ) (1, 2) = (2, 1) = (p + 1 + A)/9.

Proof: Only four of the numbers (i, j), $0 \le i, j < 3$, are independent. Linear properties (2) and (4) are used to establish the equalities (0, 1) = (1, 0), (0, 1) = (2, 2), (0, 2) = (2, 0), and (1, 2) = (2, 1) in the matrix above. This leads to a system of equations with only three equations (rows/columns sums):

$$\begin{aligned} (0,0) + (0,1) + (0,2) &= k - 1, \\ (0,1) + (0,2) + (1,2) &= k, \\ (0,2) + (1,2) + (0,1) &= k. \end{aligned}$$

The other equations needed to complete the calculations, and other details of the proof appears in [1, Dickson, p. 397], also [1, Evans, et al, p. 71. ]. Some of these numbers will be used later on.

## Case n = 4

Let p = kn + 1 = 4k + 1 be a prime and g be a primitive root modulo p. The prime p splits in the quadratic field $\mathbf{Q}(e^{i\pi/2})$ in four different ways. But the quadratic partition specified by the relations below gives a unique one.

( 1 ) $p = A^2 + B^2$,
( 2 ) $A \equiv -(2 \mid p) \bmod 4$, and
( 3 ) $B \equiv Ag^{(p-1)/4} \bmod p$.

Line (2) above specifies a solution (A, ±B) of equation (1). And condition (3) determine a unique B or −B. The parity of k determines one of two cyclotomic matrices of type (k, 4) possible. Using properties (2) and (4) to remove any dependency, the matrices reduce to

$$C_k(4) = \begin{bmatrix} (0,0) & (0,1) & (0,2) & (0,3) \\ (0,1) & (0,3) & (1,2) & (1,2) \\ (0,2) & (1,2) & (0,2) & (1,2) \\ (0,3) & (1,2) & (1,2) & (0,1) \end{bmatrix} \quad \text{or} \quad C_k(4) = \begin{bmatrix} (0,0) & (0,1) & (0,2) & (0,3) \\ (1,0) & (1,0) & (?1,2) & (0,1) \\ (0,0) & (1,0) & (0,0) & (1,0) \\ (1,0) & (0,3) & (0,1) & (1,0) \end{bmatrix}$$

for k even or odd respectively.

***Lemma 4.19.*** The cyclotomic numbers of type (k, 4) are given by either

( 1 ) 16(0, 0) = p – 11 + 6A,
( 2 ) (0, 1) = (1, 0) = (3, 3) = (p – 3 – 2A + 4B)/16,
( 3 ) (0, 2) = (2, 0) = (2, 2) = (p – 3 – 2A)/16,
( 4 ) (0, 3) = (3, 0) = (1, 1) = (p – 3 – 2A – 4B)/16,
( 5 ) (1, 2) = (2, 1) = (1, 3) = (2, 3) = (3, 1) = (3, 2) = (p + 1 + 2A )/16,

if k is even, or

( 6 ) 16(2, 2) = p – 7 + 2A,
( 7 ) (0, 1) = (1, 3) = (3, 2) = (p + 1 + 2A + 4B)/16,
( 8 ) (0, 2) = (0, 0) = (2, 0) = (p + 1 – 6A)/16,
( 9 ) (0, 3) = (1, 2) = (3, 1) = (p + 1 + 2A – 4B)/16,
( 10 ) (1, 0) = (1, 1) = (2, 1) = (2, 3) = (3, 0) = (3, 3) = (p – 3 – 2A)/16,

if k is odd.

Proof: Only five of the numbers (i, j), $0 \le i, j < 4$, are independent. Linear properties (2) and (4) are used to establish the equalities or dependencies among the numbers (i, j) in the matrix above. This leads to a system of equations with only three equations (rows/columns sums):

$$(0, 0) + (0, 1) + (0, 2) + (0, 3) = k - 1,$$
$$(0, 1) + (0, 3) + 2(1, 2) = k,$$

$$2(0, 2) + 2(2, 1) = k.$$

The other equations used to complete the calculations, and other details of the proof appears in [1, Dickson, p. 397], also [1, Evans, et al, p. 74. ].

## Case n = (p – 1)/2

Let p = 2n + 1 be a prime and g be a primitive root modulo p. The quadratic partition of p is not required, and has no part in the calculations of the cyclotomic numbers. This case turns out be a lot easier.

***Lemma 4.20.*** Let $K_i = \{ g^i, g^{n+i} \}$ be the nth power residues partition of $\mathbf{F}_p^*$, and $\eta_i = \psi(g^i) + \psi(g^{n+i})$, where $0 \le i < n$, and $\psi(x) = e^{i2\pi x/p}$. Then the cyclotomic numbers of type (2, n) are the followings.

( 1 ) (0, j) = 1, if $g^j \equiv 2 \bmod p$, otherwise (0, j) = 0.
( 2 ) (i, j) = 1, if $(1 + g^i) \equiv j \bmod p$ or $(1 - g^i) \equiv j \bmod p$, otherwise (i, j) = 0.

Proof: Since $-1 = g^n$, the element $\eta_i = \psi(g^i) + \psi(-g^i)$, for all $i \ne 0$, and it follows that

$$\begin{aligned}\eta_0\eta_i &= (\psi(1) + \psi(-1))(\psi(g^i) + \psi(-g^i)) = \psi(1+g^i) + \psi(-1-g^i) + \psi(1-g^i) + \psi(-1+g^i)\\ &= \eta_a + \eta_b,\\ &= (i, 0)\eta_0 + (i, 1)\eta_1 + \cdots + (i, n-1)\eta_{n-1}\\ &= (i, a)\eta_a + (i, b)\eta_b\end{aligned}$$

where $g^a = 1 + g^i \bmod p$, and $g^b = 1 - g^i \bmod p$. And last for i = 0, $\eta_0\eta_0 = \eta_a + 2 = (0, a)\eta_a + 2$, where $g^a = 2 \bmod p$. Note that $\varepsilon_0 = 2$, and $\varepsilon_i = 0$ for $i \ne 0$. ■

The cyclotomic numbers of type type (3, n) for the periods $\eta_i = \psi(g^i) + \psi(g^{n+i}) + \psi(g^{2n+i})$ derived from the partitions $K_i = \{ g^i, g^{n+i} , g^{2n+i} \}$, p = 3n + 1, exhibit similar pattern. More generally, for p = kn + 1, the periods of type (k, n) are derived from the cosets $K_i = \{ g^i, g^{n+i} , ..., g^{(k-1)n+i} \}$.

**Remark:** It is clear that the complexity of multiplication in $\mathbf{Q}(\psi(g)+\psi(g^n))$ is $w(T_0) = 3n - 2$. If the periods $\eta_i = \psi(g^i) + \psi(g^{n+i})$ of type (2, n) form a normal basis in characteristic 2, then the multiplication matrix has the minimal complexity possible: 2n – 1. And the periods of type (4, n) has the minimal complexity possible: 4n – 3, etc. Similarly in characteristic 3, the multiplication matrix with respect to a period normal basis of type (3, n) has the minimal complexity possible: 3n – 2. A type (6, n) has the minimal complexity possible: 6n – 5, etc.

## 4.7 Extension to the Ring $\mathbf{Z}_r$

This section considers a generalization of the gaussian periods to the residue numbers system $\mathbf{Z}_r$, $r = p^v$, $p$ prime and $v \geq 1$. This generalization is precisely the gaussian period if $v = 1$. Only the quadratic case will be given here. Put $K_0 = < g^2 >$, and $K_1 = gK_0$, where $g$ is primitive in $\mathbf{Z}_r$. The two $K_0$ and $K_1$ together with the set $M = \{ 0, p, 2p, ..., (p^{v-1} - 1)p \}$, which is a maximal ideal, forms a partition of $\mathbf{Z}_r$.

***Lemma 4.21.*** The quadratic cyclotomic numbers over the ring $\mathbf{Z}_r$ are as follows.

(1) If $p = 4c + 1$, then $(0, 1)_v = (1, 0)_v = (1, 1)_v = p^{v-1}(p - 1)/4$, and $(0, 0)_v = p^{v-1}(p - 5)/4$.

(2) If $p = 4c + 3$, then $(0, 0)_v = (1, 0)_v = (1, 1)_v = p^{v-1}(p - 3)/4$, and $(0, 1)_v = p^{v-1}(p + 1)/4$.

Other advanced details and applications of this generalization are covered in [1, Ding and Helleseth, 1999].

## 4.8 Extension of the Finite Field $\mathbf{F}_p$

Let $r = kn + 1 = p^v$ be a prime power. The $n$th power residues partitions or cosets of a $v$th degree extension $\mathbf{E}_v$ of the finite field $\mathbf{F}_p$ is utilized to augment those of $\mathbf{F}_p$ and produce other varieties of periods. This is one of the many possible generalization of the cyclotomic periods. See [3, Gurak], [ , Rawarte]. Let $K_0 \cup K_1 \cup \cdots \cup K_{n-1}$ be a partition of $\mathbf{F}^*_{p^v}$, and $Tr : \mathbf{F}^*_{p^v} \to \mathbf{F}_p$ be the absolute trace. The cyclotomic periods are defined by

$$\eta_i = \sum_{x \in K_i}^{n-1} \omega^{Tr(g^i x^k)},$$

where $0 \leq i < n$.

**Uniform Cyclotomic Numbers**

***Definition 1.*** The cyclotomic numbers are called *uniform* if the following holds.
( 1 ) $(0, 1) = (0, i) = (i, 0) = (i, i)$ for $i \neq 0$.
( 2 ) $(1, 2) = (i, j)$ for $i \neq j$, and $i, j \neq 0$.

***Lemma 2.*** If the cyclotomic numbers are uniform, then either $p = 2$ or $k$ is even.

Proof: On the contrary if both $p$ and $k$ are odd, then the integer $n$ is even, and by definition of uniformity it becomes

(0, 0) = ($n$/2, $n$/2) = (0, 1).

Further, for $n > 2$, by uniformity, and property (3),

(1, 2) = ($n$/2 – 1, $n$/2) = (0, $n$ – 1) = (0, 1).

This implies that all the cyclotomic numbers in any row or any column are the same, for instance,

$(0, j) + (1, j) + (2, j) + \cdots + (n-1, j) = (0, 0)n$,

in contradiction of the basic property of the cyclotomic numbers ($i$, $j$) given by

$(0, j) + (1, j) + (2, j) + \cdots + (n-1, j) = k - \delta$, ???

where $\delta = 0, 1$. ■

***Lemma 3.*** Let $r = kn + 1 = p^{2v}$ be a prime power and let $n$ □ 3 be a divisor of $p^v + 1$. Then
( 1 ) The cyclotomic numbers are uniform.
( 2 ) The periods of type ($k$, $n$) are linearly dependent.

Proof of (2): If $n$ divides $p^v + 1$, the cyclotomic numbers are uniform since $k$ is even, which imply that the periods are uniform. Consequently there are linear dependencies among the periods, see [1, Baumert *et al*] for more details.

The corresponding $n \times n$ cyclotomic matrices $C_k(n)$ = ( ($i$, $j$) ) has at most 3 independent entries (0, 0), (0, 1), and (1, 2).

## 4.9 Exponential Sums Properties of the Periods

Let $\eta_{n-1}, \ldots, \eta_1, \eta_0$ be the periods of degree $n$, and define the g*amma* and *beta* polynomials by $G(x) = \eta_{n-1}x^{n-1} + \cdots + \eta_1 x + \eta_0$, and $B(s,t) = G(s)G(t)/G(st)$. The restriction of the gamma beta polynomials to the unit disk are called the gamma and beta sums respectively, (or the *resolvent* of the periods), see [2, Lehmers 1967, 1968].

**Properties of the Beta and Gamma Sums**
The beta and gamma sums satisfy the following properties.

( 1 ) $G(\theta^s) = -1$ if $n$ divides $s$, and $G(1) + G(\theta) + G(\theta^2) + \cdots + G(\theta^{n-1}) = 0$.

( 2 ) $G(\theta^s)\, G(\theta^{-s}) = (-1)^{sk}p$ if $n$ does not divide $s$.

( 3 ) $B(s,t) = -\frac{G(\theta^s)G(\theta^t)}{G(\theta^{s+t})}$,

where and $\theta = e^{i2\pi/n}$ is an $n$th primitive root of unity.

( 4 ) $G(e^{i2\square s/n})$ and $\eta_t = \frac{1}{n}\sum_{s=0}^{n-1} G(e^{i2\pi s/n})e^{-i2\pi st/n}$ is a discrete Fourier transform pair.

Property (2) stems from the fact that $G(e^{i2\pi s/n}) = G(\chi_t)$ is a gaussian sum, where $\chi_t(x) = e^{i2\pi t log(x)/p} = \xi^{t log(x)}$ is a multiplicative character on $\mathbf{F}_p$ of order $ord(\chi) = d \mid n$. and $\xi$ be an nth root of unity, possibly $\xi = e^{i2\pi/n}$.

$$G(\chi_s) = \sum_{t=0}^{n-1} \eta_t \xi^{st} = \sum_{t=0}^{n-1} \eta_t e^{i2\pi st/n} \quad \text{and} \quad \eta_t = \frac{1}{n}\sum_{s=0}^{n-1} G(e^{i2\pi s/n})e^{-i2\pi st/n}$$

$$\eta_i = \sum_{x \in \mathbf{F}_{p^a}}^{n-1} \omega^{Tr(g^i x^k)}$$

a gauss sum and the resolvent sum.

(1) $\eta_0 + \eta_1 + \eta_2 + \cdots + \eta_{v-1} = -1$,

(2) $R_0 + R_1 + R_2 + \cdots + R_{n-1} = 0$,

(3) $R_i = 1 + n\eta_i$,

$$G(\chi_s) = \sum_{t=0}^{n-1} \eta_t \xi \omega^{st} \quad , \quad \text{and} \quad \eta_i = \frac{1}{n}\sum_{s=0}^{n-1} G(\chi_s)\xi^{-st} \quad ,$$

(4) If integer n is invertible in $\mathbf{F}_p$, then
are a discrete Fourier transform pair.

Let η0, η1, ..., ηn-1 be periods of type (k, n) and let ηd,0, ηd,1, ..., ηd,n-1 be the extended periods of type (d, k, n). Define the gamma polynomial

$H(x) = H_1(x) = \eta_{n-1}x^{n-1} + \eta_{n-2}x^{n-2} + \cdots + \eta_1 x + \eta_0$,
and let

$$H_d(x) \equiv \sum_{t=0}^{r^d-1} x^t e^{i2\pi Tr(\gamma^t)/p} \mod (x^{r-1} - 1),$$

where $\gamma$ is primitive in $\mathbf{F}_d$, $d = [\mathbf{F}_d : \mathbf{F}_r]$, and $r = kn + 1 = p^v$. Another round of modular reduction gives

$$H_d(x) \equiv \sum_{i=0}^{n-1} \eta_{d,i} x^i \mod (x^n - 1),$$

where

$$\eta_{d,i} = \sum_{j \equiv i \bmod n} e^{i2\pi Tr(\gamma^t)/p}$$

are the periods of type (d, k, n). These polynomials will provides a formula for computing the periods $\eta_{d,0}, \eta_{d,1}, ..., \eta_{d,n-1}$ of type (d, k, n) in terms of the periods $\eta_0, \eta_1, ..., \eta_{n-1}$ of type (k, n).

$f(x) = (-H_1(x))^d = -(\eta_{d,n-1}x^{n-1} + \eta_{d,n-2}x^{n-2} + \cdots + \eta_{d,1}x + \eta_{d,0})$,

for all $x \in \Omega = \{ 1 \neq \xi : \xi^n = 1 \}$. The function f(x) is a periodic function, and it has a Fourier representation with coefficients given by

$$\eta_{d,t} = \frac{-1}{n}\sum_{s=0}^{n-1}\left( -H(s)\right)^d \xi^{-st},$$

for t = 0, 1, 2, ..., n – 1.

$$G(\chi_s) = \sum_{x=0}^{r^d-2} e^{i2\pi x/(r-1)} e^{i2\pi Tr(\xi^x)/p} = (-1)^{d+1}\left(\sum_{x=0}^{r-2} e^{i2\pi x/(r-1)} e^{i2\pi \tau^x)/p}\right)^d,$$

(Hasse- Davenport Relation) Let $r = p^v = kn + 1$, $< \tau > = \mathbf{F}_r^*$, and $< \xi > = \mathbf{L}^*$, $r^d = \#\mathbf{L}$ .
where.

Other identities and related topics are treated in McEnliece and Rumsey [2].

**Hyperperiods**
The generalization of the periods to hyperperiods replaces the inner term with an exponential sum

$$\eta_s^{(k)} = \sum_{s \in K_t}\sum_{x=0}^{p-1} e^{-i2\pi(x+sx^k)/p}$$

see [Lehmers 1967, 1968].

# Chapter 5

# Period Polynomials

# 5.1 Definition of Period Polynomials

Let p = kn + 1 be prime, and let g be a primitive root modulo p. The periods $\eta_i$ and reduced periods $g_i$ of type (k, n) are the incomplete exponential sums

$$\eta_d = \sum_{x=0}^{p-1} e^{i2\pi g^{jn+d}/p}, \text{ and } g_d = \sum_{x=0}^{p-1} e^{i2\pi g^d x^n/p},$$

for d = 0, 1, …, n – 1.

In general, if n > 2 and $p \neq 1 \bmod n$, then the map $x \to x^n$ is one-to-one, and the periods reduce to the trivial period $\eta_0 = \eta_1 = \cdots = \eta_{n-1}$. Nontrivial periods of type (k,n) are constructible if and only if prime p = kn + 1.

***Definition 5.1.*** The period polynomial of type (k, n) is defined by

$$\begin{aligned}\Psi_p(x) &= (x - \eta_0)(x - \eta_1) \cdots (x - \eta_{n-1}) \\ &= x^n + c_1x^{n-1} + c_2x^{n-2} + \cdots + c_{n-1}x + c_n.\end{aligned}$$

The associated reduced period polynomial of type (k, n) is defined by

$$\begin{aligned}\Theta_p(x) &= (x - g_0)(x - g_1) \cdots (x - g_{n-1}) \\ &= x^n + d_1x^{n-1} + d_2x^{n-2} + \cdots + d_{n-1}x + d_n.\end{aligned}$$

From the defining sums for $\eta_i$ and $g_i$, it is clear that $g_0 = 1 + n\eta_0$, $g_1 = 1 + n\eta_1$, ..., $g_{n-1} = 1 + n\eta_{n-1}$. Accordingly, the period and the reduced period polynomials are equivalent up to a linear transformation. The linking formulae are

$$\Psi_p(x) = \prod_{i=0}^{n-1}(x - \eta_i) = n^{-n}\Theta_p(nx+1),$$

and

$$\Theta_p(x) = \prod_{i=0}^{n-1}(x - \theta_i) = n^n\Psi_p((x-1)/n).$$

The Galois group Gal($\mathbf{Q}(\eta_0),\mathbf{Q}$) = { $\tau^i : 0 \le i < n$ } of the extension $\mathbf{Q}(\eta_0)$ of the rational numbers $\mathbf{Q}$ is a cyclic group of order n, and the generating map $\tau$ is defined by $\tau^i(\eta_j) = \eta_{i+j}$, this is exponentiation by g. Consequently

$$\tau(\Psi_p(x)) = \prod_{i=0}^{n-1}(x - \tau(\eta_i)) = \Psi_p(x).$$

In other words, both the period and the reduced period polynomials have integer coefficients.

The traces of $\eta_i$ and $g_i$ are

$$Tr(\eta_i) = \tau^0(\eta_i) + \tau^1(\eta_i) + \cdots + \tau^{n-1}(\eta_i) = \eta_0 + \eta_1 + \cdots + \eta_{n-1} = -1,$$

and

$$Tr(g_i) = 0.$$

These data immediately yield the first coefficients $c_1 = 1$ and $d_1 = 0$ of $\Psi_p(x)$ and $\Theta_p(x)$. Thus for any p, the coefficients have the pattern

$$\Psi_p(x) = x^n + x^{n-1} + c_2x^{n-2} + \cdots + c_{n-1}x + c_n,$$

and

$$\Theta_p(x) = x^n + d_{n-2}x^{n-2} + \cdots + d_{n-1}x + d_n.$$

The other coefficients $c_2, c_3, ..., c_n$ of $\Psi_p(x)$ are more difficult to calculate. Only in a few cases are the period polynomials completely determined. For example, it is quite easy to determine the period polynomials of maximal degree $n = p - 1 = deg(\Psi_p)$.

The coefficients of $\Psi_p(x)$ of low degree $n = deg(\Psi_p)$ are usually determined in terms of the quadratic partition of the prime p. The complexity of the quadratic partition of a prime p increases as n increases. For n = 5 it involves a Diophantine equation in 4 vaiables, see [ , Lehmer].

Two different approaches can be used to investigate the period polynomials.

( 1 ) $k = (p - 1)/n$ fixed and n as a function of the prime p.
( 2 ) $n = (p - 1)/k$ fixed and k as a function of the prime p.

If the parameter k = 2, 3, 4, 5, ... is small, and constant, then the periods are short. The number of terms k remains constant as p increases, but the degree $n = deg(\Psi_p)$ of period polynomial increases as a function of p. On the other hand, if the parameter n = 2, 3, 4, ... is small, and constant, then the periods are long, (have k terms which increases as p increases), but the degree $n = deg(\Psi_p)$ of period polynomial remains constant independent of p.

**Collection of Period Polynomials**
For every fixed n, the infinite sequence of primes $p = kn + 1$ determines an infinite sequence of period polynomials of degree $n = deg(\psi_p)$. The collection of period polynomials

$$\Psi_p(x) = x^7 + x^6 + c_2x^5 + \cdots + c_6x + c_7$$

is discussed in [1, Thaine].

## 5.2 Discriminant and Factorization of Period Polynomials

The discriminant of period polynomial provides various information about the factorization of $\Psi_p(x)$ over $\mathbf{F}_q$ and reciduacity of the primes q others than p.

***Definition 5.2.*** The discriminant of the period polynomial of type (k, n) is given by

$$D(\Psi_p) = \prod_{0 \le i < j < n} \left(\eta_i - \eta_j\right)^2 .$$

This integer also has the alternate form $D(\Psi_p) = P_1 P_2 \cdots P_{n-1}$, where the individual terms are

$$P_j = \prod_{0=i}^{n-1} (\eta_i - \eta_{i+j})^2 ,$$

$j = 1, 2, \ldots, n - 1$. Each $P_i$ is in fact an integer which is divisible by p. Thus the discriminant is itself divisible by $p^{n-1}$, see [E. Lehmer, 1987].

***Lemma 5.3.*** For every prime p, $D(\Theta_p) = n^{n(n-1)} D(\Psi_p)$.

Proof: Use $g_i = 1 + n\eta_i$. ■

***Lemma 5.4.*** (1) A prime $q \ne p$ not a divisor of $D(\Psi_p)$ is an nth power residue modulo p if and only if the equation $\Psi_p(x) \equiv 0 \bmod q$ is solvable.
(2) Every prime q divisor of $D(\Psi_p)$ is an nth power residue modulo p.

Proof: More details, see [E. Lehmer, 1968], [KS Williams, 1976]. ■

***Theorem 5.5.*** *(Discriminant Theorem)* If the polynomial $f(x) \in \mathbf{F}_q[x]$ is monic of degree n and nonzero discriminant D(f), then the quadratic symbol of $\mathbf{F}_q$ satisfies the relation

$$\left(\frac{D(f)}{q}\right) = (-1)^{n-s}$$

where s is the number of irreducible factors of f(x) in $\mathbf{F}_q[x]$.

Proof: Confer [1, Stepanov, p. 34.]. ■

***Theorem 5.6.*** Let $n$ be squarefree and let $q$ be a prime, $gcd(p, q) = 1$. Let $K_0$ ba a subgroup of index $k$ in the group of units of $\mathbf{Z}_n$. Let $l$ denotes the smallest integer for which $q^l \in K_0$. Then $l = lcm(d_1, d_2, \ldots, d_r)$, where $\Psi_n(x) = f_1(x) f_2(x) \cdots f_r(x) \in \mathbf{F}_q[x]$, $f_i(x)$, and $deg(f_i(x) = d_i$. In particular, $\Psi_n(x)$ splits into linear factors over $\mathbf{F}_q$ if $q \in K_0$.

Proof: See [Evans, p. 1077, 1989]. ■

***Theorem 5.7.*** (*Kummer-Dedekind*) Same assumption as above, if $f(x)$ is the minimal polynomial of $\alpha \in \mathbf{Q}(\eta)$ over the rational numbers, then $f(x)$ splits into irreducible factors over the rational numbers.

Proof: See [Januz, p. 32]. ■

## 5.3 Period Polynomials of Low Degrees

Periods polynomials of low degrees correspond to long periods. The number $k = (p - 1)/n$ of terms in a long period increases as a function of $p$, and the degree of the period polynomial is a small constant $n$. These polynomials are investigated here.

### Case n = 2

Let $p = 2k + 1$ be prime. The quadratic period and reduced period polynomials are the simplest case.

***Lemma 5.8.*** The quadratic period and reduced period polynomials are given by

$$\Psi_p(x) = x^2 + x + (1 - (-1)^{(p-1)/2}p)/4,$$

and

$$\Theta_p(x) = x^2 - (-1)^{(p-1)/2}p.$$

Proof: Compute the minimal polynomial $\Theta_p(x) = (x - g_0)(x - g_1)$ of the reduced periods, and then use the transformation formula $x \to (2x + 1)/4$ to obtain $\Psi_p(x)$. The reduced periods of type (k, 2) are the quadratic exponential sums

$$g_0 = \sum_{x=0}^{p-1} e^{i2\pi x^2/p} = \sum_{x=0}^{p-1} \left(\frac{x}{p}\right) e^{i2\pi x/p} \quad \text{and} \quad g_1 = \sum_{x=0}^{p-1} e^{i2\pi a x^2/p} = \sum_{x=0}^{p-1} \left(\frac{ax}{p}\right) e^{i2\pi x/p}$$

where $a$ is a quadratic nonresidue modulo p, and $( x \mid p )$ is the quadratic symbol. ■

**The Discriminant of $\Psi_p(x) = x^2 + x + c_2$**

The discriminant of the quadratic period polynomial is dependent on the congruence class of the prime p. The discriminant is given by

$$D(\Psi_p) = (\eta_0 - \eta_1)^2 = (-1)^{(p-1)/2}p.$$

**The Power Sums**

***Lemma 5.9.*** The power sums of the quadratic periods are given by
( 1 ) $S_1 = \eta_0 + \eta_1 = -1$,
( 2 ) $S_2 = \eta_0^2 + \eta_1^2 = (1-(-1)^{(p-1)/2}p)/2$,
( 3 ) $S_k = \eta_0^k + \eta_1^k = S_{k-1} + \left(\frac{1-(-1)^{(p-1)/2}p}{4}\right)S_{k-2}, \quad k > 2.$

**Irreducibility Conditions**

***Lemma 5.10.*** $\Psi_p(x) = x^2 + x + (1 - (-1)^{(p-1)/2}p)/4 \in \mathbf{Z}[x]$ is irreducible for all prime p.

***Lemma 5.11.*** The polynomial $\Psi_p(x) = x^2 + x + (1 - (-1)^{(p-1)/2}p)/4$ is irreducible over $\mathbf{F}_q$ if and only if the integer $(-1)^{(p-1)/2}p$ is a quadratic nonresidue modulo $q \neq 2^v$.

***Lemma 5.12.*** (1) $\Psi_p(x) = x^2 + x + (1 - (-1)^{(p-1)/2}p)/4 \in \mathbf{F}_2[x]$ is primitive if and only if $(1 - (-1)^{(p-1)/2}p)/4 \equiv 1 \bmod 2$, or equivalently $p \equiv \pm 3 \bmod 8$.
(2) $\Psi_p(x) = x^2 + x + (1 - (-1)^{(p-1)/2}p)/4 \in \mathbf{F}_3[x]$ is primitive if and only if $(1 - (-1)^{(p-1)/2}p)/4 \equiv -1 \bmod 3$ or equivalently $p \equiv 7 \bmod 12$.

**Table of Period Polynomials of Type (k, 2), p < 100.**

| | |
|---|---|
| $\Psi_3(x) = x^2 + x + 1$ | $\Psi_{43}(x) = x^2 + x + 11$ |
| $\Psi_5(x) = x^2 + x - 1$ | $\Psi_{47}(x) = x^2 + x + 12$ |
| $\Psi_7(x) = x^2 + x + 2$ | $\Psi_{53}(x) = x^2 + x - 13$ |
| $\Psi_{11}(x) = x^2 + x + 3$ | $\Psi_{59}(x) = x^2 + x + 15$ |
| $\Psi_{13}(x) = x^2 + x - 3$ | $\Psi_{61}(x) = x^2 + x - 15$ |
| $\Psi_{17}(x) = x^2 + x - 4$ | $\Psi_{67}(x) = x^2 + x + 17$ |
| $\Psi_{19}(x) = x^2 + x + 5$ | $\Psi_{71}(x) = x^2 + x + 18$ |
| $\Psi_{23}(x) = x^2 + x + 6$ | $\Psi_{73}(x) = x^2 + x - 18$ |
| $\Psi_{29}(x) = x^2 + x - 7$ | $\Psi_{79}(x) = x^2 + x + 20$ |
| $\Psi_{31}(x) = x^2 + x + 8$ | $\Psi_{83}(x) = x^2 + x + 21$ |
| $\Psi_{37}(x) = x^2 + x - 9$ | $\Psi_{89}(x) = x^2 + x - 22$ |
| $\Psi_{41}(x) = x^2 + x - 10$ | $\Psi_{97}(x) = x^2 + x - 24$ |

**Irreducible Quadratic Polynomials Algorithm**

The structure of the quadratic periods immediately leads to a deterministic algorithm for irreducible quadratic polynomials in $\mathbf{F}_q[x]$. The algorithm described here determines an irreducible quadratic polynomials in $\mathbf{F}_q[x]$ in fewer than $\log(q)^2$ trials.

*Quadratic Irreducible Polynomial Algorithm*

*Step 1.* If the prime $q \equiv 3 \bmod 4$, go to step (3). Otherwise, for the prime $q \equiv 1 \bmod 4$, choose an arbitrary prime p, and compute the quadratic symbol ( q | p ).

*Step 2*. If the quadratic symbol( q | p ) ≠ 1, then

$$\Psi_p(x) = x^2 + x + (1 - (-1)^{(p-1)/2}p)/4$$

is an irreducible quadratic polynomial in $\mathbf{F}_q[x]$, otherwise repeat step 1.
*Step 3*. Choose an arbitrary prime p ≡ 1 mod 4, and compute the quadratic symbol ( q | p ).
*Step 4*. If the quadratic symbol ( q | p ) ≠ 1, then $\Psi_p(x) = x^2 + x + (1 - (-1)^{(p-1)/2}p)/4 \in \mathbf{F}_q[x]$ is irreducible. Otherwise repeat step 3.

## Case n = 3

For the primes p = 3k + 2, the map $x \rightarrow x^3$ is one-to-one, and the cubic periods $\eta_0 = \eta_1 = \eta_3$ are trivial. In light of this put p = 3k + 1. The primes p split in the quadratic field $\mathbf{Q}(\omega)$, $\omega = (-1+\sqrt{-3})/2$, as $p = (a + \omega b)(a + \omega^2 b) = a^2 - ab + b^2$.

The cubic period and reduced period polynomials are given in terms of the parameters A and B in the quadratic partition of the primes $4p = A^2 + 27B^2$, A ≡ 1 mod 3, (or $p = A^2 + 3B^2$, A ≡ 1 mod 3).

***Lemma 5.13.*** The cubic period and reduced period polynomials are given by

$$\Psi_p(x) = x^3 + x^2 - \frac{p-1}{3}x - \frac{p(A+3)-1}{27}$$

and

$$\Theta_p(x) = x^3 - 3px - Ap.$$

Proof: To determine $\Psi_p(x)$, compute the minimal polynomial $\Theta_p(x) = (x - g_0)(x - g_1)(x - g_2)$ of the reduced periods, and then use the linear transformation formula $x \rightarrow (3x + 1)/27$ to obtain $\Psi_p(x)$. The reduced periods of type (k, 3) are the cubic exponential sums

$$g_0 = \sum_{x=0}^{p-1} e^{i2\pi a x^3/p}, \quad g_1 = \sum_{x=0}^{p-1} e^{i2\pi a x^3/p}, \quad g_2 = \sum_{x=0}^{p-1} e^{i2\pi a^2 x^3/p}$$

where *a* is a cubic nonresidue modulo p. To verify that $g_0$, $g_1$, $g_2$ are the roots of the polynomial $\Theta_p(x)$, rewrite each $g_i$ in terms of gaussian sums as

$$\begin{aligned}\sum_{x=0}^{p-1} e^{i2\pi a^v x^3/p} &= \sum_{x=0}^{p-1} \left(1 + \chi(ax) + \chi^2(a^2 x)\right) e^{i2\pi x/p} \\ &= G_v(\chi) + G_v(\chi^2),\end{aligned}$$

where $\chi$ is the cubic residue symbol, and v = 0, 1, 2. Now it is quickly resolved that the three real numbers $g_0, g_1, g_2 \in R$ are roots of $\Theta_p(x)$. In addition, these numbers satisfy the inequality $|G_v(\chi) + G_v(\chi^2)| \le 2p^{1/2}\cos(\rho)$, $0 \le \rho < \pi$. ■

***Lemma 5.14.*** The power sums of the cubic periods are given by

( 1 ) $S_1 = \eta_0 + \eta_1 + \eta_2 = -1$,
( 2 ) $S_2 = \eta_0^2 + \eta_1^2 + \eta_2^2 = (2p+1)/3$,
( 3 ) $S_3 = \eta_0^3 + \eta_1^3 + \eta_2^3 = [p(A-6)-1]/9$,
( 4 ) $S_k = \eta_0^k + \eta_1^k + \eta_2^k = -3pS_{k-2} - ApS_{k-3}$, $k > 3$.

Proof: The first sum is just the trace of $\eta_i$. The other two sums are computed with the recurring relation

$$c_0 S_k + c_1 S_{k-1} + \cdots + c_{k-1} S_1 + k c_k = 0$$

where $c_0 = 1$, $c_1 = 1$, $c_2$ , ..., $c_n$ are the coefficients of $\Psi_p(x)$, $k \leq 3$. ■

The previous recursive formula can be used to compute $Tr(\alpha^k) \in \mathbf{F}_q$, $\alpha$ a root of an irreducible polynomial $f(x)$ over $\mathbf{F}_q$. The sequence { $Tr(\alpha^k) : k \geq 0$ } is a recurring sequence of period $\leq q^n - 1$. Algorithms of this type are of interest in cryptography, see [ Lenstra ] LNCS 1270, 1716, 1718.

**The Discriminant of $\Psi_p(x) = x^3 + x^2 + c_2x + c_3$**

***Lemma 5.15.*** Let $4p = A^2 + 27B^2$, and $A \equiv 1 \bmod 3$, p prime. Then the discriminants of $\Psi_p(x)$ and $\Theta_p(x)$ are given by

$$D(\Psi_p) = (\eta_0 - \eta_1)^2(\eta_0 - \eta_2)^2(\eta_1 - \eta_2)^2 = p^2B^2,$$

and

$$D(\Theta_p) = (g_0 - g_1)^2(g_0 - g_2)^2(g_1 - g_2)^2 = 3^6p^2B^2.$$

Proof: Apply the discriminant formula $D(f) = 4a^3 + 27b^2$ of the polynomial $f(x) = x^3 + ax + b$, to $\Theta_p(x) = x^3 - 3px - Ap \in \mathbf{Z}[x]$, and then trace it back to $\Psi_p(x)$. ■

***Corollary 5.16.*** A prime $q \neq p$ is a cubic residue modulo p if and only if $B \equiv 0 \bmod q$.

**Irreducibility Conditions**

***Lemma 5.17.*** $\Psi_p(x) = x^3 + x^2 + c_2x + c_3 \in \mathbf{Z}[x]$ is irreducible for all prime p.

Proof: Since p divides the coefficients $d_1 = 0$, $d_2 = -3p$, and $d_3 = -Ap$ of $\Theta_p(x) = x^3 - 3px - Ap \in \mathbf{Z}[x]$, but $p^2$ does not divide $-Ap$, it follows that it is irreducible. Likewise, any linear transformation of $\Theta_p(x)$ like $\Psi_p(x)$ is also irreducible. ■

***Lemma 5.18.*** $\Psi_p(x) = x^3 + x^2 + c_2x + c_3$ either splits completely as $\Psi_p(x) = (x - \eta_0)(x - \eta_1)(x - \eta_2)$ or is irreducible over $\mathbf{F}_q$, q a prime power.

Proof: Apply the Discriminant Theorem with the parameters n = 3, and $D(\Psi_p) = (pB)^2$. Specifically

$$\left(\frac{D(\Psi_p)}{q}\right) = (-1)^{n-s} = (-1)^{3-s} = 1.$$

From this it is clear that s = 1 or s = 3, and this is precisely the number of irreducible factors that $\Psi_p(x)$ has over $\mathbf{F}_q$. ■

The quadratic and cubic period polynomials were determined by [1, Gauss, Articles. 356 and 358.].

**Table of Period and Reduced Period Polynomials of Type (k, 3), p < 100**

| | |
|---|---|
| $\Psi_7(x) = x^3 + x^2 - 2x - 1$ | $\Theta_7(x) = x^3 - 3\cdot7x - 1\cdot7$ |
| $\Psi_{13}(x) = x^3 + x^2 - 4x + 1$ | $\Theta_{13}(x) = x^3 - 3\cdot13x - 15\cdot13$ |
| $\Psi_{19}(x) = x^3 + x^2 - 6x - 7$ | $\Theta_{19}(x) = x^3 - 3\cdot19x - 7\cdot19$ |
| $\Psi_{31}(x) = x^3 + x^2 - 10x - 11$ | $\Theta_{31}(x) = x^3 - 3\cdot31x - 4\cdot31$ |
| $\Psi_{37}(x) = x^3 + x^2 - 12x + 11$ | $\Theta_{37}(x) = x^3 - 3\cdot37x - 11\cdot37$ |
| $\Psi_{37}(x) = x^3 + x^2 - 12x + 11$ | $\Theta_{43}(x) = x^3 - 3\cdot43x - A11\cdot43$ |
| $\Psi_{37}(x) = x^3 + x^2 - 12x + 11$ | $\Theta_{37}(x) = x^3 - 3\cdot61x - 11\cdot61$ |
| $\Psi_{37}(x) = x^3 + x^2 - 12x + 11$ | $\Theta_{37}(x) = x^3 - 3\cdot67x - 11\cdot67$ |
| $\Psi_{37}(x) = x^3 + x^2 - 12x + 11$ | $\Theta_{37}(x) = x^3 - 3\cdot73x - 11\cdot73$ |
| $\Psi_{37}(x) = x^3 + x^2 - 12x + 11$ | $\Theta_{37}(x) = x^3 - 3\cdot83x - 11\cdot83$ |
| $\Psi_{37}(x) = x^3 + x^2 - 12x + 11$ | $\Theta_{37}(x) = x^3 - 3\cdot97x - 11\cdot97$ |

## Case *n* = 4

In the quartic case the prime 4*k* + 1 split completely as product of two distinct primes $p = (a + ib)(a - ib)$ in the quadratic numbers field $\mathbf{Q}(i)$. The conditions
(1) $a \equiv 1 \bmod 4$, and
(2) $b \equiv ag^{(p-1)/4} \bmod 4$,

uniquely determine the parameters *a* and *b* in the quadratic partition $p = a^2 + 4b^2$.

***Lemma 5.19.*** The quartic period and reduced period polynomials are given by

$$\Psi_p(x) = x^4 + x^3 - \frac{3p-3}{8}x^2 - \frac{(2a+3)p-1}{16}x - \frac{(4a^2+8a-p+7)p-1}{256}$$

and

$$\Theta_4(x) = (x^2 - p)^2 - 4p(x - a)^2,$$

if the parameter $k$ in $p = kn + 1 = 4k + 1$ is even. Otherwise

$$\Psi_p(x) = x^4 + x^3 - \frac{p-3}{8}x^2 - \frac{p-(a+1)2p-1}{16}x - \frac{(4a^2-8a-9p-2)p-1}{256}$$

and

$$\Theta_p(x) = (x^2 + 3p)^2 - 4p(x - a)^2,$$

if the parameter $k$ in $p = kn + 1 = 4k + 1$ is odd.

The reduced periods of type $(k, 4)$ are the quartic exponential sums

$$g_0 = \sum_{x=0}^{p-1} e^{i2\pi x^4/p}, \quad g_1 = \sum_{x=0}^{p-1} e^{i2\pi a x^4/p}, \quad g_2 = \sum_{x=0}^{p-1} e^{i2\pi a^2 x^4/p}, \quad g_3 = \sum_{x=0}^{p-1} e^{i2\pi a^3 x^4/p}$$

where $a$ is a quartic nonresidue modulo $p$. The four roots $\eta_0, \eta_1, \eta_2, \eta_3$ of $\Psi_p(x) = x^4 + x^3 + c_2x^2 + c_1x + c_0 \in \mathbf{Z}[x]$ are

$$\sqrt{p} \pm \sqrt{2(-1)^{(p^2-1)/8}\left(p \pm a\sqrt{p}\right)}.$$

The value of the quadratic symbol $( 2 \mid p ) = (-1)^{(p^2-1)/8}$ determine the nature of the roots. The roots are distinct and real if $p \equiv 1 \bmod 8$, i.e., $\eta_0, \eta_1, \eta_2, \eta_3 \in \mathbf{R}$. Otherwise, the roots are distinct and nonreal if $p \equiv 5 \bmod 8$, i.e., $\eta_0, \eta_1, \eta_2, \eta_3 \in \mathbf{C}-\mathbf{R}$.

### The Discriminant of $\Psi_p(x) = x^4 + x^3 + c_2x^2 + c_3x + c_4$

***Lemma 5.20.*** The discriminant $D(\Psi_p)$ of $\Psi_p(x)$ is given by

$$D(\Psi_p) = \prod_{0\le i<j<4}(\eta_i - \eta_j)^2 = 2^{14}b^2p^3\left[a^2 + \left(1-2(-1)^{(p^2-1)/8}\right)p^2\right]^2.$$

### Irreducibility Conditions

***Lemma 5.21.*** $\Psi_p(x) = x^4 + x^3 + c_2x^2 + c_3x + c_4 \in \mathbf{Z}[x]$ is irreducible for all prime p.

Proof: Since $p = a^2 + b^2$, and $\Theta_4(x) = (x^2 - p)^2 - 4p(x - a)^2 = x^4 - 6px^2 + 8apx + b^2p$ or $(x^2 + 3p)^2 - 4p(x - a)^2 \in \mathbf{Z}[x]$, it readily follows that $p \mid d_1, d_2, d_3, d_4$, but $p^2$ does not divide $d_4$, (these are the coefficients of $\Theta_4(x)$). This confirms the irreducibility claim, and the same applies to $\Psi_p(x)$. ■

***Lemma 5.22.*** The polynomial $\Psi_p(x) = x^4 + x^3 + c_2x^2 + c_3x + c_4$ is irreducible over $\mathbf{F}_q$ if and only if the quadratic symbol $( p \mid q ) = -1$. Otherwise $( p \mid q ) = 1$, and it has either two quadratic or four linear factors.

Proof: Use the discriminant theorem. ■

## 5.4 Period Polynomials of High Degrees

The period polynomials of high degrees $n = (p - 1)/k$ correspond to short periods with k constant and n variable. The integer k is the number of terms in each period and n is the degree of the periods.

The calculations of the coefficients of periods polynomials of highest degrees $n = p - 1$ or $(p - 1)/2$ do not require information about the partition of the primes $p = kn + 1$.

### Case n = p – 1

On the upper end of the degree scale, there are the period polynomials of maximal degree $n = p - 1$. Since the cosets are $\{ y^{jn} \} = \{ 1 \}$, $\{ y^{jn+1} = y \}$, ..., $\{ y^{jn+n-2} \} = \{ y^{n-2} \}$, the periods of degree n are just the nontrivial pth roots of unity $\eta_0 = \omega$, $\eta_1 = \omega^2$, ..., $\eta_{n-1} = \omega^{n-1}$, so

$$\Psi_p(x) = \Phi_p(x) = x^n + x^{n-1} + \cdots + x + 1.$$

The period and cyclotomic polynomials of degree $n = p - 1$ coincide.

**The Discriminant of $\Psi_p(x) = \Phi_p(x)$**

***Lemma 5.23.*** The discriminant $D(\Phi_p)$ of $\Phi_p(x)$ is given by p???

**Irreducibility Conditions**

As a polynomial with integers coefficients, this polynomial is always irreducible over the rational numbers $\mathbf{Q} \subset R$, but as a polynomial with coefficients in the finite field $\mathbf{F}_q$ its irreducibility is a function of the order of q modulo n; also note that it cannot be primitive.

***Lemma 5.24.*** $\Psi_p(x) = x^n + x^{n-1} + \cdots + x + 1 \in \mathbf{Z}[x]$ is irreducible.

***Lemma 5.25.*** The polynomial $\Psi_p(x)$ is irreducible over $\mathbf{F}_q$ if and only if the integer q is of order $p - 1$ modulo p.

Proof: Let ω be a pth root of unity in the extension $\mathbf{F}_{q^n}$ of $\mathbf{F}_q$, and let d be the smallest integer such that $q^d \equiv 1 \bmod p$, $p = n + 1$. Then each root of $\Psi_p(x)$ is repeated $\varphi(p)/d = n/d$ times in the conjugates sequence

$$\omega, \quad \omega^q, \quad \omega^{q^2}, \quad ...., \quad \omega^{q^{d-1}}, \omega^{q^d} = \omega, \quad \omega^{q^{d+1}} = \omega^q, \quad ..., \quad \omega^{q^{n-1}} = \omega^{q^{d-1}}.$$

Hence the minimal polynomial of $\omega$ is a divisor of the polynomial

$$\prod_{i=0}^{n-1}\left(x-\omega^{q^i}\right)=\left(\prod_{i=0}^{d-1}\left(x-\omega^{q^i}\right)\right)^{n/d}$$

which is irreducible if q has order $p-1=n$ modulo p. Conversely, if $\Psi_p(x) \in \mathbf{F}_q[x]$ is irreducible then the conjugates sequence consists of n distinct elements. This implies that the exponents sequence $1, q, q^2, ..., q^{n-1}$ is a permutation of $0, 1, 2, ..., n-1$. Hence q has order $p-1=n$ modulo p. ■

In general the polynomial

$$f(x)=\prod_{i=0}^{n-1}\left(x-\omega^{q^i}\right)\in \mathbf{F}_q[x]$$

is always a power of the minimal polynomial $m_\alpha(x) \in \mathbf{F}_q[x]$ of the element $\alpha \in \mathbf{F}_{q^n}$. The exponent is determined by the order of q modulo p. Moreover, if the integer q does not have order p - 1 modulo p, then the polynomial f(x) splits into one irreducible factor over $\mathbf{F}_q$ of degree d raised to the power (or multiplicity) $\varphi(p)/d = (p-1)/d$.

**Case n = (p – 1)/2**

Let $p = 2n + 1$ be prime, g be a primitive root modulo p, and $\omega$ be a pth root of unity. In this case, the cosets are $K_0 = \{ -1, 1 \}$, $K_1 = \{ -g, g \}$, $K_2 = \{ -g^2, g^2 \}$, ..., $K_{n-1} = \{ -g^{n-1}, g^{n-1} \}$, and the periods are

$$\omega+\omega^{-1},\quad \omega^{g}+\omega^{-g},\quad \omega^{g^2}+\omega^{-g^2},\quad ...,\quad \omega^{g^{n-1}}+\omega^{-g^{n-1}}$$

The period polynomial $\Psi_p(x) = x^n + x^{n-1} + c_2x^{n-2} + \cdots + c_{n-1}x + c_n \in$ of degree $n = (p-1)/2$, which is the minimal polynomial of $\eta_0 = \omega + \omega^{-1}$, can be computed utilizing several different methods. Two of these methods are

*(1) The Recursive Formula Method, and*
*(2) Nonrecursive Method, (using power sums).*

**The Recursive Formula Method**

The recursive formula used to compute $\Psi_p(x)$ over the integers $\mathbf{Z}$ is

$$f_v(x) = xf_{v-1}(x) - f_{v-2}(x),$$

with initial conditions $f_0(x) = 1$, and $f_1(x) = x + 1$, for all $v \geq 2$. The recursion ends at the nth step with the polynomial $\Psi_p(x) = f_n(x)$.

The recursion is derived from the identity

$$1+2\sum_{i=1}^{n}\cos(ix)=0, \quad \text{and} \quad 1+2\sum_{i=1}^{n}T_i(x/2),$$

where $\eta_0 = 2\cos(ix)$ is a complex number, the polynomial $T_n(x) = \cos(n\arccos(x))$ satisfies the recursive relation $T_n(x) = 2xT_{n-1}(x) - T_{n-2}(x)$, with initial conditions $T_0(x) = 1$, and $T_1(x) = x$, see [1, Rybowicz].

An effective algorithm for computing the coefficients of polynomials defined by difference equations of second order, or recursive formula, is outlined in [1, Bini and Pan, p. 66.]. In the case of $f_n(x)$, this is written as

$$\begin{bmatrix} f_n(x) \\ f_{n-1}(x) \end{bmatrix} = \begin{bmatrix} x & -1 \\ 1 & 0 \end{bmatrix} \begin{bmatrix} f_{n-1}(x) \\ f_{n-2}(x) \end{bmatrix}.$$

**Nonrecursive Method**

A nonrecursive formula for computing these polynomials specifies the coefficients as certain binomial coefficients. The exact expression, cf. [1, Gauss, Article. 337], [{2, Gurak], and [1, Lehmer], has the shape

***Lemma 5.26.*** Let p = 2n + 1 be a prime. Then

$$\Psi_p(x)=\sum_{k=0}^{n}(-1)^{[k/2]}\binom{n-[(k+1)/2]}{[k/2]}x^k$$

where the bracket [ ? ] is the largest integer function and the braces ( ? ) is the binomial symbol.

Proof: Let $x^i + x^{-i} = \omega^i + \omega^{-i}$, $0 \le i \le n$, and let

$$g_n(x)=x^n+x^{-n}=\sum_{k=0}^{n/2}\frac{(-1)^k n}{n-k}\binom{n-k}{k}\left(x+x^{-1}\right)^{n-2k}$$

Since each $\omega^i + \omega^{-i} \neq 2$ is a distinct root of $\Psi_p(x)$ and $\omega^{jn} + \omega^{-jn} = \omega^{j(n+1)} + \omega^{-j(n+1)}$, it follows that each $\omega^i + \omega^{-i}$ is also a root of the polynomial

$$g_{n+1}(x) - g_n(x) = x^{n+1} + x^{-(n+1)} - (x^n + x^{-n}) = (x-2)\Psi_p(x)$$

of degree n + 1. Hence after simplification, this becomes

$$\Psi_p(x) = \sum_{k=0}^{n/2} \frac{(-1)^k n}{n-k} \binom{n-k}{k} x^{n-2k} + \sum_{k=0}^{(n-1)/2} \frac{(-1)^k n}{n-k} \binom{n-1-k}{k} x^{n-1-2k}$$
$$= \sum_{k=0}^{n/2} (-1)^k \left[ \binom{n-k}{k} x + \binom{n-1-k}{k} \right] x^{n-k}.$$

The identity $-\omega^{-1}(-1 - \omega^2) = \omega + \omega^{-1}$ permits a calculation of the norm of the element $\omega + \omega^{-1}$ up to a sign:

$$N(\omega + \omega^{-1}) = \prod_{v=0}^{2n} (\omega^v + \omega^{-v})$$
$$= \prod_{v=0}^{2n} (-\omega^{-v})(-1 + \omega^{2v}) = \Phi_p(-1) = 1.$$

Thus the norm $N(\omega + \omega^{-1}) = \pm 1$.

**The Discriminant of $\Psi_p(x)$**

***Lemma 5.27.*** The discriminant of $\Psi_p(x)$ is given by $D(\Psi_p) = p^{(p-3)/2}$.

**Irreducibility Condition**

Unlike the previous case this polynomial is not always irreducible over the rational numbers $\mathbf{Q} \subset$ R, but its irreducibility in both over rational numbers $\mathbf{Q}$ and over the finite field $\mathbf{F}_q$ is a function of the q and n etc..

***Lemma 5.28.*** $\Psi_p(x) = x^n + x^{n-1} + c_2x^{n-2} + \cdots + c_{n-1}x + c_n \in \mathbf{Z}[x]$ is irreducible if 2n + 1 is prime. Proof: Cf [Lehmer 1930].

***Lemma 5.29.*** Let p = 2n + 1. Then $\Psi_p(x) = x^n + x^{n-1} + \cdots \pm 1$ is irreducible over $\mathbf{F}_q$ if and only if the integer q is of order n or 2n modulo p.

For the important case $\Psi_p(x) = x^n + x^{n-1} + \cdots + 1 \in \mathbf{F}_2[x]$ there are a few sequences of primes 2n + 1 that satisfy the condition in the previous lemma. These are sequences or clusters of primes:

( 1 ) If p = 2n + 1, n = 4c + 3 prime, then 2 has order n modulo p.
( 2 ) If p = 2n + 1, n = 4c + 1 prime, then 2 is a primitive root modulo p.
( 3 ) If p = 4k + 1, k prime, then 2 is a primitive root modulo p. These primes occur in clusters.

The polynomial $\Psi_p(x)$ can also be primitive but only over the finite fields $\mathbf{F}_2$ and $\mathbf{F}_3$. Recall that a *normal polynomial* is a polynomial with linearly independent roots over the coefficients field.

***Theorem 5.30.*** If both n and 2n + 1 are primes, then the period polynomials $\Psi_p(x) = x^n + x^{n-1} + \cdots + 1 \in \mathbf{F}_2[x]$ is a primitive normal polynomial.

This theorem is a direct consequence of the *primitive optimal normal basis theorem* for n degree extension $\mathbf{F}_{q^n}$ over $\mathbf{F}_2$.

**Table of Period Polynomials of Type (2, n).**

$\Psi_0(x) = 1$
$\Psi_3(x) = x + 1$
$\Psi_5(x) = x^2 + x - 1$
$\Psi_7(x) = x^3 + x^2 - 2x - 1$
$\Psi_9(x) = x^4 + x^3 - 3x^2 - 2x + 1^*$
$\Psi_{11}(x) = x^5 + x^4 - 4x^3 - 3x^2 + 3x + 1$
$\Psi_{13}(x) = x^6 + x^5 - 5x^4 - 4x^3 + 6x^2 + 3x - 1$
$\Psi_{15}(x) = x^7 + x^6 - 6x^5 - 5x^4 + 10x^3 + 6x^2 - 4x - 1^*$
$\Psi_{17}(x) = x^8 + x^7 - 7x^6 - 6x^5 + 15x^4 + 10x^3 - 10x^2 - 4x + 1$
$\Psi_{19}(x) = x^9 + x^8 - 8x^7 - 7x^6 + 21x^5 + 15x^4 - 20x^3 - 10x^2 + 5x + 1$
$\Psi_{21}(x) = x^{10} + x^9 - 9x^8 - 8x^7 + 28x^6 + 21x^5 - 35x^4 - 20x^3 + 15x^2 + 5x - 1^*$
$\Psi_{23}(x) = x^{11} + x^{10} - 10x^9 - 9x^8 + 36x^7 + 28x^6 - 41x^5 - 35x^4 + 35x^3 + 20x^2 - 6x - 1$
$\Psi_{25}(x) = x^{12} + x^{11} - 11x^{10} - 10x^9 + 45x^8 + 36x^7 - 63x^6 - 41x^5 + 70x^4 + 21x^3 - 26x^2 - 7x - 1^*$

The entries with a star do not correspond to primes p = 2n + 1. These were generated with the recursive formula given above.

**Table of Primitive Period Polynomials of Type (2, n) over $\mathbf{F}_2$**

$\Psi_3(x) = x + 1$
$\Psi_5(x) = x^2 + x + 1$
$\Psi_7(x) = x^3 + x^2 + 1$
$\Psi_{11}(x) = x^5 + x^4 + x^2 + x + 1$
$\Psi_{13}(x) = x^6 + x^5 + x^4 + x + 1$
$\Psi_{19}(x) = x^9 + x^8 + x^6 + x^5 + x^4 + x + 1$
$\Psi_{23}(x) = x^{11} + x^{10} - 10x^9 - 9x^8 + 36x^7 + 28x^6 - 41x^5 - 35x^4 + 35x^3 + 20x^2 - 6x - 1$

There are also cases like $\Psi_{19}(x)$, $\Psi_{37}(x)$, ... in $\mathbf{F}_2[x]$, (for n = 9, 18, ... and p = 2n + 1 = 19, 37, ... ), which are primitives, but the theory has not been worked out. A larger table of the primitive polynomials $\Psi_p(x)$ for $11 \le p \le 201$, (or n < 100), appears in [1, Rybowicz ].

The norm $N(\alpha) = (-1)^n f(0)$ of any root $\alpha$ of a primitive polynomial $f(x) = x^n + a_{n-1}x^{n-1} + \cdots + a_1x + a_0 \in \mathbf{F}_q[x]$ must be a primitive element in $\mathbf{F}_q$. Thus $\Psi_p(x) = x^n + x^{n-1} + c_{n-2}x^{n-2} + \cdots + c_1x + c_0$ $\in \mathbf{F}_3[x]$ is a potential primitive polynomial if and only if n = 4c +1 or 4c + 2, use the nonrecursive formula to verify that $c_0 = \pm 1$. This leads to primes of the form p = 2n + 1 = 8c +3 or 8c + 5.

**Irreducible Period Polynomials of Type (2, n) in $\mathbf{F}_3$[x]**

$\Psi_5(x) = x^2 + x - 1$, primitive.
$\Psi_7(x) = x^3 + x^2 - 2x - 1$, irreducible but not primitive.
$\Psi_{11}(x) = x^5 + x^4 - x^3 + 1$, primitive.
$\Psi_{17}(x) = x^8 + x^7 - x^6 + x^3 - x^2 - x + 1$, irreducible but not primitive.
$\Psi_{19}(x) = x^9 + x^8 - x^7 - x^6 + x^3 + x^2 + 2x + 1$, primitive.
$\Psi_{21}(x) = x^{10} + x^9 - 9x^8 - 8x^7 + 28x^6 + x^5 - 35x^4 - 20x^3 + 15x^2 + 5x - 1^*$.?
$\Psi_{23}(x) = x^{11} + x^{10} - 10x^9 - 9x^8 + 36x^7 + 28x^6 - 41x^5 - 35x^4 + 35x^3 + 20x^2 - 6x - 1$?.
$\Psi_{25}(x) = x^{12} + x^{11} - 11x^{10} - 10x^9 + 45x^8 + 36x^7 - 63x^6 - 41x^5 + 70x^4 + 21x^3 - 26x^2 - 7x - 1^*$.?
$\Psi_{31}(x) = x^{15} + x^{14} - \cdots - 1$, irreducible but not primitive.
$\Psi_{127}(x) = x^{63} + x^{62} - \cdots - 1$, irreducible but not primitive.
$\Psi_p(x) = x^n + x^{n-1} - \cdots - 1$, $p = 2n + 1 = 2^r - 1$, is irreducible if 3 has order n or 2n modulo p, but not primitive since the norm $N(\omega+\omega^{-1}) = (-1)^n\Psi_p(0) = 1$ is not primitive in $\mathbf{F}_3$.
$\Psi_{257}(x) = x^{128} + x^{127} - \cdots + 1$ is irreducible since 3 has order 256 modulo 257, but not primitive since the norm $N(\omega+\omega^{-1}) = (-1)^n\Psi_p(0) = 1$ is not primitive in $\mathbf{F}_3$.
$\Psi_{65535}(x) = x^{32767} + x^{32767} - \cdots + 1$ is irreducible since 3 has order 65534 modulo 65535, but not primitive since the norm $N(\omega+\omega^{-1}) = (-1)^n\Psi_p(0) = 1$ is not primitive in $\mathbf{F}_3$.

***Lemma 5.31.*** Let $\omega \in GF(q^{kn})$ be an element of order ord($\omega$) = kn + 1, where kn + 1 is a prime and k is an even integer. Then the degree $d = \deg(\Psi_d(x))$ of the minimal polynomial $\Psi_d(x)$ of $\omega + \omega^{-1}$ is the smallest integer such that $q^d \equiv \pm 1 \bmod (kn + 1)$.

Proof: $GF(q^d)$ is the smallest subfield of $GF(q^{kn})$ containing $\omega + \omega^{-1}$, so

$$f(x) = \prod_{i=0}^{kn-1}\left(x - (\omega + \omega^{-1})^{q^i}\right) = \left(\prod_{i=0}^{d-1}\left(x - (\omega + \omega^{-1})^{q^i}\right)\right)^{kn/d} = \Psi_p(x)^{kn/d}.$$

The same apply to a linear combination $a\omega + b\omega^{-1}$, $0 \neq a, b \in \mathbf{F}_q$. The instance of kn + 1 = 2n + 1 prime and $2^n \equiv \pm 1 \bmod (2n + 1)$ is a special case with $f(x) = \Psi_p(x)^2$.

## 5.5 Coefficients Calculations Via Power Sums Method

The linking formulae and Newton's identities are some of the basic tools used to compute the coefficients of period polynomials.

The dth power sums $S_d$ of the root of the period polynomials, which appear in Newton's identity, are defined by

$$S_d = \eta_0^d + \eta_1^d + \eta_2^d + \cdots + \eta_{n-1}^d .$$

The dth sum $S_d$ is computed recursively from a given list $S_0, S_1, \ldots, S_{d-1}$. Several of the sums $S_d$ for small parameter k in the prime p = kn + 1 have been determined. This section introduces the methods utilized to compute these integers.

***Lemma 5.32.*** (*Lehmers 1983*) Let p = 2n + 1 be a prime, $\eta_0 = \omega + \omega^{-1}$, and let $d \geq 1$. Then

$$S_d = -2^{d-1} + \frac{p}{2} \sum_{2u \equiv d \bmod p} \binom{d}{u} .$$

Proof: As the index runs from i = 1 to p – 1, each term in the list $(\omega + \omega^{-1})^v, (\omega^2 + \omega^{-2})^v, \ldots, (\omega^{p-1} + \omega^{-(p-1)})^v$, with v fixed, is repeated twice. Thus

$$\begin{aligned} 2S_d &= \sum_{i=1}^{p-1} (\omega^i + \omega^{-i})^d \\ &= \sum_{u=0}^{d} \binom{d}{u} \sum_{i=1}^{p-1} \omega^{i(2u-d)} \\ &= \sum_{u=0}^{d} \binom{d}{u} \left(p\delta_{d,2u} - 1\right), \end{aligned}$$

where $\delta_{i,j}$ is the unit impulse function. The claim follows from this. ■

The pattern of the even and odd power sums are simply

$$S_d = -2^{d-1} + \frac{p}{2}\binom{d}{d/2} \quad \text{and} \quad S_d = -2^{d-1} + \frac{p}{2}\binom{d}{u} .$$

with 2u ≡ d mod p, respectively. For example, $S_1 = -1$, $S_2 = p - 2$, $S_3 = -4$, $S_4 = 3p - 8$, ... These values immediately lead to a recursive determination of the coefficients of $\Psi_p(x)$. The first few are

$c_0 = 1,$
$c_1 = -S_1 = 1,$
$c_2 = -(c_0S_2 + c_1S_1)/2 = (3 - p)/2,$
$c_3 = -(c_0S_3 + c_1S_2 + c_2S_1)/3 = (5 - p)/2,$
$c_4 = -(c_0S_4 + c_1S_3 + c_2S_2 + c_3S_1)/4 = (p^2 - 12p + 35)/8,$
… … …

$$c_k = \frac{-1}{k} \sum_{v=1}^{k} S_v c_{k-v}, \qquad k < \sqrt{p} .$$

## Case n = (p – 1)/3

In the next result for p = 3n + 1, the period polynomials $\Psi_p(x)$ are polynomials of degree $\deg(\Psi_p(x)) = n$. The degree n = (p – 1)/3 increases as p increases, but the number of terms in each period is three independently of p. The first one is

$$\eta_0 = \omega + \omega^a + \omega^{a^2} \ ,$$

where $1 \neq a$, $a^3 = 1$. There are two possibilities: $a = (-1+\sqrt{-3})/2$ in $\mathbf{F}_p$.

The power sums for small d is as follows.

***Lemma 5.33.*** (*Lehmers 1983*) If $d < \sqrt{p}$, then

$$S_d = \begin{cases} -3^{d-1} + \dfrac{p}{3}\dbinom{3s}{s}\dbinom{2s}{s} & \text{if } d = 3s, \\ -3^{d-1} & \text{otherwise.} \end{cases}$$

Proof: See [1, Lehmer (Short period), p. 749.]. ■

These values immediately lead to a recursive determination of the coefficients of $\Psi_p(x)$. The first few are

$c_0 = 1$,
$c_1 = -S_1 = 1$,
$c_2 = -(c_0S_2 + c_1S_1)/2 = 2$
$c_3 = -(c_0S_3 + c_1S_2 + c_2S_1)/3 = (14 - 2p)/3$,
$c_4 = -(c_0S_4 + c_1S_3 + c_2S_2 + c_3S_1)/4 = (35 - 2p)/3$,
$c_5 = -(c_0S_5 + c_1S_4 + c_2S_3 + c_3S_2 + c_4S_1)/5 = (35 - 2p)/3$,
$c_6 = -(c_0S_6 + c_1S_5 + c_2S_4 + c_3S_3 + c_4S_2 + c_5S_1)/6 = (2p^2 - 73p + 728)/9$,
$c_7 = -(c_0S_7 + c_1S_6 + c_2S_5 + c_3S_4 + c_4S_3 + c_5S_2 + c_6S_1)/7 = (2p^2 - 115p + 1976)/9$,
$c_8 = -(c_0S_8 + c_1S_7 + c_2S_6 + c_3S_5 + c_4S_4 + c_5S_3 + c_6S_2 + c_7S_1)/8 = (4p^2 - 272p + 5434)/9$,
$c_9 = -(c_0S_9 + c_1S_8 + c_2S_7 + c_3S_6 + c_4S_5 + c_5S_4 + c_6S_3 + c_7S_2 + c_8S_1)/9 = (4p^3 - 354p^2 + 11298p - 135850)/81$, …

**Table of Period Polynomials of Type (3, n).**

A small table of the polynomials $\Psi_p(x)$ is provided here. The parameter *a* is the smallest value of $a = (-1+\sqrt{-3})/2$ in $\mathbf{F}_p$.

$\Psi_7(x) = x^2 + x + 2$, a = 2,
$\Psi_{13}(x) = x^4 + x^3 + 2x^2 - 4x + 3$, a = 3,

$\Psi_{19}(x) = x^6 + x^5 + 2x^4 - 8x^3 - x^2 + 5x + 7$, a = 7,
$\Psi_{31}(x) = x^{10} + x^9 + 2x^8 - 16x^7 - 9x^6 - 11x^5 + 43x^4 + 6x^3 + 63x^2 + 20x + 25$, a = 5,
$\Psi_{37}(x) = x^{12} + x^{11} + 2x^{10} - 20x^9 - 13x^8 - 19x^7 + 85x^6 + 51x^5 + 94x^4 - 2x^3 - 13x^2 - 77x + 47$, a = 10,
$\psi_{43}(x) = x^{14} + x^{13} + 2x^{12} - 24x^{11} - 17x^{10} - 27x^9 + 143x^8 + 81x^7 + 83x^6 - 209x^5 + 163x^4 + 88x^3 + 235x^2 - 168x + 79$, a = 6.

## Case n = (p – 1)/4

Similarly, for p = 4n + 1, and $\Psi_p(x)$ is the minimal polynomial of the period

$$\eta_0 = \omega + \omega^a + \omega^{a^2} + \omega^{a^3} = \omega + \omega^{-1} + \omega^a + \omega^{-a}$$

where $1 \neq a$, $a^4 = 1$. There is one possibility: $a = (-1)^{1/2}$ in $\mathbf{F}_p$.

***Lemma 5.34.*** (*Lehmers 1983*) If $d < \sqrt{p}$, then

$$S_d = \begin{cases} -4^{d-1} + \frac{p}{4}\binom{d}{d/2}^2 & \text{if } d = 2s, \\ -4^{d-1} & \text{otherwise.} \end{cases}$$

Proof: See [1, Lehmer (Short period), p. 754.]. ■

These values immediately lead to a recursive determination of the coefficients of $\Psi_p(x)$. The first few are

$c_0 = 1$,
$c_1 = -S_1 = 1$,
$c_2 = -(c_0S_2 + c_1S_1)/2 = -(p - 5)/2$
$c_3 = -(c_0S_3 + c_1S_2 + c_2S_1)/3 = -(p - 15)/2$,
$c_4 = -(c_0S_4 + c_1S_3 + c_2S_2 + c_3S_1)/4 = (p^2 - 28p + 195)/8$,
$c_5 = -(c_0S_5 + c_1S_4 + c_2S_3 + c_3S_2 + c_4S_1)/5 = (p^2 - 48p + 663)/8$,
$c_6 = -(c_0S_6 + c_1S_5 + c_2S_4 + c_3S_3 + c_4S_2 + c_5S_1)/6 = -(p^3 - 69p^2 + 1655p - 13923)/48$,
$c_7 = -(c_0S_7 + c_1S_6 + c_2S_5 + c_3S_4 + c_4S_3 + c_5S_2 + c_6S_1)/7 = -(p^3 - 99p^2 + 3599p - 49725)/48, \ldots$

**Table of Period Polynomials of Degree n and Type (4, n).**
A small table of the polynomials $\Psi_p(x)$ is provided here, $a = (-1)^{1/2}$.

$\Psi_5(x) = x + 1$,
$\Psi_{13}(x) = x^3 + x^2 - 4x + 1$,
$\Psi_{17}(x) = x^4 + x^3 - 6x^2 - x + 1$,
$\Psi_{29}(x) = x^7 + x^6 - 12x^5 - 7x^4 + 28x^3 + 14x^2 - 9x + 1$,
$\Psi_{37}(x) = x^9 + x^8 - 16x^7 - 11x^6 + 66x^5 + 32x^4 - 73x^3 - 7x^2 + 7x + 1$,
$\Psi_{41}(x) = x^{10} + x^9 - 18x^8 - 13x^7 + 91x^6 + 47x^5 - 143x^4 - 7x^3 + 72x^2 - 23x + 1$,

Techniques for computing the period polynomials with respect to composite parameters kn + 1 also appear in { 1, Lehmers]. More materials about the factorizations of period polynomials and related questions, see [1, 2, Gurak], [1, Gupta and Zagier], [2, Lehmer ], and [1, Meyers], etc. The period polynomials are very closely related to the Dickson polynomials of the first and second kinds:

$$D_n(x,a) = \sum_{k=0}^{n/2} \frac{n}{n-k}\binom{n-k}{k}(-a)^k x^{n-2k},$$

and

$$E_n(x,a) = \sum_{k=0}^{n/2} \binom{n-k}{k}(-a)^k x^{n-2k}.$$

The recursive relations for these polynomials are

$D_{n+2}(x) = xD_{n+1}(x) - D_n(x)$, with initial conditions $D_0(x) = 2$, and $D_1(x) = x$, and

$E_{n+2}(x) = xE_{n+1}(x) - E_n(x)$, with initial conditions $D_0(x) = 1$, and $D_1(x) = 2x$. The exact expression is

$$\Psi_p(x) = \frac{D_{n+1}(x,1) - D_n(x,1)}{x-2}.$$

## 5.6 Sequences of Period Polynomials

The *n*th power residuacity of *q* modulo *p* is employed to determine the irreducibility status of the *n*th period polynomials, and generate sequence of irreducible polynomials.

***Theorem 5.35.*** Let $p = kn + 1 = 2^u k + 1$, and q be primes, and suppose that q is a quadratic nonresidue modulo p. then the period polynomial

$$\Psi_p(x) = x^n + x^{n-1} + c_2x^{n-2} + \cdots + c_{n-2}x^2 + c_{n-1}x + c_n \in \mathbf{F}_q$$

is irreducible for all n = = $2^u$.

Proof: The hypothesis ( q | p ) = –1 implies that $q^{(p-1)/n} \neq 1$ mod p, so q is an nth power nonresidue mod p. The latter in turn implies that the periods $\eta_0, \eta_1, \ldots, \eta_{n-1}$ are linearly

independent over $\mathbf{F}_q$. ■

***Example 5.36.*** Take q = 3, and $p = 3\cdot 2^u a + 1$. Then q = 3 is a quadratic nonresidue modulo p. Accordingly the period polynomial

$$\Psi_p(x) = x^n + x^{n-1} + c_2x^{n-2} + \cdots + c_{n-2}x^2 + c_{n-1}x + c_n \in \mathbf{F}_3$$

is irreducible for all $n = 2^u$.

$\Psi_p(x) = x^8 + x^7 + c_2x^6 + c_3x^5 + c_4x^4 + c_5x^3 + c_6x^2 + c_7x + c_8{}^*$

is irreducible if $( p \mid q ) = -1$. Otherwise $( p \mid q ) = 1$ and it has 1, 2, 4, or 8 factors.

$\Psi_{25}(x) = x^{12} + x^{11} - 11x^{10} - 10x^9 + 45x^8 + 36x^7 - 63x^6 - 41x^5 + 70x^4 + 21x^3 - 26x^2 - 7x - 1^*$

The discriminant of $\Psi_{pq}(x)$ is given in [1, Brillhart].

# Chapter 6

# Periods Normal Bases

# 6.1 Definitions and Existence

Let $r \in \mathbb{N}$ be an integer such that $\gcd(r, q) = 1$, and put $\varphi(r) = kn$, where $\varphi$ is the totient function. Let $\omega$ be an rth root of unity in the rth cyclotomic field extension $\mathbf{F}_q(\omega)$ of $\mathbf{F}_q$. The finite field $\mathbf{F}_q(\omega)$ is an extension of $\mathbf{F}_q$ of degree $m = [\,\mathbf{F}_q(\omega) : \mathbf{F}_q\,]$, where $m \mid kn$ is the smallest integer such that $q^m - 1 \equiv 0 \bmod r$. The constraint $\gcd(r, q) = 1$ ensures the existence of nontrivial rth roots of unity in characteristic char($\mathbf{F}_q$).

***Definition 6.1.*** Let $K_0, K_1, ..., K_{n-1}$ be subsets of $\mathbf{Z}_r^*$. The quasi-periods are defined by

$$\eta_0 = \sum_{x \in K_0} \omega^x,\ \eta_1 = \sum_{x \in K_1} \omega^x,\ ...,\ \eta_{n-1} = \sum_{x \in K_{n-1}} \omega^x .$$

The set of elements $\{ \eta_0, \eta_1, ..., \eta_{n-1} \}$ generated by an arbitrary list of subsets $K_0, K_1, ..., K_{n-1}$ is a potential basis of $\mathbf{F}_{q^n}$ over $\mathbf{F}_q$. Normal bases require conjugate elements:

$$\eta_i,\ \eta_{i+1} = \eta_i^q,\ \eta_{i+2} = \eta_i^{q^2},\ ...,\ \eta_{i+n-1} = \eta_i^{q^{n-1}} .$$

A simple method for constructing conjugate elements is described here.

***Lemma 6.2.*** Let K be a subset of $\mathbf{Z}_r^*$, and suppose that $K_0 = K$, $K_1 = qK$, $K_2 = q^2K$, ..., $K_{n-1} = q^{n-1}K$ are distinct subsets of $\mathbf{Z}_r^*$. Then the quasi-periods $\eta_0, \eta_1, ..., \eta_{n-1}$ are conjugates over $\mathbf{F}_q$.

Proof: Take an automorphism $\sigma^j$ of $\mathbf{F}_{q^n}$. Then

$$\sigma^j(\eta_i) = \eta_i^{q^j} = \left( \sum_{x \in K_i} \omega^x \right)^{q^j} = \eta_{i+j} \qquad \blacksquare$$

A set of conjugates elements might fail to be a normal basis simply because either the trace is identically zero or the trace is not an element of $\mathbf{F}_q$.

***Lemma 6.3.*** Let $K \cup qK \cup q^2K \cup \cdots \cup q^{n-1}K = \mathbf{Z}_r^*$, where $q^iK \neq q^jK$ for $i \neq j$, and let n be a divisor of $\varphi(r) = kn$, r squarefree. Then

( 1 ) The quasi-periods $\eta_0, \eta_1, ..., \eta_{n-1}$ are conjugates over $\mathbf{F}_q$.
( 2 ) $Tr(\eta_i) = \eta_0 + \eta_1 + \cdots + \eta_{n-1} = \pm 1$.

Proof: ( 2 ) Since the subsets are disjoint, and r is squarefree,

$$Tr(\eta_i) = \sum_{\gcd(x,r)=1} \omega^x = \mu(r) = \pm 1,$$

where $\mu$ is the Mobius function. ■

In general a union $K_0 \cup K_1 \cup K_2 \cup \cdots \cup K_{n-1} = \mathbf{Z}_r^*$ of nondisjoint subsets, perhaps all of the subsets of the same cardinality does not satisfy the nonzero trace property. Likewise nonsquarefree integers r yield periods of zero traces.

The construction of period normal bases requires well-structured partitions of the multiplicative group $\mathbf{Z}_r^*$. The best-known construction method of period normal bases utilizes cosets partitions of $\mathbf{Z}_r^*$. The cosets of nth power residues and non-residues of the multiplicative group $\mathbf{Z}_r^*$, (also the multiplicative groups of extensions of $\mathbf{Z}_r$) have the structure sought after. The partition of nth power residues and nonresidues consists of n distinct cosets $K_0 \cup K_1 \cup K_2 \cup \cdots \cup K_{n-1}$ of $\mathbf{Z}_r^*$, each of cardinality k. The group $\mathbf{Z}_r^*$ has one or more subgroups of order n, depending on r, and each subgroup induces a coset partition of $\mathbf{Z}_r^*$.

***Definition 6.4.*** If the integer r is a prime then the subset $K \subset \mathbf{Z}_r^*$ of nth power residues is a unique subgroup, and the induced periods $\eta_0, \eta_1, \ldots, \eta_{n-1}$ are called *prime gaussian periods*, otherwise *nonprime gaussian periods*.

For applications to bases of finite field extension $\mathbf{F}_{q^n}$ of $\mathbf{F}_q$, it is convenient to take g = q if the integer q is primitive in $\mathbf{F}_p$, (or a nth power nonresidue). Under this condition the gaussian periods of type (k, n) are written in the more convenient form

$$\eta_0 = \sum_{j=0}^{k-1} \omega^{q^{jn}},\ \eta_1 = \sum_{j=0}^{k-1} \omega^{q^{jn+1}},\ \ldots,\ \eta_{n-1} = \sum_{j=0}^{k-1} \omega^{q^{jn+n-1}} .$$

## 6.2 Existence of Period Normal Bases

One of the main problem in the theory of normal bases of finite fields is to determine under what condition a list of conjugate elements $\eta_0, \eta_1, \ldots, \eta_{n-1}$ forms a normal basis of $\mathbf{F}_{q^n}$ over $\mathbf{F}_q$. Specifically, is the trivial solution $\mathbf{a} = (a_{n-1},\ldots,a_1,a_0) = \mathbf{0}$ the only solution of the equation

$$a_0\eta_0 + a_1\eta_1 + \cdots + a_{n-1}\eta_{n-1} = 0 .$$

The focus here is on the identification of period normal bases.

**Some Necessary Properties of Normal Bases**
(1) The elements $\eta_0, \eta_1, \ldots, \eta_{n-1}$ are conjugates.
(2) The trace $Tr(\eta_0) = \cdots = Tr(\eta_{n-1}) = \eta_0 + \eta_1 + \cdots + \eta_{n-1} \neq 0$ is an element in $\mathbf{F}_q$.
(3) $a_0\eta_0 + a_1\eta_1 + \cdots + a_{n-1}\eta_{n-1} = 0 \quad \Leftrightarrow \quad (a_{n-1},\ldots,a_1,a_0) = (0,\ldots,0,0)$.

Property (3) is sufficient for a basis, but not for a normal basis. A general normal basis test automatically verifies all theses properties. Several results or tests are stated below.

***Theorem 6.5.*** Let K be an arbitrary subset of $\mathbf{Z}_r^*$. Then the quasi-periods derived from this set is a normal basis of $\mathbf{F}_{q^n}$ over $\mathbf{F}_q$ if the following conditions hold.

(1) The union $\mathbf{Z}_r^* = K \cup qK \cup q^2K \cup \cdots \cup q^{n-1}K$ is a disjoint union, i.e., $q^iK \neq q^jK$ for $i \neq j$, where n is a divisor of $\varphi(r) = kn$.
(2) For all irreducible factor a(x) of $x^n - 1 = a(x)b(x) \in \mathbf{F}_q[x]$,

$$\frac{x^n - 1}{a(x)} \circ \eta = \sum_{i=0}^{r} b_i \eta^{q^i} \neq 0,$$

where n = de, and $b(x) = b_d x^d + b_{d-1}x^{d-1} + \cdots + b_1 x + b_0$.

The simplest case involves the two-factor factorization $x^n - 1 = (x - 1)f(x)$, f(x) irreducible in $\mathbf{F}_q[x]$. Under this setting, the system of inequalities in the above theorem has two lines
(i) $Tr(\eta_i) \neq 0$, and (ii) $\eta_i^q - \eta_i \neq 0$. In other simple cases the polynomial $x^n - 1$ has only a few irreducible factors.

Another simple period normal basis test is covered below. This test indirectly verifies all the normality conditions using only integers arithmetic operations.

***Theorem 6.6.*** Suppose r = kn + 1 is a prime, $r \nmid q$, and $d = ord_r(q)$. Then a gaussian period of type (k, n) is a normal element in $\mathbf{F}_{q^n}$ over $\mathbf{F}_q$ if and only if gcd(kn/d, n) = 1.

The integer e = kn/d is the index of q modulo r. This stipulates that n be a divisor of the order of q modulo r, so $d = ord_r(q) \geq n$, and q generates a cyclic subgroup of order $\#\langle q \rangle = an$, $a \geq 1$. A proof of the linear independence, (from the point of view of group theory), appears in [1, Menezes et al., p.101.].

For a given prime r and prime power q, gcd(r, q) = 1, this result provides a fast algorithm for identifying a normal gaussian period of type (k, n) in an extension $\mathbf{F}_{q^n}$ of $\mathbf{F}_q$.

A similar result but unconditional and based on the matrix determinant test det(N) ≠ 0 of the regular matrix representation N attached to the subset $\{ \eta_0, \eta_1, ..., \eta_{n-1} \}$ is illustrated here.

***Theorem 6.7.*** Let r = kn + 1 be a prime, and let q be an nth power nonresidue modulo r. Then the gaussian periods $\eta_i$ of degree n constitute a normal basis of $\mathbf{F}_{q^n}$ over $\mathbf{F}_q$.

Proof: Let g be primitive modulo r, and let $K = \{ g^{jn} : 0 \leq j < k \}$. Then the n cosets satisfy the inequality $q^iK \neq q^jK$ for $i \neq j$, $0 \leq i\ j < n$, ( note that $q = g^{jn+i}$ ). This condition implies the existence of n distinct periods $\eta_0, \eta_1, ..., \eta_{n-1}$ of type (k, n). Further, the matrix N is a symmetric circulant matrix, that is,

$$N = \left(\eta_i^{q^j}\right) = \left(\eta_{i+j}\right).$$

Since the determinant of a circular matrix is the product of the discrete Fourier transform of the generating vector ($\eta_0, \eta_1, ..., \eta_{n-1}$), the determinant of the circulant matrix N is

$$\det(N) = \prod_{s=0}^{n-1} \sum_{t=0}^{n-1} \eta_t \xi^{st} = \prod_{s=0}^{n-1} G(\chi_s) \neq 0,$$

where $\xi$ is a primitive nth root of unity in $\mathbf{F}_{q^{kn}}$. In the third term above the inner sum is rephrased as a gaussian sum $G(\chi_s)$, where the character $\chi_s(g) = \eta^s$. The sum $G(\chi_s)$ has the value $G(\chi_0) = -1$, and $|G(\chi_s)|^2 = p$, for $0 < s < n$. Hence $|\det(N)|^2 = p^{n-1}$, or equivalently $\det(N) = \pm p^{(n-1)/2}$. This proves the linear independence of the subset of periods. ■

The proof of this last result combine ideas from algebraic number theory and linear algebra, see the chapter on Structured Matrices for more material on this.

Observe that this result also provides a fast two steps algorithm for identifying the period normal bases using only integers arithmetics:

( 1 ) Primality testing of kn + 1, and
( 2 ) Determination of the residuacity of q modulo kn + 1.

The nth power nonresidue test is verified with the inequality

$q^{(r-1)/n} \not\equiv 1 \bmod r.$

The existence of the prime kn + 1 is taken up in the following result; reference [Gao et al., p. 318, 1995] is an alternative source.

***Theorem 6.8.*** (*Wassermann 1993*) Let $q = p^v$, p prime, and $n \in \mathbb{N}$. Suppose that both

(1) gcd(v, n) = 1, and
(2) $2p \nmid n$ if $p \equiv 1 \bmod 4$ or $4p \nmid n$ if $p \equiv 2, 3 \bmod 4$.

Then assuming the extended Riemann hypothesis, there is an integer $k < cn^3(\log(np))^2$ such that $r = kn + 1$ is a prime and $\gcd(kn/d, n) = 1$, where c is a constant independent of v, n, p, and $d = \text{ord}_r(q)$ is the order of q modulo r.

***Example 6.9.*** Let the integer r = 2n + 1 be a prime and let q = 2 be of order n or 2n modulo 2n + 1, then the type (2, n) periods in the rth cyclotomic field $\mathbf{F}_2(\omega)$ extension of $\mathbf{F}_2$ are given by

$$\eta_i = \sum_{j=0}^{k-1} \omega^{2^{jn+i}} = \omega^{2^i} + \omega^{-2^i},$$

i = 0, 1, 2, …, n–1. The degree of the extension is $[\mathbf{F}_2(\omega) : \mathbf{F}_2] = n$ or $2n$ depending on $n = 4c + 3$ or $n = 4c + 1$, (also $4c + 2$). This scheme of producing periods works whenever order $\operatorname{ord}_r(q) = 2n$, or q generates the nth power residues modulo r. But can fail in other cases, for example, let $r = 2n + 1$ be prime, and 2 be of order $\operatorname{ord}_r(q) = n$ modulo r. Since 2 does not generate the nth power residues $K = \{ -1, 1 \}$ modulo r, the previous method for constructing the type (n, 2) periods in $\mathbf{F}_2(\omega) = \mathbf{F}_{q^{2n}}$ over $\mathbf{F}_2$ cannot be used. Instead the periods

$$\eta_0 = \sum_{x\in K_0}\omega^x = \omega + \omega,\ \eta_1 = \sum_{x\in K_1}\omega^x = \omega^2 + \omega^{-2}, ..., \eta_{n-1} = \sum_{x\in K_{n-1}}\omega^x = \omega^{2n-1} + \omega^{-2^{n-1}}$$

are derived directly from the cosets $K_0 = \{ -1, 1 \}$, $K_1 = \{ -2, 2 \}$, $K_2 = \{ -2^2, 2^2 \}$, ..., $K_{n-1} = \{ -2^{n-1}, 2^{n-1} \}$. The resulting set of periods $\{ \eta_0, \eta_1, ..., \eta_{n-1} \}$ constitutes a normal basis of $\mathbf{F}_{q^n}$ over $\mathbf{F}_2$ of degree $n = [\mathbf{F}_{q^n} : \mathbf{F}_2]$. The field $\mathbf{F}_2(\omega+\omega^{-1}) = \mathbf{F}_{q^n}$ is the maximal subfield of $\mathbf{F}_2(\omega)$. The diagram below shows the subgroup/subfield correspondence; the group of automorphisms of $\mathbf{F}_2(\omega)$ is $\operatorname{Gal}(\mathbf{F}_{q^{2n}}/\mathbf{F}_2) = \{ \sigma^i : 0 \le i < 2n \}$, and the subgroups are of the form $\Sigma_d = \{ \sigma^{di} : 0 \le i < 2n/d \}$, $d \mid 2n$.

$$\begin{array}{ccccc} \mathbf{F}_2 & \xrightarrow{n} & \mathbf{F}_2(\omega+\omega^{-1}) & \xrightarrow{2} & \mathbf{F}_2(\omega) \\ | & & | & & | \\ \Sigma_1 & \longleftarrow & \Sigma_n & \longleftarrow & \Sigma_{2n} \end{array}$$

***Example 6.10.*** The 29th cyclotomic field $\mathbf{F}_2(\omega)$ extension of $\mathbf{F}_2$, ω a primitive 29th root of unity. The value of the totient function is $\varphi(29) = 2^2\cdot 7$, so there are intermediate subfields of degree 2, 4, 7 and 14. Each of these subfields has a period normal basis. The generating periods are

$$\alpha_i = \sum_{j=0}^{13}\omega^{2^{2j+i}},\quad \beta_i = \sum_{j=0}^{6}\omega^{2^{4j+i}},\quad \eta_i = \sum_{j=0}^{3}\omega^{2^{7j+i}},\quad \theta_i = \omega^{2^i} + \omega^{-2^i},$$

where the index i = 0,1, …. The subfield $\mathbf{F}_2(\omega+\omega^{-1})$ is the maximal subfield of $\mathbf{F}_2(\omega)$ of degree $14 = [\mathbf{F}_2(\omega+\omega^{-1}) : \mathbf{F}_2]$. The next subfield in this tower is $\mathbf{F}_2(\eta)$ of degree $7 = [\mathbf{F}_2(\eta) : \mathbf{F}_2]$.
The map $\sigma(\omega) = \omega^2$ generates the automorphisms group $\operatorname{Gal}(\mathbf{F}_2(\omega)/\mathbf{F}_2) = \{ \sigma^i : 0 \le i < 28 \}$of $\mathbf{F}_2(\omega)$. The subgroups are $\Sigma_d = \{ \sigma^{id} : 0 \le i < 28/d \}$, for d = 2, 4, 7, 14.

The subgroups and subfields correspondences are the followings:

$\sigma^2(\alpha_i) = \alpha_i$ ⇔ the subgroup $\Sigma_2 = \{ \sigma^{2i} : 0 < i < 14 \}$ fixes $\mathbf{F}_2(\alpha_i)$,

$\sigma^7(\beta_i) = \beta_i$ ⇔ the subgroup $\Sigma_4 = \{ \sigma^{4i} : 0 < i < 7 \}$ fixes $\mathbf{F}_2(\beta_i)$,

$\sigma^7(\eta_i) = \eta_i$ $\Leftrightarrow$ the subgroup $\Sigma_7 = \{ \sigma^{7i} : 0 < i < 4 \}$ fixes $\mathbf{F}_2(\eta)$, and

$\sigma^{14}(\omega + \omega^{-1}) = \omega + \omega^{-1}$ $\Leftrightarrow$ the subgroup $\Sigma_{14} = \{ 1, \sigma^{14} \}$ fixes the subfield $\mathbf{F}_2(\omega+\omega^{-1})$.

Some of the fields and groups diagrams are illustrated below.

$$\begin{array}{ccccccc}
\mathbf{F}_2 & \xrightarrow{7} & \mathbf{F}_2(\eta) & \xrightarrow{2} & \mathbf{F}_2(\omega+\omega^{-1}) & \xrightarrow{2} & \mathbf{F}_2(\omega) \\
| & & | & & | & & | \\
\Sigma_1 & \longleftarrow & \Sigma_7 & \longleftarrow & \Sigma_2 & \longleftarrow & \Sigma_{28}
\end{array}$$

and

$$\begin{array}{ccccccc}
\mathbf{F}_2 & \xrightarrow{2} & \mathbf{F}_2(\alpha) & \xrightarrow{2} & \mathbf{F}_2(\beta) & \xrightarrow{7} & \mathbf{F}_2(\omega) \\
| & & | & & | & & | \\
\Sigma_1 & \longleftarrow & \Sigma_2 & \longleftarrow & \Sigma_4 & \longleftarrow & \Sigma_{28}
\end{array}$$

**Other Related Period Normal Bases**

Let $\mathbf{F}_{q^d}$ be an extension of $\mathbf{F}_q$ and a subfield of $\mathbf{F}_{q^n}$, n = de. The relative trace $Tr_{n:d} : \mathbf{F}_{q^n} \to \mathbf{F}_{q^d}$ is defined by

$$Tr_{n:d}(\alpha) = \alpha + \alpha^{q^d} + \alpha^{q^{2d}} + \dots + \alpha^{q^{(e-1)d}}.$$

***Lemma 6.11.*** *(Period projection lemma)* Suppose $\{ \eta_0, \eta_1, ..., \eta_{n-1} \}$ is a period normal basis of $\mathbf{F}_{q^n}$ over $\mathbf{F}_q$, and let $\gamma_i = Tr_{n:d}(\eta_i) = \eta_i + \eta_{i+d} + \eta_{i+2d} + \eta_{i+3d} + \cdots + \eta_{i+n-d}$. Then the subset of elements $\{ \gamma_0, \gamma_1, ..., \gamma_{d-1} \}$ is a normal basis of $\mathbf{F}_{q^d}$ over $\mathbf{F}_q$.

Proof: It is easy to verify that the elements $\gamma_0, \gamma_1, ..., \gamma_{d-1}$ are distinct and conjugates, moreover $Tr_{d:1}(\gamma_i) = \eta_0 + \eta_1 + \cdots + \eta_{n-1} = -1$. To confirm the linear independence of these elements, assume that $a_0\gamma_0 + a_1\gamma_1 + \cdots + a_{d-1}\gamma_{d-1} = 0$ for some nonzero vector $\mathbf{a} = (a_0,a_1,..., a_{d-1}) \neq (0,0,...,0)$ and use the fact that $\{ \eta_0, \eta_1, ..., \eta_{n-1} \}$ is a basis to reach a contradiction. ■

As a demonstration of the period projection lemma, let r = 2n +1 be a prime, let the integer q be of order 2n modulo r, and $1 \neq \omega \in \mathbf{F}_{q^n}$ but $\omega^r = 1$, r minimal. Then the followings hold.

(1) The element ω generates a period normal basis of $\mathbf{F}_{q^{2n}}$ over $\mathbf{F}_q$.

(2) The relative trace $Tr_{2n:n}(\omega) = \omega + \omega^{-1}$ (since $q^n \equiv -1 \bmod r$) of the element ω generates a period normal basis of $\mathbf{F}_{q^n}$ over $\mathbf{F}_q$. Here $2 = [\, \mathbf{F}_q(\omega) : \mathbf{F}_q(\omega+\omega^{-1}) \,]$, and the relative trace $Tr_{2n:n}$ :

$\mathbf{F}_{q^{2n}} \to \mathbf{F}_{q^n}$.

## 6.3 Polynomial Representations of Periods

The normal elements in finite fields are represented in the multiplicative group of the polynomials algebra $\mathbf{F}_q[x]/(x^n - 1)$ as invertible polynomials.

***Definition 6.12.*** The polynomial representation of an element $\alpha \in \mathbf{F}_{q^n}$ is defined by

$$c_\alpha(x) = \sum_{i=0}^{n-1} Tr(\alpha^{1+q^i})\, x^i .$$

The significance of this correspondence is that the inverses of the polynomial representations of normal elements generate dual normal elements. The same map is also utilized to produce cyclic codes in this algebra.

More generally, for a fixed element $\xi \in \mathbf{F}_{q^n}$, the function $\kappa : \mathbf{F}_{q^n} \to \mathbf{F}_q[x]/(x^n - 1)$ defined by

$$\alpha \to \kappa(\alpha) = \sum_{i=0}^{n-1} Tr\left(\alpha\, \xi^{q^i}\right) x^i$$

establishes a similar correspondence between elements and polynomial representations or n-tuples. The latter map is also used to reproduce the coefficients vectors $(Tr(\alpha\xi_{n-1}), ..., Tr(\alpha\xi_1)Tr(\alpha\xi_0))$, $\alpha \in \mathbf{F}_{q^n}$ of the code words of cyclic codes.

***Lemma 6.13.*** If the element $\eta_0$ is normal over $\mathbf{F}_q$, then the map

$\eta_0 \to c(x) = Tr(\eta_0\eta_{n-1})x^{n-1} + Tr(\eta_0\eta_{n-2})x^{n-2} + \cdots + Tr(\eta_0\eta_1)x + Tr(\eta_0\eta_0)$

is into and *sn*-to-one, $s \geq 1$. Moreover, $c(x)$ is invertible in $\mathbf{F}_q[x]/(x^n - 1)$.

The representations of arbitrary normal elements are not necessarily easy to determine, but the representations of period normal elements are not difficult to determined.

***Lemma 6.14.*** Let $r = kn + 1$ be prime, and suppose that $r \equiv 1 \bmod 2n$. Then
( 1 ) The polynomial representation of the period $\eta_0$ of type (k, n) and its inverse are given by

$c(x) = -k(x^{n-1} + x^{n-2} + \cdots + x) + (n - 1)k + 1$,
and

$$c^{-1}(x) = \frac{k}{kn+1}\left( x^{n-1} + x^{n-1} + \cdots + x + \frac{k+1}{k} \right).$$

( 2 ) If q | k, then the representations of all the periods are mapped to c(x) = 1, viz,

$\{ \eta_0, \eta_1, \ldots, \eta_{n-1} \} \rightarrow c(x) = 1$.

Proof: ( 1 ) Let $c(x) = c_{n-1}x^{n-1} + c_{n-2}x^{n-2} + \cdots + c_1x + c_0$ be the polynomial representation of normal element $\eta_0 \in \mathbf{F}_{q^n}$ over $\mathbf{F}_q$. Then the ith coefficient $c_i$ of c(x) is given by

$$Tr(\eta_0\eta_i) = Tr\left( k\varepsilon_i + \sum_{j=0}^{n-1}(i,j)\eta_j \right) = \begin{cases} kn - \sum_{j=0}^{n-1}(i,j) & \text{if } -1 \in K_i, \\ -\sum_{j=0}^{n-1}(i,j) & \text{if } -1 \notin K_i. \end{cases}$$

The evaluation of $c_i = Tr(\eta_0\eta_i)$ utilizes the linear expansion of the pairwise products $\eta_0\eta_j$, see the chapter on Periods. The linearity of the trace function, and $Tr(\eta_i) = -1$. Moreover, applying the properties of the cyclotomic numbers (i, j), this becomes

$$c_i = Tr(\eta_0\eta_i) = \begin{cases} k(n-1) & \text{if } -1 \in K_i, \\ -k & \text{if } -1 \notin K_i. \end{cases}$$

By hypothesis $r \equiv 1 \bmod 2n$, so it follows that $-1 \in K_0$, so the coefficients for the representation c(x) of $\eta_0$ are: $c_0 = Tr(\eta_0\eta_0) = (n-1)k + 1$, and $c_i = Tr(\eta_0\eta_i) = -k$, $0 < i < n$. To verify that the given polynomials are inverses, compute $c(x)c^{-1}(x) \bmod (x^n - 1)$. ■

The constraint $-1 \in K_0$ holds whenever $-1$ is an nth power residue modulo r, or equivalently if and only if the equation $x^n = -1$ is solvable modulo r. This equation is always solvable if the parameter n in r = kn + 1 is odd or $r \equiv 1 \bmod 2n$. . For the case $-1 \in K_j$ for some $j \neq 0$, $c_j = Tr(\eta_0\eta_0) = (n-1)k + 1$, and $c_i = Tr(\eta_0\eta_i) = -k$ for all $i \neq j$.

It also should be observed that these representations of the periods work whether or not the periods form normal bases of the finite field $\mathbf{F}_{q^n}$ over $\mathbf{F}_q$.

***Example 6.15.*** ( 1 ) For r = 2k + 1 prime, k even, the polynomial representation of the quadratic period $\eta_0$ and its inverse are given by

$$c(x) = -kx + k + 1 \quad \text{and} \quad c(x)^{-1} = \frac{k}{2k+1}x + \frac{k+1}{2k+1}.$$

And for k odd, the polynomial representation its inverse are given by

$$c(x) = (k+1)x - k \quad \text{and} \quad c(x)^{-1} = \frac{k+1}{2k+1}x + \frac{k}{2k+1}.$$

( 2 ) For r = 3k + 1, k even only, the polynomial representation of the cubic period $\eta_0$ its inverse are given by

$$c(x) = -k(x^2+x)+2k+1 \text{ and } c(x)^{-1} = \frac{k(x^2+x)}{3k+1} + \frac{k+1}{3k+1}.$$

( 3 ) For r = 4k + 1, k even, the polynomial representation of the quartic period is

$$c(x) = -k(x^3+x^2+x) + 3k+1$$

and its inverse is

$$c(x)^{-1} = \frac{k(x^3+x^2+x)}{4k+1} + \frac{k+1}{4k+1} .$$

And for $r = 4k + 1$, $k$ odd, $c(x) = c_3 x^3 + c_2 x^2 + c_1 x - k$, where $c_i = k(n-1)$ if $-1 \in K_i$, otherwise $c_i = -k$.

## 6.4 Dual Period Normal Bases

This section considers a general method of determining dual bases of period normal bases, and continues to the self-dual periods normal bases. The general construction of period dual normal bases uses the inverses of the polynomial representations of the period $\eta_0$ in the polynomials algebra $\mathbf{F}_q[x]/(x^n-1)$

***Lemma 6.16.*** (*Gao, von zur Gathen, Penario 1995*) Let $-1 \in K_r$, $0 \le r < n-1$, and let $\delta = (\eta_{-r} - k)/(kn+1)$. Then

$$\delta_0 = \delta,\ \delta_1 = \delta^q,\ \delta_2 = \delta^{q^2}\ \ldots,\ \delta_{n-1} = \delta^{q^{n-1}}$$

is the dual basis of $\eta_0, \eta_1, \ldots, \eta_{n-1}$.

Proof: Compute the dual bases relation $\mathrm{Tr}(\delta_i\eta_j)$, specifically

$$Tr\left(\delta^{q^i}\eta^{q^j}\right) = Tr\left(\frac{\eta^{q^{i-r}} - k}{kn+1}\eta^{q^j}\right) = \frac{1}{kn+1}\left[Tr\left(\eta^{q^{i-r}}\eta^{q^j}\right) - k\,Tr\left(\eta^{q^j}\right)\right]$$
$$= \frac{1}{kn+1}\left[Tr\left(\eta^{q^{i-r}}\eta^{q^j}\right) + k\right]$$

since the trace $\mathrm{Tr}(\eta_i) = -1$. Now replace the linear expansion of the pairwise terms $\eta_i\eta_j$ in $Tr(\eta_i\eta_j)$, use the linearity of the trace function, and $Tr(x) = nx$ for $x \in \mathbf{F}_q$, to obtain (for k even)

$$Tr\left(\eta\eta^{q^{j-i+r}}\right) = Tr\left(k\varepsilon_i + \sum_{d=0}^{n-1}(j-i+r, d)\eta^{q^d}\right) = \begin{cases} kn - \sum_{d=0}^{n-1}(j-i+r, d) & \text{if } -1 \in K_{j-i+r}, \\ -\sum_{d=0}^{n-1}(j-i+r, d) & \text{if } -1 \notin K_{j-i+r}. \end{cases}$$

Using properties of the cyclotomic numbers (i, j), this reduces to

$$Tr\left(\eta\,\eta^{q^{j-i+r}}\right) = Tr\left(\eta^{q^{i-r}}\,\eta^{q^j}\right) = \begin{cases} k(n-1)+1 & \text{if } -1 \in K_{j-i+r}, \\ -k & \text{if } -1 \notin K_{j-i+r}. \end{cases}$$

Since k even mean r = 0, substituting this into the last equation returns

$$Tr\left(\delta^{q^i}\eta^{q^j}\right) = \begin{cases} 1 & \text{if } i = j, \\ 0 & \text{if } i \neq j. \end{cases}$$

This complete the verification of the dual bases equation $Tr(\delta_i\eta_j) = \delta_{i,j}$. ■

Recall that a self-dual normal basis of $\mathbf{F}_{q^n}$ over $\mathbf{F}_q$ exists if and only if one of the following conditions holds.

( 1 ) $n = 2c + 1$ and $q = p^v$ or
( 2 ) $n \neq 4c$ and $q = 2^v$.

***Lemma 6.17.*** (*Gao, von zur Gathen, Penario* 1995) Let n > 2. A normal basis { $\eta_0, \eta_1, ..., \eta_{n-1}$ } of $\mathbf{F}_{q^n}$ over $\mathbf{F}_q$ consisting of the type (k, n) periods is a self-dual basis if and only if one of the following conditions are met.

( 1 ) $n \not\equiv 0 \bmod 4$, $q = 2^v$, and k is even.
( 2 ) $n \not\equiv 0 \bmod 2$, $q = p^v$, p > 2 prime, and k is even, and divisible by p.

Proof: (The original proof is not complete) The condition on the degree n of the extension is required for the existence of self-dual normal bases of $\mathbf{F}_{q^n}$ over $\mathbf{F}_q$. The condition k is even $\Rightarrow$ $-1 \in K_0$, i.e., $-1$ is an nth power residue in $\mathbf{F}_r$, $r = kn + 1$, so that $\delta = (\eta_0 - k)/(kn + 1)$. The divisibility condition $p \mid k$ is precisely the constraint needed to force the dual element to $\delta = (\eta_r - k)/(kn + 1) = (\eta_0 - k)/(kn + 1) = \eta_0 = \eta$. ■

In the first case, the condition k = even is sufficient to yield a self-dual period normal basis in characteristic 2, the divisibility statement is already met.

The 2-term periods normal bases are the easiest to manage in calculations, and quite often these are self-dual period normal bases.

***Lemma 6.18.*** Let r = 2n + 1 be a prime, and suppose that { $\eta_0, \eta_1, ..., \eta_{n-1}$ } is a period normal basis of $\mathbf{F}_{q^n}$ over $\mathbf{F}_q$. Then the type (2, n) periods constitute a self-dual basis if and only if conditions (1) or (2) are satisfied.

( 1 ) $n \neq 4c$, and $q = 2^v$ has order $ord_r(q) = n$ or $2n$ modulo r
( 2 ) $n = 2c + 1$, and $q = p^v$ has order $ord_r(q) = 2n$ modulo r $q = p^v$, $p > 2$ prime.

Proof: ( 1 ) The parameters $n \neq 4c$, $q = 2^v$ signal the possibility of a self-dual normal basis. Since q has order n or 2n modulo r = 2n + 1, the n periods

$$\eta_0 = \omega + \omega^{-1}, \quad \eta_1 = \omega^q + \omega^{-q}, \quad \eta_2 = \omega^{q^2} + \omega^{-q^2}, \quad ..., \eta_{n-1} = \omega^{q^{n-2}} + \omega^{-q^{n-2}}$$

form a normal basis of $\mathbf{F}_{q^n}$ over $\mathbf{F}_q$. To verify the dual property, consider

$$\eta_i \eta_j = \left(\omega^{q^i} + \omega^{-q^i}\right)\left(\omega^{q^j} + \omega^{-q^j}\right) = \left(\omega^{q^a} + \omega^{-q^a}\right) + \left(\omega^{q^b} + \omega^{-q^b}\right).$$

Now $q^i + q^j = q^a$, and $q^i - q^j = q^b$ in $\mathbf{Z}_r$, $0 \leq a, b < r$, whenever $q$ generates the subset of quadratic residues, confer [1, Mullin et al.]. Combine these to obtain

$$Tr(\eta_i \eta_j) = \begin{cases} 1 & \text{if } i = j \\ 0 & \text{if } i \neq j \end{cases}$$

as claimed. ■

***Example 6.19.*** Generators of self-dual period normal bases of $\mathbf{F}_{2^n}$ over $\mathbf{F}_2$, for 1 < n < 11.

n = 2, r = 2·2 + 1 = 5, $\eta_0 = \omega + \omega^{-1}$, where $1 \neq \omega \in \mathbf{F}_{2^4}$ and $\omega^5 = 1$.

n = 3, r = 2·3 + 1 = 7, $\eta_0 = \omega + \omega^{-1}$, where $1 \neq \omega \in \mathbf{F}_{2^6}$ and $\omega^7 = 1$.

n = 4, none.

n = 5, r = 2·5 + 1 = 11, $\eta_0 = \omega + \omega^{-1}$, where $1 \neq \omega \in \mathbf{F}_{2^{10}}$ and $\omega^{11} = 1$.

n = 6, r = 2·6 + 1 = 13, $\eta_0 = \omega + \omega^{-1}$, where $1 \neq \omega \in \mathbf{F}_{2^{12}}$ and $\omega^{13} = 1$.

n = 7, r = 4·7 + 1 = 29, $\eta_0 = \omega + \omega^{2^7} + \omega^{2^{14}} + \omega^{2^{21}}$, where $1 \neq \omega \in \mathbf{F}_{2^{28}}$ and $\omega^{29} = 1$.

n = 8, none.

n = 9, r = 2·9 + 1 = 19, $\eta_0 = \omega + \omega^{-1}$, where $1 \neq \omega \in \mathbf{F}_{2^{18}}$ and $\omega^{19} = 1$.

n = 10, r = 4·10 + 1 = 41, $\eta_0 = \omega + \omega^{2^{10}} + \omega^{2^{20}} + \omega^{2^{30}}$, where $1 \neq \omega \in \mathbf{F}_{2^{40}}$ and $\omega^{41} = 1$.

# CHAPTER 7

# Period Normal Bases For Extensions of Low Degrees

## 7.1 Quadratic Extensions

Several results concerning period normal bases for extensions of low degrees will be considered here. The existence results for extensions of low degree are quite often easy to verify and are subject to a variety of techniques not applicable to the general cases. The simplest period normal bases are those for quadratic, cubic, and quartic extensions of $\mathbf{F}_q$.

The group of automorphisms of the quadratic extension $\mathbf{F}_{q^2}$ of $\mathbf{F}_q$ is the set Gal($\mathbf{F}_{q^2}/\mathbf{F}_q$) = { 1, σ } of linear maps $\sigma : \mathbf{F}_{q^2} \rightarrow \mathbf{F}_q$. The nontrivial automorphism is defined by $\sigma(x) = x^q$. Sometimes is convenient to use the conjugation form $\sigma(a+b\theta) = a - b\theta$, where $a, b \in \mathbf{F}_q$, and some $\theta \notin \mathbf{F}_q$.

***Lemma 1.*** Let $p = 2k + 1$ be a fixed prime, and suppose that q is a nonquadratic residue modulo p, this requires $q \not\equiv 1 \bmod p$. Then the pair $\{ \eta_0, \eta_1 \}$ is a period normal basis of the quadratic extension $\mathbf{F}_{q^2}$ over $\mathbf{F}_q$ .

Proof: Let $\omega \neq 1$ be a pth root of unity in some extension $\mathbf{F}_{q^d}$ of $\mathbf{F}_q$, $q \neq 2^v$, let g be a generator of the multiplicative group of $\mathbf{F}_p$, and write $q = g^{2a+1}$. Then

$$\eta_0 = \sum_{x \in K_0} \omega^x = \frac{1}{2}\sum_{x=1}^{p-1} e^{i2\pi x^2/p}, \quad \eta_1 = \sum_{x \in K_1} \omega^x = \frac{1}{2}\sum_{x=1}^{p-1} e^{i2\pi g x^2/p},$$

where $K_0 = \{ 1, g^2, g^4, ..., g^{p-3} \}$ and $K_1 = \{ q, qg^2, qg^4, ..., qg^{p-3} \}$. Since $q \not\equiv 1 \bmod p$, the element ω is not in $\mathbf{F}_q$. Similarly $\eta_0, \eta_1 \notin \mathbf{F}_q$, this follows from $\eta_i^q \neq \eta_i$ .The linear independence follows from the matrix determinant test for basis. Specifically, using the identity $\eta_0 + \eta_1 = -1$, this is

$$\det \begin{bmatrix} \eta_0 & \eta_1 \\ \eta_1 & \eta_0 \end{bmatrix} = \eta_1 - \eta_0 \neq 0 \ .$$

Thus the set $\{ \eta_0, \eta_1 \}$ constitutes a normal basis of $\mathbf{F}_{q^2}$ over $\mathbf{F}_q$. ■

***Lemma 2.*** Let p be a prime and let $q \neq 2^v$ a prime power. Then the quadratic periods

$$\eta_0 = \frac{-1 + (-1)^{(p-1)/4}\sqrt{p}}{2}, \quad \eta_1 = \frac{-1 - (-1)^{(p-1)/4}\sqrt{p}}{2}$$

form a basis of $\mathbf{F}_{q^2}$ over $\mathbf{F}_q$ if and only if the quadratic symbol $((-1)^{(p-1)/2}p \mid q) = -1$.

Proof: Compute the discriminant and the roots of the period polynomial $\psi_p(x) = x^2 + x + (1 - (-1)^{(p-1)/2}p)/4$. ■

**Remark:** If the gcd(n, q) ≠ 1, the identity

$$\eta_s = \sum_{x \in K_s} \omega^x = \frac{1}{n} \sum_{x=1}^{p-1} e^{i2\pi x^n / p},$$

s = 0, 1, … n − 1, is not defined in characteristic char($\mathbf{F}_q$), so it cannot be used in calculations. However, the periods can still exist, in fact this is often the case for the quadratic periods in characteristic char($\mathbf{F}_q$) = 2.

From Lemmas 1 and 2 it is clear that quadratic periods are of the form a + bθ, with a, b ∈ $\mathbf{F}_q$, but θ ∉ $\mathbf{F}_q$. It is also clear that there is a close link between the existence of quadratic period normal bases and the reciduacity of the integers p and q. The reciprocity law

$$\left(\frac{p}{q}\right)\left(\frac{q}{p}\right) = (-1)^{(p-1)(q-1)/4}$$

is one of the main tools used for determining the existence of various types of quadratic period normal bases.

***Theorem 3.*** A pair of primes p and q determines a period normal basis of $\mathbf{F}_{q^2}$ over $\mathbf{F}_q$ if and only if either $((-1)^{(p-1)/2}p \mid q) = -1$ or $(q \mid p) = -1$ holds.

Proof: Easy. ■

The suborder of an element $\alpha \in \mathbf{F}_{q^2} - \mathbf{F}_q$ is the smallest integer d such that $\alpha^d \in \mathbf{F}_q$. The test $\alpha^d \notin \mathbf{F}_q$ for all proper divisors d of q + 1 implies a suborder of subord(α) = q + 1.

***Lemma 4.*** The period normal basis { $\eta_0$, $\eta_1$ } is a primitive normal basis of $\mathbf{F}_{q^2}$ over $\mathbf{F}_q$ if and only if the integer $(1 - (-1)^{(p-1)/2}p)/4$ is a primitive root modulo q, and $\eta_i^d \notin \mathbf{F}_q$ for all $d \mid q + 1$.
Proof: The period polynomial $\psi_p(x) = x^2 + x + (1 - (-1)^{(p-1)/2}p)/4$ is a primitive polynomial in $\mathbf{F}_q[x]$ if and only if the norm $N(\eta_0) = (-1)^2\psi_p(0) = (1 - (-1)^{(p-1)/2}p)/4$ is a primitive root modulo q, and the element $\eta_0$ has suborder of subord($\eta_0$) = q + 1. ■

## 7.2 Dual Period Normal Bases of Quadratic Extensions

Let $a(x) = a_1x + a_0 = kx + (1 - k)$, where $a_i = \mathrm{Tr}(\eta_0\eta_i)$, and $b(x) = b_1x + b_0$. The polynomial a(x)

is the representation of the normal element $\eta_0$ in the polynomials algebra $\mathbf{F}_q[x]/(x^2 - 1)$, and b(x) is the inverse of a(x), it est, $a(x)b(x) \equiv 1 \bmod (x^2 - 1)$.

***Lemma 5.*** If the pair $\{ \eta_0, \eta_1 \}$ is a period normal basis of $\mathbf{F}_{q^2}$ over $\mathbf{F}_q$, then

$$\delta_0 = \frac{k+1}{2k+1}\eta_0 + \frac{k}{2k+1}\eta_1, \quad \delta_1 = \frac{k}{2k+1}\eta_0 + \frac{k+1}{2k+1}\eta_1$$

or

$$\delta_0 = \frac{k}{2k+1}\eta_0 + \frac{k+1}{2k+1}\eta_1, \quad \delta_1 = \frac{k+1}{2k+1}\eta_0 + \frac{k}{2k+1}\eta_1$$

is the dual basis depending on p = 2k + 1 = 4m +1 or 4m + 3.

Proof: Find the coefficients of b(x), and compute $b(x)\circ\eta_0 = b_0\eta_0 + b_1\eta_1$. ■

The self-dual normal basis theorem rules out the existence of self-dual normal bases for quadratic extensions, (more generally n = even). The result below gives a different way of confirming this for quadratic period normal bases. By definition a basis $\{ \eta_0, \eta_i \}$ is a self-dual basis if and only if the trace matrix

$$\left(Tr(\eta_i\eta_j)\right) = \begin{bmatrix} Tr(\eta_0\eta_0) & Tr(\eta_0\eta_1) \\ Tr(\eta_1\eta_0) & Tr(\eta_1\eta_1) \end{bmatrix} = \begin{bmatrix} 1 & 0 \\ 0 & 1 \end{bmatrix}.$$

***Lemma 6..*** (1) There are no self-dual period normal basis of any quadratic extension $\mathbf{F}_{q^2}$ of $\mathbf{F}_q$ for all odd prime q > 2.
(2) Let $q = 2^v$. Then the quadratic extension $\mathbf{F}_{q^2}$ of $\mathbf{F}_q$ has a self-dual period normal basis if and only if $k \equiv 2 \bmod 4$, and $v$ is odd.

Proof: ( 1 ) Let p = 2k + 1, k even. The basis $\{ \eta_0, \eta_1 \}$ is a self-dual basis if and only if

$$Tr(\eta_0\eta_0) = Tr(\eta_1\eta_1) = 1 \text{ and } Tr(\eta_0\eta_1) = Tr(\eta_1\eta_0) = 0.$$

But $Tr(\eta_0\eta_0) = Tr((0,0)\eta_0+(0,1)\eta_1+k) = -(0,0) - (0,1) + 2k = k + 1$, and $Tr(\eta_1\eta_1) = Tr(\eta_0\eta_0) = k + 1$. Likewise, $Tr(\eta_0\eta_1) = Tr(\eta_1\eta_0) = -k$. These imply that the integer k must be divisible by the characteristic of $\mathbf{F}_q$, or equivalently k = 2qc, c ≥ 1. Put p = 2k + 1 = 4cq + 1. Then

$$\left(\frac{p}{q}\right)\left(\frac{q}{p}\right) = \left(\frac{1}{q}\right)\left(\frac{q}{p}\right) = \left(\frac{q}{p}\right) = (-1)^{(p-1)(q-1)/4} = 1 \ .$$

This immediately precludes the existence of quadratic period normal bases since $\eta_0 = (-1+\sqrt{p})/2,\ \eta_1 = (-1-\sqrt{p})/2 \ \in \mathbf{F}_q$. The case of k odd is similar to this.
( 2 ) If $p \neq 8a + 1$ and $q = 2^v$, $v$ is odd, then q is not a quadratic residue modulo $p$, so $\eta_0, \eta_1 \notin \mathbf{F}_q$ and $Tr(\eta_1\eta_1) = Tr(\eta_0\eta_0) = k + 1 = 1$. ■

## 7.3 Multiplications in Quadratic Extensions

The linear expansions of the pairwise products $\eta_i\eta_j$ of the periods will be employed to generate the multiplication table T = [ $T_0$ $T_1$ ] for multiplication in a quadratic extension $\mathbf{F}_{q^2}$ of $\mathbf{F}_q$ with respect to the period normal basis { $\eta_0, \eta_1$ }. The 2×2 submatrices $T_i = ( t_{i,j,k} )$ are nonsingular and $T_i$ coincides with the change of basis matrix for $\eta_j \to \eta_i\eta_j$. In the quadratic case, the submatrix $T_1$ is equal to a row permutation of $T_0$.

The pairwise products $\eta_i\eta_j$ of the quadratic periods are:

$$\eta_0\eta_0 = \left(\sum_{x\in K_0}\omega^x\right)\left(\sum_{y\in K_0}\omega^y\right) = (0,0)\eta_0 + (0,1)\eta_1 + k\delta_i,$$

$$\eta_0\eta_1 = \left(\sum_{x\in K_0}\omega^x\right)\left(\sum_{y\in K_1}\omega^y\right) = (1,0)\eta_0 + (1,1)\eta_1 + k\delta_i$$

where $K_0 = \{ x^2 : 0 \neq x \in \mathbf{F}_q \}$, $K_1 = qK_0$, and $\delta_i = 1$ if $-1 \in K_i$, otherwise $\delta_i = 0$. And the other combination is $\eta_1\eta_1 = \sigma(\eta_0\eta_0) = (0, 1)\eta_0 + (0, 0)\eta_1$. Here the action of the nontrivial automorphism of $\mathbf{F}_{q^2}$ over $\mathbf{F}_q$ is defined by $\sigma(\eta_0) = \eta_1$ and $\sigma(\eta_1) = \eta_0$.
The cyclotomic numbers (i, j) which appear in the linear expansions of $\eta_i\eta_j$ are:

(1) (0, 1) = (1, 0) = (1, 1) = k/2, and (0, 0) = (k – 2)/2, if k is even, or
(2) (0, 1) = (1, 0) = (1, 1) = (k + 1)/2, and (0, 0) = (k – 1)/2, if k is odd.

Accordingly, depending on the parameter k, there are two cases.

***Case of $-1 \in K$ or k is even***. The linear expansions of the pairwise products of the periods are:

$$\eta_0\eta_0 = (0,\,0)\,\eta_0 + (0,1)\,\eta_1 + k = \frac{k-2}{2}\eta_0 + \frac{k}{2}\eta_1 + k = -\frac{k+2}{2}\eta_0 - \frac{k}{2}\eta_1,$$
$$\eta_0\eta_1 = (1,\,0)\,\eta_0 + (1,1)\,\eta_1 = \frac{k}{2}\eta_0 + \frac{k}{2}\eta_1.$$

***Case of –1 ∉ K and k is odd.*** The linear expansions of the pairwise products of the periods are:

(1) $\eta_0\eta_0 = (0,0)\eta_0 + (0,1)\eta_1 = \frac{k-1}{2}\eta_0 + \frac{k+1}{2}\eta_1,$

(2) $\eta_0\eta_1 = (1,0)\eta_0 + (1,1)\eta_1 + k = -\frac{k+1}{2}\eta_0 - \frac{k+1}{2}\eta_1 .$

***Lemma 7.*** The matrix $T_0$ is given by either

$$T_0 = \begin{bmatrix} \frac{-k-2}{2} & \frac{-k}{2} \\ \frac{k}{2} & \frac{k}{2} \end{bmatrix} \quad \text{or} \quad T_0 = \begin{bmatrix} \frac{k-1}{2} & \frac{k+1}{2} \\ -\frac{k+1}{2} & -\frac{k+1}{2} \end{bmatrix}$$

for k even or odd respectively.

Proof: The entries in $T_0 = (t_{i,j})$ are the coefficients in linear expansions $\eta_0\eta_i = t_{i,j}\eta_0 + t_{i,j}\eta_1$. ■

Another simple and direct method of deducing this table is to calculate $\eta_0\eta_0 = t_{0,j}\,\eta_0 + t_{0,j}\eta_1$, and $\eta_0\eta_1 = t_{1,j}\,\eta_0 + t_{1,j}\eta_1$ in the field of complex numbers $\mathbb{C}$, which also yields the cyclotomic numbers (i, j).

***Lemma 8.*** Let $\{\eta_0, \eta_1\}$ be a period normal basis of $\mathbf{F}_{q^2}$ over $\mathbf{F}_q$, and let $x = x_0\eta_0 + x_1\eta_1$ and $y = y_0\eta_0 + y_1\eta_1$, $x_i, y_i \in \mathbf{F}_q$. Then the product xy in $\mathbf{F}_{q^2}$ is given by either

$$\begin{aligned} xy &= (x_0\eta_0 + x_1\eta_1)(y_0\eta_0 + y_1\eta_1) \\ &= \left(\frac{-k-2}{2}x_0y_0 + \frac{k}{2}(x_0y_1 + x_1y_0) - \frac{k}{2}x_1y_1\right)\eta_0 \\ &\qquad + \left(\frac{-k}{2}x_0y_0 + \frac{k}{2}(x_0y_1 + x_1y_0) + \frac{-k-2}{2}x_1y_1\right)\eta_1 , \end{aligned}$$

if k is even; or if k is odd the product is given by

$$\begin{aligned} xy &= (x_0\eta_0 + x_1\eta_1)(y_0\eta_0 + y_1\eta_1) \\ &= \left(\frac{k-1}{2}x_0y_0 - \frac{k+1}{2}(x_0y_1 + x_1y_0) + \frac{k+1}{2}x_1y_1\right)\eta_0 \\ &\qquad + \left(\frac{k+1}{2}x_0y_0 - \frac{k+1}{2}(x_0y_1 + x_1y_0) + \frac{k-1}{2}x_1y_1\right)\eta_1 . \end{aligned}$$

Let the weight of the matrix $T_i$ be defined by $w(T_i) = \#\{ t_{i,j,k} \neq 0 \}$. An *optimal* normal basis has a complexity of C = $w(T_i) = 2n - 1$. All optimal normal bases are generated by normal period, confer [1, Mullin].

The optimal normal basis theorem claims the existence of an optimal normal basis of $\mathbf{F}_{q^n}$ over $\mathbf{F}_q$ of degree $n = [\mathbf{F}_{q^n} : \mathbf{F}_q]$, for the following parameters.

(1) The integer n + 1 is prime and q has order n modulo n + 1.
(2) The integer 2n + 1 is prime and q = 2 has order n or 2n modulo 2n + 1.
Case (1) above only includes only even n, for example, quadratic extension $\mathbf{F}_{q^2}$ of $\mathbf{F}_q$ for which $q \equiv 2 \bmod 3$, but not $q \equiv 1 \bmod 3$. Here is a proof for case (1) with n = 2. There is no claim about the optimal normal basis of $\mathbf{F}_{q^2}$ of $\mathbf{F}_q$ in characteristic char($\mathbf{F}_q$) = 3. The quadratic extension $\mathbf{F}_{q^2}$ of $\mathbf{F}_q$, $q \equiv 2 \bmod 3$, has the normal basis $\{ \eta, \eta^q \}$, where $\eta$ is a root of the period polynomial $\Psi_p(x) = x^2 + x + (1 - (-1)^{(p-1)/4}p)/4$, for example $\Psi_3(x) = x^2 + x + 1$.

***Theorem 9.*** Let q be a quadratic nonresidue modulo p. Then
( 1 ) The quadratic extension $\mathbf{F}_{q^2}$ of $\mathbf{F}_q$ has an optimal normal basis for all odd primes q = 3m + 2.
( 2 ) If $q = 3^v$, then the quadratic extension $\mathbf{F}_{q^2}$ of $\mathbf{F}_q$ has an optimal normal basis if either v = odd and $p \equiv 1, 5 \bmod 12$, or v = even and $p \equiv 7, 11 \bmod 12$.

Proof: ( 1 ) Let $q \neq 3^v$. First of all a quadratic extension $\mathbf{F}_{q^2}$ of $\mathbf{F}_q$ has a period normal basis if and only if q is a quadratic nonresidue modulo p, so the quadratic symbol ( q | p ) = –1. Second a normal basis is an optimal normal basis if and only if the submatrix $T_0$ of the multiplication table $T = [\, T_0\ T_1\,]$ of $\mathbf{F}_{q^2}$ over $\mathbf{F}_q$ has the discrete weight $w(T_0) = 2n - 1 = 3$. By inspection it is easy to find that the submatrix

$$T_0 = \begin{bmatrix} \frac{-k-2}{2} & \frac{-k}{2} \\ \frac{k}{2} & \frac{k}{2} \end{bmatrix} = \begin{bmatrix} 0 & -1 \\ 1 & 1 \end{bmatrix} \quad \text{or} \quad T_0 = \begin{bmatrix} \frac{k-1}{2} & \frac{k+1}{2} \\ -\frac{k+1}{2} & -\frac{k+1}{2} \end{bmatrix} = \begin{bmatrix} 0 & 1 \\ -1 & -1 \end{bmatrix}$$

has a complexity of $w(T_0) = 3$ if and only if $k = 2(cq - 1)$ or $2cq + 1$. This requires primes of the form $p = 2k + 1 = 4cq - 3$ or $4cq + 3$. Furthermore, the reciprocity law, gives

$$\left(\frac{p}{q}\right)\left(\frac{q}{p}\right) = \left(\frac{-3}{q}\right)(-1) = (-1)^{(p-1)(q-1)/4} = 1$$

for even $k = 2(cq - 1)$, which implies that q = 12a + 5 or 12a + 11 (since ( –3 | q ) = –1), or

$$\left(\frac{p}{q}\right)\left(\frac{q}{p}\right)=\left(\frac{3}{q}\right)(-1)=(-1)^{(p-1)(q-1)/4}=(-1)^{(q-1)/2}$$

for odd k = 2cq + 1, which implies that q = 12a + 5 or 12a + 11 (since ( 3 | q ) = –1). The case of q = 3c + 1 fails because ( ±3 | q ) = 1. The proof for q = $3^v$ uses the same technique. ■

Note: Use the relation

$$\left(\frac{3}{q}\right)=\begin{cases}-1 & \text{if } q\equiv \pm 5 \bmod 12\\ 1 & \text{if } q\equiv \pm 1 \bmod 12\end{cases}$$

to evaluate the previous quadratic symbols equations, and determine the forms of the prime q.

For instance, in characteristic q = 3, 5, 11, 17, 23, ..., the sequence

$$\mathbf{F}_{3^2}, \mathbf{F}_{5^2}, \mathbf{F}_{11^2}, \mathbf{F}_{17^2}, \mathbf{F}_{23^2}, \mathbf{F}_{29^2}, \ldots$$

of quadratic extensions have optimal normal bases over $\mathbf{F}_q$.

An optimal normal basis multiplier requires n(2n – 1) multiplications (AND cells if $\mathbf{F}_q = \mathbf{F}_2$). For example, the quadratic multiplier uses 6 multiplication cells. The maximal number of multiplications in $\mathbf{F}_q$ needed to compute the product xy in $\mathbf{F}_{q^2}$ is $n^3 = 8$. But a proper selection of the integer k in the prime kn +1 = 2k + 1 can reduce this number to the minimal n(2n – 1) = 6. And regrouping of the terms involving the product $\eta_0\eta_1$ and $\eta_1\eta_0$ further reduces it to 4.

***Example 10.*** Construct a normal basis of the quadratic extension $\mathbf{F}_9$ of $\mathbf{F}_3$ of low complexity.

Solution: Select a prime p = 2k + 1 which reduces the weight of the matrix $T_0$. This in turn reduces the numbers of terms $x_iy_j$ in the product to a minimal. The primes p = 2k + 1 = 3c + 2 are suitable because q is a quadratic nonresidue modulo p, and there exits a period normal basis $\eta_0=(-1+\sqrt{p})/2,\ \eta_1=(-1-\sqrt{p})/2 \in \mathbf{F}_{q^2}$ of complexity 2n – 1 = 3. Thus to compute the product, (using 2·8 + 1 = 17, with k = 8),

$$xy = (x_0y_1 + x_1y_0 - x_1y_1)\eta_0 + (-x_0y_0 + x_0y_1 + x_1y_0)\eta_1$$
$$= [x_0y_1 + x_1(y_0 - y_1)]\eta_0 + [x_1y_0 - x_0(y_0 - y_1)]\eta_1,$$

xy in $\mathbf{F}_9$, a total of n(2n – 1) = 6 multiplications in $\mathbf{F}_3$ are needed. But in practice this is simpler since $x_i$, $y_j \in \mathbf{F}_3$, and regrouping the terms reduces it to only four multiplications. Note that the period polynomial $\Psi_{17}(x) = x^2 + x - 1 \in \mathbf{F}_3[x]$ is irreducible.

## 7.4 Criteria For Cubic Nonresidues

The partition of the primes p = 3k + 1 is used throughout this section. These primes factor as p = (a + bω)(a – bω) in the quadratic domain ℤ[ω], where $\omega^2 + \omega + 1 = 0$. Note that $p = 3k + 1 = a^2 - ab + b^2$. In addition there is a pair of unique integers A and B such that $4p = A^2 + 27B^2$.

Information about the cubic reciduacity of the primes q modulo p is essential in the construction of cubic period bases. A few criteria for selecting pair of primes p, q which are parameters of cubic period normal bases are given here. These criteria are of general interest in various area of mathematics.

***Theorem 11.*** Let p = 3k + 1 be a prime, and $4p = A^2 + 27B^2$, A ≡ 1 mod 3. Then

( 1 ) $2^{(p-1)/3} \equiv \dfrac{A+9B}{A-9B} \bmod p$, if A ≡ B mod 4.

( 2 ) $3^{(p-1)/3} \equiv \dfrac{A+9B}{A-9B} \bmod p$, if B ≡ –1 mod 3.

Proof: See [1, Williams 1975.]. ■

***Theorem 12.*** Let p = 3k + 1, and $4p = A^2 + 27B^2$, A ≡ 1 mod 3. Then

( 1 ) q = 2 is a cubic residue modulo p if and only if 2 | A.
( 2 ) q = 3 is a cubic residue modulo p if and only if 9 | 3B.
( 3 ) q = 5 is a cubic residue modulo p if and only if 5 | A and 5 | 3B.
( 4 ) q = 7 is a cubic residue modulo p if and only if 7 | A and 7 | 3B.

Proof: Confer [1, Berndt et al, pp. 213-215.]. ■

These last results are useful tools in the implementation of search algorithms for sequences of primes for which q = 2, 3, 5, and 7 are cubic nonresidues.

***Example 13.*** ( 1 ) The integer q = 2 is a cubic nonresidue for all primes for all primes p = 3k + 1 for which $4p = (6a + 1)^2 + 27(6b + 1)^2$, a, b ≥ 1. For instance, p = 7, 19, ….

( 2 ) The integer q = 5 is a cubic nonresidue for all primes p = 3k + 1 for which $4p = (30a + 1)^2 + 27(30b + 1)^2$, a, b ≥ 1, et cetera..

***Theorem 14.*** Let $4p = A^2 + 27B^2$, A ≡ 1 mod 3, and let q ≠ 2, 3, p so that q ≡ ± 1 mod 6. Then q is a cubic nonresidue modulo p if and only if

$$\left(\frac{A+iB\sqrt{3}}{A-iB\sqrt{3}}\right)^{(q-\varepsilon)/6} \neq \left(\frac{p}{q}\right) \bmod q$$

where ε = ±1, and q ≡ ε mod 6.

For a fixed q ≠ 2, 3, p, the congruence parametizes a sequence of primes p, (or for a fixed p, the congruence parametizes a sequence of primes q ≠ 2, 3, p). In turns, this realizes a search algorithm for the primes p such that q is a cubic nonresidue modulo p. For example, let q = 7, then any prime p = 3k + 1 > 7 which satisfies the congruence

$$\frac{A+2B}{A-2B} \neq \left(\frac{p}{7}\right) \bmod 7,$$

is a cubic nonresidue modulo p.

## 7.5 Period Normal Bases of Cubic Extensions

The next simple period normal bases are those for cubic extensions $\mathbf{F}_{q^3}$ of $\mathbf{F}_q$. The group of automorphisms of the cubic extension $\mathbf{F}_{q^3}$ of $\mathbf{F}_q$ is the set Gal($\mathbf{F}_{q^3}/\mathbf{F}_q$) = { 1, σ, $\sigma^2$ } of linear maps $\sigma : \mathbf{F}_{q^3} \rightarrow \mathbf{F}_q$. The nontrivial automorphism is defined by $\sigma(x) = x^q$.

***Theorem 15.*** Let p = 3k + 1 be a prime, and $4p = A^2 + 27B^2$, A ≡ 1 mod 3, and let ω be a cube root of unity. Then the three conjugates of

$$\eta_i = \frac{-1}{3} + \frac{p^{1/3}}{3}\left[\omega\left(\frac{A+iB\sqrt{3}}{2}\right)^{1/3} + \omega^2\left(\frac{A-iB\sqrt{3}}{2}\right)^{1/3}\right]$$

form the set of cubic periods.

Proof: Determine the three roots $g_0$, $g_1$, and $g_2$ of the reduced period polynomial $\Theta_p(x) = x^3 - 3px - Ap$, and use the linear relation

$$\eta_i = \frac{-1+g_i}{3}.$$

The permutation π(0) = i, π(1) = j, and π(2) = k of { 0, 1, 2 } which identifies the index of the three periods $\eta_i$, $\eta_j$, $\eta_k$ is unknown. ■

***Theorem 16.*** Let p = 3k + 1 be a fixed prime, $4p = A^2 + 27B^2$, A ≡ 1 mod 3, and let $q \neq 3^v$ be a cubic nonresidue modulo p, this requires q ≢ 1 mod p. Then

(1) The triple { $\eta_0$, $\eta_1$, $\eta_2$ } is a period normal basis of the cubic extension $\mathbf{F}_{q^3}$ over $\mathbf{F}_q$.

(2) If the integer [p(A + 3) – 1]/27 is a primitive root modulo q, and $\eta_i^d \notin \mathbf{F}_q$ for all divisors d of $q^2 + q + 1$ then { $\eta_0$, $\eta_1$, $\eta_2$ } is a primitive period normal basis of $\mathbf{F}_{q^3}$ over $\mathbf{F}_q$.

Proof: (1) Let $\omega \neq 1$ be a pth root of unity in some extension $\mathbf{F}_{q^d}$ of $\mathbf{F}_q$, let g be a generator of the multiplicative group of $\mathbf{F}_p$, and write $q = g^{3a+1}$. Then

$$\eta_0 = \sum_{x\in K_0}\omega^x = \frac{1}{3}\sum_{x=1}^{p-1} e^{i2\pi x^3/p}, \quad \eta_1 = \sum_{x\in K_1}\omega^x = \frac{1}{3}\sum_{x=1}^{p-1} e^{i2\pi g x^3/p}, \quad \eta_2 = \sum_{x\in K_1}\omega^x = \frac{1}{3}\sum_{x=1}^{p-1} e^{i2\pi g^2 x^3/p}$$

where $K_0 = \{ 1, g^3, g^6, ..., g^{3(k-1)} \}$, $K_1 = \{ q, qg^3, qg^6, ..., qg^{3(k-1)} \}$, and $K_2 = \{ q^2, q^2g^3, q^2g^6, ..., q^2g^{3(k-1)} \}$. Since $q \not\equiv 1 \bmod p$, the element $\omega$ is not in $\mathbf{F}_q$, and similarly $\eta_0, \eta_1, \eta_2 \notin \mathbf{F}_q$. Moreover q is a noncubic residue the triple are conjugates and their traces are $Tr(\eta_i) = \eta_0 + \eta_1 + \eta_2 = -1$. To confirm the linear independence, the matrix determinant test requires that

$$\det\begin{bmatrix} \eta_0 & \eta_1 & \eta_2 \\ \eta_1 & \eta_2 & \eta_0 \\ \eta_2 & \eta_0 & \eta_1 \end{bmatrix} = 3\,\eta_0\,\eta_1\,\eta_2 - \left(\eta_0^3 + \eta_1^3 + \eta_2^3\right) \neq 0 .$$

From the coefficients of the cubic period polynomial $_p(x) = x^3 + x^2 - \frac{p-1}{9}x - \frac{p(A+3)-1}{27}$ it readily follows that

$$\eta_0\eta_1\eta_2 = \frac{p(A+3)-1}{27}, \text{ and } \eta_0^3 + \eta_1^3 + \eta_2^3 = \frac{p(A-6)-1}{9}.$$

Thus det(N) = p, and the set { $\eta_0$, $\eta_1$, $\eta_2$ } is a basis of $\mathbf{F}_{q^3}$ over $\mathbf{F}_q$.

(2) The polynomial $\psi_p(x) \in \mathbf{F}_q[x]$ is primitive if and only if the integer $N(\eta_i) = (-1)^3\psi_p(0)$ = [p(A + 3) – 1]/27 is a primitive root modulo q, ( $N(\eta_i)$ is the norm of $\eta_i$ ), and the element $\eta_i$ is of suborder $q^2 + q + 1$. ■

***Example 17.*** (1 ) The cubic periods { $\eta_0$, $\eta_1$, $\eta_2$ } is a basis of $\mathbf{F}_{2^3}$ over $\mathbf{F}_2$ for all primes p = 3k + 1, such that $4p = (6a+1)^2 + 27(6b+1)^2$, $A = 6a + 1 \equiv 1 \bmod 3$, since the period polynomial is $\psi_p(x) = x^3 + x^2 - 1 \in \mathbf{F}_2[x]$ is irreducible.
( 2 ) Let p = 7, A = 1, and $\eta_0, = \omega + \omega^{-1}$, $\omega^7 \neq 1$, (note $K_0 = \{ -1, 1 \}$ is the set of cubic residues in $\mathbf{F}_7$). The periods form a basis of $\mathbf{F}_{q^3}$ over $\mathbf{F}_q$ for all $q \not\equiv 2, 3 \bmod 7$ not a cubic residue. For these q, the period polynomial $\psi_p(x) = x^3 + x^2 - 2\cdot 3^{-1}x - 1 \in \mathbf{F}_q[x]$ is irreducible.

# 7.6 Dual Period Normal Bases of Cubic Extensions

Let $a(x) = a_2x^2 + a_1x + a_0,= -k(x^2 + x) + 2k + 1$, and $b(x) = b_2x^2 + b_1x + b_0$, where $a_i = Tr(\eta_0\eta_i)$, and $a(x)b(x) \equiv 1 \bmod x^3 - 1$. If the pair of elements $\{ \eta_0, \eta_1, \eta_2 \}$ is a normal basis of $\mathbf{F}_q(\eta_0)$ over $\mathbf{F}_q$, then the linear combination

$b_0\eta_0 + b_1\eta_1 + b_2\eta_2$

is the dual element of $\eta_0$. For the cubic period of type (k, 3), 3k + 1 prime.

***Lemma 18.*** Let $p = 3k + 1$, and suppose that $\{ \eta_0, \eta_1, \eta_2 \}$ is a the period normal basis of a cubic extension $\mathbf{F}_{q^3}$ over $\mathbf{F}_q$. Then

$$\delta_0 = \frac{k}{p}(\eta_2 + \eta_1) + \frac{2k+1}{p}\eta_0,\ \ \delta_1 = \frac{k}{p}(\eta_0 + \eta_2) + \frac{2k+1}{p}\eta_1,\ \ \delta_2 = \frac{k}{p}(\eta_0 + \eta_1) + \frac{2k+1}{p}\eta_2$$

is the dual basis. Moreover if either $q = 2^v$, or $k = aq$, then this is a self-dual basis.

Proof: Direct calculations. ■

***Lemma 19.*** Every cubic extension $\mathbf{F}_{q^3}$ of $\mathbf{F}_q$ has a self-dual period normal basis.

Proof: The self-dual basis the trace matrix equation is

$$\left(Tr(\eta_i\eta_j)\right) = \begin{bmatrix} Tr(\eta_0\eta_0) & Tr(\eta_0\eta_1) & Tr(\eta_0\eta_2) \\ Tr(\eta_1\eta_0) & Tr(\eta_1\eta_1) & Tr(\eta_1\eta_2) \\ Tr(\eta_2\eta_0) & Tr(\eta_2\eta_1) & Tr(\eta_2\eta_2 \end{bmatrix} = \begin{bmatrix} 3k+1 & -k & -k \\ -k & 3k+1 & -k \\ -k & -k & 3k+1 \end{bmatrix}.$$

Now take $k = 2qc$, $c \geq 1$ and select a prime $p = 6qc + 1$ such that q is a cubic nonresidue modulo p. ■

# 7.7 Multiplications in Cubic Extensions

The action of the generator of the group of automorphisms $\sigma : \mathbf{F}_{q^3} \to \mathbf{F}_q$ is defined by $\sigma(\eta_0) = \eta_1$, $\sigma(\eta_1) = \eta_2$, and $\sigma(\eta_2) = \eta_0$. This is utilized to determine all the products $\eta_i\eta_j$ and to complete the multiplication table for cubic extensions.

The relevant pairwise products for this case are:

$$(1)\ \eta_0\eta_0 = \left(\sum_{x\in K_0}\omega^x\right)\left(\sum_{y\in K_0}\omega^y\right) = (0,0)\eta_0 + (0,1)\eta_1 + (0,1)\eta_2 + k\,,$$

$$(2)\ \eta_0\eta_1 = \left(\sum_{x\in K_0}\omega^x\right)\left(\sum_{y\in K_1}\omega^y\right) = (1,0)\eta_0 + (1,1)\eta_1 + (1,1)\eta_2\,,$$

$$(3)\ \eta_0\eta_2 = \left(\sum_{x\in K_0}\omega^x\right)\left(\sum_{y\in K_2}\omega^y\right) = (2,0)\eta_0 + (2,1)\eta_1 + (2,1)\eta_2\,,$$

where $K_0$, $K_1$, and $K_2$ are the cosets of cubic residues and nonresidues in $\mathbf{F}_p$. The other pairwise products are computed via the automorphism:

$\eta_1\eta_1 = \sigma(\eta_0\eta_0) = (0, 2)\eta_0 + (0, 0)\eta_1 + (0, 1)\eta_2$,

$\eta_1\eta_2 = \sigma(\eta_0\eta_1) = (1, 2)\eta_0 + (1, 0)\eta_1 + (1, 1)\eta_2$, and

$\eta_2\eta_2 = \sigma(\eta_1\eta_1) = (0, 1)\eta_0 + (0, 2)\eta_1 + (0, 0)\eta_2$.

For the prime $p = 3k + 1$, $4p = A^2 + 27B^2$, $A \equiv 1 \bmod 3$, $3B \equiv (g^{(p-1)/3} - g^{2(p-1)/3}) \bmod 3$, and g primitive modulo p, the cyclotomic numbers (i, j) which appear in the linear expansions are

(1) $(0, 0) = (p - 8 + A)/9$

(2) $(0, 1) = (1, 0) = (2, 2) = (2p - 4 - A + 9B)/18$,

(3) $(0, 2) = (2, 0) = (1, 1) = (2p - 4 - A - 9B)/18$, and

(4) $(1, 2) = (2, 1) = (p + 1 + A)/9$.

The expressions for the products of the periods are

$$(1)\ \eta_0\eta_0 = \frac{A-2p-5}{9}\eta_0 + \frac{2-4p-A+9B}{18}\eta_1 + \frac{2-4p-A-9B}{18}\eta_2\ ,$$

$$(2)\ \eta_0\eta_1 = \frac{2p-4-A+9B}{18}\eta_0 + \frac{2p-4-A-9B}{18}\eta_1 + \frac{p+1+A}{9}\eta_2\,,$$

and

$$(3)\ \eta_0\eta_2 = \frac{2p-4-A-9B}{18}\eta_0 + \frac{p+1+A}{9}\eta_1 + \frac{2p-4-A+9B}{18}\eta_2\,.$$

The constraint $3B \equiv (g^{(p-1)/3} - g^{2(p-1)/3}) \bmod 3$ identifies an unique choice from the two possibilities ±B. These linear expansions generate the multiplication table for multiplication in a cubic extension $\mathbf{F}_{q^3}$ of $\mathbf{F}_q$ with respect to the period normal basis { $\eta_0$, $\eta_1$ $\eta_2$ }. The change of basis matrix $\eta_i \rightarrow \eta_0\eta_i$ coincides with the 3×3 submatrix $T_0$ of the multiplication table T = [$T_0$ $T_1$ $T_2$ ] of the cubic extension.

***Lemma 20.*** The multiplication matrix of $\mathbf{F}_{q^3}$ over $\mathbf{F}_q$ is generated by

$$T_0 = \begin{bmatrix} t_{0,0} & t_{0,1} & t_{0,2} \\ t_{1,0} & t_{1,1} & t_{1,2} \\ t_{2,0} & t_{2,1} & t_{2,2} \end{bmatrix} = \frac{1}{18}\begin{bmatrix} 2(A-2p-5) & 2-4p-A+9B & 2-4p-A-9B \\ 2p-4-A+9B & 2p-4-A-9B & 2(p+1+A) \\ 2p-4-A-9B & 2(p+1+A) & 2p-4-A+9B \end{bmatrix}.$$

Proof: These are the coefficients in the linear expansions $\eta_0\eta_i = t_{i,0}\eta_0 + t_{i,1}\eta_1 + t_{i,2}\eta_2$. ■

***Lemma 21.*** Let { $\eta_0$, $\eta_1$, $\eta_2$} be a period normal basis of $\mathbf{F}_{q^3}$ over $\mathbf{F}_q$. Then the product of the two elements $x = x_0\eta_0 + x_1\eta_1 + x_2\eta_2$ and $y = y_0\eta_0 + y_1\eta_1 + y_2\eta_2$, $x_i, y_i \in \mathbf{F}_q$, in $\mathbf{F}_{q^3}$ is given by

$xy = c_0\eta_0 + c_1\eta_1 + c_2\eta_2$,

where

$c_0 = x_0y_0t_{0,0} + x_2y_2t_{0,1} + x_1y_1t_{0,2} + (x_0y_1 + x_1y_0)t_{1,0} + (x_1y_2 + x_2y_1)t_{1,2}+ (x_0y_2 + x_2y_0)t_{2,0}$,

$c_1 = x_0y_0t_{0,1} + x_2y_2t_{0,2} + x_1y_1t_{0,0} + (x_0y_1 + x_1y_0)t_{1,1} + (x_1y_2 + x_2y_1)t_{1,0}+ (x_0y_2 + x_2y_0)t_{2,1}$,

$c_2 = x_0y_0t_{0,2} + x_2y_2t_{0,0} + x_1y_1t_{0,1} + (x_0y_1 + x_1y_0)t_{1,2} + (x_1y_2 + x_2y_1)t_{1,1}+ (x_0y_2 + x_2y_0)t_{2,2}$,

and $c_i \in \mathbf{F}_q$.

Proof: Direct calculations. ■

***Theorem 22.*** The cubic extension $\mathbf{F}_{q^3}$ of $\mathbf{F}_q$ has an optimal normal basis for all even primes power $q = 2^v$ iff gcd(v, 3) = 1.

A proper selection of the integer k in the prime kn + 1 = 3k + 1 can reduce the number of multiplications in $\mathbf{F}_q$ needed to compute the product xy in $\mathbf{F}_{q^3}$ to a minimal of 3n − 2 = 7 multiplications in characteristic char($\mathbf{F}_q$) ≠ 2 or to a minimal of 2n − 1 = 5 multiplications in characteristic char($\mathbf{F}_q$) = 2 with $q = 2^v$, gcd(v, 3) = 1.

***Example 23.*** Construct a normal basis of the cubic extension $\mathbf{F}_{5^3}$ of $\mathbf{F}_5$ of low complexity. Solution: Select a prime p = 3k + 1 which reduces the numbers of terms $x_iy_j$ in the product to a minimal. Any prime in the sequences of primes p = 6k + 1, such that $4p = (30a+1)^2$ +

$27(\pm(30b+1))^2$ is suitable because q = 5 is a cubic nonresidues and the multiplication submatrix is

$$T_0 = \begin{bmatrix} t_{0,0} & t_{0,1} & t_{0,2} \\ t_{1,0} & t_{1,1} & t_{1,2} \\ t_{2,0} & t_{2,1} & t_{2,2} \end{bmatrix} = 2\begin{bmatrix} -1 & -2-B & -2+B \\ -1-B & -1+B & 3 \\ -1+B & 3 & -1-B \end{bmatrix}.$$

There are two matrices depending on B =1 or B = –1.These are

$$T_0 = \begin{bmatrix} t_{0,0} & t_{0,1} & t_{0,2} \\ t_{1,0} & t_{1,1} & t_{1,2} \\ t_{2,0} & t_{2,1} & t_{2,2} \end{bmatrix} = \begin{bmatrix} -2 & -1 & -2 \\ 1 & 0 & 1 \\ 0 & 1 & 1 \end{bmatrix} \text{ or } \begin{bmatrix} -2 & -2 & -1 \\ 0 & 1 & 1 \\ 1 & 1 & 0 \end{bmatrix}.$$

The period polynomial $\psi_p(x) = x^3 + x^2 + x - 1 \in \mathbf{F}_5[x]$ is irreducible and invariant with respect to these primes.

The product $xy = c_0\eta_0 + c_1\eta_1 + c_2\eta_2$, where

$c_0 = x_0y_0t_{0,0} + x_2y_2t_{0,1} + x_1y_1t_{0,2} + (x_0y_1 + x_1y_0)t_{1,0} + (x_1y_2 + x_2y_1)t_{1,2} + (x_0y_2 + x_2y_0)t_{2,0}$,

$c_1 = x_0y_0t_{0,1} + x_2y_2t_{0,2} + x_1y_1t_{0,0} + (x_0y_1 + x_1y_0)t_{1,1} + (x_1y_2 + x_2y_1)t_{1,0} + (x_0y_2 + x_2y_0)t_{2,1}$,

$c_2 = x_0y_0t_{0,2} + x_2y_2t_{0,0} + x_1y_1t_{0,1} + (x_0y_1 + x_1y_0)t_{1,2} + (x_1y_2 + x_2y_1)t_{1,1} + (x_0y_2 + x_2y_0)t_{2,2}$,

is computed with just $n(3n - 2) + n^2 = 30$ multiplications in $\mathbf{F}_q$. But regrouping the terms reduces it to only 22 multiplications. Further reduction is possible depending on the characteristic char($\mathbf{F}_q$).

**Note**: The quadratic term $n^2$ in the total number of multiplications $n(3n - 2) + n^2$ for bases of low complexity $3n - 2$ accounts for the multiplications needed to obtain $x_iy_j$.

# CHAPTER 8

# Asymptotic Proofs For Primitive Polynomials

# On the Coefficients of Polynomials in Finite Fields

N. A. Carella, January 2000.

**ABSTRACT.** In this paper we will demonstrate the existence of primitive polynomials of degree $n$ with $k$ prescribed consecutive coefficients in the finite field $\mathbf{F}_q$ for all sufficiently large odd $q = p^v$ such that $k < p$, and $k < n(1/2 - \varepsilon) - 1$.

Key words: Finite fields, coefficients, primitive and primitive normal polynomials.

Mathematics Subject Classifications: 11T06.

## 1 Introduction

Several authors have proved the existence of primitive polynomials $f(x) = x^n + c_1x^{n-1} + c_2x^{n-2} + \cdots + c_{n-1}x + c_n \in \mathbf{F}_q[x]$ of prescribed traces $c_1 \in \mathbf{F}_q$, see [2] and [7]. There are exceptions, for instance, the quadratic primitive polynomials must be trinomials. The existence of primitive polynomials with one prescribed coefficient $c_i \in \mathbf{F}_q$, $0 < i < n$, with many exceptions, has been conjectured in [5], this conjecture is valid for irreducible polynomials, see [3]; a limited version of this conjecture has recently been proved in [13]. In addition, a formula that enumerates the number of irreducible polynomials of degree $n$ in $\mathbf{F}_2[x]$ with $c_1 = c_{n-1} = 1$, and $c_n = 1$ appears in [9], and [4] has proven the existence of primitive polynomials with two prescribed consecutive coefficients $c_1, c_2 \in \mathbf{F}_q$ for all sufficiently large odd $q$, see Theorem 7. In this paper we will prove the existence of primitive polynomials with $k$ prescribed consecutive coefficients $c_1, c_2, ..., c_k \in \mathbf{F}_q$ provided certain restrictions are met. The precise statement is given below.

***Theorem 1.*** (*Coefficients Theorem*) Let $q = p^v$, and $\varepsilon > 0$. If $k < p$, and $k < n(1/2 - \varepsilon) - 1$, then there exist a constant $q_\varepsilon$ and a primitive polynomials $f(x) = x^n + c_1x^{n-1} + c_2x^{n-2} + \cdots + c_{n-1}x + c_n \in \mathbf{F}_q[x]$ with $k$ prescribed consecutive coefficients $c_1, c_2, ..., c_k \in \mathbf{F}_q$ for all odd $q \geq q_\varepsilon$.

Since the reciprocal polynomial $f^*(x) = x^nf(1/x) = x^n + a_{n-1}x^{n-1} + a_{n-2}x^{n-2} + \cdots + a_1x + a_0$, $a_i = c_{n-i} / c_n$, of a primitive polynomial $f(x)$ is also a primitive polynomial, the k prescribed consecutive coefficients can also be chosen in reversed order.

In section 2 we will post a few results needed in the proof, and section 3 covers the proof of the Coefficients Theorem.

## 2 Auxiliary Results

Let $\mathbf{F}_{q^n}$ be a finite field extension of $\mathbf{F}_q$ of degree $n = [\mathbf{F}_{q^n} : \mathbf{F}_q]$, and let $z_1, z_2, ..., z_n \in \mathbf{F}_{q^n}$ be the roots of the polynomial

$$f(x) = (x - z_1)(x - z_2) \cdots (x - z_n) = x^n + c_1 x^{n-1} + c_2 x^{n-2} + \cdots + c_{n-1}x + c_n \in \mathbf{F}_q[x]. \quad (1)$$

The coefficients $c_i$ of $f(x)$ are functions of $z_1, z_2, ..., z_n$ given by the formulae

$$c_1 = -\sum_{i=1}^{n} z_i, \quad c_2 = \sum_{1 \le v_1 < v_2 \le n} z_{v_1} z_{v_2}, \quad ..., \quad c_i = (-1)^i \sum_{1 \le v_1 < v_2 < \cdots < v_i \le n} z_{v_1} z_{v_2} \cdots z_{v_i}, \quad c_n = (-1)^n z_1 z_2 \cdots z_n.$$
(2)

These functions are related to the symmetric functions $\sigma_i$ by the relation $c_i = (-1)^i \sigma_i$. The associated power sums are defined by

$$w_i = \sum_{j=1}^{n} z_j^i = z_1^i + z_2^i + \cdots + z_n^i, \quad i \ge 1. \quad (3)$$

Due to the structure of the roots of polynomials in characteristic $\neq 0$ the power sums have the alternate form

$$w_i = \sum_{j=1}^{n} z_j^i = \sum_{j=1}^{n} \alpha^{iq^j} = Tr(\alpha^i), \quad i \ge 1, \quad (4)$$

where $Tr(x)$ is the trace of $x$, and $\alpha \in \mathbf{F}_{q^n}$ is an element of degree $n$. These formulae allow one to compute or solve a problem about the coefficients of polynomials via power functions. Fast algorithms for computing the coefficients $\{ c_i \}$ from the power functions $\{ w_i \}$ (or conversely) in $O(\log n)$ steps are discussed in [11].

Substituting the power sum $w_i = w_i(\alpha) = Tr(\alpha^i)$ in Newton's identities yields the corresponding *coefficients identities* for finite fields. The coefficients identities are given by

$$(1) \; Tr(\alpha^k) = (-1)^{k+1} k c_k + \sum_{i=1}^{k-i} (-1)^i c_i Tr(\alpha^{k-i}), \qquad \text{for } 1 \le k \le n, \quad (5)$$

$$(2) \; Tr(\alpha^k) = \sum_{i=1}^{n-1} (-1)^i c_i Tr(\alpha^{k-i}), \qquad \text{for } n < k.$$

Since the reduction of formula (5-1) modulo $p$ loses some information, it is convenient to consider formula (5-1) as two separate cases: $n < p$, and $p \le n$.

**Case of $n < p$.** In this case the reduction of formula (5-1) modulo is the same as in characteristic 0, so it works as usual – all the coefficients can be computed. In addition, we can use identity (5-1) to rewrite the coefficient functions

$$c_i = c_i(w_1,w_2,...,w_i) = c_i(Tr(\alpha),Tr(\alpha^2),...,Tr(\alpha^i)) \tag{6}$$

as function of the power functions $w_1 = Tr(\alpha)$, $w_2 = Tr(\alpha^2)$, ..., $w_i = Tr(\alpha^i)$, or equivalently as function of the traces of the elements $\alpha$, $\alpha^2$, ..., $\alpha^i$. A few of the coefficients are given below.

(1) $c_1 = -w_1$, (7)

(2) $c_2 = \frac{1}{2!}\left(w_1^2 - w_2\right)$,

(3) $c_3 = \frac{1}{3!}\left(w_1^3 - w_1 w_2 + w_3\right)$,

(4) $c_4 = \frac{1}{4!}\left(2w_1^4 + 4w_1^2 w_2 + 4w_1 w_3 - 3w_2^2 - 6w_4\right)$, ….,

($k$) $c_k = \frac{1}{k!}\sum_{i} a_i w_1^{i_1} w_2^{i_2} \cdots w_k^{i_k}$ ,

where $1 \le k \le n$, and $\boldsymbol{i} = (i_1, i_2, \ldots, i_k)$.

In characteristic $p > 2$, a direct calculation of the second coefficient (7-2) using (2), namely,

$$c_2 = \sum_{0 \le i < j < n} \alpha^{q^i + q^j} = \frac{1}{2}\left(Tr(\alpha)^2 - Tr(\alpha^2)\right), \tag{8}$$

can be accomplished in a few steps, (the right hand side of this formula is well known in matrix theory). But a direct calculation of any of the other coefficients $c_3$, $c_4$, ..., $c_{n-1}$, appears to be a lengthy and times consuming task. Even the third coefficient

$$c_3 = \sum_{0 \le i < j < k < n} \alpha^{q^i + q^j + q^k} = \frac{1}{6}\left(2Tr(\alpha^3) - Tr(\alpha)Tr(\alpha^2) - Tr(\alpha)^3\right), \tag{9}$$

in characteristic $p > 3$, is difficult to calculate directly. However, using identity (5-1), this task becomes a simple algebraic manipulation. In summary, if the coefficients $c_1$, $c_2$, ..., $c_k$ are prescribed, then so are the powers functions $w_1$, $w_2$, ..., $w_k$, and conversely.

**Case of $p \le k$.** If $q = p^v$ and $p \le n$, the first $p - 1$ coefficients $c_1$, $c_2$, ..., $c_{p-1}$ of any polynomial of degree $n$ in $\mathbf{F}_q[x]$ can be computed using formula (5-1). But the remaining coefficients $c_p$, $c_{p+1}$,

..., $c_n$ might not be computable since the reduction of identity (5-1) modulo $p$ loses information. This case also has many exceptions, and will not be considered here.

Several other results which will be used later on are also included here, the proofs can be found in [1] or similar sources.

***Lemma 2.*** Let $N >>1$, and let $p$ runs through the prime divisors of $N$. Then

$$\frac{e^{-\gamma}}{\log\log N}\left(1+O(1/\log\log N)\right) \le \frac{\varphi(N)}{N} = \prod_{p|N}\left(1-1/p\right) \tag{10}$$

The constant $e^{-\gamma} = .5615...$ is important whenever $N$ is a small integer, but for sufficiently large $N$, we can use the simpler asymptotic formula $1/\log\log N \le \varphi(N)/N$. It is clear that the inequalities $0 < \varphi(N)/N < 1$ holds, and $\varphi(N)/N \le 1/2$ if $N$ is even. If $N = q^n - 1$, q a prime power, then the ratio $P_1 = \varphi(N)/N$ is interpreted as the probability of primitive elements in $\mathbf{F}_{q^n}$.

The omega function is defined by $\omega(N) = \#\{$ distinct prime divisors of $N$ $\}$. The inequality $2^{\omega(N)} < 2^{2\log N/\log\log N}$ is quite useful in analysis, but we will utilize the following.

***Lemma 3.*** If $\varepsilon > 0$, and $N \ge 1$, then exists a constant $C_\varepsilon > 0$ such that $2^{\omega(N)} < C_\varepsilon N^\varepsilon$.

The construction of the existence equation (15) of primitive polynomials with prescribed consecutive coefficients is based on the two exponential sums

$$\text{(1)}\quad \frac{\varphi(q^n-1)}{q^n-1}\sum_{d|q^n-1}\frac{\mu(d)}{\varphi(d)}\sum_{ord(\chi)=d}\chi(\xi) \quad \text{and} \tag{11}$$

$$\text{(2)}\quad \frac{1}{q}\sum_{\alpha\in\mathbf{F}_{q^n}}\sum_{x\in\mathbf{F}_q}\psi(x(Tr(f(\alpha))-a))$$

where $\chi$ is a multiplicative character of order $d = ord(\chi)$, $d \mid q^n - 1$, on the finite field $\mathbf{F}_{q^n}$ of $q^n$ elements, and $\psi \ne 1$ is a nontrivial additive character on $\mathbf{F}_q$ respectively. The first sum is the characteristic function of primitive elements in $\mathbf{F}_{q^n}$, and the second sum gives the cardinality of the solution set $\{ \alpha \in \mathbf{F}_{q^n} : Tr(f(\alpha)) - a = 0 \}$ of the equation $Tr(f(x)) - a = 0$ in $\mathbf{F}_{q^n}$, where $f(x)$ is a function on $\mathbf{F}_{q^n}$, and $a \in \mathbf{F}_q$ is a constant. The product of these two sums gives the number of primitive elements solutions of the equation $Tr(f(x)) - a = 0$ in $\mathbf{F}_{q^n}$.

***Lemma 4.*** *([14])* Let $\chi$, $\psi \ne 1$ be a pair of nontrivial multiplicative and additive characters, and let $f(x), g(x) \in \mathbf{F}_q[x]$ be polynomials of degrees $k = \deg(f)$ and $m = \deg(g)$ respectively. Then

$$\left| \sum_{x \in \mathbf{F}_{q^n}} \chi(f(x))\psi(Tr(g(x))) \right| \leq (k+m-1)q^{n/2}. \tag{12}$$

## 3 Proof of the Coefficients Theorem

The total number $\varphi(q^n-1)/n$ of primitive polynomials of degree $n$ in $\mathbf{F}_q[x]$ is an upper limit of the cardinality of any proper subset of primitive polynomials. Specifically, we are interested in the subset $\Im$ = { primitive polynomials of degree $n$ with $k$ prescribed coefficients } of cardinality $\#\Im \geq q^k$. Clearly the absolute maximal number $k$ of independent coefficients satisfies the inequalities

$$q^k \leq \frac{\varphi(q^n-1)}{n} \leq \frac{q^n-1}{n}P_1 \tag{13}$$

This count also includes $c_n$ as an independent coefficient. Consequently, $k < n - \log n/\log q$, which means that at least one coefficient among $c_1, c_2, ..., c_n$ remains dependent on the other coefficients. For small $q < n$, there are fewer linearly independent coefficients since

$$w_q = Tr(\alpha^q) = Tr(\alpha) = w_1, \quad w_{2q} = Tr(\alpha^{2q}) = Tr(\alpha^2) = w_2, \text{ et cetera.}$$

***Example 5.*** (1) The quadratic primitive polynomials $x^2 + c_1x + c_2 \in \mathbf{F}_q[x]$ can be selected with prescribed traces $0 \neq c_1 \in \mathbf{F}_q$, or norm $c_2$ for all $q$. But since $k \leq n - 1 = 1$, the two coefficients $0 \neq c_1$ and $c_2$ cannot be selected independently, in fact, $\varphi(q^2-1)/2 < (q-1)\varphi(q-1)$, for all $q > 3$.

(2) The cubic primitive polynomials $x^3 + c_1x^2 + c_2x + c_3 \in \mathbf{F}_q[x]$ can be selected with at most two independent coefficients since $k \leq n - 1 = 2$, in fact $\varphi(q^3-1)/3 < q^2\varphi(q-1)$ for all $q \geq 2$. For instance, $x^3 + c_1x^2 + c_2x + c_3 = x^3 + x^2 + 1$ or $x^3 + x + 1 \in \mathbf{F}_2[x]$ is a primitive polynomial with one independent coefficient $c_1$ or $c_2$. In characteristic $p > 2$, some cubic primitive polynomials can be selected with two independent coefficients $c_1$ and $c_3$, for instance, $x^3 + c_1x^2 + c_2x + c_3 \in \mathbf{F}_5[x]$.

The key idea in the proof of Theorem 1 is to show that the system of polynomials equations

$$Tr(x) = a_1,\ Tr(x^2) = a_2,\ ...,\ Tr(x^k) = a_k, \tag{14}$$

$a_i \in \mathbf{F}_q$ constants, has at least one primitive element solution $x = \xi \in \mathbf{F}_{q^n}$. This in turns implies the existence of at least one primitive polynomial $f(x) = x^n + c_1x^{n-1} + c_2x^{n-1} + \cdots + c_{n-1}x + c_n$ with $k$ prescribed consecutive coefficients. The constants $a_1, a_2, ..., a_k$ indirectly prescribes the $k$ coefficients, viz,

$$c_1 = -Tr(\xi),\ c_2 = c_2(Tr(\xi), Tr(\xi^2)),\ ...,\ c_k = c_k(Tr(\xi), Tr(\xi^2), ..., Tr(\xi^k)).$$

This system of polynomials equations, which is a direct consequence of (4), generalize the existence equation (15) used in the investigation of the distribution of the values of $c_1, c_2, ..., c_k$ to all the coefficients possible. The analysis of the case $k = 1$ appears in [2], and [7], et cetera, and the more recent analysis of the case $k = 2$ is given in [4].

The technique we use to estimate the number of primitive solutions of this system of equations has some similarity to a method used to estimate the covering radii of algebraic codes, see [6], and [10].

The number of primitive polynomials with k prescribed consecutive coefficients is given by (existence equation for $k < p$)

$$N(q^n, a_i) = \sum_{0 \neq \xi \in \mathbf{F}_{q^n}} \left( \frac{\varphi(q^n - 1)}{q^n - 1} \sum_{d \mid q^n - 1} \frac{\mu(d)}{\varphi(d)} \sum_{ord(\chi) = d} \chi(\xi) \right) \left( \prod_{i=1}^{k} \frac{1}{q} \sum_{x_i \in \mathbf{F}_q} \psi(x_i (Tr(f(\xi^i)) - a^i)) \right) \quad (15)$$

This is a $k$-fold product of formula (11-1) and (11-2). The usual strategy of dealing with this type of sum is to decompose it into several sums, and then compute a lower estimate, consult [4], [7], [8], [10], and [12] for background details and other references on this type of analysis. To obtain a suitable decomposition of equation (15), rewrite it as

$$N(q^n, a_1, ..., a_k) = q^k P_1 - \sum_{d \mid q^n - 1} \frac{\mu(d)}{\varphi(d)} \sum_{ord(\chi) = d} \sum_{x_1, x_2, ..., x_k \in \mathbf{F}_q} \sum_{0 \neq \xi \in \mathbf{F}_{q^n}} \chi(\xi) \psi(Tr(x_1 \xi + x_2 \xi^2 + \cdots + x_k \xi^k) - a_1 x_1 + a_2 x_2 + \cdots + a_k x_k))$$

The decomposition that employed is given by

$$N(q^n, a_1, ..., a_k) = q^k P_1 (S_0 + S_1 + S_2 + S_3) \quad (17)$$

The first sum is $S_0 = q^n - 1$, the second sum is

$$S_1 = \sum_{(x_1, x_2, ..., x_k) \neq (0,0,...,0),\ 0 \neq \xi \in \mathbf{F}_{q^n}} \sum \psi(Tr(x_1 \xi + x_2 \xi^2 + \cdots + x_k \xi^k) - a_1 x_1 + a_2 x_2 + \cdots + a_k x_k))$$

(18)

the third sum is $S_2 = 0$, and the fourth sum is

$$S_3 = \sum_{1 \neq d \mid q^n - 1} \frac{\mu(d)}{\varphi(d)} \sum_{ord(\chi) = d,} \sum_{(x_1, x_2, ..., x_k) \neq (0,0,...0),} \sum_{0 \neq \xi \in \mathbf{F}_{q^n}} \chi(\xi) \psi(Tr(x_1 \xi + x_2 \xi^2 + \cdots + x_k \xi^k) - a_1 x_1 + a_2 x_2 + \cdots + a_k x_k))$$

These sums are determined by the restrictions of the inner term of (16) to the subsets of vectors

$V_0 = \{ (d=1,0,0,...,0,\xi \neq 0) \}$,
$V_1 = \{ (d=1,x_1,x_2,...,x_k,\xi \neq 0) : (x_1,x_2,...,x_k) \neq (0,...,0),\ x_i \in \mathbf{F}_q \}$,
$V_2 = \{ (d \neq 1,0,0,...,0,\xi \neq 0) : d \mid q^n - 1 \}$, and
$V_3 = \{ (d \neq 1,x_1,x_2,...,x_k,\xi \neq 0) : d \mid q^n - 1,\ (x_1,x_2,...,x_k) \neq (0,...,0),\ x_i \in \mathbf{F}_q \}$ respectively.

To compute a lower estimate of $N(q^n,a_1,...,a_k)$, based on this decomposition, oneneeds he upper estimates of the sums $S_1$ and $S_3$.

***Lemma 6.*** (1) $| S_1 | \leq (q^k - 1)(k - 1)q^{n/2}$, and (2) $| S_3 | \leq (2^\omega - 1)(q^k - 1)(kq^{n/2})$, where $\omega = \omega(q^n-1)$ is the number of distinct prime divisors of $q^n - 1$.

*Proof*: Take the absolute value of equations (18) and (19) and apply Lemma 5. ■

***Theorem 1.*** (*Coefficients Theorem*) Let $q = p^v$, and $\varepsilon > 0$. If $k < p$, and $k < n(1/2 - \varepsilon) - 1$, then there exist a constant $q_\varepsilon$ and a primitive polynomials $f(x) = x^n + c_1x^{n-1} + c_2x^{n-2} + \cdots + c_{n-1}x + c_n \in \mathbf{F}_q[x]$ with $k$ prescribed consecutive coefficients $c_1, c_2, ..., c_k \in \mathbf{F}_q$ for all odd $q \geq q_\varepsilon$.

*Proof* : The condition $k < p$ allows us to use (15). Replace $S_0 = q^n - 1$, $S_2 = 0$, and the minimal estimate $S_1, S_3 \geq -(2^\omega - 1)\cdot q^k\cdot(kq^{n/2})$ in (17) to obtain

$$N(q^n, a_1,..., a_k) \geq q^k P_1\left(q^n - 1 - 2(2^{\omega(q^n-1)} - 1)kq^{n/2+k}\right) \tag{20}$$

Substitute $2^{\omega(N)} = C_\varepsilon N^\varepsilon$, (see Lemma 3), and rearrange to obtain

$$N(q^n, a_1,..., a_k) \geq q^{k-k} P_1\left(1 - 2C_\varepsilon q^{k+1-n(1/2-\varepsilon)}\right) \tag{21}$$

Since $0 < P_1 < 1$, and by assumption $k < n(1/2 - \varepsilon) - 1$, it immediately follows that there exists a constant $q_\varepsilon$ such that $N(q^n,a_1,...,a_k) \geq 1$ for all odd $q \geq q_\varepsilon$. ■

In the (asymptotic) analysis over large finite fields $\mathbf{F}_q$ given above, the estimate $| S_i | \leq (2^\omega - 1)(kq^{n/2+k})$ is adequate. However, to investigate the existence of primitive polynomials of small degrees with $k$ prescribed coefficients over small finite fields $\mathbf{F}_q$, detailed analysis of all the individual sums $S_i$ can improve the upper estimates $| S_i | \leq (2^\omega - 1)(kq^{n/2+k})$, this in turns reduces the number of cases to be checked. The detailed analysis on the case $k = 2$ appears in [4], which concludes in the following, (valid for odd $q$ only).

***Theorem 7.*** (*[4]*) Suppose $n \geq 7$, then there exists a primitive polynomial in $\mathbf{F}_q[x]$ of degree $n$ with the first and second coefficients prescribed in advance.

**Case of k ≤ 4.** This case corresponds to primitive polynomials with $k \leq 4$ prescribed consecutive coefficients $c_1, c_2, c_3, c_4 \in \mathbf{F}_q$. Since this case requires $4 < p$, let $q = p^v$, with $p \geq 5$, and let $\varepsilon = \log 2/\log 29 = 0.2058468324604...$ . It is not difficult to verify that this choice of constant $\varepsilon$ leads to $2^{\omega(N)} < C_\varepsilon N^\varepsilon = 10N^{.205847}$. Moreover, one has

$$q^{-(\varepsilon \log N + \log C_\varepsilon)/\log q} < \left(\frac{1}{2}\right)^{\omega(N)} \frac{\varphi(N)}{N} \tag{22}$$

for all $q > 1$ and $N > 0$. Now equation (20) becomes

$$N(q^n, a_1, \ldots, a_k) \geq q^{n(1-\varepsilon)-k-\log C_\varepsilon / \log q}\left(1 - q^{k+1-n(1/2-\varepsilon)+(\log 2+\log C_\varepsilon)/\log q}\right) \tag{23}$$

Using equation (23) and the stated parameters we determine the following two cases.

(1) There are primitive polynomials of degree $n > 10$ with $k = 3$ prescribed consecutive coefficients $c_1, c_2, c_3 \in \mathbf{F}_q$ for all $q = p^v$, with $p \geq 5$.

(2) There are primitive polynomials of degree $n > 11$ with $k = 4$ prescribed consecutive coefficients $c_1, c_2, c_3, c_4 \in \mathbf{F}_q$ for all $q = p^v$, with $p \geq 5$.

These estimates assume that $q = p = 5$, but for larger $q = p^v >> 5$, the parameters $\varepsilon$ and $C_\varepsilon$ can be readjusted to reduce the minimal $n$.

January 2000, Rev. 4.

**On the Coefficients of Primitive Normal Polynomials**

N. A. Carella, October 2005.

***Abstract***: The previous paper of the Winter 1998/1999 proved the existence of primitive polynomials, and primitive normal polynomials of degree $n$ with $k$ prescribed coefficients in the finite field $\mathbf{F}_q$ for all sufficiently large powers $q = p^v$ such that $k < p$, and $k < n(1 - \varepsilon) - 1$, $\varepsilon > 0$. This paper presents a longer version of the result on primitive normal polynomials.

Keywords: Polynomials in finite fields, primitive polynomials, primitive normal polynomials.

Mathematics subject classification: 11T06.

## 1 Introduction

Let $q = p^v$ be a sufficiently large prime power, and let $k$ and $n$ be a pair of integers such that $k < p$, and $k < n(1 - \varepsilon) - 1$, $\varepsilon > 0$. This work is a longer version of the analysis on the distribution of the coefficients of primitive polynomials and primitive normal polynomials in finite fields started in [4]. The previous results proved the existence of primitive polynomials and to primitive normal polynomials $f(x) = x^n + a_1x^{n-1} + \cdots + a_{n-1}x + a_n \in \mathbf{F}_q[x]$ with $k$ prescribed consecutive coefficients $a_1 \neq 0, a_2, \ldots, a_k \in \mathbf{F}_q$ in finite fields of odd characteristic $p$. A *primitive polynomial* has roots of multiplicative order $q^n - 1$; and a *normal polynomial* has roots of additive order $x^n - 1$, which is the same as having nonzero trace and linearly independent roots.

***Theorem* 1.** (Extended Coefficient Theorem) Let $q = p^v$, $p$ prime, and let $k$ and $n$ be a pair of integers such that $k < p$, and $k < n(1 - \varepsilon) - 1$, $\varepsilon > 0$. Then there exists a primitive normal polynomial $f(x) \in \mathbf{F}_q[x]$ with $k$ prescribed coefficients for all odd $q \geq q_0$.

The case $k = 1$ was considered in [3], it showed that the existence of primitive normal polynomials of prescribed trace $a_1 \neq 0$ for all pairs $(n, q)$ but a finite numbers of exceptions.

A *primitive normal basis* is a basis $\{\eta, \eta^q, \eta^{q^2}, \ldots, \eta^{q^{n-1}}\}$ of the vector space $\mathbf{F}_{q^n}$ over $\mathbf{F}_q$ generated by a primitive normal element $\eta \in \mathbf{F}_{q^n}$. The asymptotic proof of the *Primitive Normal Basis Theorem* was first established by both [5], and [8]. And the complete version for all pair $n$, $q$ was established by [15].

***Theorem* 2.** (*Primitive Normal Basis Theorem*) Let $\mathbf{F}_{q^n}$ be an $n$-degree extension of $\mathbf{F}_q$. Then $\mathbf{F}_{q^n}$ has a primitive normal basis over $\mathbf{F}_q$.

The next result is a refinement of the Theorem 2, it calls for primitive normal elements of arbitrary traces.

***Theorem* 3.** (*Primitive Normal Basis Theorem Of Arbitrary Trace*) For every $a \neq 0$ in $\mathbf{F}_q$, there exists a primitive normal element in $\mathbf{F}_{q^n}$ of trace $a$.

The refinement was proposed as a conjecture in [17], and its first asymptotic proof was completed in [3]. About two years later it was extended to all pairs $n$, $q$ in [7]. A much more general extension remains an open problem.

***Conjecture* 4.** (*Morgan-Mullen* 1996) Let $q$ be a prime power and let $n > 3$ be an integer. Then there exists a completely normal primitive polynomial of degree $n$ over $\mathbf{F}_q$.

## 2 Auxiliary Results

This section introduces several concepts used to study the distribution of elements and the coefficients of polynomials in finite fields.

*Characteristic Functions*

A characteristic function encapsulates certain properties of a subset of elements of $\mathbf{F}_{q^n}$. It effectively filters out those elements that do not satisfy the constraints. The equation of a characteristic function is of the form

$$C(\alpha) = \begin{cases} 1 & \text{if } \alpha \text{ satisfies the properties,} \\ 0 & \text{otherwise,} \end{cases} \tag{1}$$

for all elements $\alpha \in \mathbf{F}_{q^n}$.

*The Characteristic Function of Primitive Elements*

Let $\chi$ be a multiplicative character of order $d = ord(\chi)$, $d \mid q^n - 1$, on $\mathbf{F}_{q^n}$, see [16]. The characteristic function of primitive elements in $\mathbf{F}_{q^n}$ is defined by

$$C_P(\alpha) = \frac{\varphi(q^n - 1)}{q^n - 1} \sum_{d \mid q^n - 1} \frac{\mu(d)}{\varphi(d)} \sum_{ord(\chi) = d} \chi(\alpha), \tag{2}$$

where $\alpha \in \mathbf{F}_{q^n}$, and the arithmetic functions $\mu$ and $\varphi$ are the Mobius and Euler functions on the ring of integers $\mathbb{Z}$ respectively, see [1], and [18]. A product version of this formula has the shape

$$C_P(\alpha) = \frac{\varphi(q^n - 1)}{q^n - 1} \prod_{p \mid q^n - 1} \left( 1 - \frac{1}{p-1} \sum_{ord(\chi) = p} \chi(\alpha) \right) , \tag{3}$$

where $p$ runs through the prime divisors of $q^n - 1$. The transformation used to obtain this form is straightforward and appears throughout the literature, other forms of this characteristic function are also possible.

The function $C_P(\alpha)$ is one of the basic tools used in the investigation of the distribution of primitive elements in finite fields. Typical applications are illustrated in [14], [20], [21], etc.

*The Characteristic Function of Normal Elements*

Let $\psi$ be an additive character of order $f(x) = Ord(\psi) \mid x^n - 1$ on $\mathbf{F}_{q^n}$, see [5], and [15]. The characteristic function of normal elements in $\mathbf{F}_{q^n}$ is defined by

$$C_N(\alpha) = \frac{\Phi(x^n - 1)}{q^n} \sum_{d(x) | x^n - 1} \frac{M(d(x))}{\Phi(d(x))} \sum_{Ord(\psi) = d(x)} \psi([(x^n - 1)/d(x)] \circ \alpha), \tag{4}$$

where $\alpha \in \mathbf{F}_{q^n}$, and the arithmetic functions M and $\Omega$ are the Mobius and Euler functions on the ring of polynomials $\mathbf{F}_q[x]$ respectively, see [9], [15], etc, for more details. The product version of this formula has the form

$$C_N(\alpha) = \frac{\Phi(x^n - 1)}{q^n} \prod_{f(x) \mid x^n - 1} \left(1 - \frac{1}{q^{\deg(f(x))} - 1} \sum_{Ord(\psi) = f(x)} \psi([x^n - 1)/f(x)] \circ \alpha)\right), \tag{5}$$

where $f(x)$ runs through the irreducible factors of $x^n - 1$, see [15], etc.

***Example* 5.** For the parameter $n = p^u$, $q = p^v$, the polynomial $x^n - 1 = (x - 1)^n \in \mathbf{F}_q[x]$ and the expression $[(x^n - 1)/g(x)] \circ \alpha = Tr(\alpha)$ is the trace $Tr : \mathbf{F}_{q^n} \rightarrow \mathbf{F}_q$, (since $f(x) = x - 1$ is the only irreducible factor of $x^n - 1$). Thus the characteristic function of normal elements in $\mathbf{F}_{q^n}$ is

$$C_N(\alpha) = (1 - 1/q)\left(1 - \frac{1}{q-1} \sum_{Ord(\psi) = x-1} \psi(Tr(\alpha))\right) = \begin{cases} 0 & \text{if } Tr(\alpha) = 0, \\ 1 & \text{if } Tr(\alpha) \neq 0. \end{cases} \tag{6}$$

*The Characteristic Function of Primitive Normal Elements*

The product of the characteristic function of primitive elements (2) and the characteristic function of normal elements (4) in $\mathbf{F}_{q^n}$ realizes the characteristic function of primitive normal elements. Specifically

$$C_{PN}(\alpha) = \frac{\varphi(q^n - 1)}{q^n - 1} \frac{\Phi(x^n - 1)}{q^n} \sum_{d \mid q^n - 1} \frac{\mu(d)}{\varphi(d)} \sum_{f(x) | x^n - 1} \frac{M(f(x))}{\Phi(f(x))} \sum_{Ord(\psi) = f(x)} \sum_{ord(\chi) = d} \chi(\alpha)\psi(\beta), \tag{7}$$

where $\alpha \in \mathbf{F}_{q^n}$, and $\beta = [(x^n - 1)/f(x)]\circ\alpha$.

The function $C_{PN}(\alpha)$ is one of the basic tools used in the investigation of the distribution of primitive normal elements in finite fields. Typical applications are illustrated in [17], [4], etc, and in this paper.

*The Characteristic Function of Completely Normal Elements*

An element $\eta \in \mathbf{F}_{q^n}$ is *completely normal* if and only if $\eta$ is a normal element in $\mathbf{F}_{q^n}$ over $\mathbf{F}_{q^d}$ for all $d \mid n$. The characteristic function of completely normal elements in $\mathbf{F}_{q^n}$ is constructed from the product of the characteristic functions of normal elements in $\mathbf{F}_{q^n}$ over $\mathbf{F}_{q^d}$ for all $d \mid n$. Here $\mathbf{F}_{q^n}$ is an extension of $\mathbf{F}_{q^d}$ of degree $e = [\mathbf{F}_{q^n} : \mathbf{F}_{q^d}]$, with $n = de$. Thus it follows that the characteristic function of completely normal elements is the product of the individual functions:

$$\begin{aligned} C_{CN}(\alpha) &= \prod_{e|n}\left(\frac{\Phi(x^e-1)}{q^e}\sum_{f(x)|x^e-1}\frac{M(f(x))}{\Phi(f(x))}\sum_{Ord(\psi)=f(x)}\psi([x^e-1)/f(x)]\circ\alpha)\right) \\ &= \prod_{e|n}\prod_{f(x)|x^e-1}\frac{\Phi(x^e-1)}{q^e}\left(1-\frac{1}{q^{\deg(f(x))}-1}\sum_{Ord(\psi)=f(x)}\psi([x^e-1)/f(x)]\circ\alpha)\right), \end{aligned} \tag{8}$$

where $\alpha \in \mathbf{F}_{q^n}$, and $f(x)$ runs through the irreducible factors of $x^e - 1 \in \mathbf{F}_{q^d}[x]$.

***Example* 6.** For the parameter $n = p^u$, $q = p^v$, the polynomial $x^n - 1 = (x - 1)^n \in \mathbf{F}_{q^{p^i}}[x]$ and the expression $[(x^n - 1)/f(x)]\circ\alpha = Tr_i(\alpha)$ is the trace $Tr_i : \mathbf{F}_{q^{p^u}} \to \mathbf{F}_{q^{p^i}}$, $0 \le i < u$, (since $f(x) = x - 1$ is the only irreducible factor of $x^n - 1$). Thus the characteristic function of completely normal elements in $\mathbf{F}_{q^n}$ is

$$C_{CN}(\alpha) = \prod_{i=0}^{u-1}\left(1-1/q^{p^i}\right)\left(1-\frac{1}{q^{p^i}-1}\sum_{Ord(\psi)=x-1}\psi(Tr(\alpha))\right). \tag{9}$$

*Some Probabilities Formulae*

The probabilities $P_1$ = P(ord($\alpha$)=$q^n$−1) and $P_2$ = P(Ord($\alpha$)=$x^n$−1) respectively of primitive elements and normal elements $\alpha$ in a finite field extension $\mathbf{F}_{q^n}$ of $\mathbf{F}_q$ are given by

$$P_1 = \frac{\varphi(q^n-1)}{q^n-1} = \prod_{p|q^n-1}(1-1/p) \ge \frac{e^{-\gamma}}{\log(q^n-1)} \tag{10}$$

where $p$ ranges over the prime divisors of $q^n - 1$, and

$$P_2 = \frac{\Phi(x^n - 1)}{q^n} = \prod_{f(x) | x^n - 1} (1 - 1/p^{\deg(f)}) \geq (1 - 1/q)^n, \tag{11}$$

where $f(x)$ ranges over the irreducible divisors of $x^n - 1$ respectively. The degree $\deg(f) = d$ of each irreducible factor $f(x)$ is a divisor of $n$.

The distribution of primitive normal elements is more intricate than either the distribution of primitive elements or the distribution of normal elements. An exact closed form formula for the number of primitive normal bases of $\mathbf{F}_{q^n}$ over $\mathbf{F}_q$ appears to be unknown. However there are asymptotic approximations. For example, the probabilities $P_3$ = P(ord(α)=$q^n$–1 and Ord(α)=$x^n$–1) of primitive normal elements is approximated by

$$P_3 \approx \frac{\varphi(q^n - 1)}{q^n - 1} \frac{\Phi(x^n - 1)}{q^n} \geq \frac{e^{-\gamma}(1 - 1/q)^n}{\log(q^n - 1)}. \tag{12}$$

It is clear that if $q$ is sufficiently large, $P_1$ and $P_3$ are essentially the same. In fact the probabilities $P_1$ and $P_2$ are asymptotically independent. This means that the cardinalities of the sets of primitive polynomials and primitive normal polynomials, namely,

$$\frac{\varphi(q^n - 1)}{n} \sim \frac{e^{-\gamma}(q^n - 1)}{n\log(q^n - 1)} \quad \text{and} \quad \frac{\varphi(q^n - 1)\Phi(x^n - 1)}{n} \sim \frac{e^{-\gamma}(q - 1)^n}{n\log(q^n - 1)}, \tag{13}$$

are very close. Thus the results for the coefficients of primitive polynomials and primitive normal polynomials are about the same. Another approximation due to [5] for the number of primitive normal elements is

$$PN_n(q) = \frac{\varphi(q^n - 1)\Phi(x^n - 1)}{q^n} + O(q^{(.5+\varepsilon)n}), \tag{14}$$

for all ε > 0.

*Estimate of exponential sums*

A few estimates are required in order to derive nontrivial results on the distribution of elements in finite fields.

***Theorem* 7.** If $\chi \neq 1$ is a nontrivial multiplicative character on $\mathbf{F}_{q^n}$, then

$$\left| \sum_{Tr(\xi) \neq 0} \chi(\xi) \right| \leq q^{(n-1)/2}. \tag{15}$$

***Theorem* 8.** Let ψ and χ be a pair of nontrivial additive and multiplicative characters on $\mathbf{F}_{q^n}$, and let $f(x)$, $g(x) \in \mathbf{F}_q[x]$ be polynomials of degrees $k$ and $m$ respectively. Then

$$\left| \sum_{x \in \mathbf{F}_{q^n}} \psi(f(x))\chi(g(x))) \right| \le (k+m-1)q^{n/2}, \tag{16}$$

where $f(x)$, $g(x)$ are not $q$ powers and $k + m$ is the number of distinct roots of $f$ and $g$ in the splitting field, see [22], [24], and [25].

*The point counting function*

The point counting function

$$S_a(\xi) = \frac{1}{q} \sum_{\xi \in F_{q^n}} \sum_{x \in F_q} \psi(x(f(\xi) - a)) \tag{17}$$

enumerates the cardinality of the solution set { $\xi \in \mathbf{F}_{q^n}$ : $\mathrm{Tr}(f(\xi)) - a = 0$ } of the equation $\mathrm{Tr}(f(\xi)) - a = 0$ in $\mathbf{F}_{q^n}$, where $f(x)$ is a function on $\mathbf{F}_{q^n}$, and $a \in \mathbf{F}_q$ is a constant.

*Newton identities in finite fields*

The coefficients of the polynomial $f(x) = x^n + a_1x^{n-1} + \cdots + a_{n-1}x + a_n$ are given by

$$a_1 = -\sum_{1 \le i \le n} z_i,\ a_2 = \sum_{1 \le i_1 < i_2 \le n} z_{i_1}z_{i_2},\ \ldots,\ a_i = (-1)^i \sum_{1 \le i_1 < \cdots < i_i \le n} z_{i_1} \cdots z_{i_i},\ \ldots,\ a_n = (-1)^n \prod_{1 \le i \le n} z_i, \tag{18}$$

where $z_1, z_2, \ldots, z_n$ are its roots. The symmetric functions $\sigma_i(z_1,\ldots,z_n)$ and the coefficients $a_i$ are equal up to a sign, that is, $a_i = (-1)^i\sigma_i(z_1,\ldots,z_n)$. The associated power sums are defined by

$$w_i = \sum_{1 \le j \le n} z_j^i = z_1^i + z_2^i + \cdots + z_n^i, \quad i \ge 1. \tag{19}$$

The cyclic structure of the roots of polynomials in cyclic extensions can be utilized to transform the power sums into different forms. Specifically, in finite fields the power sums become

$$w_i = \sum_{1 \le j \le n} z_j^i = \sum_{1 \le j \le n} \alpha^{iq^j} = Tr(\alpha^i), \quad i \ge 1, \tag{20}$$

where $z_i = \alpha^{q^j}$ some $j = 0, 1, \ldots, n - 1$, and $f(\alpha) = 0$. Replacing the power sums in Newton identities (1707) yields the corresponding identities for finite fields:

(1) $$Tr(\alpha^k) = -ka_k - \sum_{i=1}^{k-1} a_i Tr(\alpha^{k-i}), \quad 1 \le k \le n, \tag{21}$$

(2) $Tr(\alpha^k) = -\sum_{i=1}^{k-1} a_i Tr(\alpha^{k-i}), \quad n < k.$

These formulae are employed to solve problems about the coefficients of polynomials via the power sums. Fast algorithms for computing the coefficients from the power sums and conversely in $O(n\log(n))$ operations are discussed in [2]. A few of the coefficients are given here in terms of the power sums $w_i = \mathrm{Tr}(\alpha^i)$.

(1) $a_1 = -w_1 = -Tr(\alpha) = -(\alpha + \alpha^q + \alpha^{q^2} + \cdots + \alpha^{q^{n-1}})$,

(22)

(2) $a_2 = (2!)^{-1}(w_1^2 - w_2)$,

(3) $a_3 = (3!)^{-1}(-w_1^3 - w_1 w_2 + 2w_3)$,

(4) $a_4 = (4!)^{-1}(2w_1^4 + 4w_1^2 w_2 - 3w_2^2 + 4w_1 w_3 - 6w_4)$,

($k$) $a_k = (k!)^{-1}\sum_e a_e w_1^{e_1} w_2^{e_2} \cdots w_k^{e_k}, \quad e = (e_1,...,e_k), e_i \geq 0, \; k \leq n$.

These formula are defined in characteristic $p > k$. A direct calculation of the second coefficient

$$a_2 = \sum_{0 \leq i < j < n} \alpha^{q^i + q^j} = (2!)^{-1}(Tr(\alpha)^2 - Tr(\alpha^2)) \tag{23}$$

in characteristic $p > 2$ can be accomplished in a few lines. But direct calculations of the other coefficients appear to be a lengthy and difficult task. Even the third coefficient

$$a_3 = \sum_{0 \leq i < j < k < n} \alpha^{q^i + q^j + q^k} = (3!)^{-1}(2Tr(\alpha)^3 - Tr(\alpha)Tr(\alpha^2) - Tr(\alpha^3)) \tag{24}$$

in characteristic $p > 3$ is difficult to compute directly. But using Newton identities in finite fields (21) this task is a simple algebraic manipulation.

*Coefficients system of equations*

Let $c_1, c_2, ..., c_k \in \mathbf{F}_q$ be constants. The strategy of the analysis is to show that the system of equations

$$Tr(x) = c_1,\ Tr(x^2) = c_2,\ ...,\ Tr(x^k) = c_k \tag{25}$$

has at least one primitive normal element solution $x = \alpha \in \mathbf{F}_{q^n}$. This in turn implies the existence of at least one primitive normal polynomial $f(x)$ with $k$ prescribed consecutive coefficients. The constants $c_1, c_2, ..., c_k$ indirectly prescribes the $k$ coefficients, namely,

$$a_1 = -Tr(\alpha),\ a_2 = a_2(Tr(\alpha),Tr(\alpha^2)),\ ...,\ a_k = a_k(Tr(\alpha),\mathrm{Tr}(\alpha^2),...,Tr(\alpha^k)). \tag{26}$$

The derivation of (24) from (23), which is a direct consequence of Newton identities in finite fields. The application of (23), (24), and its use in the exponential sum (17) are the key ideas in the proofs of several problems on the existence and distribution of the coefficients of certain polynomials over finite fields. I observed the relationship between the trace function and Newton identities in the early 1990's while investigating the trace representations of sequences, and applied it to the theory of coefficients of primitive polynomials in [4].

**3 The Extended Coefficient Theorem**

There are various ways of extending the result on primitive polynomials of degree *n* with *k* prescribed coefficients in the finite field $\mathbf{F}_q$ given in [4]. The refinement given here extends it to primitive normal polynomials. The constraint of linearly independent roots on the roots of polynomials restricts the set of primitive normal polynomials to a proper subset of set primitive polynomials. There is one exception; the set of quadratic primitive polynomials and the set of quadratic primitive normal polynomials coincide.

The first line of the proof given below, using a reductio ad absurdum argument, combines the various characteristic functions and the point counting function:

$$\sum_{Tr(\alpha)\neq 0,\ \alpha\in\mathbf{F}_{q^n}} C_N(\alpha)C_P(\alpha)\prod_{i=1}^{k} S_{c_i}(\alpha)=0\,. \tag{27}$$

This claims that there is no primitive normal element solution of equation (23). The usual and standard strategy of dealing with this type of equation is to decompose it into several sums, and then compute a lower estimate, which contradicts the equality, consult [3], [4], [5], [8], [13], [14], [15], [16], [20], and [21] for background details and other references on this type of analysis.

*Proof of Theorem* 1: The number of primitive normal polynomials with *k* prescribed consecutive coefficients $a_1 \neq 0, a_2, \ldots, a_k \in \mathbf{F}_q$, $q$ odd, is given by

$$N^*(n,q,c_1,\ldots,c_k)=\sum_{Tr(\xi)\neq 0,\xi\in\mathbf{F}_{q^n}}\left(\frac{\varphi(q^n-1)}{q^n-1}\sum_{d\mid q^n-1}\frac{\mu(d)}{\varphi(d)}\sum_{ord(\chi)=d}\chi(\xi)\right)$$
$$\times\left(\frac{\Phi(x^n-1)}{q^n-1}\sum_{f\mid x^n-1}\frac{M(f)}{\Phi(f)}\sum_{Ord(\psi)=f}\psi(Tr(\xi))\right)\left(\prod_{i=1}^{k}\frac{1}{q}\sum_{x_i\in\mathbf{F}_q}\psi(x_i(Tr(\xi)-c_i))\right), \tag{28}$$

One possible approach to deal with (26) is to decompose it according to the order $\mathrm{Ord}(\psi) = 1$ or $\mathrm{Ord}(\psi) \neq 1$ of the additive character $\psi$. To accomplish this, rewrite it as

(29)

$$N^*(n,q,c_1,...,c_k) = \frac{P_1P_2}{q^k} \sum_{f|x^n-1} \frac{M(f)}{\Phi(f)} \sum_{Ord(\psi)=f} \sum_{d|q^n-1} \frac{\mu(d)}{\varphi(d)}$$

$$\times \sum_{ord(\chi)=d} \sum_{x_i\in\mathbf{F}_q} \sum_{Tr(\xi)\neq 0,\xi\in\mathbf{F}_{q^n}} \chi(\xi)\psi(Tr(\xi + \sum_{i=1}^{k} x_i\xi^i) - \sum_{i=1}^{k} c_i x_i).$$

Now separate (27) into two terms corresponding to Ord(ψ) = 1 or Ord(ψ) ≠ 1, and simplify to obtain

(30)

$$N^*(n,q,c_1,...,c_k) = P_1P_2\left( q^{n-1}(q-1) - \sum_{Tr(\xi)\neq 0,\xi\in\mathbf{F}_{q^n}} \chi(\xi) \right) - P_2N(n,q,c_1,...,c_k)$$

$$= P_1P_2\left( q^{n-1}(q-1) - \frac{1}{P_1}N(n,q,c_1,...,c_k) - \sum_{Tr(\xi)\neq 0,\xi\in\mathbf{F}_{q^n}} \chi(\xi) \right),$$

where $q^{n-1}(q-1)$ is the number of elements in $\mathbf{F}_{q^n}$ of nonzero traces, and $N(n,q,c_1,\ldots, c_k)$ denotes the total number of primitive polynomials $f(x) = x^n + a_1x^{n-1} + \cdots + a_{n-1}x + a_n \in \mathbf{F}_q[x]$ with $k$ prescribed consecutive coefficients $a_1 \neq 0, a_2, \ldots, a_k \in \mathbf{F}_q$.

Replacing the estimates $q^k \leq N(n,q,c_1,\ldots, c_k) \leq q^{n-s}$, $\prod_{p|q^n-1}(1-1/p)^{-1} = P_1^{-1} < c_\varepsilon q^{n\varepsilon}$, where $s > 0$ is a small integer, $\varepsilon > 0$, and $c_\varepsilon$ is a constant, and the exponential sum estimate yields

$$N^*(n,q,c_1,...,c_k) \geq P_1P_2q^{n-1}\left( q - 1 - \frac{1}{q^{(n-1)/2}} - \frac{c_\varepsilon}{q^{s-\varepsilon n-1}} \right). \qquad (31)$$

Since the product of the probabilities is in the range $0 < P_1P_2 < 1$, it readily follows that if $k < p$, then there exists a constant $q_1$ such that $N^*(n,q,c_1,\ldots, c_k) \geq 1$ for all $q \geq q_1$, and $s < \varepsilon n + 1$. ■

This approach leads to a relatively easy proof of the Theorem 1, it avoids the need to deal with the estimates of various exponential sums, and the function $\Omega(x^n-1)$ which enumerates the distinct irreducible factors of $x^n - 1 \in \mathbf{F}_q[x]$. The case $k = 1$ reduces to the primitive normal basis theorem of arbitrary trace. In fact this is a simpler proof than the first proof given in [3] using a variation of this technique.

**T**